%% file: main.tex
\documentclass[a4paper,11pt,leqno]{amsart}

\input{Preamble/Preamble}

\usepackage{subfiles}

\title{Path decompositions of tournaments}

\begin{document}
\renewcommand{\onlyinsubfile}[1]{}
\renewcommand{\notinsubfile}[1]{#1}

\begin{abstract}
	\subfile{Abstract}
\end{abstract}
\maketitle

\section{Introduction}\label{sec:intro}
\subfile{Introduction}

\section{Proof overview}\label{sec:sketch}
\subfile{Sketch}

\section{Notation}\label{sec:notation}
\subfile{Notation}

\section{Preliminaries}\label{sec:preliminaries}

In this section, we introduce some tools which will be used throughout the rest of the paper.

\subsection{Robust outexpanders}\label{sec:rob}
\subfile{Robust_Outexpanders}

\subsection{Probabilistic estimates}
\subfile{Probabilistic_Estimates}

\subsection{Some tools for finding matchings}
\subfile{Matchings}

\subsection{Some properties of the excess function}
\subfile{Excess}

\section{Exceptional tournaments}\label{sec:annoyingT}
\subfile{Exceptional_Tournaments}

\section{Deriving Theorem \ref{thm:even} and Corollary \ref{cor:approx} from Theorem \ref{thm:main}}\label{sec:cor}
\subfile{Corollaries}

\section{Approximate decomposition of robust outexpanders}\label{sec:approxdecomp}
\subfile{Approximate_Decomposition}

\section{Good partial path decompositions and absorbing edges}\label{sec:good}
\subfile{Good_Partial_Decompositions}

\section{Constructing layouts in general tournaments}\label{sec:constructinglayouts}
\subfile{Constructing_Layouts_Definitions_Statements}

\section{Deriving Theorem \ref{thm:main}}\label{sec:main}
\subfile{Proof_of_Main_Theorem}

\section{Auxiliary excess function}\label{sec:texc}
\subfile{Auxiliary_Excess_Function}

\section{The cleaning step: proof of Lemma \ref{lm:cleaning}}\label{sec:cleaning}
\subfile{Cleaning}

\section{Constructing layouts: proof of Lemma \ref{lm:layouts}}\label{sec:layouts}
\subfile{Constructing_Layouts_Proof}

\section{Concluding remarks}\label{sec:conclusion}
\subfile{Concluding_Remarks}

\section*{\NEWTWO{Acknowledgements}}
\NEWTWO{We thank the referee for helpful suggestions.}

\bibliographystyle{abbrv}
\bibliography{Bibliography/Bibliography}

\end{document}

%% file: Preamble/Preamble.tex
\usepackage[shortlabels]{enumitem}
\usepackage{xr-hyper, zref}
\usepackage[hypertexnames=false]{hyperref}
\usepackage{etoolbox}
\usepackage{bm}
\usepackage{amsthm,amssymb,amsmath}
\usepackage[labelformat=simple]{subcaption}

\usepackage{tikz,graphicx,placeins}
\usetikzlibrary{bending,decorations.pathmorphing}
\usetikzlibrary{shapes.geometric}
\usetikzlibrary{arrows}
\usetikzlibrary{arrows.meta}
\tikzset{>={Latex[width=2mm,length=2mm]}}
\usepackage{psfrag,mathrsfs}
\usepackage{mathtools}
\usepackage[foot]{amsaddr}
\usepackage[british]{babel}
\usepackage[babel]{microtype}
\usepackage[margin=1in]{geometry}
\usepackage[capitalise]{cleveref}
\usepackage[textsize=footnotesize]{todonotes}
\usepackage{environ}
\usepackage{float}
\usepackage{bigfoot}
\usepackage[noadjust]{cite}
\usepackage{dsfont}
\usepackage{moreenum}
\usepackage{eqparbox}

\date{}
\address[A. Gir\~{a}o]{\NEW{Mathematical Institute, University of Oxford, Oxford, OX2 6GG, United Kingdom}}
\address[B. Granet]{\NEWTWO{Institut f\"ur Informatik, Universit\"at Heidelberg, 69120 Heidelberg, Deutschland}}
\address[D. K\"{u}hn, A. Lo, and D. Osthus]{School of Mathematics, University of Birmingham, Edgbaston, Birmingham, B15 2TT, United Kingdom}

\author[Ant\'{o}nio~Gir\~{a}o]{Ant\'{o}nio Gir\~{a}o}
\email{\NEW{girao@maths.ox.ac.uk}}

\author[Bertille~Granet]{Bertille Granet}
\email{\NEWTWO{granet@informatik.uni-heidelberg.de}}

\author[Daniela~K\"{u}hn]{Daniela K\"{u}hn}
\email{d.kuhn@bham.ac.uk}

\author[Allan~Lo]{Allan Lo}
\email{s.a.lo@bham.ac.uk}

\author[Deryk~Osthus]{Deryk Osthus}
\email{d.osthus@bham.ac.uk}

\subjclass[2010]{\NEW{05C20, 05C38, 05B40, 05D40}}

\thanks{This project has received partial funding from the European Research 
 	Council (ERC) under the European Union's Horizon 2020 research and innovation programme (grant agreement no. 786198, D.~K\"{u}hn and D.~Osthus).
 	The research leading to these results was also partially supported by the EPSRC, grant nos. EP/N019504/1 (A.~Gir\~{a}o and D.~K\"{u}hn) and EP/S00100X/1 (D.~Osthus).}

\newcommand{\COMMENT}[1]{}
\newcommand{\APPENDIX}[1]{}
\newcommand{\NOAPPENDIX}[1]{#1}
\renewcommand{\APPENDIX}[1]{#1}\renewcommand{\NOAPPENDIX}[1]{}% comment out for the appendix-free version

\newcommand{\NEW}[1]{#1}
\newcommand{\OLD}[1]{}
\newcommand{\NEWTWO}[1]{#1}
\newcommand{\OLDTWO}[1]{}
\newcommand{\onlyinsubfile}[1]{#1}
\newcommand{\notinsubfile}[1]{}

\hypersetup{colorlinks=true,
	citecolor=blue,
	filecolor=blue,
	linkcolor=blue,
}

\SetEnumitemKey{longlabel}{wide=0.5cm, leftmargin=*, align=right}
\setlist[enumerate]{itemsep=3pt, topsep=5pt,leftmargin=1.2cm}
\setlist[itemize]{itemsep=3pt, topsep=2pt}

\newlist{steps}{enumerate}{1}
\setlist[steps,1]{
	label=\textbf{Step \arabic*:},
	ref=\arabic*, % if to be used with \cref, don't provide label string or parentheses
	wide,
	parsep=0pt,
	itemsep=10pt,
	topsep=10pt}
\Crefname{stepsi}{Step}{Steps}
\crefname{stepsi}{Step}{Steps}

\newlist{case}{enumerate}{1}
\setlist[case,1]{
	label=\textbf{Case \arabic*:},
	ref=\arabic*, % if to be used with \cref, don't provide label string or parentheses
	wide,
	parsep=0pt,
	itemsep=10pt,
	topsep=10pt}
\Crefname{casei}{Case}{Cases}
\crefname{casei}{Case}{Cases}

\AtBeginEnvironment{thm}{\setlist[enumerate,1]{label=\upshape(\alph*)}}
\AtBeginEnvironment{lm}{\setlist[enumerate,1]{label=\upshape(\alph*)}}
\AtBeginEnvironment{prop}{\setlist[enumerate,1]{label=\upshape(\alph*)}}
\AtBeginEnvironment{cor}{\setlist[enumerate,1]{label=\upshape(\alph*)}}
\AtBeginEnvironment{claim}{\setlist[enumerate,1]{label=\upshape(\alph*)}}
\AtBeginEnvironment{fact}{\setlist[enumerate,1]{label=\upshape(\alph*)}}

\setlist[itemize,1]{label=--}

\creflabelformat{equation}{#2(#1)#3}
\crefdefaultlabelformat{#2\textup{#1}#3}

\Crefname{enumi}{}{}
\Crefname{thm}{Theorem}{Theorems}
\Crefname{lm}{Lemma}{Lemmas}
\Crefname{cor}{Corollary}{Corollaries}
\Crefname{prop}{Proposition}{Propositions}
\Crefname{claim}{Claim}{Claims}
\Crefname{equation}{}{}
\Crefname{conjecture}{Conjecture}{Conjectures}
\Crefname{figure}{Figure}{Figures}
\Crefname{fact}{Fact}{Facts}

\newtheorem{definition}{Definition}[section]
\newtheorem{claim}{Claim}
\newtheorem{prop}[definition]{Proposition}
\newtheorem{thm}[definition]{Theorem}
\newtheorem{cor}[definition]{Corollary}
\newtheorem{lm}[definition]{Lemma}
\newtheorem{fact}[definition]{Fact}
\newtheorem{conjecture}[definition]{Conjecture}

\NewEnviron{property}[1]{\begin{equation}
	\tag{\(#1\)}
	\parbox{0.9\linewidth}{{\itshape\BODY}}
	\end{equation}}

\newenvironment{proofclaim}[1][Proof of Claim]{\begin{proof}[#1]}{\renewcommand{\qedsymbol}{\rotatebox[origin=c]{45}{\(\square\)}}\end{proof}}

\newcounter{claimnumber}
\AtBeginEnvironment{proof}{\setcounter{claim}{0}}
\AtBeginEnvironment{proofclaim}{\setcounter{claimnumber}{\value{claim}}}
\AtEndEnvironment{proofclaim}{\setcounter{claim}{\value{claimnumber}}}

\numberwithin{equation}{section}

\newcommand{\eps}{\varepsilon}
\renewcommand{\epsilon}{\varepsilon}

\DeclareMathOperator{\Var}{Var}
\DeclareMathOperator{\Bin}{Bin}
\DeclareMathOperator{\HGeom}{Hyp}
\DeclareMathOperator{\exc}{ex}
\DeclareMathOperator{\texc}{\widetilde{ex}}
\DeclareMathOperator{\hexc}{\widehat{ex}}

\DeclareMathOperator{\pn}{pn}

\newcommand{\cA}{\mathcal{A}}
\newcommand{\cB}{\mathcal{B}}
\newcommand{\cC}{\mathcal{C}}
\newcommand{\cD}{\mathcal{D}}

\newcommand{\cF}{\mathcal{F}}

\newcommand{\cH}{\mathcal{H}}

\newcommand{\cP}{\mathcal{P}}
\newcommand{\cQ}{\mathcal{Q}}

\newcommand{\cT}{\mathcal{T}}

\newcommand{\tD}{\widetilde{D}}
\newcommand{\tE}{\widetilde{E}}
\newcommand{\tF}{\widetilde{F}}

\newcommand{\tL}{\widetilde{L}}

\newcommand{\tU}{\widetilde{U}}
\newcommand{\tV}{\widetilde{V}}
\newcommand{\tW}{\widetilde{W}}

\newcommand{\hF}{\widehat{F}}

\newcommand{\hL}{\widehat{L}}

\newcommand{\hQ}{\widehat{Q}}

\newcommand{\hU}{\widehat{U}}

\newcommand{\tm}{\widetilde{m}}
\newcommand{\tn}{\widetilde{n}}

%% file: Abstract.tex
	In 1976, Alspach, Mason, and Pullman conjectured that any tournament~$T$ of even order can be decomposed into exactly~$\exc(T)$ paths, where $\exc(T)\coloneqq \frac{1}{2}\sum_{v\in V(T)}|d_T^+(v)-d_T^-(v)|$.
	We prove this conjecture for all sufficiently large tournaments. 
	We also prove an asymptotically optimal result for tournaments of odd order.

%% file: Introduction.tex
	\onlyinsubfile{\section{Introduction}}

Path and cycle decomposition problems have a long history.
For example, the Walecki construction~\cite{lucas1883recreationsII}, which goes back to the 19$^{\rm th}$ century, gives a decomposition of the complete graph of odd order into Hamilton cycles (see also~\cite{alspach2008wonderful}).
A version of this for (regular) directed tournaments was conjectured by Kelly in 1968 and proved for large tournaments by \NEW{K\"{u}hn and Osthus}\OLD{in} \cite{kuhn2013hamilton}.
Beautiful open problems in the area include the Erd\H{o}s--Gallai conjecture which asks for a decomposition of any graph into linearly many cycles and edges. The best bounds for this are due to Conlon, Fox, and Sudakov~\cite{conlon2014cycle}.
Another famous example is the linear arboricity conjecture, which asks for a decomposition of a~$d$-regular graph into~$\left\lceil\frac{d+1}{2}\right\rceil$ linear forests. The latter was resolved asymptotically by Alon~\cite{alon1988linear} and the best current bounds are due to Lang and Postle~\cite{lang2020improved}.

\subsection{Background}

The problem of decomposing digraphs into paths was first explored by Alspach and Pullman~\cite{alspach1974path}, who provided\OLD{sharp} bounds for the minimum number of paths needed in path decompositions of digraphs.
(Throughout this paper, in a digraph, for any two vertices~$u\neq v$, we allow a directed edge~$uv$ from~$u$ to~$v$ as well as a directed edge~$vu$ from~$v$ to~$u$, whereas in an oriented graph we allow at most one directed edge between any two distinct vertices.)
Given a digraph~$D$, define the \emph{path number of~$D$}, denoted by~$\pn(D)$, as the minimum integer~$k$ such that~$D$ can be decomposed into~$k$ paths.
Alspach and Pullman~\cite{alspach1974path} proved that, for any oriented graph~$D$ on~$n$ vertices,~$\pn(D)\leq \frac{n^2}{4}$, with equality holding for transitive tournaments. O'Brien~\cite{obrien1977upper} showed that the same bound holds for digraphs on at least~$4$ vertices.

The path number of digraphs can be bounded below by the following quantity.
Let~$D$ be a digraph and~$v\in V(D)$. Define the \emph{excess} at~$v$ as $\exc_D(v)\coloneqq d_D^+(v)-d_D^-(v)$. Let $\exc_D^+(v)\coloneqq \max\{0,\exc_D(v)\}$ and $\exc_D^-(v)\coloneqq \max\{0,-\exc_D(v)\}$ be the \emph{positive excess} and \emph{negative excess} at~$v$, respectively.
Then, as observed in~\cite{alspach1974path}, if $d_D^+(v)>d_D^-(v)$, a path decomposition of~$D$ contains at most~$d_D^-(v)$ paths which have~$v$ as an internal vertex, and thus at least $d_D^+(v)-d_D^-(v)=\exc_D^+(v)$ paths starting at~$v$. Similarly, a path decomposition will contain at least~$\exc_D^-(v)$ paths ending at~$v$.
Thus, the \emph{excess} of~$D$, defined as
\begin{equation}\label{eq:exc}
	\exc(D)\coloneqq \sum_{v\in V(D)}\exc_D^+(v)=\sum_{v\in V(D)}\exc_D^-(v)=\frac{1}{2}\sum_{v\in V(D)}|\exc_D(v)|,
\end{equation}
provides a natural lower bound for the path number of~$D$, i.e.\ any digraph~$D$ satisfies 
\begin{equation}\label{thm:pn>exc}
	\pn(D)\geq \exc(D).
\end{equation}
It was shown in~\cite{alspach1974path} that equality is satisfied for acyclic digraphs.
A digraph satisfying equality in \cref{thm:pn>exc} is called \emph{consistent}.
Clearly, not all digraphs are consistent (e.g. regular digraphs have excess~$0$).
However, Alspach, Mason, and Pullman~\cite{alspach1976path} conjectured in 1976 that tournaments of even order are consistent.

\begin{conjecture}[Alspach, Mason, and Pullman {\NEW{\cite{alspach1976path}}}]\label{conj:even}
	Let~$n\in \mathbb{N}$ be even. Then, any tournament~$T$ on~$n$ vertices satisfies~$\pn(T)= \exc(T)$.
\end{conjecture}

This conjecture is discussed also e.g.\ in the Handbook of Combinatorics~\cite{bondy1995basic}.%
\OLD{Moreover, it is listed as being of ``high importance" in the list of problems on the ``open problem garden" website.} 

Note that the results of Alspach and Pullman~\cite{alspach1974path} mentioned above imply that \cref{conj:even} holds for tournaments of excess~$\frac{n^2}{4}$.
Moreover, as observed by Lo, Patel, Skokan, and Talbot~\cite{lo2020decomposing}, \cref{conj:even}
for tournaments of excess~$\frac{n}{2}$ is equivalent to 
Kelly's conjecture on Hamilton decompositions of regular tournaments.
Recently, \cref{conj:even} was verified in~\cite{lo2020decomposing} for sufficiently large tournaments of sufficiently large excess. Moreover, they extended this result to tournaments of odd order~$n$ whose excess is at least~$n^{2-\frac{1}{18}}$.

\begin{thm}[{\cite{lo2020decomposing}}]\label{thm:evenlarge}
	The following hold.
	\begin{enumerate}
		\item There exists $C\in \mathbb{N}$ such that, for any tournament~$T$ of even order~$n$, if~$\exc(T)\geq Cn$, then~$\pn(T)=\exc(T)$.\label{thm:evenlarge-main}
		\item There exists~$n_0\in \mathbb{N}$ such that, for any~$n\geq n_0$, if~$T$ is a tournament on~$n$ vertices satisfying~$\exc(T)\geq n^{2-\frac{1}{18}}$, then~$\pn(T)=\exc(T)$.\label{thm:even-verylarge}
	\end{enumerate}
\end{thm}

\subsection{New results}

Building on the results and methods of~\cite{lo2020decomposing,kuhn2013hamilton}, we prove \cref{conj:even} for large tournaments.

\begin{thm}\label{thm:even}
	There exists~$n_0\in \mathbb{N}$ such that, for any even~$n\geq n_0$, any tournament~$T$ on~$n$ vertices satisfies~$\pn(T)=\exc(T)$.
\end{thm}

In fact, our methods are more general and allow us \NEW{to} determine the path number of most tournaments of odd order, whose behaviour turns out to be more complex.
As mentioned above, not every digraph is consistent.

Let~$D$ be a digraph. Let~$\Delta^0(D)$ denote the largest semidegree of~$D$, that is $\Delta^0(D) \coloneqq \max\{ d^+(v), d^-(v) \mid v \in V(D) \}$.
Note that~$\Delta^0(D)$ is a natural lower bound for~$\pn(D)$ as every vertex~$v \in V(D)$ must be in at least $\max\{ d^+(v), d^-(v) \}$ paths.
This leads to the \NEW{notion}\OLD{notation} of the \emph{modified excess} of a
digraph~$D$, which is defined as \[\texc(D)\coloneqq \max\{\exc(D), \Delta^0(D)\}.\] This provides a natural lower bound for the path number of any digraph~$D$.

\begin{fact}\label{prop:pn>texc}
	Any digraph~$D$ satisfies~$\pn(D)\geq \texc(D)$.
\end{fact}

\NEWTWO{(Note that one can easily verify that any tournament $T$ of even order satisfies $\texc(T) = \exc(T)$ (see e.g.\ \cref{prop:even}),	so \cref{prop:pn>texc} is consistent with \cref{conj:even}.)} 

Observe that, by \cref{thm:evenlarge}\cref{thm:even-verylarge}, equality holds for large tournaments of excess at least~$n^{2-\frac{1}{18}}$.
However, note that equality does not hold for regular digraphs.
(Here a digraph is~$r$-regular if for every vertex, both its in- and outdegree equal~$r$.)
Indeed, by considering the number of edges, one can show that any path decomposition of an~$r$-regular digraph will contain at least~$r+1$ paths. Thus, any \NEW{$r$-}regular digraph satisfies 
\begin{equation}\label{eq:pnreg}
	\pn(D)\geq \NEW{r+1=}\texc(D)+1.
\end{equation}
Alspach, Mason, and Pullman~\cite{alspach1976path} conjectured that equality holds \NEWTWO{in \cref{eq:pnreg} whenever $D$ is a regular tournament}\OLDTWO{for regular tournaments}.
We verify this conjecture for sufficiently large tournaments.

\begin{thm}\label{thm:reg}
	There exists~$n_0\in \mathbb{N}$ such that any regular tournament~$T$ on~$n\geq n_0$ vertices satisfies $\pn(T)=\frac{n+1}{2}=\texc(T)+1$.
\end{thm}

In fact, our argument also applies to regular oriented graphs of large \NEW{enough} degree. 

\begin{thm}\label{thm:oriented3/8}
	For any~$\varepsilon>0$, there exists~$n_0\in \mathbb{N}$ such that, if~$D$ is an~$r$-regular oriented graph on~$n\geq n_0$ vertices satisfying $r\geq \left(\frac{3}{8}+\varepsilon\right)n$, then $\pn(D)=r+1=\texc(D)+1$.
\end{thm}

More generally, we will see that \cref{thm:oriented3/8} can be extended to regular digraphs of linear degree which are ``robust outexpanders'' (see \cref{cor:reg}).

There also exist non-regular tournaments for which equality does not hold in \cref{prop:pn>texc}.
Indeed, let~$\cT_{\rm apex}$ be the set of tournaments~$T$ on~$n\geq 5$ vertices for which there exists a partition $V(T)=V_0\cup\{v_+\}\cup \{v_-\}$ such that~$T[V_0]$ is a regular tournament on~$n-2$ vertices (and so~$n$ is odd), $N_T^+(v_+)= V_0=N_T^-(v_-)$, $N_T^-(v_+)=\{v_-\}$, and $N_T^+(v_-)=\{v_+\}$.
The tournaments in~$\cT_{\rm apex}$ are called \emph{apex tournaments}.
We show that any sufficiently large tournament~$T\in \cT_{\rm apex}$ satisfies $\pn(T)=\texc(T)+1$ (see \cref{thm:annoyingT}).
Denote by~$\cT_{\rm reg}$ the class of regular tournaments and let $\cT_{\rm excep}\coloneqq \cT_{\rm apex}\cup\cT_{\rm reg}$. The tournaments in~$\cT_{\rm excep}$ are called \emph{exceptional}.
We conjecture that the tournaments in~$\cT_{\rm excep}$ are the only ones which do not satisfy equality in \cref{prop:pn>texc}.

\begin{conjecture}\label{conj:all}
	There exists~$n_0\in \mathbb{N}$ such that any tournament $T\notin \cT_{\rm excep}$ on~$n\geq n_0$ vertices satisfies~$\pn(T)= \texc(T)$.
\end{conjecture}

We prove an approximate version of this conjecture (see \cref{cor:approx}). Moreover, in \cref{thm:main}, we prove \cref{conj:all} exactly unless~$n$ is odd and~$T$ is extremely close to being a regular tournament\OLD{(in the sense that the number of vertices of nonzero excess is~$o(n)$, the excess of each vertex is~$o(n)$, and the total excess is~$\frac{n}{2}\pm o(n)$)}.

\begin{thm}\label{thm:main}
	For all~$\beta>0$, there exists~$n_0\in \mathbb{N}$ such that the following holds. If~$T$ is a tournament on~$n\geq n_0$ vertices such that~$T\notin \cT_{\rm excep}$ and
	\begin{enumerate}
		\item $\texc(T)\geq \frac{n}{2}+\beta n$, or\label{thm:main-largeexc}
		\item $N^+(T), N^-(T)\geq \beta n$, where $N^+(T)\coloneqq |\{v\in V(T)\mid \exc_T^+ (v)>0\}|+\texc(T)-\exc(T)$ and $N^-(T)\coloneqq |\{v\in V(T)\mid \exc_T^- (v)>0\}|+\texc(T)-\exc(T)$,\label{thm:main-largeU}
	\end{enumerate}
	then~$\pn(T)=\texc(T)$.
\end{thm} 

In \NEW{\cref{sec:cor}}\OLD{\cref{sec:main}}, we will derive \cref{thm:even} (i.e.\ the exact solution when~$n$ is even) from \cref{thm:main}.
This will make use of the fact that~$\texc(T)=\exc(T)$ for~$n$ even (see \cref{prop:even}).
We will also derive an approximate version of \cref{conj:all} from \cref{thm:main}.

\begin{cor}\label{cor:approx}
	For all~$\beta>0$, there exists~$n_0\in \mathbb{N}$ such that, for any tournament~$T$ on~$n\geq n_0$ vertices, $\pn(T)\leq \texc(T)+\beta n$.
\end{cor}

Note that \cref{thm:main}\cref{thm:main-largeU} corresponds to the case where linearly many different vertices can be used as endpoints of paths in a \NEW{path}\OLD{n optimal} decomposition \NEW{of size $\texc(T)$}.
Indeed, let~$T$ be a tournament and~$\cP$ be a path decomposition of~$T$. Then, as mentioned above, each~$v\in V(T)$ must be the starting point of at least~$\exc_T^+(v)$ paths in~$\cP$.
Thus, for any tournament~$T$,~$N^+(T)$ is the maximum number of distinct vertices which can be a starting point of a path in a \NEW{path} decomposition of~$T$ of size~$\texc(T)$ and similarly for~$N^-(T)$ and the ending points of paths.

	\COMMENT{Moreover, note that, for a tournament~$T$ of odd order,~$\exc(T)$ is even ($n$ odd implies~$n-1$ even so~$\exc_T(v)$ is even for all~$v$) and, by \cref{thm:pn>exc} and a result of~\cite{alspach1974path} mentioned above,~$\exc(T) \leq pn(T)\leq \frac{n^2}{4}$.
	Thus, \cref{thm:main,thm:annoyingT} verify the following conjecture of~\cite{alspach1976path} for large~$n$. (Note that $\Delta^0(T),\delta^0(T)+1\leq n-1$ and~$\pn(T)\in \mathbb{N}$.)
	\begin{conjecture}[Alspach, Pullman, and Mason]\label{conj:odd}
		Let~$n\in \mathbb{N}$ be odd and~$T$ be a tournament on~$n$ vertices. Then,~$\pn(T)=k$ for some integer~$\frac{n+1}{2}\leq k\leq n-1$ or even integer~$n-1\leq k\leq \frac{n^2}{4}$.
	\end{conjecture}
	}

One can show that almost all large tournaments satisfy $\exc(T)=n^{\frac{3}{2}+o(1)}$%
.
Indeed, consider a tournament~$T$ on~$n$ vertices, where the orientation of each edge is chosen uniformly at random,
independently of all other orientations. For the upper bound on~$\exc(T)$, one can simply apply a Chernoff bound to show that for a given vertex~$v$ and~$\eps>0$, we have
$\exc^+_T(v)\le n^{\frac{1}{2}+\eps}$ with probability $1-o\left(\frac{1}{n}\right)$. The result follows by a union bound over all vertices. 
For the lower bound, let~$X$ denote the number of vertices~$v$ with $d^-_T(v) \in \left[\frac{n}{2}-2\sqrt{n},\frac{n}{2}-\sqrt{n}\right]$.
Then it is easy to see that, for large enough~$n$, we have $\mathbb{E}[X] \ge \frac{n}{10^4}$, say. Moreover, Chebyshev's inequality can be used to show that, with probability~$1-o(1)$, we have~$X \ge \frac{n}{2\cdot 10^4}$, again with room to spare. 	
Thus, \cref{thm:main} implies the following.

\begin{cor}\label{cor:almostall}
	As $n\rightarrow \infty$, the proportion of tournaments~$T$ on~$n$ vertices satisfying~$\pn(T)=\texc(T)$ tends to~$1$.
\end{cor}

Note that the case when~$n$ is even already follows from \cref{thm:evenlarge}\cref{thm:evenlarge-main}. \Cref{cor:almostall} is an analogue of a result of \NEW{K\"{u}hn and Osthus} \cite{kuhn2014hamilton}, which states that almost all sufficiently large tournaments $T$ contain~$\delta^0(T)\coloneqq \min\{d_T^+(v),d_T^-(v)\mid v\in V(T)\}$ edge-disjoint Hamilton cycles and which proved a conjecture of Erd\H{o}s (see~\cite{thomassen1982edge}).

\NEW{Rather than random tournaments, it is also natural to consider the following related question: for which densities $p$ is the random binomial digraph $D_{n,p}$ likely to be consistent? Very recently, significant partial results towards this question were obtained by Espuny D\'{i}az, Patel, and Stroh~\cite{diaz2021path}.}

Finally, we will see in \cref{sec:conclusion} that our methods give a short proof of (a stronger version of) a result of Osthus and Staden~\cite{osthus2013approximate}, which guarantees an approximate decomposition of regular ``robust outexpanders'' of linear degree into Hamilton cycles and which was used as a tool in the proof of Kelly's conjecture~\cite{kuhn2013hamilton}.

\subsection{Organisation of the paper}
This paper is organised as follows. In \cref{sec:sketch}, we give a proof overview of \cref{thm:main}. Notation will be introduced in \cref{sec:notation}, while tools and preliminary results will be collected in \cref{sec:preliminaries}.
\NEW{We consider exceptional tournaments in \cref{sec:annoyingT} and derive \cref{thm:even,cor:approx} from \cref{thm:main} in \cref{sec:cor}. Then, \cref{sec:good,sec:constructinglayouts,sec:main,sec:cleaning,sec:layouts,sec:approxdecomp,sec:texc} are devoted to proving \cref{thm:main}. In particular, the approximate decomposition step is carried out in \cref{sec:approxdecomp} and \cref{thm:main} is derived in \cref{sec:main}.}%
\OLD{Moreover, exceptional tournaments will be considered in \cref{sec:annoyingT}. \cref{sec:cleaning,sec:layouts,sec:approxdecomp,sec:main} will be devoted to the proof of \cref{thm:main}. \cref{thm:even,cor:approx} are derived from \cref{thm:main} in \cref{sec:main}.}
Finally, in \cref{sec:conclusion}, we discuss Hamilton decompositions of robust outexpanders and conclude with a remark about \cref{conj:all}.

\onlyinsubfile{\bibliographystyle{abbrv}
	\bibliography{Bibliography/Bibliography}}

%% file: Sketch.tex
	\onlyinsubfile{
		\setcounter{section}{1}
\section{Proof overview}}

\subsection{Robust outexpanders}\label{sec:sketch-rob}

Our proof of \cref{thm:main} will be based on the concept of robust outexpanders. Roughly speaking, a digraph~$D$ is called a \emph{robust outexpander} if, for any set~$S\subseteq V(D)$ which is neither too small nor too large, there exist significantly more than~$|S|$ vertices with many inneighbours in~$S$. (Robust outexpanders will be defined formally in \cref{sec:rob}.)
Any (almost) regular tournament is a robust outexpander and we will use that this property is inherited by random subdigraphs. The main result of~\cite{kuhn2013hamilton} states that any regular robust outexpander of linear degree has a Hamilton decomposition (see \cref{thm:robHamdecomp}). We can apply this to obtain an optimal path decomposition in the following setting.
Let~$D$ be a digraph on~$n$ vertices,~$0<\eta<1$, and suppose that $X^+\cup X^-\cup X^0$ is a partition of~$V(D)$ such that $|X^+|=|X^-|=\eta n$ and the following hold.
\begin{property}{\dagger}\label{sketch:degree}
	Each~$v\in X^0$ satisfies \NEW{$d_D^+(v)= \eta n=d_D^-(v)$}.\\
	Each~$v\in X^+$ satisfies~$d_D^+(v)=\eta n$ and~$d_D^-(v)=\eta n-1$.\\
	Each~$v\in X^-$ satisfies~$d_D^+(v)=\eta n-1$ and~\NEWTWO{$d_D^-(v)=\eta n$}.
\end{property}
Then, the digraph~$D'$ obtained from~$D$ by adding a new vertex~$v$ with~\NEW{$N_{D'}^+(v) = X^+$ and $N_{D'}^-(v) = X^-$}\OLD{$N_{D'}^\pm(v) = X^\pm$} is~$\eta n$-regular. Thus, if~$D$ is a robust outexpander, then \NEW{so is $D'$ and} there exists a decomposition of~$D'$ into Hamilton cycles. This induces a decomposition~$\cP$ of~$D$ into~$\eta n$ Hamilton paths, where each vertex in~$X^+$ is the starting point of exactly one path in~$\cP$ and each vertex in~$X^-$ is the ending point of exactly one path in~$\cP$.
This is formalised in \cref{thm:niceD}.
(A similar observation was already made and used in~\cite{lo2020decomposing}.)
Our main strategy will be to reduce our \NEWTWO{tournament}\OLDTWO{tournaments} to a digraph of the above form. This will be achieved as follows.

\subsection{Simplified approach for well behaved tournaments}\label{sec:sketch-simple}

Let~$\beta>0$ and fix additional constants such that $0<\frac{1}{n_0}\ll \varepsilon \ll \gamma\ll \eta  \ll \beta$. Let~$T$ be a tournament on~$n\geq n_0$ vertices. Note that by \cref{thm:evenlarge}, we may assume that~$\texc(T)\leq \varepsilon^2 n^2$.
Moreover, for simplicity, we first also assume that each~$v\in V(T)$ satisfies~$|\exc_T(v)|\leq \varepsilon n$ (i.e.~$T$ is almost regular),~$\texc(T)=\exc(T)$, and both \NEW{$|\{v\in V(T)\mid \exc_T^+(v)>0\}|,|\{v\in V(T)\mid \exc_T^-(v)>0\}|\geq \eta n$}\OLD{$|\{v\in V(T)\mid \exc_T^\pm(v)>0\}|\geq \eta n$}. In \cref{sec:sketch-general}, we will \NEW{briefly} explain how the argument can be generalised if any of these conditions is not satisfied. \NEW{(An in-depth discussion of these modifications can be found in \cref{sec:good,sec:constructinglayouts}.)}

\NEW{Since~$T$ is almost regular, it is a robust outexpander. Let $\Gamma$ be obtained by including each edge of $T$ with probability $\gamma$. Using Chernoff bounds, we may assume that $\Gamma$ is a robust outexpander of density \NEWTWO{almost} $\gamma$ and $D\coloneqq T\setminus \Gamma$ is almost regular.}\OLD{Firstly, since~$T$ is almost regular, it is a robust outexpander and so we can fix a random spanning subdigraph~$\Gamma\subseteq T$ of density~$\gamma$ such that~$\Gamma$ is a robust outexpander and~$T\setminus \Gamma$ is almost regular.}
The digraph~$\Gamma$ will serve two purposes. Firstly, its robust outexpansion properties will be used to construct an approximate path decomposition of~$T$. Secondly, provided few edges of~$\Gamma$ are used throughout this approximate decomposition, it will guarantee that the leftover (consisting of all of those edges of~$T$ not covered by the approximate path decomposition) is still a robust outexpander, as required to complete our decomposition of~$T$ in the way described in \cref{sec:sketch-rob}.

Fix \NEW{$X^+ \subseteq \{v\in V(T)\mid \exc_T^+(v)>0\}$ and $X^- \subseteq \{v\in V(T)\mid \exc_T^-(v)>0\}$, both}\OLD{$X^\pm \subseteq \{v\in V(T)\mid \exc_T^\pm(v)>0\}$} of size~$\eta n$ and denote $X^0\coloneqq V(T)\setminus (X^+\cup X^-)$. Our goal is then to find an approximate path decomposition~$\cP$ of~$T$ such that $|\cP|= \texc(T)-\eta n$ and such that the leftover~\NEW{$T\setminus E(\cP)$}\OLD{$D\coloneqq T\setminus\bigcup\cP$} satisfies the degree conditions in \cref{sketch:degree}. Thus, it suffices \NEWTWO{to show} that~$\cP$ satisfies the following.
\begin{enumerate}[label=\upshape(\roman*)]
	\item Each~$v\in X^+$ is the starting point of exactly~$\exc_T^+(v)-1$ paths in~$\cP$, while each~$v\in V(T)\setminus X^+$ is the starting point of exactly~$\exc_T^+(v)$ paths in~$\cP$. Similarly, each~$v\in X^-$ is the ending point of exactly~$\exc_T^-(v)-1$ paths in~$\cP$, while each~$v\in V(T)\setminus X^-$ is the ending point of exactly~$\exc_T^-(v)$ paths in~$\cP$.\label{sketch:endpoint}
	\item Each~$v\in V(T)\setminus(X^+\cup X^-)$ is the internal vertex of exactly $\frac{(n-1)-|\exc_T(v)|}{2}-\eta n$ paths in~$\cP$, while each~$v\in X^+\cup X^-$ is the internal vertex of exactly $\frac{(n-1)-|\exc_T(v)|}{2}-\eta n+1$ paths in~$\cP$.\label{sketch:internal}
\end{enumerate}
Indeed, \cref{sketch:endpoint} ensures that~$|\cP|=\exc(T)-\eta n$ and each vertex has the desired excess in~\NEW{$T\setminus E(\cP)$}\OLD{$D$}, namely \NEW{$\exc_{T\setminus E(\cP)}(v)=+ 1$ if~$v\in X^+$, $\exc_{T\setminus E(\cP)}(v)=- 1$ if~$v\in X^-$,}\OLD{$\exc_D(v)=\pm 1$ if~$v\in X^\pm$} and~\NEW{$\exc_{T\setminus E(\cP)}(v)=0$}\OLD{$\exc_D(v)=0$} otherwise. In addition, \cref{sketch:internal} ensures that the degrees in~\NEW{$T\setminus E(\cP)$}\OLD{$D$} satisfy \cref{sketch:degree}.

Recall that, by assumption,~$T$ is almost regular.
Thus, in a nutshell, \cref{sketch:endpoint,sketch:internal} state that we need to construct edge-disjoint paths with specific endpoints and such that each vertex is covered by about~$(\frac{1}{2}-\eta)n$ paths.
To ensure the latter, we will in fact approximately decompose~$T$ into about~$(\frac{1}{2}-\eta)n$ spanning sets of internally vertex-disjoint paths.
To ensure the former, we will start by constructing~$(\frac{1}{2}-\eta)n$ auxiliary digraphs on~$V(T)$ such that, for each~$v\in V(T)$, the total number of edges starting (and ending) at~$v$ is the number of paths that we want to start (and end, respectively) at~$v$.
These auxiliary digraphs will be called \emph{layouts}. These layouts are constructed in \cref{sec:layouts}.
Then, it will be enough to construct, for each layout~$L$, a spanning set~$\cP_L$ of paths, called a \emph{spanning configuration of shape~$L$}, such that each path~$P\in \cP_L$ corresponds to some edge~$e\in E(L)$ and such that the starting and ending points of~$P$ equal those of~$e$.
\NEW{Roughly speaking, a spanning configuration $\cP_L$ is a set of internally-disjoint paths and $L$ indicates the starting and ending points of these paths. (See \cref{sec:approxdecomp} for further motivation of layouts.)}

These spanning configurations will be constructed one by one as follows. \NEW{(See also \cref{fig:sketch}.)} At each stage, given a layout~$L$, fix an edge \NEW{$yz\in E(L)$}\OLD{$UV\in E(L)$}. Then, using the robust \NEWTWO{outexpansion} properties of (the remainder of)~$\Gamma$, find short internally vertex-disjoint paths with endpoints corresponding to the endpoints of the edges in \NEW{$L\setminus \{yz\}$}\OLD{$L\setminus\{uv\}$}. Denote by~$\cP_L'$ the set containing these paths. Then, it only remains to construct a path from \NEW{$y$ to~$z$}\OLD{$u$ to $v$} spanning~$V(T)\setminus V(\cP_L')$. We achieve this as follows.

Let~$D'$ and~$\Gamma'$ be obtained from (the remainders of) \NEW{$D-V(\cP_L')$}\OLD{$(T\setminus \Gamma)-V(\cP_L')$} and $\Gamma-V(\cP_L')$ by merging the vertices \NEW{$y$ and $z$}\OLD{$u$ and~$v$} into a new vertex \NEW{$v_{yz}$}\OLD{$w$} whose outneighbourhood is the outneighbourhood of \NEW{$y$}\OLD{$u$} and whose inneighbourhood is the inneighbourhood of \NEW{$z$}\OLD{$v$}. Then, observe that a Hamilton cycle of~$D'\cup \Gamma'$ corresponds to a path from \NEW{$y$ to $z$}\OLD{$u$ to~$v$} of~$T$ which spans~$V(T)\setminus V(\cP_L')$.
\NEW{To}\OLD{We} construct \NEW{such} a Hamilton cycle of~$D'\cup \Gamma'$\OLDTWO{as follows}\OLD{. 
Of course}, one can simply use the fact that~$\Gamma'$ is a robust expander to find a Hamilton cycle. However, if we proceed in this way, then the robust \NEWTWO{outexpansion} property of~$\Gamma'$ might be destroyed before constructing all the desired spanning configurations.

So instead we construct a Hamilton cycle \NEW{of $D'\cup \Gamma'$} with only few edges in~$\Gamma'$ as follows.
\NEW{As a preparatory step in advance of choosing the spanning configurations, we consider a random partition of $V(T)$ into $A_1, \dots, A_a$ each of size $\frac{n}{a}$. We choose one $A_i$ for the current layout.  We restrict ourselves to use $\Gamma'$ inside $A_i$ only. Note that $\Gamma'[A_i]$ is a robust outexpander and $D'-A_i$ is a dense almost regular digraph. The latter means that we can find a spanning linear forest $F$ in $D'-A_i$ which has few components. Since $F$ has few components, we can then greedily extend the components of $F$ to obtain a linear forest $F'\subseteq D'$ which covers all the vertices in $V(D')\setminus A_i$ and whose endpoints are all in $A_i$. Finally, we use the robust outexpansion properties of $\Gamma'[A_i]$ to close $F'$ in to a Hamilton cycle of $D'\cup \Gamma'$. None of the $A_i$ will be used too often when constructing the spanning configurations, which will mean that $\Gamma'[A_i]$ is always a robust outexpander.}%
\OLD{Using the fact that $T\setminus \Gamma$ is almost regular, we first find an almost spanning linear forest~$F$ in~$D'$ which has few components. Then, we use the robust outexpanding properties of~$\Gamma'$ to tie up~$F$ into a Hamilton cycle~$C$ of~$D'\cup \Gamma'$.}
\NEWTWO{When the desired spanning configuration is a Hamilton cycle, this approach of finding many edge disjoint spanning configurations by first finding a suitable linear forest $F$, and then tying $F$ together together via some small set $A_i$  (with varying $A_i$ in order to avoid over-using a particular set of vertices) has been used successfully in several earlier papers (e.g.\, \cite{kuhn2013hamilton,ferber2018counting}).}
This construction of spanning configurations\OLD{(respecting the layouts formed in \cref{sec:layouts})} is carried out in \cref{sec:approxdecomp}.

\NEW{We illustrate this argument with the following example. Suppose that $L$ is a layout consisting of three edges $uv,wx$, and $yz$ (\cref{fig:sketch-L}). We want to construct a spanning configuration of shape $L$, that is, a set of paths which consists of a path from $u$ to $v$, a path from $w$ to $x$, and a path from $y$ to $z$ such that these three paths are vertex-disjoint and altogether cover all the vertices of $T$. First, we use robust outexpansion to construct a short path $P_1$ from $u$ to $v$ and a short path $P_2$ from $w$ to $x$ in $\Gamma$ (\cref{fig:sketch-yz}). Denote $V'\coloneqq V(T)\setminus (V(P_1)\cup V(P_2)\cup \{y,z\})$. The goal is now to construct a path from $y$ to $z$ which covers all the vertices in $V'$. To do so, we replace $y$ and $z$ by an auxiliary vertex $v_{yz}$ whose outneighbourhood is $N^+(v_{yz})\coloneqq N_D^+(y)\cap V'$ and whose inneighbourhood is $N^-(v_{yz})\coloneqq N_D^-(z)\cap V'$ (\cref{fig:sketch-yz}) and we consider a small preselected random set of vertices $A_i\subseteq V'$. It is then enough to find a cycle on $V'\cup \{v_{yz}\}$ which uses $\Gamma$ inside $A_i$ only. Denote $V''\coloneqq (V'\cup \{v_{yz}\})\setminus A_i$. Firstly, we use almost regularity of $D$ to find a spanning linear forest on $V''$ which consists of few components (\cref{fig:sketch-F}). Secondly, we use the large degree of $D$ to extend the endpoints of the linear forest to $A_i$ (\cref{fig:sketch-FA}). Finally, we use the robust outexpansion of $\Gamma$ to close the linear forest into a cycle which covers all the vertices in $A_i$ (\cref{fig:sketch-C}). This gives a cycle on \NEWTWO{$V'\cup \{v_{yz}\}$}. Replacing the auxiliary vertex $v_{yz}$ by the original vertices $y$ and $z$, we obtain a path from $y$ to $z$ which covers all the vertices in $V'$, as desired (\cref{fig:sketch-P}).}

\subfile{Figures/Figure_Sketch.tex}

\subsection{General tournaments}\label{sec:sketch-general}

For a general tournament~$T$, we adapt the above argument as follows.
Let~$W$ be the set of vertices~$v\in V(T)$ such that $|\exc_T(v)|>\varepsilon n$. If~$W\neq \emptyset$, then~$T$ is no longer almost regular and we cannot proceed as above. However, since $\exc(T)\leq \varepsilon^2 n^2$,~$|W|$ is small. Thus, we can start with a cleaning procedure which efficiently decreases the excess and degree at~$W$ by taking out \NEW{a} few edge-disjoint paths. The corresponding proof is deferred until \cref{sec:cleaning}, as it is quite involved and carrying out the other steps first helps to give a better picture of the overall argument. Then, we apply the above argument to (the remainder of)~$T-W$. We incorporate all remaining edges at~$W$ in the approximate decomposition by generalising the concept of a layout introduced above. \NEW{This is discussed in more detail in \cref{sec:constructinglayouts}.}

If $|\{v\in V(T)\mid \exc_T^+(v)>0\}|<\eta n$ but $\texc(T)=\exc(T)$, say, then we cannot choose $X^+\subseteq \{v\in V(T)\mid \exc_T^+(v)>0\}$ of size~$\eta n$. We circumvent this problem as follows. Select a small set of vertices~$W_A$ such that $\sum_{v\in W_A}\exc_T^+(v)\geq \eta n$ and let~$A$ be a set of~$\eta n$ edges such that the following hold. Each edge in~$A$ starts in~$W_A$ and ends in~$V(T)\setminus W_A$. Moreover, each~$v\in W_A$ is the starting point of at most~$\exc_T^+(v)$ edges in~$A$ and each~$v\in V(T)\setminus W_A$ is the ending point of at most one edge in~$A$. We will call~$A$ an \emph{absorbing set of starting edges}.
\NEW{Let $V_A$ be the set of ending points of the edges in $A$. Then, $V_A\subseteq V(T)\setminus W_A$. Observe that any path which starts in $V_A$ and is disjoint from $W_A$ can be extended to a path starting in $W_A$ using an edge from $A$. Thus, we can}
let the ending points of the edges in~$A$ play the role of~$X^+$ and add the vertices in~$W_A$ to~$W$ so that, at the end of the approximate decomposition, the only remaining edges at~$W_A$ are the edges in~$A$. Thus, in the final decomposition step, we can use the edges in~$A$ to extend the paths starting at~$X^+$ into paths starting in~$W_A$. \NEW{(See \cref{sec:A} for details.)} If $|\{v\in V(T)\mid \exc_T^-(v)>0\}|<\eta n$, then we proceed analogously.

If $\texc(T)>\exc(T)$, then not all paths will ``correspond" to some excess. \NEW{To be able to adopt a unified approach}\OLD{For simplicity}, we will choose which additional endpoints to use at the beginning and artificially add excess to those vertices. This then enables us to proceed as if~$\exc(T)=\texc(T)$.
More precisely, we will choose a set $U^*\subseteq \{v\in V(T)\mid \exc_T(v)=0\}$ of size~$\texc(T)-\exc(T)$ and we will treat the vertices in~$U^*$ in the same way as we treat those with \NEW{$\exc_T^+(v)=1$ and $\exc_T^-(v)=1$}\OLD{$\exc_T^\pm(v)=1$}.
Note that selecting additional endpoints in this way maximises the number of distinct endpoints, which will enable us to choose \NEW{$X^+\subseteq \{v\in V(T)\mid \exc_T^+(v)>0\}\cup U^*$ and/or $X^-\subseteq \{v\in V(T)\mid \exc_T^-(v)>0\}\cup U^*$ when $N^+(T)=|\{v\in V(T)\mid \exc_T^+(v)>0\}|+\texc(T)-\exc(T)\geq \eta n$ and/or $N^-(T)=|\{v\in V(T)\mid \exc_T^-(v)>0\}|+\texc(T)-\exc(T)\geq \eta n$}\OLD{$X^\pm\subseteq \{v\in V(T)\mid \exc_T^\pm(v)>0\}\cup U^*$ when $N^\pm(T)=|\{v\in V(T)\mid \exc_T^\pm(v)>0\}|+\texc(T)-\exc(T)\geq \eta n$}, and use absorbing edges otherwise, i.e.\ if condition \cref{thm:main-largeU} fails in \cref{thm:main}. More details of this approach are given in \cref{sec:A}.

\onlyinsubfile{\bibliographystyle{abbrv}
	\bibliography{Bibliography/Bibliography}}

%% file: Figures/Figure_Sketch.tex
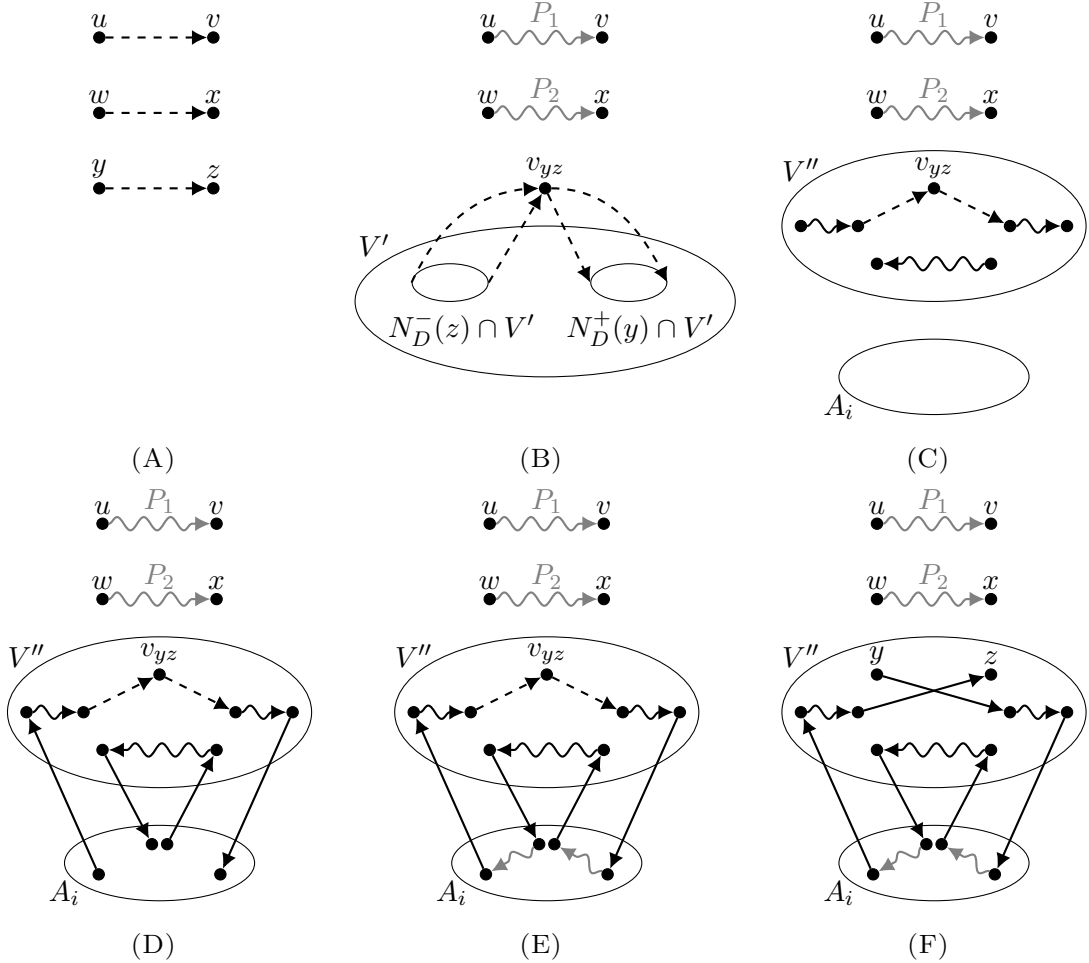
\begin{figure}[htb]
	\centering
	\begin{subfigure}{0.32\textwidth}
	\centering
	\begin{tikzpicture}
		%%vertices%%
		\draw node[circle, draw=black,fill=black, inner sep=1.5pt](u) at (-.75,0) {};
		\node[above] at (-.75,0) {$u$};
		\node[circle, draw=black,fill=black, inner sep=1.5pt](v) at (0.75,0) {};
		\node[above] at (0.75,0) {$v$};
		\draw node[circle, draw=black,fill=black, inner sep=1.5pt](w) at (-.75,-1) {};
		\node[above] at (-.75,-1) {$w$};
		\node[circle, draw=black,fill=black, inner sep=1.5pt](x) at (0.75,-1) {};
		\node[above] at (0.75,-1) {$x$};
		\draw node[circle, draw=black,fill=black, inner sep=1.5pt](y) at (-.75,-2) {};
		\node[above] at (-.75,-2) {$y$};
		\node[circle, draw=black,fill=black, inner sep=1.5pt](z) at (0.75,-2) {};
		\node[above] at (0.75,-2) {$z$};
		\node[white] at (0,-4.9) {$A_i$};
		%%edges%%
		\draw[->,thick,dashed] (u)-- node[above,white] {$P_1$} (v);
		\draw[->,thick,dashed] (w)-- (x);
		\draw[->,thick,dashed] (y)-- (z);
	\end{tikzpicture}
	\caption{\label{fig:sketch-L}}
	\end{subfigure}%
%	\hspace{0.01\textwidth}
	\begin{subfigure}{0.32\textwidth}
		\centering
		\begin{tikzpicture}
			%%vertices%%
			\draw node[circle, draw=black,fill=black, inner sep=1.5pt](u) at (-.75,0) {};
			\node[above] at (-.75,0) {$u$};
			\node[circle, draw=black,fill=black, inner sep=1.5pt](v) at (0.75,0) {};
			\node[above] at (0.75,0) {$v$};
			\draw node[circle, draw=black,fill=black, inner sep=1.5pt](w) at (-.75,-1) {};
			\node[above] at (-.75,-1) {$w$};
			\node[circle, draw=black,fill=black, inner sep=1.5pt](x) at (0.75,-1) {};
			\node[above] at (0.75,-1) {$x$};
			\draw node[circle, draw=black,fill=black, inner sep=1.5pt](yz) at (0,-2) {};
			\node[above] at (0,-2) {$v_{yz}$};
			\draw (0,-3.5) ellipse (2.5cm and 1cm);
			\node at (-2.25,-2.75) {$V'$};
			\node[white] at (0,-4.9) {$A_i$};
			%%edges%%
			\draw[->,thick,decorate,decoration={snake,post length=2pt},gray] (u)-- node[above] {$P_1$} (v);
			\draw[->,thick,decorate,decoration={snake,post length=2pt},gray] (w)-- node[above] {$P_2$} (x);
			\node[ellipse,draw,minimum width = 1cm,minimum height = 0.5cm] (e-) at (-1.25,-3.25) {};
			\node[below] at (-1.1,-3.5) {$N_D^-(z)\cap V'$};
			\draw[->,thick,dashed] (e-.west)to [out=65,in=180] (yz);
			\draw[->,thick,dashed] (e-.east)--(yz);
			\node[ellipse,draw,minimum width = 1cm,minimum height = 0.5cm] (e+) at (1.1,-3.25) {};
			\node[below] at (1.25,-3.5) {$N_D^+(y)\cap V'$};
			\draw[->,thick,dashed] (yz)--(e+.west);
			\draw[->,thick,dashed] (yz)to [out=0,in=115] (e+.east);
		\end{tikzpicture}
		\caption{\label{fig:sketch-yz}}
	\end{subfigure}%
%	\hspace{0.01\textwidth}
	\begin{subfigure}{0.32\textwidth}
		\centering
		\begin{tikzpicture}
			%%vertices%%
			\draw node[circle, draw=black,fill=black, inner sep=1.5pt](u) at (-.75,0) {};
			\node[above] at (-.75,0) {$u$};
			\node[circle, draw=black,fill=black, inner sep=1.5pt](v) at (0.75,0) {};
			\node[above] at (0.75,0) {$v$};
			\draw node[circle, draw=black,fill=black, inner sep=1.5pt](w) at (-.75,-1) {};
			\node[above] at (-.75,-1) {$w$};
			\node[circle, draw=black,fill=black, inner sep=1.5pt](x) at (0.75,-1) {};
			\node[above] at (0.75,-1) {$x$};
			\draw node[circle, draw=black,fill=black, inner sep=1.5pt](yz) at (0,-2) {};
			\node[above] at (0,-2) {$v_{yz}$};
			\draw node[circle, draw=black,fill=black, inner sep=1.5pt](1) at (-1.75,-2.5) {};
			\draw node[circle, draw=black,fill=black, inner sep=1.5pt](2) at (-1,-2.5) {};
			\draw node[circle, draw=black,fill=black, inner sep=1.5pt](3) at (1,-2.5) {};
			\draw node[circle, draw=black,fill=black, inner sep=1.5pt](4) at (1.75,-2.5) {};
			\draw node[circle, draw=black,fill=black, inner sep=1.5pt](5) at (-0.75,-3) {};
			\draw node[circle, draw=black,fill=black, inner sep=1.5pt](6) at (0.75,-3) {};
			\draw (0,-2.5) ellipse (2cm and 1cm);
			\node at (-1.75,-1.75) {$V''$};
			\draw (0,-4.5) ellipse (1.25cm and 0.5cm);
			\node at (-1.25,-4.9) {$A_i$};
			%%edges%%
			\draw[->,thick,decorate,decoration={snake,post length=2pt},gray] (u)-- node[above] {$P_1$} (v);
			\draw[->,thick,decorate,decoration={snake,post length=2pt},gray] (w)-- node[above] {$P_2$} (x);
			\draw[->,thick,decorate,decoration={snake,post length=2pt}] (1)--(2);
			\draw[->,thick,dashed] (2)--(yz);
			\draw[->,thick,dashed] (yz)--(3);
			\draw[->,thick,decorate,decoration={snake,post length=2pt}] (3)--(4);
			\draw[->,thick,decorate,decoration={snake,post length=2pt}] (6)--(5);
		\end{tikzpicture}
		\caption{\label{fig:sketch-F}}
	\end{subfigure}

	\begin{subfigure}{0.32\textwidth}
		\centering
		\begin{tikzpicture}
			%%vertices%%
			\draw node[circle, draw=black,fill=black, inner sep=1.5pt](u) at (-.75,0) {};
			\node[above] at (-.75,0) {$u$};
			\node[circle, draw=black,fill=black, inner sep=1.5pt](v) at (0.75,0) {};
			\node[above] at (0.75,0) {$v$};
			\draw node[circle, draw=black,fill=black, inner sep=1.5pt](w) at (-.75,-1) {};
			\node[above] at (-.75,-1) {$w$};
			\node[circle, draw=black,fill=black, inner sep=1.5pt](x) at (0.75,-1) {};
			\node[above] at (0.75,-1) {$x$};
			\draw node[circle, draw=black,fill=black, inner sep=1.5pt](yz) at (0,-2) {};
			\node[above] at (0,-2) {$v_{yz}$};
			\draw node[circle, draw=black,fill=black, inner sep=1.5pt](1) at (-1.75,-2.5) {};
			\draw node[circle, draw=black,fill=black, inner sep=1.5pt](2) at (-1,-2.5) {};
			\draw node[circle, draw=black,fill=black, inner sep=1.5pt](3) at (1,-2.5) {};
			\draw node[circle, draw=black,fill=black, inner sep=1.5pt](4) at (1.75,-2.5) {};
			\draw node[circle, draw=black,fill=black, inner sep=1.5pt](5) at (-0.75,-3) {};
			\draw node[circle, draw=black,fill=black, inner sep=1.5pt](6) at (0.75,-3) {};
			\draw node[circle, draw=black,fill=black, inner sep=1.5pt](7) at (-0.1,-4.25) {};
			\draw node[circle, draw=black,fill=black, inner sep=1.5pt](8) at (0.1,-4.25) {};
			\draw node[circle, draw=black,fill=black, inner sep=1.5pt](9) at (-0.8,-4.65) {};
			\draw node[circle, draw=black,fill=black, inner sep=1.5pt](10) at (0.8,-4.65) {};
			\draw (0,-2.5) ellipse (2cm and 1cm);
			\node at (-1.75,-1.75) {$V''$};
			\draw (0,-4.5) ellipse (1.25cm and 0.5cm);
			\node at (-1.25,-4.9) {$A_i$};
			%%edges%%
			\draw[->,thick,decorate,decoration={snake,post length=2pt},gray] (u)-- node[above] {$P_1$} (v);
			\draw[->,thick,decorate,decoration={snake,post length=2pt},gray] (w)-- node[above] {$P_2$} (x);
			\draw[->,thick,decorate,decoration={snake,post length=2pt}] (1)--(2);
			\draw[->,thick,dashed] (2)--(yz);
			\draw[->,thick,dashed] (yz)--(3);
			\draw[->,thick,decorate,decoration={snake,post length=2pt}] (3)--(4);
			\draw[->,thick,decorate,decoration={snake,post length=2pt}] (6)--(5);
			\draw[->,thick] (4)--(10);
			\draw[->,thick] (9)--(1);
			\draw[->,thick] (5)--(7);
			\draw[->,thick] (8)--(6);
		\end{tikzpicture}
		\caption{\label{fig:sketch-FA}}
	\end{subfigure}%
	%\hspace{0.04\textwidth}
	\begin{subfigure}{0.32\textwidth}
		\centering
		\begin{tikzpicture}
			%%vertices%%
			\draw node[circle, draw=black,fill=black, inner sep=1.5pt](u) at (-.75,0) {};
			\node[above] at (-.75,0) {$u$};
			\node[circle, draw=black,fill=black, inner sep=1.5pt](v) at (0.75,0) {};
			\node[above] at (0.75,0) {$v$};
			\draw node[circle, draw=black,fill=black, inner sep=1.5pt](w) at (-.75,-1) {};
			\node[above] at (-.75,-1) {$w$};
			\node[circle, draw=black,fill=black, inner sep=1.5pt](x) at (0.75,-1) {};
			\node[above] at (0.75,-1) {$x$};
			\draw node[circle, draw=black,fill=black, inner sep=1.5pt](yz) at (0,-2) {};
			\node[above] at (0,-2) {$v_{yz}$};
			\draw node[circle, draw=black,fill=black, inner sep=1.5pt](1) at (-1.75,-2.5) {};
			\draw node[circle, draw=black,fill=black, inner sep=1.5pt](2) at (-1,-2.5) {};
			\draw node[circle, draw=black,fill=black, inner sep=1.5pt](3) at (1,-2.5) {};
			\draw node[circle, draw=black,fill=black, inner sep=1.5pt](4) at (1.75,-2.5) {};
			\draw node[circle, draw=black,fill=black, inner sep=1.5pt](5) at (-0.75,-3) {};
			\draw node[circle, draw=black,fill=black, inner sep=1.5pt](6) at (0.75,-3) {};
			\draw node[circle, draw=black,fill=black, inner sep=1.5pt](7) at (-0.1,-4.25) {};
			\draw node[circle, draw=black,fill=black, inner sep=1.5pt](8) at (0.1,-4.25) {};
			\draw node[circle, draw=black,fill=black, inner sep=1.5pt](9) at (-0.8,-4.65) {};
			\draw node[circle, draw=black,fill=black, inner sep=1.5pt](10) at (0.8,-4.65) {};
			\draw (0,-2.5) ellipse (2cm and 1cm);
			\node at (-1.75,-1.75) {$V''$};
			\draw (0,-4.5) ellipse (1.25cm and 0.5cm);
			\node at (-1.25,-4.9) {$A_i$};
			%%edges%%
			\draw[->,thick,decorate,decoration={snake,post length=2pt},gray] (u)-- node[above] {$P_1$} (v);
			\draw[->,thick,decorate,decoration={snake,post length=2pt},gray] (w)-- node[above] {$P_2$} (x);
			\draw[->,thick,decorate,decoration={snake,post length=2pt}] (1)--(2);
			\draw[->,thick,dashed] (2)--(yz);
			\draw[->,thick,dashed] (yz)--(3);
			\draw[->,thick,decorate,decoration={snake,post length=2pt}] (3)--(4);
			\draw[->,thick,decorate,decoration={snake,post length=2pt}] (6)--(5);
			\draw[->,thick] (4)--(10);
			\draw[->,thick] (9)--(1);
			\draw[->,thick] (5)--(7);
			\draw[->,thick] (8)--(6);
			\draw[->,thick,decorate,decoration={snake,post length=3pt},gray] (7)--(9);
			\draw[->,thick,decorate,decoration={snake,post length=3pt},gray] (10)--(8);
		\end{tikzpicture}
		\caption{\label{fig:sketch-C}}
	\end{subfigure}%
	%\hspace{0.04\textwidth}
	\begin{subfigure}{0.32\textwidth}
		\centering
		\begin{tikzpicture}
			%%vertices%%
			\draw node[circle, draw=black,fill=black, inner sep=1.5pt](u) at (-.75,0) {};
			\node[above] at (-.75,0) {$u$};
			\node[circle, draw=black,fill=black, inner sep=1.5pt](v) at (0.75,0) {};
			\node[above] at (0.75,0) {$v$};
			\draw node[circle, draw=black,fill=black, inner sep=1.5pt](w) at (-.75,-1) {};
			\node[above] at (-.75,-1) {$w$};
			\node[circle, draw=black,fill=black, inner sep=1.5pt](x) at (0.75,-1) {};
			\node[above] at (0.75,-1) {$x$};
			\draw node[circle, draw=black,fill=black, inner sep=1.5pt](y) at (-.75,-2) {};
			\node[above] at (-.75,-2) {$y$};
			\node[circle, draw=black,fill=black, inner sep=1.5pt](z) at (0.75,-2) {};
			\node[above] at (0.75,-2) {$z$};
			\draw node[circle, draw=black,fill=black, inner sep=1.5pt](1) at (-1.75,-2.5) {};
			\draw node[circle, draw=black,fill=black, inner sep=1.5pt](2) at (-1,-2.5) {};
			\draw node[circle, draw=black,fill=black, inner sep=1.5pt](3) at (1,-2.5) {};
			\draw node[circle, draw=black,fill=black, inner sep=1.5pt](4) at (1.75,-2.5) {};
			\draw node[circle, draw=black,fill=black, inner sep=1.5pt](5) at (-0.75,-3) {};
			\draw node[circle, draw=black,fill=black, inner sep=1.5pt](6) at (0.75,-3) {};
			\draw node[circle, draw=black,fill=black, inner sep=1.5pt](7) at (-0.1,-4.25) {};
			\draw node[circle, draw=black,fill=black, inner sep=1.5pt](8) at (0.1,-4.25) {};
			\draw node[circle, draw=black,fill=black, inner sep=1.5pt](9) at (-0.8,-4.65) {};
			\draw node[circle, draw=black,fill=black, inner sep=1.5pt](10) at (0.8,-4.65) {};
			\draw (0,-2.5) ellipse (2cm and 1cm);
			\node at (-1.75,-1.75) {$V''$};
			\draw (0,-4.5) ellipse (1.25cm and 0.5cm);
			\node at (-1.25,-4.9) {$A_i$};
			%%edges%%
			\draw[->,thick,decorate,decoration={snake,post length=2pt},gray] (u)-- node[above] {$P_1$} (v);
			\draw[->,thick,decorate,decoration={snake,post length=2pt},gray] (w)-- node[above] {$P_2$} (x);
			\draw[->,thick,decorate,decoration={snake,post length=2pt}] (1)--(2);
			\draw[->,thick] (2)--(z);
			\draw[->,thick] (y)--(3);
			\draw[->,thick,decorate,decoration={snake,post length=2pt}] (3)--(4);
			\draw[->,thick,decorate,decoration={snake,post length=2pt}] (6)--(5);
			\draw[->,thick] (4)--(10);
			\draw[->,thick] (9)--(1);
			\draw[->,thick] (5)--(7);
			\draw[->,thick] (8)--(6);
			\draw[->,thick,decorate,decoration={snake,post length=3pt},gray] (7)--(9);
			\draw[->,thick,decorate,decoration={snake,post length=3pt},gray] (10)--(8);
		\end{tikzpicture}
		\caption{\label{fig:sketch-P}}
	\end{subfigure}
\caption{\NEW{Constructing a spanning set of vertex-disjoint paths in $D\cup \Gamma$ with prescribed endpoints and few edges of $\Gamma$. Dashed edges represent auxiliary edges, full black edges represent edges of $D$, and grey edges represent edges of $\Gamma$. Wavy black edges represent paths in $D$ and wavy grey edges represent paths in $\Gamma$.}\label{fig:sketch}}	
	\end{figure}

%% file: Notation.tex
	\onlyinsubfile{
		\setcounter{section}{2}
\section{Notation}}

In this section, we collect the notation that will be used throughout this paper. The non-standard pieces of notation will be recalled to the reader when first needed.

\subsection{Hierarchies}
We denote by~$\mathbb{N}$ the set of natural numbers (including~$0$).
Let~$a,b,c\in \mathbb{R}$. We write~$a=b\pm c$ if~$b-c\leq a\leq b+c$. For simplicity, we use hierarchies instead of explicitly calculating the values of constants for which our statements hold. More precisely, if we write~$0<a\ll b\ll c\leq 1$ in a statement, we mean that there exist non-decreasing functions $f\colon (0,1]\longrightarrow (0,1]$ and $g\colon (0,1]\longrightarrow (0,1]$ such that the statement holds for all~$0<a,b,c\leq 1$ satisfying~$b\leq f(c)$ and~$a\leq g(b)$. Hierarchies with more constants are defined in a similar way.
We assume large numbers to be integers and omit floors and ceilings, provided this does not affect the argument.

\subsection{\texorpdfstring{$\pm$}{Plus/minus}-notation}

In general, a statement $\cC^\pm$ will mean that both statements $\cC^+$ and $\cC^-$ hold simultaneously.
If used in the form that $\cC^\pm$ is the statement ``$\cA^\pm$ implies $\cB^\pm$'', the convention means that ``$\cA^+$ implies $\cB^+$'' and ``$\cA^-$ implies $\cB^-$''. Similarly, the statement ``$\cA^\pm$ implies $\cB^\mp$'' means that ``$\cA^+$ implies $\cB^-$'' and ``$\cA^-$ implies $\cB^+$''.

\subsection{Graphs and digraphs}

A \emph{digraph}~$D$ is a directed graph without loops which contains, for any distinct vertices~$u$ and~$v$ of~$D$, at most two edges between~$u$ and~$v$, at most one in each direction. A digraph~$D$ is called an \emph{oriented graph} if it contains, for any distinct vertices~$u$ and~$v$ of~$D$, at most one edge between~$u$ and~$v$; that is,~$D$ can be obtained by orienting the edges of an undirected graph.

Let~$G$ be a (di)graph. We denote by~$V(G)$ and~$E(G)$ the vertex and edge sets of~$G$, respectively. We say~$G$ is \emph{non-empty} if~$E(G)\neq \emptyset$.
Let~$u,v\in V(G)$ be distinct. If~$G$ is undirected, then we write~$uv$ for an edge between~$u$ and~$v$. If~$G$ is directed, then we write~$uv$ for an edge directed from~$u$ to~$v$, where~$u$ and~$v$ are called the \emph{starting} and \emph{ending points} of the edge~$uv$, respectively.
Let~$A,B\subseteq V(G)$ be disjoint.
Denote $E_A(G)\coloneqq \{e\in E(G)\mid V(e)\cap A\neq \emptyset\}$. Moreover, we write~$G[A,B]$ for the undirected graph with vertex set~$A\cup B$ and edge set $\{ab\in E(G)\mid a\in A, b\in B\}$ and $e(A,B)\coloneqq|E(G[A,B])|$.

Given~$S\subseteq V(G)$, we write~$G[S]$ for the sub(di)graph of~$G$ induced on~$S$ and~$G-S$ for the (di)graph obtained from~$G$ by deleting all vertices in~$S$.
Given~$E\subseteq E(G)$, we write~$G\setminus E$ for the (di)graph obtained from~$G$ by deleting all edges in~$E$. Similarly, given a sub(di)graph~$H\subseteq G$, we write~$G\setminus H\coloneqq G\setminus E(H)$. If~$F$ is a set of non-edges of~$G$, then we write~$G\cup F$ for the (di)graph obtained by adding all edges in~$F$. Given a (di)graph~$H$, if~$G$ and~$H$ are edge-disjoint, then we write~$G\cup H$ for the (di)graph with vertex set~$V(G)\cup V(H)$ and edge set~$E(G)\cup E(H)$.

\subsection{Degrees and neighbourhood}

Assume~$G$ is an undirected graph. For any~$v\in V(G)$, we write~$N_G(v)$ for the \emph{neighbourhood} of~$v$ in~$G$ and~$d_G(v)$ for the \emph{degree} of~$v$ in~$G$. Given~$S\subseteq V(G)$, we denote $N_G(S)\coloneqq \bigcup_{v\in S}N_G(v)$.

Let~$D$ be a digraph and $v\in V(D)$. We write~$N_D^+(v)$ and~$N_D^-(v)$ for the \emph{outneighbourhood} and \emph{inneighbourhood} of~$v$ in~$D$, respectively, and define the \emph{neighbourhood} of~$v$ in~$D$ as $N_D(v)\coloneqq N_D^+(v)\cup N_D^-(v)$. 
We denote by~$d_D^+(v)$ and~$d_D^-(v)$ the \emph{outdegree} and \emph{indegree} of~$v$ in~$D$, respectively, and define the \emph{degree} of~$v$ in~$D$ as $d_D(v)\coloneqq d_D^+(v)+d_D^-(v)$.
Denote $d_D^{\min}(v)\coloneqq \min\{d_D^+(v), d_D^-(v)\}$ and $d_D^{\max}(v)\coloneqq \max\{d_D^+(v), d_D^-(v)\}$. If $d_D^+(v)\neq d_D^-(v)$, then define 
\begin{equation}\label{eq:Nmax}
	N_D^{\min}(v)\coloneqq 
\begin{cases}
	N_D^+(v)& \text{if } d_D^{\min}=d_D^+(v),\\
	N_D^-(v)& \text{if } d_D^{\min}=d_D^-(v),\\
\end{cases}
\quad
\text{and}
\quad
N_D^{\max}(v)\coloneqq 
\begin{cases}
	N_D^+(v)& \text{if } d_D^{\max}=d_D^+(v),\\
	N_D^-(v)& \text{if } d_D^{\max}=d_D^-(v).\\
\end{cases}
\end{equation}
The \emph{minimum semidegree} of~$D$ is defined as $\delta^0(D)\coloneqq \min\{d_D^{\min}(v)\mid v\in V(D)\}$ and, similarly, $\Delta^0(D)\coloneqq \max\{d_D^{\max}(v)\mid v\in V(D)\}$ is called the \emph{maximum semidegree} of~$D$. Define the \emph{minimum degree} and \emph{maximum degree} of~$D$ by $\delta(D)\coloneqq \min\{d_D(v)\mid v\in V(D)\}$ and $\Delta(D)\coloneqq \max\{d_D(v)\mid v\in V(D)\}$, respectively.
Given~$S\subseteq V(D)$, we denote $N_D^\pm(S)\coloneqq \bigcup_{v\in S}N_D^\pm(v)$ and $N_D(S)\coloneqq \bigcup_{v\in S}N_D(v)$.

Let~$D$ be a digraph on~$n$ vertices. We say~$D$ is~\emph{$r$-regular} if, for any~$v\in V(D)$, $d_D^+(v)=d_D^-(v)=r$. We say~$D$ is \emph{regular} if it is~$r$-regular for some~$r\in \mathbb{N}$.
Let~$\varepsilon, \delta>0$. We say~$D$ is \emph{$(\delta, \varepsilon)$-almost regular} if, for each~$v\in V(D)$, both $d_D^+(v)=(\delta\pm \varepsilon)n$ and $d_D^-(v)=(\delta\pm \varepsilon)n$.

\subsection{Multidigraphs}

Let~$A$ and~$B$ be multisets. The \emph{support} of~$A$ is the set $S(A)\coloneqq \{a\mid a\in A\}$. For each~$a\in S(A)$, we denote by~$\mu_A(a)$ the \emph{multiplicity} of~$a$ in~$A$. For any~$a\notin S(A)$, we define~$\mu_A(a)\coloneqq 0$.
We write~$A\cup B$ for the multiset with support $S(A\cup B)\coloneqq S(A)\cup S(B)$ and such that, for each $a\in S(A\cup B)$, $\mu_{A\cup B}(a)\coloneqq \mu_A(a)+\mu_B(b)$.
We denote by~$A\setminus B$ the multiset with support $S(A\setminus B)\coloneqq \{a\in S(A)\mid \mu_A(a)>\mu_B(a)\}$ and such that, for each $a\in S(A\setminus B)$, $\mu_{A\setminus B}(a)\coloneqq \mu_A(a)-\mu_B(a)$.
We say~$A$ is a \emph{submultiset} of~$B$, denoted~$A\subseteq B$, if~$S(A)\subseteq S(B)$ and, for each~$a\in S(A)$,~$\mu_A(a)\leq \mu_B(a)$.

By a \emph{multidigraph}, we mean a directed graph where we allow multiple edges but no loops.
All the notation and definitions introduced thus far extend naturally to multidigraphs, with unions/differences of edge sets now interpreted as multiset unions/differences.
In a multidigraph, two instances of an edge are considered to be distinct.
In particular, given a multidigraph~$D$, we say~$D_1, D_2\subseteq D$ are \emph{edge-disjoint} submultidigraphs of~$D$ if, for any~$e\in E(D)$, $\mu_{E(D_1)}(e)+\mu_{E(D_2)}(e)\leq \mu_{E(D)}$.

\subsection{Paths}

In this paper, all paths and cycles are directed, with edges consistently oriented. The \emph{length} of a path~$P$, denoted by~$e(P)$, is the number of edges it contains. A path on one vertex, i.e.\ a path of length~$0$ is called \emph{trivial}.
Let~$P=v_1v_2\dots v_\ell$ be a path. We say~$v_1$ is the \emph{starting point} of~$P$ and~$v_\ell$ is the \emph{ending point} of~$P$. We say~$v$ is an \emph{endpoint} of a path~$P$ if~$v$ is the starting or ending point of~$P$. We say $v_2, \dots, v_{\ell-1}$ are \emph{internal vertices} of~$P$.
We write $V^+(P)=\{v_1\}$, $V^-(P)=\{v_\ell\}$, and $V^0(P)=\{v_2, \dots, v_{\ell-1}\}$.
We say that a path~$P$ is a \emph{$(u,v)$-path} if~$V^+(P)=\{u\}$ and~$V^-(P)=\{v\}$.
Given~$1\leq i<j\leq \ell$, we denote $v_iPv_j\coloneqq v_iv_{i+1}\dots v_j$.
A \emph{linear forest} is a set of pairwise vertex-disjoint paths.

Similarly, given a (multi)set~$\cP$ of paths, we write~$V^+(\cP)$ for the set of vertices which are the starting point of a path in~$\cP$. Similarly, we write~$V^-(\cP)$ for the set of vertices which are the ending point of a path in~$\cP$ and~$V^0(\cP)$ for the set of vertices which are an internal vertex of a path in~$\cP$. (Note that~$V^\pm(\cP)$ and~$V^0(\cP)$ are sets and not multisets.)

Given \NEW{a}\OLD{an (auxiliary)} directed edge~$xy$ and a path~$P$, we say~$P$ has \emph{shape}~$xy$ if~$P$ is an~$(x,y)$-path.
Similarly, let~$E$ be a (multi)set of (auxiliary) directed edges and~$\cP$ be a (multi)set of paths. We say~$\cP$ has \emph{shape}~$E$ if there exists a bijection $\phi:E\longrightarrow \cP$ such that, for each~$xy\in E$,~$\phi(xy)$ is an~$(x,y)$-path.

For convenience, a (multi)set~$\cP$ of paths will sometimes be viewed as the (multi)digraph consisting of their union. In particular, given a (multi)set~$\cP$ of paths, we write~$V(\cP)$ for the set of vertices of~$\cP$ and~$E(\cP)$ for the (multi)set of edges of~$\cP$, i.e.~$V(\cP)$ is the set $\bigcup_{P\in \cP}V(P)$ and~$E(\cP)$ is the (multi)set $\bigcup_{P\in \cP}E(P)$. (Note that~$V(\cP)$ is a set and not a multiset.) For any~$v\in V(\cP)$, we write~$d_{\cP}^{\pm}(v)$ and~$\exc_{\cP}^{\pm}(v)$ for the in/outdegree and positive/negative excess of~$v$ in~$\cP$ when viewed as a multidigraph, i.e.\ $d_{\cP}^{\pm}(v)\coloneqq d_{\bigcup \cP}^{\pm}(v)$ and $\exc_{\cP}^{\pm}(v)\coloneqq \exc_{\bigcup\cP}^{\pm}(v)$.
We define~$d_\cP(v)$ and~$\exc_\cP(v)$ similarly. \NEW{For any digraph $D$, we} denote \NEW{$D\setminus \cP\coloneqq D\setminus E(\cP)$}\OLD{$D\setminus \cP\coloneqq D\setminus\bigcup\cP$}.

\subsection{Subdivisions and contractions}

Let~$D$ and~$D'$ be digraphs and~$uv\in E(D)$.
We say~$D'$ is obtained from~$D$ by \emph{subdividing}~$uv$, if $V(D')=V(D)\cup \{w\}$, for some $w\notin V(D)$, and $E(D')=(E(D)\setminus \{uv\})\cup \{uw,wv\}$.
We say~$D'$ is a \emph{subdivision} of~$D$ if~$D'$ is obtained by successively subdividing some edges of~$D$.
Let~$P$ be a~$(u,v)$-path satisfying $V^0(P)\cap V(D)=\emptyset$.
We say~$D'$ is obtained from~$D$ by \emph{subdividing~$uv$ into~$P$}, if $V(D')=V(D)\cup V^0(P)$ and $E(D')=(E(D)\setminus \{uv\})\cup E(P)$.
Similarly, given an induced~$(u,v)$-path~$P\subseteq D$, we say~$D'$ is obtained from~$D$ by \emph{contracting the path~$P$ into an edge~$uv$} if $V(D')=V(D)\setminus V^0(P)$ and $E(D')=(E(D)\setminus E(P))\cup \{uv\}$.

\subsection{Decompositions}

Let~$D$ be a digraph. A \emph{decomposition} of~$D$ is set~$\cD$ of non-empty%
	\COMMENT{At least one edge.}
edge-disjoint subdigraphs of~$D$ such that every edge of~$D$ is in one of these subdigraphs.
A \emph{(Hamilton) path decomposition} of~$D$ is a decomposition~$\cP$ of~$D$ such that each subdigraph~$P\in \cP$ is a (Hamilton) path of~$D$.
Similarly, a \emph{(Hamilton) cycle decomposition} of~$D$ is a decomposition~$\cC$ of~$D$ such that each subdigraph~$C\in \cC$ is a (Hamilton) cycle of~$D$.
By a \emph{Hamilton decomposition} of~$D$, we mean a Hamilton cycle decomposition of~$D$.

%% file: Robust_Outexpanders.tex
	\onlyinsubfile{
		\setcounter{section}{4}
\subsection{Robust outexpanders}}

Let~$D$ be a digraph on~$n$ vertices. Given~$S\subseteq V(D)$, the \emph{$\nu$-robust outneighbourhood of~$S$} is the set $RN_{\nu, D}^+(S)\coloneqq \{v\in V(D)\mid |N_D^-(v)\cap S|\geq \nu n\}$.
We say that~$D$ is a \emph{robust~$(\nu, \tau)$-outexpander} if, for any~$S\subseteq V(D)$ satisfying $\tau n\leq |S|\leq (1-\tau)n$, $|RN_{\nu, D}^+(S)|\geq |S|+\nu n$.

In this section we state some useful properties of robust outexpanders. First, observe that the next \lcnamecref{fact:robparameters} follows immediately from the definition.

\begin{fact}\label{fact:robparameters}
	Let~$D$ be a robust~$(\nu,\tau)$-outexpander. Then, for any~$\nu'\leq \nu$ and~$\tau'\geq \tau$,~$D$ is a robust~$(\nu',\tau')$-outexpander.
\end{fact}

The following \lcnamecref{lm:verticesedgesremovalrob} states that robust outexpansion is preserved when few edges are removed and/or few vertices are removed and/or added. \NEW{This follows immediately from the definition and so we omit details. (Note that a similar result was already observed e.g. in \cite[Lemma 4.8]{kuhn2015robust}.)}

\begin{lm}\label{lm:verticesedgesremovalrob}
	Let $0<\varepsilon\leq \nu\leq \tau\leq 1$.
	Let~$D$ be a robust~$(\nu, \tau)$-outexpander on~$n$ vertices. 
	\begin{enumerate}
		\item If~$D'$ is obtained from~$D$ by removing at most~$\varepsilon n$ \NEW{inedges and at most $\varepsilon n$ outedges}\OLD{in/out edges} at each vertex, then~$D'$ is a robust~$(\nu-\varepsilon, \tau)$-outexpander.\label{lm:verticesedgesremovalrob-edges}
		\item \NEW{Suppose that $\tau\geq (1+2\tau)\varepsilon$.} If~$D'$ is obtained from~$D$ by adding or removing at most~$\varepsilon n$ vertices, then~$D'$ is a robust~$(\nu-\varepsilon, 2\tau)$-outexpander.\label{lm:verticesedgesremovalrob-vertices}
	\end{enumerate}
\end{lm}

One can easily show that the~$\tau$-parameter of robust outexpansion can be decreased when the minimum semidegree is large. This will enable us to state some results of~\cite{kuhn2013hamilton,kuhn2015robust} with slightly adjusted parameters.

\begin{lm}\label{lm:robparameters}
	Let $0<\frac{1}{n}\ll \nu\ll \tau\leq \frac{\delta}{2}\leq 1$. Assume~$D$ is a robust~$(\nu, \frac{\delta}{2})$-outexpander on~$n$ vertices satisfying~$\delta^0(D)\geq \delta n$. Then, $D$ is a robust~$(\nu, \tau)$-outexpander. 
\end{lm}

\begin{proof}
	\NEW{Let $S\subseteq V(D)$ satisfy $\tau n\leq |S|\leq (1-\tau)n$ and denote $T\coloneqq RN_{\nu, D}^+(S)$.
	We need to show that $|T|\geq |S|+\nu n$.
	If $\frac{\delta n}{2}\leq |S|\leq (1-\frac{\delta}{2})n$, then we are done by assumption. If $|S|>(1-\frac{\delta}{2})n$, then each $v\in V(D)$ satisfies $|N_D^\pm(v)\cap S|\geq \delta n- |V(D)\setminus S|\geq \frac{\delta n}{2}\geq \nu n$ and so $T=V(D)$.}
	
	\NEW{We may therefore assume that $|S|<\frac{\delta n}{2}$. Note that the number of edges of $D$ which start in $S$ is $\sum_{v\in S}d_D^+(v)$. By definition of $T$, we have \[|S|\delta n\stackrel{\delta^0(D)\geq \delta n}{\leq} \sum_{v\in S}d_D^+(v)\leq |T||S|+(n-|T|)\nu n\leq |T||S|+\nu n^2.\]
	Therefore, $|T|\geq \frac{|S|\delta n-\nu n^2}{|S|}\geq \delta n-\frac{\nu n}{\tau}\geq \frac{\delta n}{2}+\nu n\geq |S|+\nu n$,
	as desired.}				
\end{proof}

The next result states that oriented graphs of sufficiently large minimum semidegree are robust outexpanders.

\begin{lm}[{\cite[Lemma 13.1]{kuhn2013hamilton}}]\label{lm:3/8rob}
	Let $0<\frac{1}{n}\ll \nu\ll \tau\leq \varepsilon\leq 1$. Let~$D$ be \NEW{an} oriented graph on~$n$ vertices with $\delta^0(D)\geq (\frac{3}{8}+\varepsilon)n$. Then, $D$ is a robust~$(\nu, \tau)$-outexpander. 
\end{lm}

	\COMMENT{$\delta^0(D)\geq (\frac{3}{8}+\varepsilon)n$ implies $\delta^+(D)+\delta^-(D)+\delta(D)\geq 4(\frac{3}{8}+\varepsilon)n$ so can apply~\cite[Lemma 13.1]{kuhn2013hamilton} with $4\varepsilon$ playing the role of $\varepsilon$, i.e.\ only need $\tau\leq 2\varepsilon$ here.}

The next \lcnamecref{lm:robshortpath} follows easily from the definition of robust outexpansion and states that robust outexpanders of linear minimum semidegree have small
diameter.

\begin{lm}[{\cite[Lemma 6.6]{kuhn2015robust}}]\label{lm:robshortpath}
	Let $0<\frac{1}{n}\ll \nu\ll \tau\leq \frac{\delta}{2}\leq 1$. Let~$D$ be a robust~$(\nu, \tau)$-outexpander on~$n$ vertices with $\delta^0(D)\geq \delta n$. Then, for any~$x,y\in V(D)$,~$D$ contains an~$(x,y)$-path of length at most~$\nu^{-1}$.
\end{lm}

\NEW{One can iteratively apply \cref{lm:robshortpath} to obtain a small set of short internally vertex-disjoint paths with prescribed endpoints. After each application of \cref{lm:robshortpath}, one can check that the remaining digraph is still a robust outexpander by applying \cref{lm:verticesedgesremovalrob}\cref{lm:verticesedgesremovalrob-vertices}.}

\begin{cor}\label{cor:robshortpaths}
	Let $0<\frac{1}{n}\ll \varepsilon\ll\nu\ll \tau\leq \frac{\delta}{2}\leq 1$. Let~$D$ be a robust~$(\nu, \tau)$-outexpander on~$n$ vertices. Suppose that $\delta^0(D)\geq \delta n$ and let $S\subseteq V(D)$ be such that $|S|\leq \varepsilon n$. Let $k\leq \nu^3 n$ and $x_1, \dots, x_k, x_1', \dots, x_k'$ be (not necessarily distinct) vertices of~$D$. Let $X\coloneqq \{x_1, \dots, x_k,x_1', \dots, x_k'\}$. Then, there exist internally vertex-disjoint paths $P_1, \dots, P_k\subseteq D$ such that, for each~$i\in [k]$,~$P_i$ is an~$(x_i, x_i')$-path of length at most~$2\nu^{-1}$ and $V^0(P_i)\subseteq V(D)\setminus (X\cup S)$.
\end{cor}

\COMMENT{\begin{proof}
	\NEW{Assume inductively that, for some $0\leq m\leq k$, we have constructed internally vertex-disjoint paths $P_1, \dots, P_m\subseteq D$ such that, for each $i\in [m]$, $P_i$ is an $(x_i, x_i')$-path of length at most~$2\nu^{-1}$ and $V^0(P_i)\subseteq V(D)\setminus (X\cup S)$.
	If $m=k$, we are done. We may therefore assume that $m<k$. We construct $P_{m+1}$ as follows.
	Let $S'\coloneqq (S\cup \bigcup_{i\in [m]}V(P_i)\cup X)\setminus \{x_{m+1}, x_{m+1}'\}$. 
	Then, $|S'|\leq \varepsilon n+\nu^3n\cdot (2\nu^{-1}+1)\leq 3\nu^2n$ and so \cref{lm:verticesedgesremovalrob}\cref{lm:verticesedgesremovalrob-vertices} implies that $D'\coloneqq D-S'$ is a robust $(\frac{\nu}{2}, 2\tau)$-outexpander. Let $P_{m+1}$ be the path obtained by applying \cref{lm:robshortpath} with $D', \frac{\delta}{2}, \frac{\nu}{2}, 2\tau, x_{m+1}$, and $x_{m+1}'$ playing the roles of $D, \delta, \nu, \tau, x$, and $y$. Then, $P_{m+1}$ is an $(x_{m+1},x_{m+1}')$ of length at most $2\nu^{-1}$ which satisfies $V^0(P_{m+1})\subseteq V(D)\setminus S'\subseteq V(D)\setminus (S\cup X)$. Moreover, our choice of $S'$ implies that $P_{m+1}$ is internally vertex-disjoint with each of $P_1, \dots, P_m$.}	
\end{proof}}

We will use the fact that robust outexpanders of linear minimum degree contain Hamilton paths from any fixed
vertex~$x$ to any vertex~$y\neq x$. 
This immediately follows \NEW{by identifying $x$ and~$y$ to a single vertex~$z$ whose outneighbourhood is that of~$x$ and whose inneighbourhood is that of~$y$. The resulting digraph is a robust outexpander which contains a Hamilton cycle.}\OLD{from the result that such digraphs contain a Hamilton cycle (by contracting~$x$ and~$y$ to a single vertex~$z$ whose outneighbourhood is that of~$x$ and whose inneighbourhood is that of~$y$).} The Hamiltonicity of such digraphs was first proven in~\cite{keevash2009exact,kuhn2010hamiltonian}.

\begin{lm}[\NEW{\cite[Corollary 6.9]{kuhn2015robust}}]\label{lm:robHampath}
	Let $0<\frac{1}{n}\ll \nu\ll \tau\leq \frac{\delta}{2}\leq 1$. Let~$D$ be a robust~$(\nu, \tau)$-outexpander on~$n$ vertices with $\delta^0(D)\geq \delta n$. Then, for any distinct $x,y\in V(D)$,~$D$ contains a Hamilton~$(x,y)$-path.
\end{lm}

\NEW{Using similar arguments as in \cref{cor:robshortpaths}, we can iteratively apply \cref{lm:robshortpath} to tie together a small set of paths into a short path (\cref{cor:robpaths}\cref{cor:robpaths-short}). By replacing the last iteration of \cref{lm:robshortpath} with an application of \cref{lm:robHampath}, we can tie together a small set of paths into a Hamilton path (\cref{cor:robpaths}\cref{cor:robpaths-long}).
Similarly, we can tie together a small set of paths into a Hamilton cycle (\cref{cor:robpaths}\cref{cor:robcycles}) by first using \cref{lm:robshortpath} to tie these paths into a short path $P$ and then tie together the endpoints of $P$ using \cref{lm:robHampath}.}\OLD{Using \cref{lm:verticesedgesremovalrob,cor:robshortpaths,lm:robparameters,lm:robHampath}, one can prove the following \lcnamecref{cor:robpaths}.}

\begin{cor}\label{cor:robpaths}
	Let $0<\frac{1}{n}\ll\nu\ll \tau\leq \frac{\delta}{2}\leq 1$ and $k\leq \nu^3 n$. Let~$D$ be a digraph and $P_1, \dots, P_k\subseteq D$ be vertex-disjoint paths. For each~$i\in [k]$, denote by~$v_i^+$ and~$v_i^-$ the starting and ending points of~$P_i$, respectively.
	Let $V'\coloneqq V(D)\setminus \bigcup_{i\in [k]}V(P_i)$ and $S\subseteq V'$.
	Suppose that $D'\coloneqq D[V'\setminus S]$ is a robust~$(\nu,\tau)$-outexpander on~$n$ vertices satisfying $\delta^0(D')\geq \delta n$. Assume that for each~$i\in [k-1]$, $|N_D^+(v_i^-)\cap (V'\setminus S)|\geq 2k$ and $|N_D^-(v_{i+1}^+)\cap (V'\setminus S)|\geq 2k$. Then, the following hold.
	\begin{enumerate}
		\item There exists a~$(v_1^+,v_k^-)$-path~$Q\subseteq D\NEW{-S}$ of length at most $2\nu^{-1}k+\sum_{i\in[k]}e(P_i)$ such that \NEW{$Q$ contains $\bigcup_{i\in [k]}P_i$}\OLD{for each~$i\in [k]$,~$P_i\subseteq Q$ and $V(Q)\setminus \bigcup_{i\in [k]} V(P_i)\subseteq V'\setminus S$}.\label{cor:robpaths-short}
		\item There exists a~$(v_1^+,v_k^-)$-Hamilton path~$Q'$ of~$D-S$ \NEW{which contains $\bigcup_{i\in [k]}P_i$}\OLD{such that, for each~$i\in [k]$,~$P_i\subseteq Q'$}.\label{cor:robpaths-long}
		\item There exists a Hamilton cycle~$C$ of~$D-S$ \NEW{which contains $\bigcup_{i\in [k]}P_i$}\OLD{such that, for each~$i\in [k]$,~$P_i\subseteq C$}.\label{cor:robcycles}
	\end{enumerate}
\end{cor}

\begin{proof}
	\NEW{By assumption, there exist distinct $w_1^+, \dots, w_k^+, w_1^-, \dots, w_k^-\in V'\setminus S$ such that, for each $i\in [k]$, $w_i^+\in N_D^-(v_i^+)$ and $w_i^-\in N_D^+(v_i^-)$. In particular, observe that $w_1^+v_1^+P_1v_1^-w_1^-, \dots$, $w_k^+v_k^+P_kv_k^-w_k^-$ are vertex-disjoint paths of $D-S$. 
	Apply \cref{cor:robshortpaths} with $(x_i,x_i')=(w_i^-,w_{i+1}^+)$ for each $i\in [k-1]$ to obtain vertex disjoint paths $P_1', \dots, P_{k-1}'$ such that for each $i\in [k-1]$, $P_i'$ is a $(w_i^-,w_{i+1}^+)$-path of length at most $2\nu^{-1}$. Note that \cref{cor:robpaths-short} holds by setting \[Q\coloneqq v_1^+P_1v_1^-w_1^-P_1'w_2^+v_2^+P_2v_2^-w_2^-\dots w_{k-1}^-P_{k-1}'w_k^+v_k^+P_kv_k^-.\]
	For \cref{cor:robpaths-long}, let $D''\coloneqq D'-\bigcup_{i\in [k-2]}V(P_i')$. By \cref{lm:verticesedgesremovalrob}\cref{lm:verticesedgesremovalrob-vertices} and \cref{lm:robHampath}, $D''$ contains a Hamilton $(w_{k-1}^-,w_k^+)$-path $P_{k-1}''$. Let \[Q'\coloneqq v_1^+P_1v_1^-w_1^-P_1'w_2^+v_2^+P_2v_2^-w_2^-\dots w_{k-2}^-P_{k-2}'w_{k-1}^+v_{k-1}^+P_{k-1}v_{k-1}^-w_{k-1}^-P_{k-1}''w_k^+v_k^+P_kv_k^-.\]To prove \cref{cor:robcycles}, a similar argument shows that there exists a Hamilton $(w_k^-, w_1^+)$-path $P_k'$ in $D'-\bigcup_{i\in [k-1]}V(P_i')$. Let $C\coloneqq w_1^+v_1^+Qv_k^-w_k^-P_k'w_1^+$.}
\end{proof}

The main result of~\cite{kuhn2013hamilton} states that regular robust outexpanders of linear degree can be decomposed into Hamilton cycles. Note that this implies Kelly's conjecture on Hamilton decompositions of regular tournaments. \NEW{Indeed, any regular tournament $T$ on $n$ vertices satisfies $\delta^0(T)=\frac{n-1}{2}$. Thus, \cref{lm:3/8rob} implies that any regular tournament is in fact a robust outexpander (of linear degree).}

\begin{thm}[{\cite[Theorem 1.2]{kuhn2013hamilton}}]\label{thm:robHamdecomp}
	Let $0<\frac{1}{n}\ll \nu\ll \tau \leq \frac{\delta}{2}\leq 1$ and $r\geq \delta n$. Suppose~$D$ is an~$r$-regular robust~$(\nu, \tau)$-outexpander on~$n$ vertices. Then, $D$ has a Hamilton decomposition.
\end{thm}

The following result is a consequence of \cref{thm:robHamdecomp} and will be used to complete our path decompositions (\NEW{as described in}\OLD{recall} the proof overview). In particular, this implies that any digraph~$D$ satisfying \cref{sketch:degree} from \cref{sec:sketch-rob} is consistent, i.e.~$\pn(D) = \exc(D)$. \NEW{Note that \cref{thm:niceD} is slightly more general than the argument described in \cref{sec:sketch-rob}. Indeed, since \cref{thm:robHamdecomp} holds for digraphs (rather than just oriented graphs), we can allow $X^+$ and $X^-$ from \cref{sec:sketch-rob} to ``intersect": we may have some vertices with one less inedge and one less outedge and then join these vertices to the auxiliary vertex by an inedge and an outedge. The set $X^*$ in \cref{thm:niceD} consists of these vertices.}

\begin{cor}\label{thm:niceD}
	Let $0<\frac{1}{n}\ll \nu\ll \tau \leq \frac{\delta}{2}\leq 1$ and $r\geq \delta n$. Suppose~$D$ is a robust~$(\nu, \tau)$-outexpander on~$n$ vertices with a vertex partition $V(D)=X^+\cup X^-\cup X^*\cup X^0$ such that $|X^+\cup X^*|=|X^-\cup X^*|=r$ and, for all $v\in V(D)$, the following hold.
	\begin{equation}\label{eq:niceD}
	\exc_D(v)=
	\begin{cases}
	\pm 1 & \text{if } v\in X^{\pm},\\
	0 & \text{otherwise},
	\end{cases}
	\quad \text{and} \quad
	d_D(v)=
	\begin{cases}
	2r - 1 & \text{if } v\in X^{\pm},\\
	2r -2 & \text{if } v\in X^*,\\
	2r & \text{otherwise}.
	\end{cases}
	\end{equation}
	Then,~$\pn(D)=r$. \NEW{Moreover, $D$ has a path decomposition which consists of precisely $r$ Hamilton paths with distinct starting points in $X^+\cup X^*$ and distinct ending points in $X^-\cup X^*$.}
\end{cor}

The proof is very similar to~\cite[Theorem 4.7]{lo2020decomposing}, but we include it here for completeness.

\begin{proof}
	By \cref{fact:robparameters,lm:robparameters}, we may assume that $\tau \ll \delta$.
	
	Note that $\pn(D)\geq r$. Indeed, if~$X^*=V(D)$, then~$D$ is~$(r-1)$-regular \NEW{and so \cref{eq:pnreg} implies that $\pn(D)\geq r$.} Otherwise, $\Delta^0(D)=r$ and so, by \cref{prop:pn>texc}, $\pn(D)\geq \texc(D)\geq r$. \OLD{Thus, it is enough to find a path decomposition of~$D$ of size~$r$}\NEW{Thus, it is enough to decompose $D$ into $r$ Hamilton paths with distinct starting points in $X^+\cup X^*$ and distinct ending points in $X^-\cup X^*$.}
	
	Let~$D'$ be obtained from~$D$ by adding a new vertex~$v$ with $N_{D'}^\pm(v)\coloneqq X^\pm \cup X^*$. Then, by \cref{lm:verticesedgesremovalrob}\NEW{\cref{lm:verticesedgesremovalrob-vertices}},~$D'$ is a\NEW{n}~$r$-regular robust~$(\frac{\nu}{2}, 2\tau)$-outexpander. Applying \cref{thm:robHamdecomp} with $D', \frac{\delta}{2}, \frac{\nu}{2}$, and $2\tau$ playing the roles of $D, \delta, \nu$, and $\tau$ yields a Hamilton decomposition of~$D'$ \NEW{into $r$ Hamilton cycles.} \NEW{Removing $v$ we obtain}\OLD{This induces} a path decomposition of~$D$ \NEW{which consists of precisely $r$ Hamilton paths with distinct starting points in $X^+\cup X^*$ and distinct ending points in $X^-\cup X^*$}\OLD{of size~$r$}, as desired.
\end{proof}

\onlyinsubfile{\bibliographystyle{abbrv}
	\bibliography{Bibliography/Bibliography}}

%% file: Probabilistic_Estimates.tex
	\onlyinsubfile{
		\setcounter{section}{4}
		\setcounter{subsection}{1}
		\setcounter{definition}{10}
\subsection{Probabilistic estimates}}

In this section, we introduce a Chernoff-type bound and derive several easy probabilistic lemmas which will be used in the approximate decomposition step.

Let~$X$ be a random variable. We write~$X\sim \Bin(n,p)$ if~$X$ follows  a binomial distribution with parameters~$n$ and~$p$. Let~$N,n,m\in \mathbb{N}$ be such that~$\max\{n, m\} \leq N$. Let~$\Gamma$ be a set of size~$N$ and~$\Gamma'\subseteq \Gamma$ be of size~$m$. Recall that~$X$ has a \textit{hypergeometric distribution with parameters~$N, n$, and $m$} if~$X=|\Gamma_n\cap \Gamma'|$, where~$\Gamma_n$ is a random subset of~$\Gamma$ with~$|\Gamma_n|=n$ (i.e.~$\Gamma_n$ is obtained by drawing~$n$ elements of~$\Gamma$ without replacement). We will denote this by~$X\sim \HGeom(N,n,m)$. \COMMENT{Note that if~$X\sim \HGeom(N,n,m)$ then~$\mathbb{E}[X]=\frac{nm}{N}$.}

We will use the following Chernoff-type bound.

\begin{lm}[{see e.g.~\cite[Theorems 2.1 and 2.10]{janson2011random}}]\label{lm:Chernoff}
	Assume that $X\sim \Bin(n,p)$ or $X\sim \HGeom(N,n,m)$. Then, for any~$0< \varepsilon \leq 1$, the following hold.
	\begin{enumerate}
		\item $\mathbb{P}\left[X\leq (1-\varepsilon)\mathbb{E}[X]\right] \leq \exp \left(-\frac{\varepsilon^2}{3}\mathbb{E}[X]\right)$.
		\item $\mathbb{P}\left[X\geq (1+\varepsilon)\mathbb{E}[X]\right] \leq \exp \left(-\frac{\varepsilon^2}{3}\mathbb{E}[X]\right)$.
	\end{enumerate}
\end{lm}

Using \cref{lm:Chernoff}, it is easy to see that robust outexpansion is preserved with high probability when taking random edge-slices\OLD{(see e.g.\ the proof of~\cite[Lemma 3.2(ii)]{kuhn2014hamilton})}.

\begin{lm}[{\cite[Lemma 3.2(ii)]{kuhn2014hamilton}}]\label{lm:slicerob}
	Let $0<\frac{1}{n}\ll \nu \ll \tau, \gamma\leq 1$. Let~$D$ \NEW{be} a robust~$(\nu,\tau)$-outexpander on~$n$ vertices. Suppose~$\Gamma$ is obtained from~$D$ by taking each edge independently with probability~$\gamma$. Then, with probability at least $1-\exp(-\nu^3 n^2)$,~$\Gamma$ is a robust~$(\frac{\gamma\nu}{2}, \tau)$-outexpander. 
\end{lm}

\NEW{The bound on the probability is not part of the statement in \cite{kuhn2014hamilton} but follows immediately from the proof. (The latter considers each $S\subseteq V(D)$ of size $\tau n\leq |S|\leq (1-\tau)n$ and uses \cref{lm:Chernoff} to show that for all $v\in RN_{\nu, D}^+(S)$, the probability that $|N_\Gamma^-(v)\cap S|$ is small is exponentially small in $n$.)}

\COMMENT{\begin{proof}
	\NEW{Let $S\subseteq V(D)$ satisfy $\tau n\leq |S|\leq (1-\tau)n$. Denote $T\coloneqq RN_{\nu, D}^+(S)$ and note that $|T|\geq |S|+\nu n$. For each $v\in T$, $\mathbb{E}[|N_\Gamma^-(v)\cap S|]=\gamma|N_D^-(v)\cap S|\geq \gamma\nu n$ and so \cref{lm:Chernoff} implies that
	\[\mathbb{P}\left[|N_\Gamma^-(v)\cap S|\leq \frac{\gamma \nu n}{2}\right]\leq \exp\left(-\frac{\gamma\nu n}{12}\right),\]
	Thus, for any $T'\subseteq T$ satisfying $|T'|>\frac{\nu}{2}n$,
	\[\mathbb{P}\left[\forall v\in T'\colon |N_\Gamma^-(v)\cap S|\leq \frac{\gamma \nu n}{2}\right]\leq \exp\left(-\frac{\gamma\nu^2 n^2}{24}\right).\]
	Therefore, taking a union bound over all such $T'\subseteq T$ gives that%
		%\COMMENT{$\mathbb{P}[|RN_{\frac{\nu}{2}, \Gamma}^+(S)|\leq |S|+\frac{\gamma\nu n}{2}]\leq \mathbb{P}[|T\setminus RN_{\frac{\gamma \nu}{2}, \Gamma}^+(S)|>\frac{\nu n}{2}]=\mathbb{P}[\exists T'\subseteq T: |T'|>\frac{\nu n}{2} \wedge \forall v\in V'(|N_\Gamma^-(v)\cap S|\leq \frac{\gamma \nu n}{2})]$.}
	\[\mathbb{P}\left[|RN_{\frac{\nu}{2}, \Gamma}^+(S)|\leq |S|+\frac{\gamma\nu n}{2}\right]\leq\mathbb{P}\left[|RN_{\frac{\nu}{2}, \Gamma}^+(S)|\leq |S|+\frac{\nu n}{2}\right]\leq 2^n\cdot\exp\left(-\frac{\gamma\nu^2 n^2}{24}\right).\]
	By taking a union bound over all such $S\subseteq V(D)$, we conclude that $\Gamma$ is a robust $(\frac{\gamma\nu}{2}, \tau)$-outexpander with probability at least
	\[1-2^n\cdot 2^n\cdot\exp\left(-\frac{\gamma\nu^2 n^2}{24}\right)>1-\exp(-\nu^3 n^2),\]
	as desired.}
\end{proof}}

Similarly, using \cref{lm:Chernoff}, it is easy to see that the property of being almost regular is preserved when a random edge-slice is taken.

\begin{lm}\label{lm:slicereg}
	Let $0<\frac{1}{n}\ll \varepsilon, \gamma\ll \delta\leq 1$. Let~$D$ be a~$(\delta, \varepsilon)$-almost regular digraph on~$n$ vertices. Let~$\Gamma$ be obtained from~$D$ by taking each edge independently with probability~$\frac{\gamma}{\delta}$. Then, with probability at least~$1-\frac{1}{n}$,~$\Gamma$ is~$(\gamma, \varepsilon)$-almost regular and~$D\setminus \Gamma$ is~$(\delta-\gamma, \varepsilon)$-almost regular.
\end{lm}

\COMMENT{\begin{proof}
	Denote $V\coloneqq V(D)$. Let $v\in V$. We have $\mathbb{E}[d_\Gamma^\pm(v)]=\frac{\gamma}{\delta} d_D^\pm(v)=(\gamma\pm \frac{\varepsilon\gamma}{\delta})n=(\gamma\pm \frac{\varepsilon}{2})n$ so, by \cref{lm:Chernoff},%
		%\COMMENT{$\frac{\varepsilon n}{2}\leq \frac{\varepsilon}{\gamma}\mathbb{E}[d_\Gamma^\pm(v)]$}
	\[\mathbb{P}\left[d_\Gamma^\pm(v)\neq (\gamma\pm \varepsilon)n\right]\leq \mathbb{P}\left[|d_\Gamma^\pm(v)-\mathbb{E}[d_\Gamma^\pm(v)]|>\frac{\varepsilon n}{2}\right]\leq 2\exp\left(-\frac{\varepsilon^2 n }{24\gamma}\right).\]
	Therefore, taking a union bound over all $v\in V$, $\Gamma$ is $(\gamma, \varepsilon)$-almost regular with probability at least $1-\frac{1}{n^2}$.\\
	Similarly, for any $v\in V$, $\mathbb{E}[d_{D\setminus\Gamma}^\pm(v)]=(1-\frac{\gamma}{\delta}) d_D^\pm(v)=(\delta-\gamma\pm \varepsilon(1-\frac{\gamma}{\delta}))n$ so, by \cref{lm:Chernoff},%
		%\COMMENT{$\frac{\varepsilon\gamma n}{\delta}\leq \frac{2\varepsilon\gamma}{\delta^2}\mathbb{E}[d_{D\setminus \Gamma}^\pm(v)]$}
	\[\mathbb{P}\left[d_{D\setminus\Gamma}^\pm(v)\neq (\delta-\gamma\pm \varepsilon)n\right]\leq \mathbb{P}\left[|d_{D\setminus\Gamma}^\pm(v)-\mathbb{E}[d_{D\setminus\Gamma}^\pm(v)]|>\frac{\varepsilon\gamma n}{\delta}\right]\leq 2\exp\left(-\frac{\varepsilon^2\gamma^2 n }{4\delta^3}\right).\]
	Therefore, taking a union bound over all $v\in V$, $D\setminus\Gamma$ is $(\delta-\gamma, \varepsilon)$-almost regular with probability at least $1-\frac{1}{n^2}$. This completes the proof.
\end{proof}}

Let~$D$ be a digraph on~$n$ vertices. We say~$D$ is an \emph{$(\varepsilon,p)$-robust~$(\nu, \tau)$-outexpander} if~$D$ is a robust~$(\nu, \tau)$-outexpander and, for any integer~$k\geq \varepsilon n$, if~$S\subseteq V(D)$ is a random subset of size~$k$, then~$D[S]$ is a robust~$(\nu, \tau)$-outexpander with probability at least~$1-p$.
Note that the following analogue of \cref{fact:robparameters} holds for this new notion of robust outexpansion.

\begin{fact}\label{fact:robspecialparameters}
	Let~$D$ be a\NEW{n}~$(\varepsilon,p)$-robust~$(\nu, \tau)$-outexpander. Then, for any $\varepsilon'\geq \varepsilon$, $p'\geq p$, $\nu'\leq \nu$, and $\tau'\geq \tau$,~$D$ is a\NEW{n}~$(\varepsilon',p')$-robust~$(\nu', \tau')$-outexpander.
\end{fact}

Moreover, by \cref{lm:verticesedgesremovalrob}\NEW{\cref{lm:verticesedgesremovalrob-edges}}, the following holds.

\begin{lm}\label{fact:verticesedgesremovalrob}
	Let $0<\varepsilon\leq \nu\leq \tau\leq 1$.
	Let~$D$ be a~$(\varepsilon, p)$-robust~$(\nu,\tau)$-outexpander on~$n$ vertices. 
	If~$D'$ is obtained from~$D$ by removing at most~$\varepsilon n$ \NEW{inedges and at most $\varepsilon n$ outedges}\OLD{in/out edges} at each vertex, then~$D'$ is a~$(\sqrt{\varepsilon}, p)$-robust~$(\nu-\sqrt{\varepsilon}, \tau)$-outexpander.
\end{lm}

\begin{proof}
	\NEW{Let $S\subseteq V(D)$ satisfy $|S|\geq \sqrt{\varepsilon}n$ and suppose that $D[S]$ is a robust $(\nu,\tau)$-outexpander. By assumption, $D'[S]$ is obtained from $D[S]$ by removing at most $\varepsilon n\leq \sqrt{\varepsilon} |S|$ inedges and at most $\sqrt{\varepsilon} |S|$ outedges at each vertex. Thus, \cref{lm:verticesedgesremovalrob}\NEW{\cref{lm:verticesedgesremovalrob-edges}} implies that $D'[S]$ is a robust $(\nu-\sqrt{\varepsilon},\tau)$-outexpander.}
\end{proof}

We will see in the concluding remarks that any robust outexpander is in fact~$(\varepsilon,p)$-robust (for some suitable parameters).
However, our method for showing this requires the regularity lemma and so, for brevity, we will not prove this result. In this paper, we work with almost regular tournaments. Thus, it will be enough to use the next \lcnamecref{lm:Gamma}, which shows that~$(\varepsilon,p)$-robustness is easily inherited from almost regular robust outexpanders of sufficiently large minimum semidegree.

\begin{lm}\label{lm:Gamma}
	Let $0<\frac{1}{n}\ll \varepsilon \ll \nu \ll \tau \ll \gamma\ll\frac{3}{7}\leq \delta\leq 1$. Let~$D$ be a~$(\delta, \varepsilon)$-almost regular oriented graph on~$n$ vertices. Then, there exists a~$(\gamma, \varepsilon)$-almost regular spanning subdigraph~$\Gamma$ of~$D$ which is an~$(\varepsilon,n^{-3})$-robust%
		\COMMENT{We get $p$ exponentially small but $n^{-3}$ is already more than enough for our purposes.}
	$(\nu, \tau)$-outexpander and such that~$D\setminus\Gamma$ is~$(\delta-\gamma, \varepsilon)$-almost regular.
\end{lm}

\begin{proof}
	Let~$\Gamma$ be obtained from~$D$ by taking each edge independently with probability~$\frac{\gamma}{\delta}$.
	By \cref{lm:slicereg}, with probability at least~$1-\frac{1}{n}$,~$\Gamma$ is~$(\gamma, \varepsilon)$-almost regular and~$D\setminus \Gamma$ is~$(\delta-\gamma, \varepsilon)$-almost regular.
	 
	By \cref{lm:3/8rob},~$D$ is a robust~$(2\gamma^{-1}\nu, \tau)$-outexpander. Therefore, by \cref{lm:slicerob},~$\Gamma$ is a robust~$(\nu, \tau)$-outexpander with probability at least~$1-\exp(-8\gamma^{-3}\nu^3n^2)$.
		
	Assume that $S\subseteq V(D)$ is such that~$|S|\geq \varepsilon n$ and~$D[S]$ is a robust~$(2\nu^{-1}\gamma, \tau)$-outexpander. Then, \cref{lm:slicerob} implies that $\Gamma[S]$ is a robust $(\nu, \tau)$-outexpander with probability at least $1-\exp(-8\gamma^{-3}\nu^3\varepsilon^2n^2)$. Therefore, the probability that~$\Gamma[S]$ is a robust~$(\nu, \tau)$-outexpander for each such~$S$ is at least $1-2^n\exp(-8\gamma^{-3}\nu^3\varepsilon^2n^2)$.
	
	Thus, by a union bound, there exists a~$(\gamma, \varepsilon)$-almost regular~$\Gamma\subseteq D$ which is a robust~$(\nu, \tau)$-outexpander and such that~$D\setminus\Gamma$ is~$(\delta-\gamma, \varepsilon)$-almost regular and, for each~$S\subseteq V(D)$ with~$|S|\geq \varepsilon n$, if~$D[S]$ is a robust~$(2\gamma^{-1}\nu, \tau)$-outexpander then~$\Gamma[S]$ is also a robust~$(\nu, \tau)$-outexpander.
	
	It now suffices to check that for any integer~$k\geq \varepsilon n$, if~$S\subseteq V(D)$ is chosen uniformly among the subsets of~$V(D)$ of size~$k$, then~$D[S]$ is a robust~$(2\gamma^{-1}\nu, \tau)$-outexpander with probability at least~$1-n^{-3}$.	
	Fix an integer~$k\geq \varepsilon n$ and let~$S\subseteq V(D)$ satisfy~$|S|=k$.
	Then, for any~$v\in V(D)$, $\mathbb{E}[d_{D[S]}^\pm(v)]=(\delta\pm\varepsilon)|S|$ and, by \cref{lm:Chernoff},
	\[\mathbb{P}\left[d_{D[S]}^\pm(v)<\left(\frac{3}{8}+\gamma\right)|S|\right]\leq \mathbb{P}\left[d_{D[S]}^\pm(v)<\frac{9}{10}\mathbb{E}[d_{D[S]}^\pm(v)]\right]\leq \exp\left(-\varepsilon^2 n\right).\]
	Therefore, by \cref{lm:3/8rob},~$D[S]$ is a robust~$(2\gamma^{-1}\nu, \tau)$-outexpander with probability at least $1-n\exp(-\varepsilon^2 n)\geq 1-n^{-3}$. 
	This completes the proof.
\end{proof}

The following result is an easy \NEW{and well-known} consequence of \cref{lm:Chernoff}.

\begin{lm}\label{lm:partition}
	Let $0<\frac{1}{n}\ll\frac{1}{k}, \varepsilon, \delta \ll 1$. Let~$D$ be a~$(\delta, \varepsilon)$-almost regular digraph on~$n$ vertices. Let $n_1, \dots, n_k\in \mathbb{N}$ be such that $\sum_{i\in [k]}n_i=n$ and, for each~$i\in [k]$, $n_i=\frac{n}{k}\pm 1$. Assume $V_1, \dots, V_k$ is a random partition of~$V(D)$ such that, for each~$i\in [k]$, $|V_i|=n_i$. Then, with probability at least~$1-n^{-1}$,
	the following holds. For each~$i\in [k]$ and~$v\in V(D)$, $|N_D^\pm(v)\cap V_i|=(\delta\pm 2\varepsilon)\frac{n}{k}$.
\end{lm}

	\COMMENT{\begin{proof}
		Let $i\in [k]$ and $v\in V(D)$. Then, $\mathbb{E}[|N_D^\pm(v)\cap V_i|]=\frac{d_D^\pm(v)n_i}{n}=(\delta\pm \varepsilon)n_i$. Therefore, by \cref{lm:Chernoff}, 
		\begin{align*}
			\mathbb{P}\left[|N_D^\pm(v)\cap V_i|\neq(\delta\pm 2\varepsilon)\frac{n}{k}\right]&\leq\mathbb{P}\left[|N_D^\pm(v)\cap V_i|\neq(\delta\pm 3\varepsilon)n_i\right]\\
			&\leq \mathbb{P}\left[|N_D^\pm(v)\cap V_i|\neq(1\pm 2\varepsilon)\mathbb{E}[|N_D^\pm(v)\cap V_i|]\right]\\
			&\leq  2\exp\left(-\frac{2\varepsilon^2\delta n}{3k}\right).
		\end{align*}
		A union bound gives that the partition $V_1, \dots, V_k$ satisfies the desired properties with probability at least 
		\[1-2kn\exp\left(-\frac{2\varepsilon^2\delta n}{3k}\right)\geq 1-\frac{1}{n}.\] This completes the proof.
	\end{proof}}

\onlyinsubfile{\bibliographystyle{abbrv}
	\bibliography{Bibliography/Bibliography}}

%% file: Matchings.tex
	\onlyinsubfile{
		\setcounter{section}{4}
		\setcounter{subsection}{2}
		\setcounter{definition}{17}
\subsection{Some tools for finding matchings}}

In this \NEW{sub}section, we record two easy consequences of Hall's theorem which will enable us to construct matchings.

\COMMENT{Hall's theorem:
	Let $G$ be a bipartite graph on vertex classes $A$ and $B$. Then, $G$ contains a matching covering $A$ if and only if, for each $A'\subseteq A$, $|A'|\leq |N_G(A')|$.}

\begin{prop}\label{cor:Hall}
	Let~$G$ be a bipartite graph on vertex classes~$A$ and~$B$ with~$|A|\leq |B|$. Suppose that, for each~$a\in A$, $d_G(a)\geq\frac{|B|}{2}$ and, for each $b\in B$, $d_G(b)\geq |A|-\frac{|B|}{2}$.
	Then,~$G$ contains a matching covering~$A$. 
\end{prop}

	\COMMENT{\begin{proof}
		Let $A'\subseteq A$. If $|A'|\leq\frac{|B|}{2}$, then $|N_D(A')|\geq |A'|$. Otherwise, $N_G(A')=B$. Indeed, assume for a contradiction that $|A'|>\frac{|B|}{2}$ and $b\in B\setminus N_G(A')$. Then, $|A|-\frac{|B|}{2}\leq d_G(b)< |A|-\frac{|B|}{2}$, a contradiction. 
	\end{proof}}

\begin{prop}\label{cor:Hallreg}
	Let $0<\frac{1}{n}\ll\varepsilon \ll \delta\leq 1$. Let~$G$ be a bipartite graph on vertex classes~$A$ and~$B$ such that $|A|,|B|=(1\pm \varepsilon)n$. Suppose that, for each~$v\in V(G)$, $d_G(v)=(\delta\pm \varepsilon) n$. Then,~$G$ contains a matching of size at least~$\left(1-\frac{3\varepsilon}{\delta}\right)n$. 
\end{prop}

\COMMENT{\begin{proof}
	Let $B'$ be a set of $\left\lceil\frac{2\varepsilon n}{\delta}\right\rceil$ new vertices. Let $G'$ be the bipartite graph on vertex classes $A$ and $B\cup B'$ with $E(G')= E(G)\cup\{ab'\mid a\in A, b'\in B'\}$. Note that, for any $a\in A$, $d_{G'}(a)\geq \left(\delta+\frac{\varepsilon}{\delta}\right)n$ and, for any $b\in B$, $d_{G'}(b)\leq (\delta+\varepsilon)n$.\\
	Let $A'\subseteq A$. If $|A'|\leq \left(\delta+\frac{\varepsilon}{\delta}\right)n$, then, clearly, $|N_{G'}(A')|\geq |A'|$. Assume $|A'|\geq\left(\delta+\frac{\varepsilon}{\delta}\right)n$. Then, $|N_{G'}(A')|= |N_G(A')|+|B'|\geq \frac{(\delta-\varepsilon)n|A'|}{(\delta+\varepsilon)n}+|B'|\geq |A'|-\frac{2|A'|\varepsilon}{\delta+\varepsilon}+\frac{2\varepsilon n}{\delta}\geq |A'|$. Therefore, by Hall's theorem, $G'$ contains a matching $M'$ covering $A$. Let $M\coloneqq M'\cap E(G)$. Then $|M|\geq |A|-\left\lceil\frac{2\varepsilon n}{\delta}\right\rceil\geq \left(1-\frac{3\varepsilon}{\delta}\right)n$.
\end{proof}}

%% file: Excess.tex
	\onlyinsubfile{
		\setcounter{section}{4}
		\setcounter{subsection}{3}
		\setcounter{definition}{19}
\subsection{Some properties of the excess function}}

\NEW{We will need the following inequalities, which hold by definition of the excess function.}

\begin{fact}\label{fact:exc}
	Let~$D$ be a digraph and~$v\in V(D)$.
	\begin{enumerate}
		\item $\texc(D)\geq\Delta^0(D)\geq \frac{\Delta(D)}{2}\geq
		\frac{\delta(D)}{2}$.\label{fact:exc-delta}
		\item $d_D^{\min}(v)\NEW{=\min\{d_D^+(v),d_D^-(v)\}}=\frac{d_D(v)-|\exc_D(v)|}{2}$.\label{fact:exc-dmin}
		\item $d_D^{\max}(v)\NEW{=\max\{d_D^+(v),d_D^-(v)\}}=\frac{d_D(v)+|\exc_D(v)|}{2}$.\label{fact:exc-dmax}
		\item~$\texc(D)\geq \Delta^0(D)\geq d_D^{\max}(v)=d_D^{\min}(v)+|\exc_D(v)|$.\label{fact:exc-excv}
	\end{enumerate}
\end{fact}

Given a digraph~$D$ and~$S\subseteq V(D)$, \NEW{define}\OLD{denote} $\exc_D^\pm(S)\coloneqq
\sum_{v\in S}\exc_D^\pm(v)$.
\NEW{The next \lcnamecref{fact:partialdecompS} follows immediately from \cref{eq:exc}.}
\OLD{Then, observe that the following holds.}

\begin{fact}\label{fact:partialdecompS}
	Let~$D$ be a digraph\OLD{,~$V\coloneqq V(D)$,} and~$S\subseteq V$. Then,
	$\exc(D)=\exc_D^\pm(V\NEW{(D)})=\exc_D^\pm (S)+\exc_D^\pm (V\NEW{(D)}\setminus S)$.
\end{fact}

\begin{fact}\label{prop:texcT}
	\NEW{Any tournament~$T$ on~$n$ vertices which is not regular satisfies~$\texc(T)\geq \Delta^0(T)\geq
		\left\lceil\frac{n}{2}\right\rceil$.}\OLD{Any tournament~$T\notin\cT_{\rm reg}$ on~$n$ vertices satisfies~$\texc(T)\geq
	\left\lceil\frac{n}{2}\right\rceil$.}
\end{fact}

%% file: Exceptional_Tournaments.tex
	\onlyinsubfile{
		\setcounter{section}{4}
\section{Exceptional tournaments}}

Recall the definition of the class $\cT_{\rm excep}=\cT_{\rm reg}\cup \cT_{\rm apex}$ of exceptional tournaments from \cref{sec:intro}.
The main purpose of this section is to prove \cref{thm:oriented3/8} as well as the following result. 

\begin{thm}\label{thm:annoyingT}
	There exists~$n_0\in \mathbb{N}$ such that any tournament~$T\in \cT_{\rm excep}$ on~$n\geq n_0$ vertices satisfies~$\pn(T)=\texc(T)+1$.
\end{thm}

By \cref{lm:3/8rob}, \cref{thm:oriented3/8} (and thus also \cref{thm:annoyingT} in the case when~$T\in \cT_{\rm reg}$) is an immediate corollary of the following result.

\begin{thm}\label{cor:reg}
	Let $0<\frac{1}{n}\ll \nu\ll \tau\leq\frac{\delta}{2} \leq 1$ and $r\geq \delta n$. Let~$D$ be a\NEW{n}~$r$-regular digraph on~$n$ vertices. Assume~$D$ is a robust~$(\nu, \tau)$-outexpander. Then, $\pn(D)=\texc(D)+1=r+1$. 
\end{thm}

\begin{proof}
	\NEW{By \cref{eq:pnreg}, we have $\texc(D)=r$ and $\pn(D)\geq r+1$.}\OLD{Clearly,~$\texc(D)= r$.}
	Let $P\coloneqq v_1\dots v_{r+1}$ be a path of~$D$. 
	Then, by \cref{lm:verticesedgesremovalrob}\NEW{\cref{lm:verticesedgesremovalrob-edges}},~$D\setminus P$ is a robust~$(\frac{\nu}{2}, \tau)$-outexpander.
	Let $X^+\coloneqq \{v_{r+1}\}$, $X^-\coloneqq \{v_1\}$, $X^*\coloneqq \{v_2, \dots, v_r\}$, and $X^0\coloneqq V(D)\setminus (X^+\cup X^-\cup X^*)$. Applying \cref{thm:niceD} with $D\setminus P$ and $\frac{\nu}{2}$ playing the roles of~$D$ and~$\nu$ completes the proof.
\end{proof}

In order to prove \cref{thm:annoyingT} for~$T\in \cT_{\rm apex}$, we need the following result.

\begin{prop}\label{prop:annoyingT}
	Any~$T\in \cT_{\rm apex}$ on~$n$ vertices satisfies \NEW{$\exc(T)=n-3$ and} $\pn(T)\geq \texc(T)+1=n-1$.
\end{prop}

\begin{proof}
	Denote by $v_\pm\in V(T)$ the unique vertices such that $v_\pm\in U^\pm(T)$ and $V_0\coloneqq V(T)\setminus \{v_+,v_-\}=U^0(T)$.
	Thus, $v^-v^+\in E(T)$\OLD{and one can easily verify that $\texc(T)=n-2$}.
	
	\begin{claim}
		\NEW{$\exc(T)=n-3$ and $\texc(T)=n-2$.}
	\end{claim}
	
	\begin{proofclaim}
		\NEW{By definition of $\cT_{\rm apex}$, we have $d_T^+(v_+)=n-2=d_T^-(v_-)$ and $d_T^-(v_+)=1=d_T^+(v_-)$. Moreover, each $v\in V_0$ satisfies $d_T^+(v)=\frac{n-1}{2}=d_T^-(v)$. Therefore, $\Delta^0(T)=n-2$ and $\exc(T)=\frac{1}{2}\sum_{v\in V(T)}|d_T^+(v)-d_T^-(v)|=n-3$.}
	\end{proofclaim}

	\NEW{It remains to} show that $\pn(T)\geq n-1$. Let~$P\subseteq T$ be a path containing the edge~$v_-v_+$. It suffices to show that $\pn(T\setminus P) \geq n-2$.
	
	Let~$v$ be the starting point of~$P$. Observe that, since $v_-v_+\in E(P)$, we have~$v\neq v_+$. If~$v=v_-$, then $\exc(D\setminus P)\geq \exc_{D\setminus P}^-(v_-)=n-2$; otherwise,~$v\in U^0(T)$ and so $\exc(D\setminus P)\geq \exc_{D\setminus P}^-(v_-)+\exc_{D\setminus P}^-(v)=(n-3)+1=n-2$.
	Thus, we have shown that $\exc(T\setminus P)\geq n-2$. By \cref{thm:pn>exc},~$T\setminus P$ cannot be decomposed into fewer than $\NEW{\texc(T\setminus \cP)\geq} n-2$ paths. Therefore, $\pn(T)\geq 1+(n-2)=\texc(T)+1$.
\end{proof}

\begin{proof}[Proof of \cref{thm:annoyingT}]
	By \cref{lm:3/8rob,cor:reg}, we may assume that~$T \in T_{\rm apex}$.
	Fix additional constants such that $0<\frac{1}{n_0}\ll \nu\ll\tau\ll 1$. Let~$T\in \cT_{\rm apex}$ be a tournament on~$n\geq n_0$ vertices. By \cref{prop:annoyingT}, $\pn(T)\geq \texc(T)+1=n-1$. Thus, it suffices to find a path decomposition of~$T$ of size~$n-1$.
	
	Let $v_\pm\in V(T)$ denote the unique vertices such that $v_\pm\in U^\pm(T)$. Let $V'\coloneqq U^0(T)$\OLD{and $W\coloneqq V(T)\setminus V'=\{v_+,v_-\}$}. 
	Let $v_1, \dots, v_{n-2}$ be an enumeration of~$V'$ \OLD{$\ell\coloneqq \frac{n-1}{2}$,} and $r\coloneqq \frac{n-3}{2}$. 
	Since~$T[V']$ is a regular tournament on~$n-2$ vertices, \cref{lm:3/8rob} implies that~$T[V']$ is a robust~$(\nu, \tau)$-outexpander.
	Thus, by \cref{lm:robHampath}, we may assume without loss of generality that \NEW{$v_1 \dots v_{r+1}$}\OLD{$v_1\dots v_\ell$} is a path in~$T[V']$. 
	
	\NEW{Define a set of $r+2$ paths in $T$ by \[\cP\coloneqq \{v_-v_+,v_+v_1\dots v_{r+1}v_-, v_+v_{r+2}v_-, \dots, v_+v_{n-2} v_-\}.\]
	We now decompose $T\setminus \cP$ into $(n-1)-(r+2)=r$ paths. Note that $d_{T\setminus \cP}^\pm(v_\pm)=d_{T\setminus \cP}(v_\pm)=r$. Thus, each path must start at $v_+$ and end at $v_-$.
	Let $A^+\coloneqq \{v_+v_i\mid 2\leq i\leq r+1\}$ and $A^-\coloneqq \{v_iv_-\mid i\in [r]\}$.
	Denote $D\coloneqq T\setminus (A^+\cup A^-\cup \cP)$. Then, $d_D^\pm(v_+)=0=d_D^\pm(v_-)$. Moreover, each $i\in [r]$ satisfies $d_D^+(v_i)=\frac{n-1}{2}-2=r-1$ and each $j\in [n-2]\setminus [r]$ satisfies $d_D^+(v_j)=\frac{n-1}{2}-1=r$.
	Similarly, each $i\in \{2,\dots, r+1\}$ satisfies $d_D^+(v_i)=\frac{n-1}{2}-2=r-1$ and each $j\in [n-2]\setminus \{2, \dots,r+1\}$ satisfies $d_D^+(v_j)=\frac{n-1}{2}-1=r$.
	Let $X^+\coloneqq \{v_{r+1}\}$, $X^-\coloneqq \{v_1\}$, $X^*\coloneqq \{v_2, \dots, v_r\}$, and $X^0\coloneqq \{v_i\mid i\in [n-2]\setminus [r+1]\}$.
	Then, $|X^+\cup X^*|=r=|X^-\cup X^*|$ and \cref{eq:niceD} holds with $D-\{v_+,v_-\}$ playing the role of $D$.
	By \cref{thm:niceD}, there exists a path decomposition $\cP'=\{P_1, \dots, P_r\}$ of $D-\{v_+,v_-\}$. For each $i\in [r]$, let $w_i^+$ and $w_i^-$ denote the starting and ending points of $P_i$. By the ``moreover part" of  \cref{thm:niceD}, we may assume that $w_1^+, \dots, w_r^+$ are distinct and $\{w_i^+\mid i\in [r]\}=X^+\cup X^*$. We may also assume that $w_1^-, \dots, w_r^-$ are distinct and $\{w_i^-\mid i\in [r]\}=X^-\cup X^*$. Thus, $\cP''\coloneqq \{v_+w_i^+P_iw_i^-v_-\mid i\in [r]\}$ is a path decomposition of $D\cup A^+\cup A^-=T\setminus \cP$. Therefore, $\cP\cup \cP''$ is a path decomposition of $T$ of size $2r+2=n-1$. I.e.\ $\pn(T)\leq n-1$, as desired.}
	\OLD{Let $\ell\coloneqq \frac{n-1}{2}$. Let $\cP\coloneqq \{v_-v_+,v_+v_1\dots v_\ell v_-, v_+v_{\ell+1}v_-, \dots, v_+v_{n-2} v_-\}$ and~$D\coloneqq T\setminus \cP$. Note that~$|\cP|=r+2$.
	Moreover, $d_D(v_\pm)=\ell-1=r$ and, for each~$i\in [n-2]$, $d_D(v_i)=n-3=2r$.
	Let $A^+\coloneqq \{v_+v_i\mid 2\leq i\leq \ell\}$ and $A^-\coloneqq \{v_iv_-\mid 1\leq i\leq \ell-1\}$. Define $X^+\coloneqq X^-\coloneqq X^*\coloneqq \emptyset$, and $X^0\coloneqq \{v_i\mid i\in [n-2]\}$. Note that $|X^\pm\cup X^*\cup A^\pm|=\ell-1=r$. Thus, we can apply \cref{prop:absorbingedges} with~$n-2$ and~$\frac{1}{4}$ playing the roles of~$n$ and~$\delta$ to obtain a path decomposition~$\cP'$ of~$D$ of size~$r$. Then, $\cP\cup \cP'$ is a path decomposition of size $r+2+r=n-1$, as desired.}
\end{proof}

\NEW{In \cref{sec:A}, we will introduce the concept of ``absorbing edges" which plays a similar role as the edge sets $A^+$ and $A^-$ above.}

We will need the following observation about tournaments in~$\cT_{\rm apex}$ \NEW{for later}.

\begin{prop}\label{fact:annoyingT}
	A tournament~$T$ satisfies $|U^+(T)|=|U^-(T)|=1$, $e(U^-(T),U^+(T))=1$ and $\texc(T)-\exc(T)<2$ if and only if $T\in \cT_{\rm apex}$.
\end{prop}

\begin{proof}
	\NEW{Suppose that $T\in \cT_{\rm apex}$. By definition and \cref{prop:annoyingT}, $|U^+(T)|=|U^-(T)|=1$, $e(U^-(T),U^+(T))=1$, $\texc(T)=n-2$, and $\exc(T)=n-3$.}
	
	\NEW{Suppose that $|U^+(T)|=|U^-(T)|=1$ and $e(U^-(T),U^+(T))=1$ and $\texc(T)-\exc(T)<2$. Let $v_\pm\in U^\pm(T)$. By \cref{eq:exc}, we have $\exc_T^+(v_+)=\exc_T^-(v_-)$ and so, as $T$ is a tournament, $d_T^+(v_+)=d_T^-(v_-)$.
	Since $e(U^-(T),U^+(T))=1$, $d_T^\mp(v_\pm)\geq 1$. On the other hand, \[2>\texc(T)-\exc(T)\geq \Delta^0(T)-\exc(T)\geq d_T^\pm(v_\pm)-\exc_T^\pm(v_\pm)\stackrel{\text{\cref{fact:exc}\cref{fact:exc-excv}}}{=}d_T^\mp(v_\pm).\] Therefore, $N^\mp(v_\pm)=\{v_\mp\}$ and $N^\pm(v_\pm)=U^0(T)$. In particular, $T-\{v_+,v_-\}=T[U^0(T)]$ is regular. Hence, $T\in \cT_{\rm apex}$, as required.}
\end{proof}

%% file: Corollaries.tex
	\onlyinsubfile{
		\setcounter{section}{5}
\section{Deriving Theorem \ref{thm:even} and Corollary \ref{cor:approx} from Theorem \ref{thm:main}}}

\NEW{In this section, we assume that \cref{thm:main} holds and derive \cref{thm:even,cor:approx}.}
\NEW{For \cref{thm:even}, we first observe that if~$n$ is even, then~$\texc(T)=\exc(T)$.}

\begin{prop}\label{prop:even}
	Let~$T$ be a tournament of even order~$n$. Then,~$\texc(T)=\exc(T)$ \NEW{and $U^0(T)=\emptyset$.}
\end{prop}

\begin{proof}
	It is easy to see that each~$v\in V(T)$ satisfies~$\exc_T(v)\neq 0$ \NEW{and so $U^0(T)=\emptyset$}.
	Let~$v\in V(T)$ be such that $d_T^{\max}(v)=\Delta^0(T)$. Thus,
	\begin{equation*}
		\exc(T)=\frac{1}{2}\sum_{u\in V(T)}|\exc_T(u)|\geq \frac{n-1+|\exc_T(v)|}{2}\stackrel{\text{\cref{fact:exc}\cref{fact:exc-dmax}}}{=}d_T^{\max}(v)=\Delta^0(T),
	\end{equation*}
	so $\texc(T)=\exc(T)$, as desired.
\end{proof}

\begin{proof}[Proof of \cref{thm:even}]
	Let $0<\frac{1}{n_0}\ll \beta\ll 1$. Let~$n\geq n_0$ be even and~$T$ be a tournament on~$n$ vertices. It is easy to see that~$T\notin\cT_{\rm excep}$. We show that one of \cref{thm:main}\cref{thm:main-largeU,thm:main-largeexc} holds.
	Suppose that \cref{thm:main}\cref{thm:main-largeU} does not hold. Without loss of generality, we may assume that~$N^+(T)\leq \beta n$. Thus, $|U^+(T)|\leq \beta n$. Since~$n$ is even, each~$v\in V(T)$ satisfies~$\exc_T(v)\neq 0$. Thus, $\texc(T)\geq \exc(T)\geq |U^-(T)|=n-|U^+(T)|\geq n-\beta n\geq \frac{n}{2}+\beta n$ and so \cref{thm:main}\cref{thm:main-largeexc} holds.
	Therefore, by \cref{thm:main,prop:even}, $\pn(T)=\texc(T)=\exc(T)$.
\end{proof}

Finally, we derive \cref{cor:approx} from \cref{thm:main}. The idea is that if none of \cref{thm:main,thm:even,thm:annoyingT} apply to~$T$ then we can transform~$T$ into a tournament~$T'$ which satisfies the conditions of \cref{thm:main} by flipping a small number of edges, and so that~$\pn(T)\sim \pn(T')$ and~$\texc(T)\sim \texc(T')$.

\begin{proof}[Proof of \cref{cor:approx}]
	We may assume without loss of generality that~$\beta \ll 1$.
	Let $0<\frac{1}{n_0}\ll \beta\ll 1$. Let~$T$ be a tournament on~$n\geq n_0$ vertices. By \cref{thm:annoyingT,thm:even}, we may assume that~$T\notin \cT_{\rm excep}$ and that~$n$ is odd.
	If $\Delta^0(T)\geq \frac{n}{2}+\frac{\beta n}{5}$, then, by \cref{thm:main} applied with~$\frac{\beta}{5}$ playing the role of~$\beta$, $\pn(T)=\texc(T)$ and we are done.
	We may therefore assume that $\Delta^0(T)< \frac{n}{2}+\frac{\beta n}{5}$.
	Let~$v\in V\NEW{(T)}$.
	Since~$T$ is not regular, we may assume without loss of generality that~$v\in U^+(T)$. Then, note that $d_T^+(v)\geq \frac{n+1}{2}$.
	Let $S\subseteq N_T^-(v)$ satisfy $|S|=\lceil\frac{n}{2}+\frac{\beta n}{5}\rceil-d_T^+(v)$ (this is possible since $d_T^-(v)=(n-1)- d_T^+(v)$). Note that $|S|\leq \lceil\frac{n}{2}+\frac{\beta n}{5}\rceil-\frac{n+1}{2}\leq \frac{\beta n}{4}$.
	
	Let~$T'$ be obtained from~$T$ by flipping the direction of all edges between~$v$ and~$S$.
	Then, observe that $\texc(T')\geq \Delta^0(T')\geq d_{T'}^+(v)=\lceil\frac{n}{2}+\frac{\beta n}{5}\rceil$ and, in particular,~$T'\notin \cT_{\rm excep}$.
	Moreover, we claim that $\texc(T')\leq \texc(T)+ 2|S|$. 
	Since $\Delta^0(T')\leq \Delta^0(T)+|S|$, it suffices to show that $\exc(T')\leq \exc(T)+2|S|$. Note that, by \cref{fact:exc}\cref{fact:exc-dmax}, $\exc_{T'}^+(v)-\exc_T^+(v)=2(d_{T'}^+(v)-d_T^+(v))=2|S|$.
	For each $\diamond\in\{+,-,0\}$, denote $S^\diamond\coloneqq S\cup U^\diamond(T)$. Then, by \cref{fact:exc}\cref{fact:exc-dmax}, for each~$u\in S^+$, $\exc_{T'}^+(u)-\exc_T^+(u)=-2$%
	\COMMENT{Note that $n$ odd and $u\in U^+(T)$ implies $\exc_T^+(u)\geq 2$ so $u$ does have non-negative excess.}
	and, for each $u\in S^-\cup S^0$, $\exc_{T'}^+(u)=0=\exc_T^+(v)$.
	Thus, $\exc(T')-\exc(T)=2|S|-2|S^+|\leq 2|S|$, as desired.
	
	By \cref{thm:main}, $\pn(T')=\texc(T')\leq \texc(T)+2|S|$ and thus, since $|S|\leq \frac{\beta n}{4}$, it suffices to show that $\pn(T)\leq \pn(T')+2|S|$.
	Let~$\cP'$ be a path decomposition of~$T'$ of size~$\pn(T')$. Let~$\cP_1'$ consist of all the paths~$P\in \cP'$ such that $E(P)\subseteq E(T)$. Let $\cP_2'\coloneqq \cP'\setminus \cP_1'$. Let~$\cP_2$ be set of paths obtained from~$\cP_2'$ by deleting all the edges in $E(T')\setminus E(T)$.
	Then, $\cP\coloneqq \cP_1'\cup \cP_2\cup (E(T)\setminus E(T'))$ is a path decomposition of~$T$. 
	Moreover, by construction, all the edges in $E(T')\setminus E(T)$ are incident to~$v$. Thus, each path in~$\cP_2'$ contains exactly one edge of $E(T')\setminus E(T)$ and so $|\cP_2\cup (E(T)\setminus E(T'))|\leq 3|\cP_2'|= |\cP_2'|+ 2|E(T')\setminus E(T)|=|\cP_2'|+2|S|$. Therefore, $\pn(T)\leq |\cP|\leq |\cP_1'|+|\cP_2'|+2|S|= \pn(T')+2|S|$. This completes the proof.
\end{proof}

%% file: Approximate_Decomposition.tex
	\onlyinsubfile{
		\setcounter{section}{6}
\section{Approximate decomposition of robust outexpanders}}

\NEW{The rest of the paper is devoted to the proof of \cref{thm:main}. We start by discussing and proving the approximate decomposition step, which is achieved via \cref{lm:approxdecomp}.}

As mentioned in the proof overview, in order to reduce the excess and the vertex degrees at the correct rate, we will approximately decompose our digraphs into sets of paths.
To do so, we will start by constructing auxiliary multidigraphs called \emph{layouts} which will prescribe the ``shape" of the structures in our approximate decomposition.

Suppose that we would like to find a Hamilton~$(v_+,v_-)$-path which contains a fixed edge~$f = u_+u_-$. 
We can view this as the task of finding two paths of shapes~$v_+u_+$ and~$u_-v_-$, respectively,  that are vertex-disjoint and cover all remaining vertices. 
(Recall from \cref{sec:notation} that, given an (auxiliary) edge~$uv$, \NEW{we say that} a~$(u,v)$-path has shape~$uv$.)
We now generalise this approach to \emph{layouts}, which will tell us the shapes of paths required, the set~$F$ of fixed edges to be included, and the vertices to be avoided by these paths.
The ``spanning'' extension of a layout will be called a \emph{spanning configuration}.
To ensure that the spanning configuration has a suitable path decomposition, we will define a layout to consist of a (multi)set of paths rather than a multiset of edges.
\NEW{The concepts of layout and spanning configuration are also illustrated in \cref{fig:layout}.}

We will be working with multidigraphs.
Let~$V$ be a vertex set. 
We say~$(L, F)$ is a \emph{layout} if the following hold.
\begin{enumerate}[label=\upshape(L\arabic*)]
	\item $L$ is a multiset consisting of paths on~$V$ and isolated vertices. \label{def:layout-L}
	\item $F \subseteq E(L)$. \label{def:layout-F}
	\item $E(L)\setminus F\neq \emptyset$.  \label{def:layout-unfixededge}
\end{enumerate}

Conditions \cref{def:layout-L,def:layout-F,def:layout-unfixededge} can be motivated as follows.%
\OLD{Firstly, the construction of the paths that we find in the approximate decomposition step will be based on the robust outexpansion property. But the robust outexpander that we work with will only span $V'=V \setminus W$ so we cannot automatically incorporate the set of edges~$F$ incident to~$W$ in the approximate decomposition step (here~$W$ will be the ``exceptional set" defined as in \cref{sec:sketch-general}).
As indicated in the above example, we will allocate the edges in~$F$ to some given paths in advance.
Replacing the edges of~$L$ which are not in~$F$ as in the above example now gives a spanning configuration (defined formally below) with endpoints induced by the endpoints of the paths in~$L$ and which contains all edges in~$F$. In particular, each path of the spanning configuration contains exactly the edges in~$F$ which were in the corresponding path in~$L$. Because we need to construct (almost) spanning structures, we need~$E(L)\setminus F$ to be non-empty, which is the reason for \cref{def:layout-unfixededge}.}
\NEW{Suppose that $(L,F)$ is a layout on $V$. As described above, our goal is to construct a spanning set of paths whose shapes correspond to those of the paths in $L$ and which contain the edges in~$F$. This will be achieved by replacing the edges in $E(L)\setminus F$ by internally vertex-disjoint paths which cover all the vertices in $V\setminus V(L)$. This motivates \cref{def:layout-unfixededge}: if $E(L)\setminus F$ was empty, then there would be no edge to replace by a path and so we would not be able cover the vertices in $V\setminus V(L)$, i.e.\ our set of paths would not be spanning.
Moreover, as we will be constructing the paths according to the shapes of the paths in $L$, we need to make sure that the fixed edge set $F$ is covered by $L$. This explains \cref{def:layout-F}.}
\NEW{Finally}\OLD{Secondly}, since we \NEW{will} have already covered some edges \NEW{before the approximate decomposition (recall from the proof overview that we will need a cleaning step)} and \NEW{since} not all vertices have the same excess, we do not actually want our structures to be completely spanning, but want them to avoid a suitable small set of vertices. This is why we allow paths of length~$0$ in \cref{def:layout-L}.

\NEW{Let $(L,F)$ be a layout on $V$.} A multidigraph~$\cH$ on~$V$ is a \emph{spanning configuration of shape~$(L,F)$} if~$\cH$ can be decomposed into internally vertex-disjoint paths $\{P_e \mid e \in E(L)\}$ such that each~$P_e$ has shape~$e$;~$P_f = f$ for all~$f \in F$; and \NEW{$\bigcup_{e \in E(L)}V^0(P_e) = V \setminus V(L)$}\OLD{$V^0(\bigcup_{e \in E(L)} \{P_e\} ) = V \setminus V(L)$} \NEW{(recall that given a path $P$, $V^0(P)$ denotes the set of internal vertices of $P$)}. (Note that the last equality implies that the isolated vertices of~$L$ remain isolated in~$\cH$.) \NEW{See \cref{fig:layout-H} for an example of a spanning configuration.}

\subfile{Figures/Figure_Layouts.tex}

\NEW{Let $(L,F)$ be a layout on $V$ and $\cH$ be a spanning configuration of shape $(L,F)$ on $V$.}
There is a natural bijection between a path~$Q$ in~$L$ and the path $P_Q \coloneqq \bigcup \{P_e \mid e \in E(Q)\}$ in~$\cH$. \NEW{(E.g. in the example presented in \cref{fig:layout}, the path $v_1e_1v_2e_2v_3$ in $L$ corresponds to the path $v_1P_{e_1}v_2P_{e_2}v_3$ in $\cH$.)}
Note that this bijection is not necessarily unique since if~$e$ has multiplicity more than~$1$ in~$L$, then there are different ways to define~$P_e$. \NEW{(E.g. in the example presented in \cref{fig:layout}, we could have exchanged $P_{e_1}$ and $P_{e_3}$.)}
A path decomposition~$\cP$ of~$\cH$ consisting of all such~$P_Q$ for all the paths~$Q\in L$ is said to be \emph{induced} by~$(L,F)$. \NEW{(E.g.\ in the example presented in \cref{fig:layout}, $\{v_1P_{e_1}v_2P_{e_2}v_3, v_1P_{e_3}v_2P_{e_4}v_{10}, v_6P_{e_5}v_7P_{e_6}v_8P_{e_7}v_9P_{e_8}v_{10}\}$ is a path decomposition of $\cH$ induced by $(L,F)$.)}
\NEW{Note that if $\cP$ is a path decomposition $\cH$ induced by $(L,F)$, then the paths in $\cP$ are non-trivial and have the same endpoints as their corresponding path in $L$.}
\OLD{(Note that each path in~$\cP$ is non-trivial.) Thus, our above example is a spanning configuration of shape $(\{v_+u_+u_-v_-\}, \{u_+u_-\})$, where $(\{v_+u_+u_-v_-\}, \{u_+u_-\})$ is a layout.}

\NEW{\begin{fact}\label{fact:layoutbijection}
	Let $(L,F)$ be a layout on $V$ and $\cH$ be a spanning configuration of shape $(L,F)$ on $V$. Let $L'$ denote the set of (non-trivial) paths contained in $L$ (i.e.\ $L'$ is obtained by deleting all the isolated vertices in $L$). For any path $Q$ in $L$, the corresponding path $P_Q\coloneqq \bigcup \{P_e \mid e \in E(Q)\}$ in $\cH$ satisfies $V^\pm(P_Q)=V^\pm(Q)$ and $V^0(Q)\subseteq V^0(P_Q)\subseteq V^0(Q)\cup (V\setminus V(L))$. Thus, if $\cP$ is a path decomposition of $\cH$ which is induced by $(L,F)$, then $V^\pm(\cP)=V^\pm(L')$ and $V^0(\cP)=V^0(L)\cup (V\setminus V(L))$.
\end{fact}}

\NEW{Let $V$ be a vertex set. Let $(L,F)$ be a layout on $V$ and $\cH$ be a spanning configuration of shape $(L,F)$ on $V$. By \cref{fact:layoutbijection}, the degree of each $v\in V$ in $\cH$ is entirely determined by the degree of $v$ in $L$. Thus, the following holds.}

\NEW{\begin{fact}\label{fact:layoutsV}
	Let~$D$ be a digraph on a vertex set~$V$ and $(L_1,F_1) \dots (L_\ell,F_\ell)$ be layouts on $V$.
	For each~$i\in [\ell]$, let~$\cH_i$ be a spanning configuration of shape~$(L_i,F_i)$. Suppose that $\cH_1, \dots, \cH_\ell$ are pairwise edge-disjoint.
	Then, for all~$v \in V$, 
	\[d^{\pm}_\cH(v)=\sum_{i\in [\ell]}d_{\cH_i}^\pm(v) 
	=\sum_{i\in [\ell]}(d_{L_i}^\pm(v) +\mathds{1}_{v\notin V(L_i)})= d^{\pm}_L(v)+ |\{i \in [\ell] \mid v \in V\setminus V(L_i)\}|.\]
\end{fact}
}

\NEW{Roughly speaking, the approximate decomposition \lcnamecref{lm:approxdecomp} says that given a dense almost regular digraph $D$ and a sparse almost regular robust outexpander $\Gamma$, we can transform a suitable set of small layouts into edge-disjoint spanning configurations of corresponding shape in $D\cup \Gamma$. Moreover, the number of layouts that we are allowed to prescribe is close to $\delta^0(D)$, in which case the the configurations form an approximate decomposition of $D\cup \Gamma$.}

\begin{lm}[Approximate decomposition lemma for robust outexpanders]\label{lm:approxdecomp}
	Let $0<\frac{1}{n}\ll \varepsilon \ll \nu\ll \tau\ll \gamma \ll\eta, \delta \leq 1$. Suppose~$\ell\in \mathbb{N}$ satisfies $\ell \leq (\delta-\eta)n$. If~$\ell\leq \varepsilon^2n$, then let~$p\leq n^{-1}$; otherwise, let~$p\leq n^{-2}$.	
	Let~$D$ and~$\Gamma$ be edge-disjoint digraphs on a common vertex set~$V$ of size~$n$. 
	Suppose that~$D$ is~$(\delta, \varepsilon)$-almost regular and~$\Gamma$ is~$(\gamma,\varepsilon)$-almost regular. Suppose further that~$\Gamma$ is an~$(\varepsilon,p)$-robust~$(\nu, \tau)$-outexpander.
	Let~$\cF$ be a multiset of directed edges on~$V$. Any edge in~$\cF$ is considered to be distinct from the edges of~$D\cup \Gamma$, even if the starting and ending points are the same (recall \cref{sec:notation}).
	Let $F_1, \dots, F_\ell$ be a partition of~$\cF$.	
	Assume that $(L_1, F_1), \dots, (L_\ell, F_\ell)$ are layouts such that~$V(L_i)\subseteq V$ for each~$i\in [\ell]$ and the following hold, where $L\coloneqq \bigcup_{i\in [\ell]}L_i$.
	\begin{enumerate}
		\item For each~$i\in [\ell]$, $|V(L_i)|\leq \varepsilon^2 n$ and $|E(L_i)|\leq \varepsilon^4 n$.\label{lm:approxdecomp-sizeLi}
		\item Moreover, for each~$v\in V$, $d_L(v)\leq \varepsilon^3n$ and there exist at most~$\varepsilon^2 n$ indices~$i\in [\ell]$ such that~$v\in V(L_i)$.\label{lm:approxdecomp-degreev}
	\end{enumerate}
	Then, there exist edge-disjoint \NEW{submultidigraphs} $\cH_1, \dots, \cH_\ell\subseteq D\cup \Gamma\cup \cF$ such that, for each~$i\in [\ell]$,~$\cH_i$ is a spanning configuration of shape~$(L_i,F_i)$ and the following hold, where $\cH\coloneqq \bigcup_{i\in [\ell]} \cH_i$, $D'\coloneqq D\setminus \cH$, and $\Gamma'\coloneqq \Gamma\setminus \cH$. 
	\begin{enumerate}[label=\rm(\roman*)]
		\item If $\ell\leq \varepsilon^2 n$, then~$\Gamma'$ is obtained from~$\Gamma$ by removing at most~$3\varepsilon^3\nu^{-4}n$ edges incident to each vertex\NEW{, i.e. $\Delta(\Gamma\setminus \Gamma')\leq 3\varepsilon^3\nu^{-4}n$}.\label{lm:approxdecomp-smallchunks}
		\item If $\ell \leq \nu^5 n$, then~$D'$ is~$(\delta-\frac{\ell}{n}, 2\varepsilon)$-almost regular and~$\Gamma'$ is~$(\gamma,2\varepsilon)$-almost regular. Moreover, $\Gamma'$ is a~$(\sqrt{\varepsilon},p)$-robust~$(\nu-\sqrt{\varepsilon}, \tau)$-outexpander.\label{lm:approxdecomp-bigchunks}
		\item $D'\cup \Gamma'$ is a robust~$(\frac{\nu}{2}, \tau)$-outexpander.\label{lm:approxdecomp-all}
	\end{enumerate} 
\end{lm}

The approximate decomposition \NEW{guaranteed by \cref{lm:approxdecomp}} is constructed in stages.
The core of the approximate decomposition occurs in \cref{lm:approxdecomp}\cref{lm:approxdecomp-smallchunks}, where a small set of layouts is converted into spanning configurations one by one \NEW{(see \cref{sec:sketch-simple})}.
Repeated applications of \cref{lm:approxdecomp}\cref{lm:approxdecomp-smallchunks} will then enable us to transform larger sets of layouts into spanning configurations (\cref{lm:approxdecomp}\cref{lm:approxdecomp-bigchunks}).
Then, one can obtain the final approximate decomposition (\cref{lm:approxdecomp}\cref{lm:approxdecomp-all}) by repeatedly applying \cref{lm:approxdecomp}\cref{lm:approxdecomp-bigchunks}, adjusting the parameters in each iteration.
This can be seen as a semirandom ``nibble" process, where the applications of \cref{lm:approxdecomp}\cref{lm:approxdecomp-smallchunks} are the ``nibbles" (which are chosen via a probabilistic argument) and the applications of \cref{lm:approxdecomp}\cref{lm:approxdecomp-bigchunks} correspond to ``bites" consisting of several ``nibbles".
We prove \cref{lm:approxdecomp-bigchunks}, \cref{lm:approxdecomp-all}, and \cref{lm:approxdecomp-smallchunks} in this order.

\begin{proof}[Proof of \cref{lm:approxdecomp}\cref{lm:approxdecomp-bigchunks}]
	Let $\ell'\coloneqq \left\lfloor \varepsilon^2n\right\rfloor$ and $k\coloneqq \left\lceil\frac{\ell}{\ell'}\right\rceil$. Note that $k\leq 2 \nu^5 \varepsilon^{-2}$.
	We now group $(L_1, F_1), \dots, (L_\ell, F_{\ell})$ into~$k$ batches, each of size at most~$\ell'$. \NEW{For each $m\in [k]$, the $m^{\rm th}$ batch will consist of $(L_i, \cF_i)$ with $(m-1)\ell'\leq i\leq \min\{m\ell', \ell\}$.}
	We aim to apply \cref{lm:approxdecomp}\cref{lm:approxdecomp-smallchunks} to each batch in turn.
	
	Assume that \NEW{we have done $m$ batches} for some~$0\leq m\leq k$. \NEW{This means that} we have constructed edge-disjoint $\cH_1, \dots, \cH_{\min\{m\ell',\ell\}}\subseteq D\cup \Gamma\cup \cF$ such that, for each $i\in [\min\{m\ell',\ell\}]$,~$\cH_i$ is a spanning configuration of shape~$(L_i,F_i)$ \NEW{satisfying $E(\cH_i)\cap E(\cF)=E(F_i)$ and} the following holds. Let $\Gamma_m\coloneqq \Gamma\setminus \bigcup_{i\in [\min\{m\ell',\ell\}]}\cH_i$. 
	Then, for each~$v\in V$,
	\begin{equation}\label{eq:approxdecomp-NGamma}
	|N_{\Gamma\setminus \Gamma_m}(v)|\leq m\cdot 25\varepsilon^3\nu^{-4}n\leq 50 \eps \nu n \le \frac{\eps n}{2}.
	\end{equation}%
	\OLD{\begin{equation*}
		|N_{\Gamma\cap \bigcup_{i\in [\min\{m\ell',\ell\}]} \cH_i}(v)|=|N_{\Gamma\setminus \Gamma_m}(v)|\leq 25\varepsilon^3\nu^{-4}mn\leq 50 \eps \nu n \le \frac{\eps n}{2}.
	\end{equation*}}%
	Let $D_m\coloneqq D\setminus \bigcup_{i\in [\min\{m\ell',\ell\}]}\cH_i$.
	Observe that, by \cref{fact:layoutsV} and \cref{lm:approxdecomp-degreev}, \NEWTWO{$\bigcup_{i \in [ \min\{m \ell', \ell\}]} \mathcal{H}_i\setminus F_i$} is $(\frac{ \min\{m \ell', \ell\}}{n}, \eps^2+\eps^3)$-almost regular. 
	Together with \cref{eq:approxdecomp-NGamma}, this implies that~$D_m$ is $(\delta-\frac{\min\{m\ell',\ell\}}{n},2\varepsilon)$-almost regular and~$\Gamma_m$ is~$(\gamma, 2\varepsilon)$-almost regular%
		\COMMENT{\label{com:stillreg}$\Gamma_m$ is obtained from $\Gamma$ by removing at most $25\varepsilon^3\nu^{-4}mn\leq 50\varepsilon\nu n\leq \varepsilon n$ in/out edges at each $v\in V$.\\
		$D_m$ is obtained from $D$ by removing at most $\min\{m\ell',\ell\}+\varepsilon^3 n$ (extra $\varepsilon^3n$ to take into account paths starting/ending at $v$ in $L$) and at least $\min\{m\ell',\ell\}-25\varepsilon^3\nu^{-4}mn-\varepsilon^2 n\geq \min\{m\ell',\ell\}-50\varepsilon\nu n-\varepsilon^2 n\geq \min\{m\ell',\ell\}-\varepsilon n$ in/out edges at each $v\in V$. (Extra $\varepsilon^2 n$ to take into account when $v\in V(L_i)$.)}.

	Moreover, by \cref{fact:verticesedgesremovalrob},~$\Gamma_m$ is a~$(\sqrt{\varepsilon},p)$-robust~$(\nu-\sqrt{\varepsilon}, \tau)$-outexpander.
	Thus, if~$m=k$, we are done.
	
	Suppose~$m<k$. We show that~$\Gamma_m$ is a~$(2\varepsilon, n^{-1})$-robust~$(\nu-\varepsilon,\tau)$-outexpander. If~$m=0$, then~$\Gamma_m=\Gamma$ and we are done. We may therefore assume that~$m\geq 1$. Then, note that~$k\geq 2$ so $\ell>\ell'=\lfloor\varepsilon^2 n\rfloor$ and, thus,~$p\leq n^{-2}$. 
	Fix an integer $k'\geq 2\varepsilon n$. Suppose~$S\subseteq V$ is a random subset of size~$k'$. We show that~$\Gamma_m[S]$ is a robust~$(\nu-\varepsilon,\tau)$-robust outexpander with probability at least~$1-n^{-1}$. 
	Let~$v\in V$. 
	If $|N_{\Gamma\setminus \Gamma_m}(v)|\leq \varepsilon^2 n$, then $|N_{\Gamma\setminus \Gamma_m}(v)\cap S|\leq \varepsilon^2 n\leq \varepsilon k'$. Suppose $|N_{\Gamma\setminus \Gamma_m}(v)|\geq \varepsilon^2 n$.
	Then, by \cref{eq:approxdecomp-NGamma}, $\mathbb{E}[|N_{\Gamma\setminus \Gamma_m}(v)\cap S|]=\frac{k'}{n}|N_{\Gamma\setminus \Gamma_m}(v)|\leq \frac{\varepsilon k'}{2}$. Thus, \cref{lm:Chernoff} implies that
	\[\mathbb{P}\left[|N_{\Gamma\setminus \Gamma_m}(v)\cap S|>\varepsilon k'\right]\leq \mathbb{P}\left[|N_{\Gamma\setminus \Gamma_m}(v)\cap S|>2\mathbb{E}[|N_{\Gamma\setminus \Gamma_m}(v)\cap S|]\right]\leq \exp\left(-\frac{2\varepsilon^3 n}{3}\right).\]
	Therefore, by a union bound, with probability at least $1- n\exp\left(-\frac{2\varepsilon^3 n}{3}\right)$, the digraph~$\Gamma_m[S]$ is obtained from~$\Gamma[S]$ by removing at most~$\varepsilon k'$ edges incident to each vertex. Our assumption on~$\Gamma$ implies that~$\Gamma[S]$ is a robust~$(\nu, \tau)$-outexpander with probability at least $1-p\geq 1-n^{-2}$. Therefore, by \cref{lm:verticesedgesremovalrob}\NEW{\cref{lm:verticesedgesremovalrob-edges}}, we conclude that~$\Gamma_m[S]$ is a robust~$(\nu-\varepsilon,\tau)$-outexpander with probability at least $1-p- n\exp\left(-\frac{2\varepsilon^3 n}{3}\right)\geq 1-n^{-1}$. Thus,~$\Gamma_m$ is a~$(2\varepsilon, n^{-1})$-robust~$(\nu-\varepsilon,\tau)$-outexpander.
		
	Let $\ell''\coloneqq \min\{\ell-m\ell', \ell'\}$ and $\cF'\coloneqq \bigcup_{i\in [\ell'']}\cF_{m\ell'+i}$. Apply \cref{lm:approxdecomp}\cref{lm:approxdecomp-smallchunks}%
		\COMMENT{Note that $\delta-\frac{m\ell'}{n}\geq \delta-\frac{\ell}{n}\geq \eta\gg \gamma$ so the hierarchy is still okay. Moreover, $\delta\geq \frac{\ell}{n}+\eta\geq \frac{m\ell'+\ell''}{n}+\eta$ so $\ell''\leq (\delta-\frac{m\ell'}{n}-\eta)n$, as needed.}
	with $D_m$, $\Gamma_m$, $\cF'$, $n^{-1}$, $\delta-\frac{m\ell'}{n}$, $\nu-\varepsilon$, $2\varepsilon$, $\ell''$, $L_{m\ell'+1}, \dots, L_{m\ell'+\ell''}$, and $F_{m\ell'+1}, \dots, F_{m\ell'+\ell''}$ playing the roles of $D$, $\Gamma$, $\cF$, $p$, $\delta$, $\nu$, $\varepsilon$, $\ell$, $L_1, \dots, L_\ell$, and $F_1, \dots, F_\ell$ to obtain edge-disjoint $\cH_{m\ell'+1}, \dots, \cH_{m\ell'+\ell''}\subseteq D_m\cup \Gamma_m\cup \cF'$ such that, for each~$i\in [\ell'']$,~$\cH_{m\ell'+i}$ is a spanning configuration of shape~$(L_{m\ell'+i},F_{m\ell'+i})$ and, for each~$v\in V$, $|N_{\Gamma_m\setminus \Gamma_{m+1}}(v)|\leq 3(2\varepsilon)^3(\nu-\varepsilon)^{-4}n\leq 25\varepsilon^3\nu^{-4}n$, where $\Gamma_{m+1}\coloneqq \Gamma_m\setminus \bigcup_{i\in [\ell'']}\cH_{m\ell'+i}$.
	In particular, \cref{eq:approxdecomp-NGamma} holds. This completes the proof. 
\end{proof}

\begin{proof}[Proof of \cref{lm:approxdecomp}\cref{lm:approxdecomp-all}]
	Let $\ell'\coloneqq \lfloor\nu^5 n\rfloor$ and $k\coloneqq \left\lceil\frac{\ell}{\ell'}\right\rceil$.
	Note that $k\leq \nu^{-5}$.
	For each~$i\in \mathbb{N}$, denote $\varepsilon_i\coloneqq 2^i\varepsilon^{\frac{1}{2^i}}$.	
	Assume inductively that, for some $0\leq m\leq k$, we have constructed edge-disjoint $\cH_1, \dots, \cH_{\min\{m\ell', \ell\}}\subseteq D\cup \Gamma\cup \cF$ such that
	\begin{itemize}[--]
		\item for each $i\in [\min\{m\ell', \ell\}]$,~$\cH_i$ is a spanning configuration of shape~$(L_i,F_i)$ \NEW{satisfying $E(\cH_i)\cap E(\cF)=E(F_i)$};
		\item for each~$i\in [m]$, $D_i\coloneqq D\setminus \bigcup_{j\in [\min\{i\ell', \ell\}]}\cH_j$ is $(\delta-\frac{\min\{i\ell', \ell\}}{n}, \varepsilon_i)$-almost regular; and 
		\item for each~$i\in [m]$, $\Gamma_i\coloneqq \Gamma\setminus \bigcup_{j\in [\min\{i\ell', \ell\}]}\cH_j$ is a $(\gamma, \varepsilon_i)$-almost regular~$(\varepsilon_i,p)$-robust~$(\nu-\varepsilon_i, \tau)$-outexpander.
	\end{itemize}

	If~$m=k$, then, since $k\leq \nu^{-5}$ and~$\varepsilon\ll\nu$,~$\Gamma_m$ is a robust~$(\frac{\nu}{2}, \tau)$-outexpander
	and so is~$D_m\cup \Gamma_m$, as desired.
	Assume~$m<k$. Let $\ell''\coloneqq \min\{\ell-m\ell', \ell'\}$ and $\cF'\coloneqq \bigcup_{i\in [\ell'']}\cF_{m\ell'+i}$. Then, apply \cref{lm:approxdecomp}\cref{lm:approxdecomp-bigchunks}%
		\COMMENT{Note that $\delta-\frac{m\ell'}{n}\geq \delta-\frac{\ell}{n}\geq \eta\gg \gamma$ so the hierarchy is still okay. Moreover, $\delta\geq \frac{\ell}{n}+\eta\geq \frac{m\ell'+\ell''}{n}+\eta$ so $\ell''\leq (\delta-\frac{m\ell'}{n}-\eta)n$, as needed.}
	with $D_m, \Gamma_m, \cF', \delta-\frac{m\ell'}{n}, \nu-\varepsilon_m, \varepsilon_m, \ell''$, $L_{m\ell'+1}, \dots, L_{m\ell'+\ell''}$, and $F_{m\ell'+1}, \dots, F_{m\ell'+\ell''}$ playing the roles of $D, \Gamma, \cF, \delta, \nu, \varepsilon, \ell$, $L_1, \dots, L_\ell$, and $F_1, \dots, F_\ell$ to obtain edge-disjoint $\cH_{m\ell'+1}, \dots, \cH_{m\ell'+\ell''}\subseteq D_m\cup \Gamma_m\cup \cF'$ such that the following hold. For each~$i\in [\ell'']$,~$\cH_{m\ell'+i}$ is a spanning configuration of shape~$(L_{m\ell'+i},F_{m\ell'+i})$. Moreover, $D_{m+1}\coloneqq D_m\setminus \bigcup_{i\in [\ell'']}\cH_{m\ell'+i}$ is $(\delta-\frac{\min\{(m+1)\ell', \ell\}}{n}, \varepsilon_{m+1})$-almost regular%
		\COMMENT{$2\varepsilon_m\leq \varepsilon_{m+1}$}
	and $\Gamma_{m+1}\coloneqq \Gamma_m\setminus \bigcup_{i\in [\ell'']}\cH_{m\ell'+i}$ is a $(\gamma, \varepsilon_{m+1})$-almost regular $(\varepsilon_{m+1},p)$-robust $(\nu-\varepsilon_{m+1}, \tau)$-outexpander%
		\COMMENT{By \cref{fact:robspecialparameters} and since $\sqrt{\varepsilon_m}\leq \varepsilon_{m+1}$ and $\varepsilon_m+\sqrt{\varepsilon_m}\leq \varepsilon_{m+1}$.}, 
	as desired.
\end{proof}

\NEW{As discussed in \cref{sec:sketch}}, the key idea in the proof of \cref{lm:approxdecomp}\cref{lm:approxdecomp-smallchunks} is how to use the robust outexpander~$\Gamma$ efficiently, i.e.\ to find the required number of spanning configurations~$\cH_i$ without using too many edges of~$\Gamma$.
We achieve this by considering a random partition $A_1, \dots, A_a$ of~$V$. To build~$\cH_i$, we find an almost cover of~$V$ in~$D$ with few long paths (which exists since~$D$ is almost regular) and tie them together into a single spanning path using only~$\Gamma[A_j]$ for a suitable~$j\in [a]$. The remainder of~$\cH_i$ is comparatively small and its construction does not affect~$\Gamma$ significantly. \NEW{(See also \cref{fig:sketch}.)}

\begin{proof}[Proof of \cref{lm:approxdecomp}\cref{lm:approxdecomp-smallchunks}]
	Let $a\coloneqq \left\lceil \varepsilon^{-1}\nu^4\right\rceil$. By \cref{lm:partition} (successively applied to~$D$ and~$\Gamma$) and since~$\Gamma$ is an~$(\varepsilon,p)$-robust~$(\nu,\tau)$-outexpander, we can fix a partition $A_1, \dots, A_a$ of~$V$ such that, for each~$i\in [a]$, the following hold%
		\COMMENT{\cref{lm:approxdecomp-Ai-Gammarob} fails with probability at most $\frac{a}{n}$. \cref{lm:approxdecomp-Ai-Gammareg} fails with probability at most $\frac{1}{n}$. \cref{lm:approxdecomp-Ai-Dreg} fails with probability at most $\frac{1}{n}$. So \cref{lm:approxdecomp-Ai-size,lm:approxdecomp-Ai'-Gammarob,lm:approxdecomp-Ai-Gammareg,lm:approxdecomp-Ai-Dreg} hold with probability at least $1-\frac{a+2}{n}>0$.}.
	\begin{enumerate}[label=\upshape(\greek*)]
		\item $|A_i|=\frac{n}{a}\pm 1=\varepsilon(\nu^{-4}\pm 1)n$\label{lm:approxdecomp-Ai-size}%
		\COMMENT{$|A_i|=\frac{n}{a}\pm 1$ so $|A_i|\leq \frac{n}{\varepsilon^{-1}\nu^4}+1\leq \varepsilon(\nu^{-4}+1)n$. Moreover, $|A_i|\geq \frac{n}{\varepsilon^{-1}\nu^4+1}-1= \varepsilon n(\frac{1}{\nu^4+\varepsilon}-\frac{1}{\varepsilon n})=\varepsilon n(\nu^{-4}-\frac{\varepsilon}{\nu^4(\nu^4+\varepsilon)}-\frac{1}{\varepsilon n})\geq \varepsilon(\nu^{-4}-1)n$.}.
		\item $\Gamma[A_i]$ is a robust~$(\nu, \tau)$-outexpander.\label{lm:approxdecomp-Ai-Gammarob}
		\item For each~$v\in V$, $|N_\Gamma^\pm(v)\cap A_i|=(\gamma\pm 2\varepsilon)\frac{n}{a}$.\label{lm:approxdecomp-Ai-Gammareg}
		\item For each~$v\in V$, $|N_D^\pm(v)\cap A_i|=(\delta\pm 2\varepsilon)\frac{n}{a}$.\label{lm:approxdecomp-Ai-Dreg}
	\end{enumerate}
	For each~$i\in [\ell]$, let~$j\in [a]$ be such that~$i\equiv j\mod a$ and define~$A_i'\coloneqq A_j\setminus V(L_i)$. Using \cref{lm:approxdecomp-sizeLi} and \cref{lm:verticesedgesremovalrob}\NEW{\cref{lm:verticesedgesremovalrob-vertices}}, it is easy to check that, for each~$i\in [\ell]$, the following hold.
	\begin{enumerate}[label=\upshape(\greek*$'$)]
		\item $|A_i'|=\varepsilon(\nu^{-4}\pm 2)n$.\label{lm:approxdecomp-Ai'-size}
		\item $\Gamma[A_i']$ and~$\Gamma-A_i'$ are both robust~$(\frac{\nu}{2}, 2\tau)$-outexpanders\label{lm:approxdecomp-Ai'-Gammarob}%
			\COMMENT{$\Gamma-A_i'$ is obtained from $\Gamma$ by removing fewer than $\sqrt{\varepsilon} n$ vertices. $\Gamma[A_i']$ is obtained from $\Gamma[A_i]$ by removing at most $\varepsilon^2 n\leq \varepsilon|A_i|$ vertices. Thus, the result follows from \cref{lm:verticesedgesremovalrob}.}.
		\item $\Gamma[A_i']$ and~$\Gamma-A_i'$ are both~$(\gamma, 3\varepsilon)$-almost regular\label{lm:approxdecomp-Ai'-Gammareg}%
		\COMMENT{\label{com:Gammareg}We first show that $\Gamma[A_i']$ is $(\gamma,3\varepsilon)$-almost regular.
			Let $v\in A_i'$. Then,
			\[|N_\Gamma^\pm(v)\cap A_i'|\geq (\gamma-2\varepsilon)\frac{n}{a}-\varepsilon^2n\geq(\gamma-2\varepsilon)(|A_i'|-1)-\varepsilon^2n \geq (\gamma-3\varepsilon)|A_i'|.\]
			Moreover, 
			\[|N_\Gamma^\pm(v)\cap A_i'|\leq (\gamma+2\varepsilon)\frac{n}{a}\leq(\gamma+2\varepsilon)(|A_i'|+1+\varepsilon^2n )\leq (\gamma+2\varepsilon)|A_i'|+\varepsilon^2 n\leq (\gamma+3\varepsilon)|A_i'|.\]
			Therefore, $\Gamma[A_i']$ is $(\gamma, 3\varepsilon)$-almost regular, as desired.\\
			Similarly, we show that $\Gamma-A_i'$ is $(\gamma, 3\varepsilon)$-almost regular.
			Let $v\in V\setminus A_i'$.
			Then,
			\begin{align*}
			|N_\Gamma^\pm(v)\setminus A_i'|&\geq (\gamma-\varepsilon)n-(\gamma+2\varepsilon)\frac{n}{a}\\
			&\geq (\gamma-\varepsilon)n-(\gamma+2\varepsilon)(|A_i'|+1+\varepsilon^2n)\\
			&\geq\gamma(n-|A_i'|)-\varepsilon(n+2|A_i'|)-2\varepsilon^2 n\\
			&\geq (\gamma-3\varepsilon)(n-|A_i'|).
			\end{align*}
			Moreover,
			\begin{align*}
			    |N_\Gamma^\pm(v)\setminus A_i'|&\leq(\gamma+\varepsilon)n-(\gamma-2\varepsilon)\frac{n}{a}+\varepsilon^2n\\
			    &\leq (\gamma+\varepsilon)n-(\gamma-2\varepsilon)(|A_i'|-1)+\varepsilon^2n\\
    			&\leq \gamma(n-|A_i'|)+\varepsilon(n+2|A_i|)+\varepsilon^2 n+1\\
    			&\leq (\gamma+3\varepsilon)(n-|A_i|).
			\end{align*}
			Therefore, $\Gamma-A_i'$ is $(\gamma, 3\varepsilon)$-almost regular, as desired.}.
		\item $D-A_i'$ is~$(\delta, 3\varepsilon)$-almost regular\label{lm:approxdecomp-Ai'-Dreg}%
		    \COMMENT{Same arguments as above.}. 
		\item For each~$v\in V\setminus A_i'$, $|N_D^\pm(v)\cap A_i'|\geq \frac{\varepsilon\delta n}{2\nu^4}$\label{lm:approxdecomp-Ai'-Gamma'reg}%
		    \COMMENT{$|N_D^\pm(v)\cap A_i'|\geq (\delta-2\varepsilon)\frac{n}{a}-\varepsilon^2 n\geq (\delta-2\varepsilon)\varepsilon(\nu^{-4}-1)n-\varepsilon^2n\geq \frac{\varepsilon \delta n}{2\nu^4}$.}.
	\end{enumerate}
	
	For each~$i\in [\ell]$, fix~$e_i\in E(L_i)\setminus F_i$ (this is possible by \cref{def:layout-unfixededge}).
	Assume inductively that for some $0\leq m\leq \ell$ we have constructed, for each~$i\in [m]$, a set of paths $\cP_i= \{P_e^i\mid e\in E(L_i)\setminus F_i\}$ in~$D\cup \Gamma$ such that $\cP_1, \dots, \cP_m$ are edge-disjoint and the following hold.
	\begin{enumerate}[label=\upshape(\Alph*)]
		\item Let $i\in [m]$. For each $e\in E(L_i)\setminus F_i$,~$P_e^i$ is a path of shape~$e$. Moreover, the paths in~$\cP_i$ are internally vertex-disjoint and $V^0(\cP_i)= V\setminus V(L_i)$. In particular,~$\cP_i\cup F_i$ is a spanning configuration of shape~$(L_i,F_i)$.\label{lm:approxdecomp-IH-shape}
		\item For each~$i\in [m]$ and $e\in E(L_i)\setminus (F_i\cup \{e_i\})$, $P_e^i\subseteq\Gamma-A_i'$ and $e(P_e^i)\leq 8\nu^{-1}$. Moreover, for each~$v\in V$, there exist at most~$\varepsilon^3 n$ indices~$i\in [m]$ such that $v\in V^0 (\cP_i\setminus \{P_{e_i}^i\})$.\label{lm:approxdecomp-IH-shortpaths}
		\item For each~$i\in [m]$, $E(P_{e_i}^i)\cap E(\Gamma)\subseteq E(\Gamma[A_i'])$.\label{lm:approxdecomp-IH-spanningpath}   
	\end{enumerate}
	Denote $D_m\coloneqq D\setminus \bigcup_{i\in [m]}E(\cP_i)$ and $\Gamma_m\coloneqq \Gamma\setminus \bigcup_{i\in [m]}E(\cP_i)$. For each~$i\in[m]$, define $\cH_i\coloneqq \cP_i\cup F_i$. Denote $\cH^m\coloneqq\bigcup_{i\in [m]} \cH_i$.
	Then, note that, for each~$v\in V$, since $d_L(v)\leq \varepsilon^3 n$, there are at most~$\varepsilon^3 n$ indices~$i\in [m]$ such that $v\in V^+(\cP_i\setminus \{P_{e_i}^i\})\cup V^-(\cP_i\setminus \{P_{e_i}^i\})$ and, by \cref{lm:approxdecomp-IH-shortpaths}, there are at most $\varepsilon^3 n$ indices $i\in [m]$ such that $v\in V^0(\cP_i\setminus \{P_{e_i}^i\})$.
	Moreover, by \cref{lm:approxdecomp-IH-spanningpath} and construction of the~$A_i'$, there are, for each~$v\in V$, at most $\left\lceil \frac{\ell}{a} \right\rceil$ indices~$i\in [m]$ such that $v\in V(E(P_{e_i}^i)\cap E(\Gamma))$. Hence, each~$v\in V$ satisfies
	\begin{equation}\label{eq:approxdecompGamma}
		|N_{\cH^m\cap \Gamma}(v)|\leq \varepsilon^3n +2\varepsilon^3 n+2\left\lceil\frac{\ell}{a}\right\rceil\leq3 \varepsilon^3 n+\frac{2\varepsilon^2n}{\varepsilon^{-1}\nu^4}+2\leq 3\varepsilon^3\nu^{-4}n.
	\end{equation}
	
	Assume $m=\ell$. Then, by \cref{lm:approxdecomp-IH-shape}, $\cH_i$ is a spanning configuration of shape $(L_i,F_i)$ for each $i\in [\ell]$. Moreover, \cref{lm:approxdecomp-smallchunks} holds by \cref{eq:approxdecompGamma} and we are done.
	
	Assume $m<\ell$. Using \cref{lm:approxdecomp-Ai'-Gamma'reg,lm:approxdecomp-Ai'-Gammareg,lm:approxdecomp-Ai'-Gammarob,lm:approxdecomp-Ai'-size,lm:approxdecomp-Ai'-Dreg}, \cref{eq:approxdecompGamma}, \cref{lm:approxdecomp-degreev}, \NEW{and \cref{lm:verticesedgesremovalrob}\cref{lm:verticesedgesremovalrob-edges}}, it is easy to check that the following hold.%
		\COMMENT{By the same arguments as above, $\Gamma_m$ is obtained from $\Gamma$ by removing at most $3\varepsilon^3\nu^{-4}n\leq 4\varepsilon^2|A_{m+1}'|$ in/out edges at each vertex.\\
		$D_m$ is obtained from $D$ by removing at most $m+\varepsilon^3n\leq 2\varepsilon^2n$ and at least $m-\varepsilon^2n-3\varepsilon^3\nu^{-4}n$ in/out edges at each vertex.}
	\begin{enumerate}[label=\upshape(\Roman*)]
		\item $\Gamma_m[A_{m+1}']$ and~$\Gamma_m-A_{m+1}'$ are robust~$(\frac{\nu}{4}, 2\tau)$-outexpanders\label{lm:approxdecomp-IS-Gammamrob}%
		    \COMMENT{By \cref{lm:verticesedgesremovalrob}.}.
		\item $\Gamma_m[A_{m+1}']$ and~$\Gamma_m-A_{m+1}'$ are both~$(\gamma, 4\varepsilon)$-almost regular.\label{lm:approxdecomp-IS-Gammamreg}
		\item $D_m-A_{m+1}'$ is~$(\delta-\frac{m}{n}, 4\varepsilon)$-almost regular.\label{lm:approxdecomp-IS-Dmreg}
		\item For each~$v\in V\setminus A_{m+1}'$, $|N_{D_m}^\pm(v)\cap A_{m+1}'|\geq \frac{\varepsilon\delta n}{3\nu^4}$. \label{lm:approxdecomp-IS-Gammamdegreeacross}
	\end{enumerate}
	
	We first construct~$P_e^{m+1}$ for each $e\in E(L_{m+1})\setminus (F_{m+1}\cup \{e_{m+1}\})$ in the following way. Let~$S$ be the set of vertices~$v\in V$ for which there exist~$\lfloor\varepsilon^3 n\rfloor$ indices~$i\in [m]$ such that $v\in V^0(\cP_i\setminus \{P_{e_i}^i\})$. Observe that, by \cref{lm:approxdecomp-sizeLi} and \cref{lm:approxdecomp-IH-shortpaths}, \NEW{$|S|\leq \frac{8\nu^{-1}\cdot \ell \cdot \varepsilon^4 n}{\lfloor\varepsilon^3 n\rfloor}\leq \varepsilon |V\setminus A_{m+1}'|$}\OLD{$|S|\leq \frac{8\nu^{-1}\cdot \ell \cdot \varepsilon^4 n}{\lceil\varepsilon^3 n\rceil}\leq \varepsilon |V\setminus A_{m+1}'|$}%
	    \COMMENT{$|S|\leq \frac{8\nu^{-1}\cdot \ell\cdot \varepsilon^4n}{\lfloor\varepsilon^3 n\rfloor}\leq 10\nu ^{-1}\varepsilon^3 n\leq \varepsilon|V\setminus A_{m+1}'|$}.
	Denote $E(L_{m+1})\setminus (F_{m+1}\cup \{e_{m+1}\})\eqqcolon \{x_1x_1', \dots, x_kx_k'\}$. Apply \cref{cor:robshortpaths} with $\Gamma_m-A_{m+1}', \frac{\nu}{4}, 2\tau, \gamma-4\varepsilon$, and $S\cup V(L_{m+1})$ playing the roles of $D, \nu, \tau, \delta$, and $S$ to obtain internally vertex-disjoint paths $P_{x_1x_1'}^{m+1}, \dots, P_{x_kx_k'}^{m+1}\subseteq \Gamma_m-A_{m+1}'$ such that, for each~$i\in [k]$,~$P_{x_ix_i'}^{m+1}$ is an~$(x_i,x_i')$-path of length at most~$8\nu^{-1}$ with $V^0(P_{x_ix_i'}^{m+1})\subseteq V\setminus (A_{m+1}'\cup S\cup V(L_{m+1}))$.
	Let $\cP_{m+1}'\coloneqq \{P_{x_ix_i'}^{m+1}\mid i\in [k]\}$.
	
	Let~$z\notin V$ be a new vertex. Let~$H$ be the digraph on vertex set $V(H)\coloneqq V\setminus (V(L_{m+1})\cup V(\cP_{m+1}'))\cup \{z\}$ defined as follows. Denote $v^+v^-\coloneqq e_{m+1}$ and recall that, by construction, $v^\pm\notin A_{m+1}'$. Then, let $N_H^\pm (z)\coloneqq N_{D_m}^\pm (v^\pm)\cap V(H)$, $H[A_{m+1}']\coloneqq \Gamma_m[A_{m+1}']$, and, for each $v\in V(H)\setminus (A_{m+1}'\cup \{z\})$, \NEW{$N_{H-\{z\}}^\pm(v)\coloneqq N_{D_m}^\pm(v)\cap V(H)$}\OLD{$N_{H-\{z\}}^\pm(v)= N_{D_m}^\pm(v)\cap V(H)$}.
	Note that, by \cref{lm:approxdecomp-IS-Dmreg,lm:approxdecomp-IS-Gammamdegreeacross,lm:approxdecomp-IS-Gammamreg,lm:approxdecomp-IS-Gammamrob}, the following hold%
		\COMMENT{$H[A_{m+1}']=\Gamma_m[A_{m+1}']$ and $H-A_{m+1}'$ is obtained from $D_m-A_{m+1}'$ by adding a vertex and deleting $|V(L_{m+1})\cup V(\cP_{m+1}')|\leq \varepsilon^2n+\varepsilon^4 n\cdot 8\nu^{-1}\leq 2\varepsilon ^2 n$ vertices.}.
	\begin{enumerate}[label=\upshape(\Roman*$'$)]
		\item $H[A_{m+1}']$ is a robust~$(\frac{\nu}{4}, 2\tau)$-outexpander.\label{lm:approxdecomp-H-rob}
		\item $H[A_{m+1}']$ is~$(\gamma, 4\varepsilon)$-almost regular.\label{lm:approxdecomp-H-Aireg}
		\item $H-A_{m+1}'$ is~$(\delta-\frac{m}{n}, 5\varepsilon)$-almost regular.\label{lm:approxdecomp-H-reg}
		\item For each~$v\in V(H)\setminus A_{m+1}'$, $|N_H^\pm(v)\cap A_{m+1}'|\geq \frac{\varepsilon\delta n}{3\nu^4}$.\label{lm:approxdecomp-H-degreeacross}
	\end{enumerate}
	Indeed, to check \cref{lm:approxdecomp-H-reg}, note that, by \cref{lm:approxdecomp-sizeLi},~$H-A_{m+1}'$ is obtained from~$D_m-A_{m+1}'$ by adding~$z$ and deleting \NEW{$|V(L_{m+1})\cup V(\cP_{m+1}')|\leq \varepsilon^2n+\varepsilon^4 n\cdot 8\nu^{-1}\leq 2\varepsilon^2 n$}\OLD{$|V(L_{m+1})\cup V(\cP_{m+1})|\leq \varepsilon^2n+\varepsilon^4 n\cdot 8\nu^{-1}\leq 2\varepsilon^2 n$} vertices.
	
	Our aim is to find a Hamilton cycle of~$H$ which contains few edges of~$\Gamma[A_{m+1}']$. 
	First, we cover~$V(H) \setminus A_{m+1}'$ with a small number of paths as follows. Let $k'\coloneqq \left\lfloor \frac{|V(H)\setminus A_{m+1}'|}{\varepsilon n}\right\rfloor$. Apply \cref{lm:partition} with $H-A_{m+1}', |V(H)\setminus A_{m+1}'|, \delta-\frac{m}{n}$, and $5\varepsilon$ playing the roles of $D, n, \delta$, and $\varepsilon$ to obtain a partition $V_1, \dots, V_{k'}$ of $V(H)\setminus A_{m+1}'$ such that, for each~$i\in [k']$, $|V_i|=(1\pm 2\varepsilon)\varepsilon n$%
		\COMMENT{$|V_i|=\frac{|V(H)\setminus A_{m+1}'|}{k'}\pm 1=\frac{|V(H)\setminus A_{m+1}'|}{\frac{|V(H)\setminus A_{m+1}'|}{\varepsilon n}\mp 1}\pm 1=\varepsilon n(1\pm \frac{\varepsilon n}{|V(H)\setminus A_{m+1}'|\mp \varepsilon n})\pm 1$.}
	and, for each~$i\in [k']$ and~$v\in V_i$, $|N_H^-(v)\cap V_{i-1}|= (\delta-\frac{m}{n}\pm 10\varepsilon)\varepsilon n$ if~$i>1$ and $|N_H^+(v)\cap V_{i+1}|= (\delta-\frac{m}{n}\pm 10\varepsilon)\varepsilon n$ if~$i<k'$.
	
	Then, for each~$i\in [k'-1]$, apply \cref{cor:Hallreg} with $H[V_i, V_{i+1}]$, $V_i$, $V_{i+1}$, $\varepsilon n$, $\delta-\frac{m}{n}$, and $10\varepsilon$ playing the roles of $G, A, B, n, \delta$, and $\varepsilon$ to obtain a matching~$M_i$ of~$H[V_i, V_{i+1}]$ of size at least~$\left(1-\frac{31\varepsilon}{\delta}\right)\varepsilon n$. For each~$i\in [k'-1]$, denote by~$\overrightarrow{M}_i$ the directed matching obtained from~$M_i$ by directing all edges from~$V_i$ to~$V_{i+1}$. Note that, by construction,~$\overrightarrow{M}_i\subseteq H$.
	Define~$F\subseteq H$ by letting $V(F)\coloneqq V(H)\setminus A_{m+1}'$ and $E(F)\coloneqq \bigcup_{i\in [k'-1]} \overrightarrow{M}_i$. Observe that~$F$ is a linear forest which spans~$V(H)\setminus A_{m+1}'$ and has~$f\leq \frac{33\varepsilon n}{\delta}$ components.
	Indeed, one can count the number of paths in~$F$ by counting the number of ending points as follows. (An isolated vertex is considered as the ending point of a trivial path of length~$0$.) Note that, for each~$i\in [k'-1]$,~$v\in V_i$ is the ending point of path in~$F$ if and only if~$v\notin V(M_i)$, while every~$v\in V_{k'}$ is the ending point of a path in~$F$.
	Moreover, for each~$i\in [k'-1]$, we have $|V_i\setminus V(M_i)|\leq |V_i|-|M_i|\leq \varepsilon n+2\varepsilon^2 n-\left(1-\frac{31\varepsilon}{\delta}\right)\varepsilon n\leq \frac{32\varepsilon^2 n}{\delta}$. Thus, since $k'-1\leq \varepsilon^{-1}-1$, we have $f\leq \frac{32\varepsilon^2 n}{\delta}(\varepsilon^{-1}-1)+|V_k|\leq \frac{33\varepsilon n}{\delta}$, as desired.
	
	Denote the components of~$F$ by $P_1, \dots, P_f$.
	We now join $P_1, \dots, P_f$ into a Hamilton cycle as follows.
	Note that, by \cref{lm:approxdecomp-Ai'-size}, $f\leq \left(\frac{\nu}{4}\right)^3 |A_{m+1}'|$%
	    \COMMENT{Indeed, $|A_{m+1}'|\geq \varepsilon(\nu^{-4}-2)n\geq \frac{\varepsilon n}{2\nu^4}$. Therefore, $\left(\frac{\nu}{4}\right)^3|A_{m+1}'|\geq \frac{\varepsilon n}{2\cdot 4^3\cdot \nu}\geq f$.}.
	For each~$i\in[f]$, denote by~$v_i^+$ and~$v_i^-$ the starting and ending points of~$P_i$. By \cref{lm:approxdecomp-H-degreeacross}, for each~$i\in[f]$, we have $|N_H^\mp(v_i^\pm)\cap A_{m+1}'|\geq 2f$.
	Apply \cref{cor:robpaths}\cref{cor:robcycles} with $H, A_{m+1}', \emptyset, f, \frac{\nu}{4}, 2\tau$, and $\gamma-\nu$ playing the roles of $D, V', S, k, \nu, \tau$, and $\delta$ to obtain a Hamilton cycle~$C$ of $H$ such that~$F\subseteq C$.
	Denote by~$u^\pm$ the (unique) vertices such that $u^\pm \in N_C^\pm(z)$, respectively. Let $P_{e_{m+1}}^{m+1}\coloneqq (C-\{z\})\cup \{v^+u^+,u^-v^-\}$. By construction, $P_{e_{m+1}}^{m+1}$ is a path of shape~$e_{m+1}$ such that $P_{e_{m+1}}^{m+1}\subseteq (D_m\cup \Gamma_m)-V(\cP_{m+1}')$ and $V^0(P_{e_{m+1}}^{m+1})=V\setminus (V(\cP_{m+1}')\cup V(L_{m+1}))$. Moreover, $P_{e_{m+1}}^{m+1}[A_{m+1}']\subseteq \Gamma_m$ and $P_{e_{m+1}}^{m+1}\setminus P_{e_{m+1}}^{m+1}[A_{m+1}']\subseteq D_m$. Let $\cP_{m+1}\coloneqq \cP_{m+1}'\cup \{P_{e_{m+1}}^{m+1}\}$. Thus, \cref{lm:approxdecomp-IH-shape,lm:approxdecomp-IH-shortpaths,lm:approxdecomp-IH-spanningpath} hold. This completes the induction step.
\end{proof}

%% file: Figures/Figure_Layouts.tex
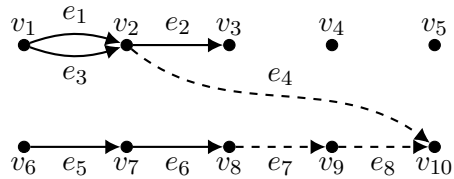
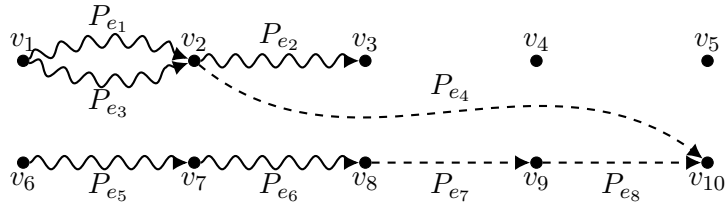
\begin{figure}[htb]
	\centering
	\begin{subfigure}{\textwidth}
	\centering
	\begin{tikzpicture}[scale=0.9]
		%%vertices%%
		\draw node[circle, draw=black,fill=black, inner sep=1.5pt](v1) at (-3,0) {};
		\node[above] at (-3,0) {$v_1$};
		\node[circle, draw=black,fill=black, inner sep=1.5pt](v2) at (-1.5,0) {};
		\node[above] at (-1.5,0) {$v_2$};
		\node[circle, draw=black,fill=black, inner sep=1.5pt](v3) at (0,0) {};
		\node[above] at (0,0) {$v_3$};
		\draw node[circle, draw=black,fill=black, inner sep=1.5pt](v9) at (1.5,0) {};
		\node[above] at (1.5,0) {$v_4$};
		\node[circle, draw=black,fill=black, inner sep=1.5pt](v10) at (3,0) {};
		\node[above] at (3,0) {$v_5$};
		\node[circle, draw=black,fill=black, inner sep=1.5pt](v4) at (-3,-1.5) {};
		\node[below] at (-3,-1.5) {$v_6$};
		\node[circle, draw=black,fill=black, inner sep=1.5pt](v5) at (-1.5,-1.5) {};
		\node[below] at (-1.5,-1.5) {$v_7$};
		\node[circle, draw=black,fill=black, inner sep=1.5pt](v6) at (0,-1.5) {};
		\node[below] at (0,-1.5) {$v_8$};
		\node[circle, draw=black,fill=black, inner sep=1.5pt](v7) at (1.5,-1.5) {};
		\node[below] at (1.5,-1.5) {$v_9$};
		\node[circle, draw=black,fill=black, inner sep=1.5pt](v8) at (3,-1.5) {};
		\node[below] at (3,-1.5) {$v_{10}$};
		%%edges%%
		\draw[->,thick] (v1) to [out=20,in=160]  node[above] {$e_1$} (v2);
		\draw[->,thick] (v2)-- node[above] {$e_2$} (v3);
		\draw[->,thick] (v1) to [out=-20,in=-160] node[below] {$e_3$} (v2);
		\draw[->,thick,dashed] (v2) to [out=-40,in=140] node[above] {$e_4$}(v8);
		\draw[->,thick] (v4)-- node[below] {$e_5$} (v5);
		\draw[->,thick] (v5)-- node[below] {$e_6$}(v6);
		\draw[->,thick,dashed] (v6)-- node[below] {$e_7$} (v7);
		\draw[->,thick,dashed] (v7) -- node[below] {$e_8$} (v8);
	\end{tikzpicture}
	\caption{\NEW{A layout $(L,F)=(\{v_1e_1v_2e_2v_3,v_1e_3v_2e_4v_{10}, v_4, v_5,v_6e_5v_7e_6v_8e_7v_9e_8v_{10}\},\{e_4,e_7,e_8\})$ on $V$.}\label{fig:layout-L}}
	\end{subfigure}
	\begin{subfigure}{\textwidth}
		\centering
		\begin{tikzpicture}[scale=0.9]
			%%vertices%%
			\draw node[circle, draw=black,fill=black, inner sep=1.5pt](v1) at (-5,0) {};
			\node[above] at (-5,0) {$v_1$};
			\node[circle, draw=black,fill=black, inner sep=1.5pt](v2) at (-2.5,0) {};
			\node[above] at (-2.5,0) {$v_2$};
			\node[circle, draw=black,fill=black, inner sep=1.5pt](v3) at (0,0) {};
			\node[above] at (0,0) {$v_3$};
			\draw node[circle, draw=black,fill=black, inner sep=1.5pt](v9) at (2.5,0) {};
			\node[above] at (2.5,0) {$v_4$};
			\node[circle, draw=black,fill=black, inner sep=1.5pt](v10) at (5,0) {};
			\node[above] at (5,0) {$v_5$};
			\node[circle, draw=black,fill=black, inner sep=1.5pt](v4) at (-5,-1.5) {};
			\node[below] at (-5,-1.5) {$v_6$};
			\node[circle, draw=black,fill=black, inner sep=1.5pt](v5) at (-2.5,-1.5) {};
			\node[below] at (-2.5,-1.5) {$v_7$};
			\node[circle, draw=black,fill=black, inner sep=1.5pt](v6) at (0,-1.5) {};
			\node[below] at (0,-1.5) {$v_8$};
			\node[circle, draw=black,fill=black, inner sep=1.5pt](v7) at (2.5,-1.5) {};
			\node[below] at (2.5,-1.5) {$v_9$};
			\node[circle, draw=black,fill=black, inner sep=1.5pt](v8) at (5,-1.5) {};
			\node[below] at (5,-1.5) {$v_{10}$};
			%%edges%%
			\draw[->,thick,decorate,decoration={snake,post length=3pt}] (v1) to [out=20,in=160]  node[above] {$P_{e_1}$} (v2);
			\draw[->,thick,decorate,decoration={snake,post length=2pt}] (v2)-- node[above] {$P_{e_2}$} (v3);
			\draw[->,thick,decorate,decoration={snake,post length=3pt}] (v1) to [out=-20,in=-160] node[below] {$P_{e_3}$} (v2);
			\draw[->,thick,dashed] (v2) to [out=-40,in=140] node[above] {$P_{e_4}$}(v8);
			\draw[->,thick,decorate,decoration={snake,post length=2pt}] (v4)-- node[below] {$P_{e_5}$} (v5);
			\draw[->,thick,decorate,decoration={snake,post length=2pt}] (v5)-- node[below] {$P_{e_6}$}(v6);
			\draw[->,thick,dashed] (v6)-- node[below] {$P_{e_7}$} (v7);
			\draw[->,thick,dashed] (v7) -- node[below] {$P_{e_8}$} (v8);
		\end{tikzpicture}
		\caption{\NEW{A spanning configuration $\cH$ of shape $(L,F)$ on $V$.
		The wavy edges represent internally vertex-disjoint paths on $V$ which altogether cover all the vertices in $V\setminus V(L)$.}\label{fig:layout-H}}
	\end{subfigure}
\caption{\NEW{A layout $(L,F)$ and a spanning configuration of shape $(L,F)$. Dashed edges represent fixed edges (i.e.\ the edges of $F$).}\label{fig:layout}}	
	\end{figure}

%% file: Good_Partial_Decompositions.tex
	\onlyinsubfile{
		\setcounter{section}{7}
\section{Good partial path decompositions and absorbing edges}}

\NEW{\Cref{lm:approxdecomp} only covers most of the edges. Moreover, we will see that we also need an extra cleaning step before being able to apply \cref{lm:approxdecomp}. This means that our path decomposition will be constructed in several stages.}
	
\NEW{Suppose that we have already constructed an intermediate set of paths $\cP$ and that we want to extend $\cP$ to a path decomposition of $T$. Then, $\texc(T\setminus \cP)$ must not be too large (for otherwise we will not have any hope of extending $\cP$ to a path decomposition of the desired size $\texc(T)$).	
This is encapsulated in the concept of a good partial path decomposition, which is defined and discussed in \cref{sec:goodpartialdecomp}.}

\NEW{Moreover, we will need to make sure that, in the last stage, the remaining digraph $D$ has a nice structure (for otherwise we may not know how to decompose $D$). In \cref{thm:niceD}, we saw an example of a digraph that we can decompose efficiently. Unfortunately, it will not always be possible to get a leftover of that form, so, in \cref{sec:A}, we will generalise \cref{thm:niceD} using the concept of absorbing edges.}

\subsection{Partial path decompositions}\label{sec:goodpartialdecomp}

Denote $U^\pm(D)\coloneqq \{v\in V(D)\mid \exc_D^\pm (v)>0\}$ and
$U^0(D)\coloneqq \{v\in V(D)\mid \exc_D (v)=0\}$.

\begin{prop}\label{prop:sizeU}
	Any oriented graph~$D$ satisfies $|U^0(D)|\geq \texc(D)-\exc(D)$.%
	\COMMENT{Note that this is false for general digraphs. E.g.\ let $D$ be obtained from a complete digraph on $n$ vertices by deleting a perfect matching. Then, $U^0(D)=\emptyset$, $\exc(D)=n$, $\texc(D)=\Delta^0(D)=2n-3$.}
\end{prop}

\begin{proof}
	Assume for a contradiction that there exists an oriented graph~$D$ such that
	$|U^0(D)|<\texc(D)-\exc(D)$. Then, note that~$\texc(D)=\Delta^0(D)$
	and let~$v\in V$ be such that $d_D^{\max}(v)=\Delta^0(D)$.
	Assume without loss of generality that~$v\in U^+(D)$.
	Then, $d_D^+(v)=\texc(D)>\exc(D)$. 
	By \cref{fact:partialdecompS}, $\exc(D)\geq \exc_D^+(v)+|U^+(D)|-1$ and so $|U^+(D)|\leq \exc(D)-\exc_D^+(v)+1$. Moreover, by assumption, we have \NEW{$|U^0(D)|< d_D^+(v)-\exc(D)$}\OLD{$|U^0(D)|\leq d_D^+(v)-\exc(D)$}. Therefore, by Facts \ref{fact:exc}\cref{fact:exc-dmin} and \ref{fact:exc}\cref{fact:exc-dmax}, we have \NEW{$\exc(D)\geq |U^-(D)|=n-|U^+(D)|-|U^0(D)|> n-(\exc(D)-\exc_D^+(v)+1)-(d_D^+(v)-\exc(D))=n-1-d_D^-(v)\geq d_D^+(v)$}\OLD{$\exc(D)\geq |U^-(D)|=n-|U^+(D)|-|U^0(D)|\geq n-(\exc(D)-\exc_D^+(v)+1)-(d_D^+(v)-\exc(D))=n-1-d_D^-(v)\geq d_D^+(v)$}, a contradiction.
\end{proof}

\NEW{Let $D$ be an oriented graph. Recall that in a path decomposition $\cP$ of $D$ of size $\texc(D)$, each $v\in V(D)$ will be the starting point of at least $\exc_D^+(v)$ paths in $\cP$ and the ending point of at least $\exc_D^-(v)$ paths in $\cP$. When $\texc(D)>\exc(D)$, there are $\texc(D)-\exc(D)$ starting (and ending) points unaccounted for. By \cref{prop:sizeU}, we can choose these endpoints (to be distinct vertices) in $U^0(D)$. Thus, our path decomposition $\cP$ will also maximise the number of distinct vertices that are an endpoint of some path in $\cP$.}
This motivates the following definition.

\begin{definition}[Partial path decomposition]
	Let~$D$ be a digraph.
	A set~$\cP$ of edge-disjoint paths of~$D$ is called a \emph{partial path
		decomposition of~$D$} if the following hold.
	\begin{enumerate}[label=\upshape(P\arabic*)]
		\item Any vertex~$v\in V(D)\setminus
		U^0(D)$ is the starting point of at most~$\exc_D^+(v)$ paths in~$\cP$ and the
		ending point of at most~$\exc_D^-(v)$ paths in~$\cP$.\label{def:partialpathdecomp-exc}
		\item Any vertex~$v\in
		U^0(D)$ is the starting point of at most one path in~$\cP$ and the ending point
		of at most one path in~$\cP$.\label{def:partialpathdecomp-U0}
		\item There are at most~$\texc(D)-\exc(D)$
		vertices~$v\in U^0(D)$ such that~$v$ is an endpoint%
		\COMMENT{i.e. the starting or ending point}
		of a path in~$\cP$, that is, $|U^0(D)\cap (V^+(\cP)\cup V^-(\cP))|\leq
		\texc(D)-\exc(D)$.\label{def:partialpathdecomp-texc}
	\end{enumerate}
\end{definition}

\NEW{By \cref{def:partialpathdecomp-texc}, we will need to construct sets of edge-disjoint paths which do not contain too many paths which start and/or end at vertices of zero excess. It will turn out to be convenient to fix in advance which zero-excess vertices will be used as endpoints. This motivates the following definition.}
\NEW{Let $D$ be a digraph and suppose that $U^*\subseteq U^0(D)$ satisfies $|U^*|\leq \texc(D)-\exc(D)$. We say that $\cP$ is a \emph{$U^*$-partial path decomposition of $D$} if $\cP$ is a partial path decomposition where $(V^+(\cP)\cup V^-(\cP))\cap U^0(D)\subseteq U^*$, i.e.\ no path in $\cP$ has an endpoint in $U^0(D)\setminus U^*$.}

\NEW{Let $D$ be a digraph.} Recall \NEW{that in \cref{thm:main} we defined}\OLD{from \cref{thm:main} that}
\begin{equation}\label{eq:N}
	N^\pm(D)= |U^\pm(D)|+\texc(D)-\exc(D).
\end{equation} 
Note that \NEW{\cref{def:partialpathdecomp-exc,def:partialpathdecomp-texc} imply that} if~$\cP$ is a partial path decomposition of~$D$, then there are at most~$N^+(D)$ distinct vertices which are the starting point of a path in~$\cP$ 
\NEW{(\cref{def:partialpathdecomp-exc} implies that the vertices in $U^-(D)$ cannot be used as starting points and \cref{def:partialpathdecomp-texc} implies that at most $\texc(D)-\exc(D)$ vertices in $U^0(D)$ may be used as starting points). Similarly, there are} at most~$N^-(D)$ distinct vertices which are the ending point of a path in~$\cP$.

\NEW{\begin{prop}\label{prop:excP}
		Let $D$ be a digraph and $\cP$ be a partial path decomposition of $D$. 
		Then, $\exc(D\setminus \cP)=\exc(D)-|\cP|+|U^0(D)\cap (V^+(\cP)\cup V^-(\cP))|\leq \texc(D)-|\cP|$.	
\end{prop}}

\begin{proof}
	\NEW{By \cref{def:partialpathdecomp-texc}, $|U^0(D)\cap (V^+(\cP)\cup V^-(\cP))|\leq \texc(D)-\exc(D)$ and so it is enough to show that $\exc(D\setminus \cP)=\exc(D)-|\cP|+|U^0(D)\cap (V^+(\cP)\cup V^-(\cP))|$. For each $v\in V(D)$, denote by $n_\cP^+(v)$ and $n_\cP^-(v)$ the number of paths on $\cP$ which start and end at $v$, respectively. Let $v\in V(D)$ and note that
	\begin{align}\label{eq:excP}
		\exc_{D\setminus \cP}(v)&=d_{D\setminus \cP}^+(v)-d_{D\setminus \cP}^-(v)=(d_D^+(v)-d_D^-(v))-(d_\cP^+(v)-d_\cP^-(v))\nonumber\\
		&=\exc_D(v)-n_\cP^+(v)+n_\cP^-(v).
	\end{align}
	Let $n^+$ be the number of paths in $\cP$ which start in $U^+(D)$. Since $\cP$ is a partial path decomposition, we have	
	\begin{align*}
		\exc(D\setminus \cP)&\stackrel{\text{\eqmakebox[excP][c]{\cref{eq:exc}}}}{=}\sum_{v\in V(D)}\exc_{D\setminus \cP}^+(v)
		\stackrel{\text{\cref{def:partialpathdecomp-exc},\cref{def:partialpathdecomp-U0}}}{=}\sum_{v\in U^+(D)}\exc_{D\setminus \cP}^+(v)+|U^0(D)\cap (V^-(\cP)\setminus V^+(\cP))|\\
		&\stackrel{\text{\eqmakebox[excP][c]{\text{\cref{def:partialpathdecomp-exc}}}}}{=}\left(\sum_{v\in U^+(D)}\exc_D^+(v)-n^+\right)+ |U^0(D)\cap (V^-(\cP)\cup V^+(\cP))|- |U^0(D)\cap V^+(\cP)|\\
		&\stackrel{\text{\eqmakebox[excP][c]{\text{\cref{def:partialpathdecomp-U0}}}}}{=}\exc(D)-|\cP|+ |U^0(D)\cap (V^-(\cP)\cup V^+(\cP))|,
	\end{align*}
	as desired.}
\end{proof}

\NEW{Let $D$ be a digraph and $\cP$ be a partial path decomposition of $D$. The next \lcnamecref{prop:excpartialdecomp} expands on \cref{prop:excP} to give further bounds on $\exc(D\setminus \cP)$ and $\texc(D\setminus \cP)$.}

\begin{prop}\label{prop:excpartialdecomp}
	Let~$D$ be a digraph and~$\cP$ be a partial path decomposition of~$D$. Then, \NEW{the following hold.}
	\begin{enumerate}
		\item \NEW{If $\Delta^0(D\setminus \cP)\leq \texc(D)-|\cP|$, then $\texc(D\setminus \cP)=\texc(D)-|\cP|$.}\OLD{$\exc(D)-|\cP|\leq \exc(D\setminus \cP)\leq \texc(D)-|\cP|\leq \texc(D\setminus \cP)\leq \texc(D)$.}\label{prop:excpartialdecomp-general}
		\item If $\texc(D)=\exc(D)$, then $\exc(D\setminus \cP)=\exc(D)-|\cP|$.\label{prop:excpartialdecomp-exc}
	\end{enumerate} 
\end{prop}

\begin{proof}
	\NEW{Note that it is enough to show that the following inequalities hold.
	\begin{equation}\label{eq:excpartialdecomp}
		\exc(D)-|\cP|\leq \exc(D\setminus \cP)\leq \texc(D)-|\cP|\leq \texc(D\setminus \cP).
	\end{equation}
	Indeed, if $\texc(D) = \exc(D)$, then the first two inequalities of \cref{eq:excpartialdecomp} are equalities implying \cref{prop:excpartialdecomp-exc}.
	If $\Delta^0(D \setminus \cP) \le \texc(D) - |\cP|$, then together with the last two inequalities of \cref{eq:excpartialdecomp}, we deduce that
	\begin{align*}
		\exc(D \setminus \cP), \Delta^0(D \setminus \cP)  \leq \texc (D) - |\cP| \leq \texc(D \setminus \cP) =  \max \{ \Delta^0(D \setminus \cP) ,  \exc(D \setminus \cP)\}.
	\end{align*}
	Thus, we must have equalities, which implies \cref{prop:excpartialdecomp-general}.}

	\NEW{First, consider the case where $\texc(D)=\exc(D)$. By \cref{def:partialpathdecomp-texc}, no path in $\cP$ has an endpoint in $U^0(D)$. Thus, \cref{prop:excP} implies that
		\begin{align*}
			\texc(D)-|\cP|=\exc(D)-|\cP|=\exc(D\setminus \cP)\leq \texc(D\setminus \cP)
		\end{align*}
		and so \cref{eq:excpartialdecomp} holds.}%
	\OLD{If $\texc(D)=\exc(D)$, then by \cref{def:partialpathdecomp-exc,def:partialpathdecomp-U0}, $\exc(D\setminus
	\cP)=\exc(D)-|\cP|$ and so, since $\Delta^0(D\setminus\cP)\leq \Delta^0(D)$, \cref{prop:excpartialdecomp-general,prop:excpartialdecomp-exc} hold.}
	We may therefore assume that $\texc(D)=\Delta^0(D)\neq \exc(D)$.\OLD{We show \cref{prop:excpartialdecomp-general} is satisfied (\cref{prop:excpartialdecomp-exc} holds vacuously).}
	Clearly, 
	\begin{equation}\label{eq:excpartialdecomp-4}
		\texc(D) - |\cP| = \Delta^0(D) - |\cP| \le \Delta^0(D \setminus \cP) \le \texc (D \setminus \cP).
	\end{equation}
	\NEW{By \cref{prop:excP}, we have
	\begin{align}\label{eq:excpartialdecomp-5}
		\exc(D)-|\cP|\leq \exc(D\setminus \cP)\leq \texc(D)-|\cP|.
	\end{align}
	Therefore, \cref{eq:excpartialdecomp} follows from \cref{eq:excpartialdecomp-4,eq:excpartialdecomp-5}.}
	\OLD{We now show that $\exc(D)-|\cP|\leq \exc(D\setminus \cP)\leq \texc(D)-|\cP|$.
		Let~$k$ be the number of paths in~$\cP$ which start in~$U^0(D)$ and let~$S$ be the set of vertices~$v\in U^0(D)$ such that no path in~$\cP$ starts at~$v$ but~$v$ is the ending point of \NEW{a} path in~$\cP$. Note that by definition of a partial path decomposition, $k+|S|\leq \texc(D)-\exc(D)$. Moreover, observe that, for each~$v\in S$, $\exc_{D\setminus \cP}^+(v)=1$ and, for each $v\in (U^0(D)\setminus S)\cup U^-(D)$, $\exc_{D\setminus \cP}^+(v)=0$.
		For each~$v\in U^+(D)$, let~$n_{\cP}^+(v)$ be the number of paths in~$\cP$ which start at~$v$. Then, by the definition of a partial path decomposition, for each~$v\in U^+(D)$, $\exc_{D\setminus \cP}^+(v)=\exc_D^+(v)-n_{\cP}^+(v)$.
		Thus,
		$\exc(D\setminus\cP)=\sum_{v\in V(D)}\exc_{D\setminus \cP}^+(v)=\exc(D)-(|\cP|-k)+|S|$ and, therefore, $\exc(D)-|\cP|\leq \exc(D\setminus\cP)\leq \exc(D)-|\cP|+(\texc(D)-\exc(D))=\texc(D)-|\cP|$, as desired.
		Finally, note that, since $\exc(D\setminus \cP)\leq \texc(D)-|\cP|\leq \texc(D)=\Delta^0(D)$ and $\Delta^0(D\setminus \cP)\leq \Delta^0(D)$, we have $\texc(D\setminus \cP)\leq \texc(D)$.}
\end{proof}

Let~$D$ be a digraph\OLD{on~$n$ vertices}. We say that a partial path
decomposition~$\cP$ of~$D$ is \emph{good} if $\texc(D\setminus \cP)= \texc(D)-|\cP|$. We say that a path decomposition~$\cP$ of~$D$ is \emph{perfect} if~$|\cP|=\texc(D)$.

\begin{fact}\label{fact:combinefull}
	Let \NEW{$k\in \mathbb{N}$}\OLD{$k\in \mathbb{N}\setminus \{0\}$} and~$D$ be a digraph. Denote~$D_0\coloneqq D$.
	Suppose that, for each~$i\in [k-1]$,~$\cP_i$ is a good partial path decomposition of~$D_{i-1}$ and~$D_i\coloneqq D_{i-1}\setminus \cP_i$. Suppose that~$\cP_k$ is a perfect path decomposition of~$D_{k-1}$. Then, $\cP\coloneqq \bigcup_{i\in [k]}\cP_i$ is a perfect path decomposition of~$D$\OLD{of size $|\cP|=\sum_{i\in [k]}|\cP_i|$}.
\end{fact}

\NEW{Let $D$ be an oriented graph on $n$ vertices. The next \lcnamecref{prop:Delta0} shows that if there is a vertex $v\in V(D)$ with $d_D^+(v)\geq \exc(D)-\varepsilon n$, then $\texc(D)\leq (1+\varepsilon)n$ (\cref{prop:Delta0}\cref{prop:Delta0-excD})}
\NEW{and most of the positive excess of $D$ is concentrated at $v$ (\cref{prop:Delta0}\cref{prop:Delta0-excv}).}
\NEW{\Cref{prop:Delta0}\cref{prop:Delta0-good} gives a sufficient condition for a small partial path decomposition to be good.}

\begin{prop}\label{prop:Delta0}
	Let $0<\frac{1}{n}\ll \eta \ll 1$. Let~$D$ be an oriented graph on~$n$ vertices satisfying $\exc(D)\geq (1-21\eta)n$.%
	\COMMENT{Could prove an analogous statement for general digraphs by requiring $\exc(D)$ to be a bit larger than $n$.}
	\NEW{Let $V^\pm\coloneqq \{v\in V(D)\mid d_D^\pm(v)\geq \texc(D)-22\eta n\}$.}
	Then, the following hold.
	\begin{enumerate}
		\item \NEW{If $V^\diamond\neq \emptyset$ for some $\diamond \in \{+,-\}$,}\OLD{For each $\diamond\in\{+,-\}$, if there exists~$v\in V$ such that $d_D^\diamond(v)\geq \texc(D)-22\eta n$.} then the following hold. \label{prop:Delta0-unique}
		\begin{enumerate}[label=\rm(\roman*), ref=\upshape(\alph{enumi}.\roman{*})]
			\item $\texc(D)\leq  (1+22\eta)n \NEW{\leq \exc(D)+43\eta n}$. \label{prop:Delta0-excD}
			\item $\exc_D^\diamond(v)\geq
			(1-86\eta)n \NEW{\geq \exc(D)-108\eta n}$ \NEW{for all $v\in V^\diamond$}.\label{prop:Delta0-excv}
		\end{enumerate}
		\item \NEW{Let $\cP$ be a partial path decomposition of~$D$ of size $|\cP|\leq 22\eta n$. Suppose that both $V^\pm\subseteq V^\pm(P)\cup V^0(P)$ for each $P\in \cP$.}\OLD{Suppose that~$\cP$ is a partial path decomposition of~$D$ such that $|\cP|\leq 22\eta n$ and, for each $\diamond\in \{+,-\}$ and~$v\in V(D)$ satisfying $d_D^\diamond(v)\geq \texc(D)-22\eta n$, $v\in V^\diamond(P)\cup V^0(P)$ for each~$P\in \cP$.} Then,~$\cP$ is good.\label{prop:Delta0-good}
	\end{enumerate}
\end{prop}

\begin{proof}
	\NEW{For \cref{prop:Delta0-unique}, we may assume that there exists $v\in V^+$ (similar arguments hold if $V^-\neq \emptyset$).
		Since $\exc(D)\geq (1-21\eta)n$, we have
		\[\texc(D)\leq d_D^+(v)+22\eta n\leq (1+22\eta)n\leq \exc(D)+43\eta n.\]
		Thus, \cref{prop:Delta0-excD} holds. By assumption, $d_D^+(v)\geq \texc(D)-22\eta n\geq \exc(D)-22\eta n\geq  \frac{n}{2}$ and so $d_D^-(v)\leq d_D^+(v)$. Thus, $\exc_D^+(v)=\exc_D(v)$ and so
		\begin{align*}
			\exc_D^+(v)&\stackrel{\text{\eqmakebox[Delta0][c]{\text{\cref{fact:exc}\cref{fact:exc-dmax}}}}}{=}2d_D^+(v)-d_D(v)\geq 2(\texc(D)-22\eta n)-n\\
			&\stackrel{\text{\eqmakebox[Delta0][c]{}}}{\geq}	(1-86\eta)n\stackrel{\text{\cref{prop:Delta0-excD}}}{\geq} \texc(D)-108\eta n\geq \exc(D)-108\eta n.
		\end{align*}
		Thus, \cref{prop:Delta0-excv} holds.}
	\OLD{For \cref{prop:Delta0-unique}, let~$\diamond\in \{+,-\}$ and assume~$v\in V$ satisfies
		$d_D^\diamond(v)\geq \texc(D)-22\eta n$. Note that, since $d_D^\diamond(v)\leq n$, $\texc(D)\leq (1+22\eta)n$.
		Moreover, $d_D^\diamond(v)=d_D^{\max}(v)$, since otherwise $\NEW{n-1\geq }d_D(v)\geq 2d_D^\diamond(v)\geq 2\texc(D)-44\eta n\geq 2\exc(D)-44\eta n>n$, a contradiction. 
		Thus, by
		\cref{fact:exc}\cref{fact:exc-dmax}, $\exc_D^\diamond(v)\geq 2\texc(D)-44\eta n-d_D(v)\geq 2\exc(D)-44\eta n-n\geq
		(1-86\eta)n$.}
	
	\NEW{For \cref{prop:Delta0-good}, let $\cP$ be a partial path decomposition of~$D$ of size $|\cP|\leq 22\eta n$. Suppose that both $V^\pm\subseteq V^\pm(P)\cup V^0(P)$ for each~$P\in \cP$.
		We need to show that $\cP$ is good, i.e.\ that $\texc(D\setminus \cP)=\max\{\exc(D\setminus \cP), \Delta^0(D\setminus \cP)\}=\texc(D)-|\cP|$. By \cref{prop:excpartialdecomp}\cref{prop:excpartialdecomp-general}, it is enough to show that $\Delta^0(D\setminus\cP)\leq \texc(D)-|\cP|$.
		Let $v\in V(D)$. We need to show that both $d_{D\setminus \cP}^\pm(v)\leq \texc(D)-|\cP|$. If both $d_D^\pm(v)\leq \texc(D)-|\cP|$, then we are done. We may therefore assume without loss of generality that $d_D^+(v)\geq \texc(D)-|\cP|\geq \texc(D)-22\eta n$. Then, $v\in V^+$ and so, by assumption, $d_P^+(v)=1$ for each $P\in \cP$. Thus,
		\[d_{D\setminus \cP}^+(v)=d_D^+(v)-d_\cP^+(v)=d_D^+(v)-|\cP|\leq \Delta^0(D)-|\cP|\leq \texc(D)-|\cP|.\]
		Moreover,
		\begin{align*}
			d_{D\setminus \cP}^-(v)&\leq d_D^-(v)=d_D(v)-d_D^+(v)\leq n-(\texc(D)-22\eta n)\\
			&\leq (1+22\eta)n-\exc(D)\leq 43\eta n
			\leq \texc(D)-|\cP|,
		\end{align*}
		as desired.}
	\OLD{Let~$\cP$ be as in \cref{prop:Delta0-good}. By \cref{prop:excpartialdecomp}, $\exc(D\setminus \cP)\leq \texc(D)-|\cP|\leq \texc(D\setminus \cP)$.
		Therefore, it only remains to show that $\Delta^0(D\setminus\cP)\leq \texc(D)-|\cP|$.	
		Let~$v\in V(D)$ and assume without loss of generality that $d_D^+(v)\geq d_D^-(v)$, i.e.\ that~$v\in U^+(D)\cup U^0(D)$. If $d_D^+(v)\geq \texc(D)-22\eta n>\frac{2n}{3}$, then~$v\in U^+(D)$ and so, by assumption and since~$\cP$ is a partial path decomposition of~$D$, $d_{D\setminus \cP}^-(v)\leq d_{D\setminus \cP}^+(v)=d_D^+(v)-|\cP|\leq \texc(D)-|\cP|$.
		If $d_D^+(v)\leq \texc(D)-22\eta n$, then $d_{D\setminus\cP}^{\max}(v)\leq d_D^+(v)\leq \texc(D)-|\cP|$. Therefore, $\Delta^0(D\setminus \cP)\leq \texc(D)-|\cP|$, as desired.}
\end{proof}

\subsection{Completing path decompositions via absorbing edges}\label{sec:A}

\NEW{As discussed in the proof overview, the goal is to complete our path decomposition by applying \cref{thm:niceD}. However, this requires linearly many vertices to serve as endpoints, which may not always be possible. The concept of absorbing edges provides an approach to overcome this issue. We motivate this concept via the following example. Let $D$ be a digraph and $w\in V(D)$ with $\exc_D^+(w)$ large. Suppose that $v\in N_D^+(w)\cap U^0(D)$. Note that in $D\setminus \{wv\}$, the excess of $v$ is now $1$ instead of $0$. Moreover, $U^+(D\setminus \{wv\})=U^+(D)\cup \{v\}$ and so the number of possible distinct starting points increases by one. A perfect path decomposition of $D\setminus \{wv\}$ must have a path $P$ starting at $v$. If $P$ does not contain $w$, then we can extend it to start at $w$ by adding the edge $wv$ and so we obtain a perfect decomposition of $D$. We can view the edge $wv$ as an absorbing starting edge which absorbs the path $P$. This motivates the following definition.}

\OLD{As discussed in the proof
overview, the goal is to complete our path decomposition by applying
\cref{thm:niceD}. But, it will not always be possible to find enough distinct
vertices to serve as endpoints. To solve this problem, we will set aside some
edges in order to extend the paths obtained from \cref{thm:niceD} to suitable
endpoints.}

\begin{definition}[Absorbing sets of edges]\label{def:A}
	Let~$D$ be a digraph. Let~$W,V'\subseteq V(D)$ be disjoint.
	\begin{itemize} 
		\item An \emph{absorbing set of~$(W, V')$-starting edges (for~$D$)} is a set~$A \subseteq E(D)$ of edges with starting point in~$W$ and ending point in~$V'$ such
		that, for each~$w\in W$, at most~$\exc_D^+(w)$ edges in~$A$ start at~$w$, and, for each~$v\in V'$, at most one edge in~$A$ ends at~$v$.
		\item An \emph{absorbing set of~$(V', W)$-ending edges (for~$D$)} is a set~$A\subseteq E(D)$ of edges with starting point in~$V'$ and ending point in~$W$ such that, for each~$w\in W$, at most~$\exc_D^-(w)$ edges in~$A$ end at~$w$, and, for each~$v\in V'$, at most one edge in~$A$ starts at~$v$.
		\item \NEW{A \emph{$(W,V')$-absorbing set (for $D$)} is the union of an absorbing set of $(W,V')$-starting edges and an absorbing set of $(V',W)$-ending edges.}
	\end{itemize}
\end{definition}

\NEW{Let~$D$ be a digraph. Let~$W,V'\subseteq V(D)$ be disjoint. Recall that an absorbing $(W,V')$-starting edge $wv$ can only absorb a path starting at $v$ that does not contain $w$. Thus, we will find a path decomposition in $D[V']$. We will find this decomposition via \cref{thm:niceD}. For this, we need to adapt the degree conditions to account for the absorbing paths a vertex is involved in. This is formalised in the following \lcnamecref{prop:absorbingedges}.}

\OLD{The following \lcnamecref{prop:absorbingedges} shows how absorbing edges are
used to the extend paths obtained from \cref{thm:niceD}.}

\begin{cor}\label{prop:absorbingedges}
	Let $0<\frac{1}{n}\ll \nu\ll\tau \leq \frac{\delta}{2}\leq 1$ and $r\geq \delta n$.
	Suppose that~$D$ is a digraph with a vertex partition $V(D)=W\cup V'$ such that
	$D[V']$ is a robust~$(\nu,\tau)$-outexpander on~$n$ vertices. Suppose that
	$A^+,A^-\subseteq E(D)$ are absorbing sets of~$(W,V')$-starting and~$(V',W)$-ending
	edges such that~$|A^\pm|\leq r$. Denote~$A\coloneqq A^+\cup A^-$. 
	Suppose furthermore that there exists a partition $V'=X^+\cup X^-\cup X^*\cup
	X^0$ such that
	\NEW{$|X^\pm\cup X^*|+ |A^\pm|=r$}\OLD{$|X^\pm\cup X^*\cup A^\pm|=r$}, $V(A^\pm)\cap V'\subseteq X^\mp \cup X^0$, and, for each~$v\in V(D)$,
	the following hold.
	\begin{equation*}
		\exc_D(v)=
		\begin{cases}
			\NEW{\exc}_A(v)& \text{if } v\in W,\\
			\pm 1 & \text{if } v\in X^{\pm},\\
			0 & \NEW{\text{if } v\in X^*\cup X^0},
		\end{cases}
		\quad \text{and} \quad
		d_D(v)=
		\begin{cases}
			d_A(v)& \text{if } v\in W,\\
			2r - 1 & \text{if } v\in X^{\pm},\\
			2r -2 & \text{if } v\in X^*,\\
			2r & \NEW{\text{if } v\in X^0}.
		\end{cases}
	\end{equation*}
	Then, $\pn(D)=r$. 
\end{cor}

\begin{proof}
	\NEW{By \cref{thm:niceD}, we may assume that $A\neq \emptyset$. Hence,}%
	\OLD{Observe that~$\pn(D)\geq r$. Indeed, if~$X^*=V'$, then~$D[V']$ is~$(r-1)$-regular and~$A^+\cup A^-=\emptyset$. Hence,~$W$ is a set of isolated vertices and, therefore,~$\pn(D)=\pn(D[V'])\geq r$. Otherwise,}
	$\pn(D)\geq \Delta^0(D)=r$. Thus, it
	suffices to find a path decomposition of~$D$ of size~$r$.
	
	Let \begin{align*}
		Y^\pm&\coloneqq (X^\pm \cup (V(A^\pm)\cap V'))\setminus (X^\mp\cup V(A^\mp))=(X^\pm \cup (V(A^\pm)\cap X^0))\setminus V(A^\mp),\\
		Y^*&\coloneqq X^*\cup
		(V(A^+)\cap V(A^-))\cup (X^+\cap V(A^-))\cup (X^-\cap V(A^+)), \text{ and}\\
		Y^0&\coloneqq X^0\setminus (V(A^+)\cup V(A^-)).
	\end{align*}
	Then, observe that $Y^+, Y^-, Y^*$, and $Y^0$ are all pairwise disjoint and form a partition of~$V'$. Moreover, $|Y^\pm\cup Y^*|=|X^\pm\cup X^*\cup (V(A^\pm)\cap V')|=r$%
	\COMMENT{$Y^\pm\cup Y^*=X^\pm\cup X^*\cup (V(A^\pm)\cap V')$, $|X^\pm\cup X^*\cup (V(A^\pm)\cap V')|=|X^\pm\cup X^*\cup A^\pm|=r$ by assumption, and $X^\pm, X^*, (V(A^\pm)\cap V')$ are all pairwise disjoint.}
	and, for each~$v\in
	V'$, the following hold.
	\begin{equation*}
		\exc_{D[V']}(v)=
		\begin{cases}
			\pm 1 & \text{if } v\in Y^{\pm},\\
			0 & \text{otherwise},
		\end{cases}
		\quad \text{and} \quad
		d_{D[V']}(v)=
		\begin{cases}
			2r - 1 & \text{if } v\in Y^{\pm},\\
			2r -2 & \text{if } v\in Y^*,\\
			2r & \text{otherwise}.
		\end{cases}
	\end{equation*}
	Thus, we can apply \cref{thm:niceD} with $D[V'], Y^\pm, Y^*$, and $Y^0$ playing the
	roles of $D, X^\pm, X^*$, and $X^0$ to obtain a path decomposition~$\cP\NEW{=\{P_1, \dots, P_r\}}$ of~$D[V']$ of size~$r$. \NEW{For each~$i\in [r]$, let $v_i^+$ and $v_i^-$ denote the starting and ending points of~$P_i$.
	By the ``moreover part" of \cref{thm:niceD}, we may assume that $v_1^+, \dots, v_r^+$ are distinct and $\{v_i^+\mid i\in [r]\}=Y^+\cup Y^*$. We may also assume that $v_1^-, \dots, v_r^-$ are distinct and $\{v_i^-\mid i\in [r]\}=Y^-\cup Y^*$.}\OLD{Note that each vertex in $Y^+\cup Y^*$ is the starting point of exactly one path in~$\cP$ and each vertex in $Y^-\cup Y^*$ is the ending point of exactly one path in~$\cP$. Indeed, by \cref{eq:degrees} and since~$|\cP|=r$, $V^\pm(\cP)\subseteq Y^\pm\cup Y^*$. Moreover, each vertex in $Y^\pm\cup Y^*$ is the starting/ending point of at most one path in~$\cP$. Thus, since $|Y^\pm\cup Y^*|=r$, each vertex in $Y^\pm\cup Y^*$ is the starting/ending point of exactly one path in~$\cP$.}
	
	\OLD{Denote $\cP\coloneqq \{P_1, \dots, P_r\}$ and, for each~$i\in [r]$, let
		$v_i^\pm$ denote the starting/ending point of~$P_i$.}
	We use~$A^+$ to absorb the paths starting at~$V(A^+)\cap V'$ as follows.
	For each~$i\in [r]$, if~$v_i^+\notin V(A^+)$, \NEW{then} let~$P_i^+\coloneqq P_i$;
	otherwise, denote by~$w_i^+v_i^+$ the unique edge in~$A^+$ which is incident to~$v_i^+$ and let $P_i^+\coloneqq w_i^+v_i^+P_iv_i^-$.
	Then, absorb the paths ending in~$V(A^-)\cap V'$ similarly. For each~$i\in
	[r]$, if~$v_i^-\notin V(A^-)$, \NEW{then} let~$P_i^-\coloneqq P_i^+$; otherwise, denote by~$v_i^-w_i^-$ the unique edge in~$A^-$ which is incident to~$v_i^-$ and let~$P_i^-$ be obtained by concatenating~$P_i^+$ and~$v_i^-w_i^-$.
	\NEW{Since $d_D(v)=d_A(v)$ for each $v\in W$, it follows that}\OLD{Then,} $\cP'\coloneqq \{P_i^-\mid i\in [r]\}$ is a path decomposition of~$D$ of size~$r$, as desired. 
\end{proof}

\NEW{Although the absorbing set is chosen at the beginning, we do not remove this set as it may affect our calculation of $\texc(D)$. Thus, we require all our partial path decompositions to avoid the edges in the absorbing set. Moreover, their endpoints should not ``overuse" the vertices in $V(A)$.}

\begin{definition}[Consistent partial path decomposition]
	Let~$D$ be a digraph, let~$W, V'\subseteq V(D)$ be disjoint, and \NEW{$A\subseteq E(D)$.}\OLD{$A^\pm\subseteq
	E(D)$.} \NEW{Note that $W$ and $V'$ do not necessarily partition $V(D)$. Suppose that $A$ is a $(W,V')$-absorbing set.}\OLD{Suppose that~$A^+$ and~$A^-$ are absorbing sets of~$(W, V')$-starting and
	$(V', W)$-ending edges for~$D$. Denote~$A\coloneqq A^+\cup A^-$.}
	Then, a partial path decomposition~$\cP$ of~$D$ is \emph{consistent} with \NEW{$A$}\OLD{$A^+$
	and~$A^-$} if~$\cP\subseteq D\setminus A$ and each~$v\in W$ is the starting point of at most~$\exc_D^+(v)-d_A^+(v)$ paths in
	$\cP$ and the ending point of at most~$\exc_D^-(v)-d_A^-(v)$ paths in~$\cP$.
\end{definition}

\begin{definition}[$(U^*,W,A)$-partial path decomposition]\label{def:U*WApartialpathdecomp}
	\NEW{Let~$D$ be a digraph, let~$W, V'\subseteq V(D)$ be disjoint, and $A\subseteq E(D)$. Suppose that $A$ is a $(W,V')$-absorbing set. Given $U^*\subseteq U^0(D)$ satisfying $|U^*|\leq \texc(D)-\exc(D)$, we say that $\cP$ is a \emph{$(U^*,W,A)$-partial path decomposition} if $\cP$ is a $U^*$-partial path decomposition which is consistent with $A$.}
\end{definition}

\NEW{Let $\cP$ be a partial path decomposition of $D$ which is consistent with $A$. By definition, $A$ is still a $(W,V')$-absorbing set for $D\setminus \cP$.}

\begin{fact}\label{fact:absorbingedges}
	Let~$D$ be a digraph and~$W,V'\subseteq V(D)$ be disjoint. Suppose that \NEW{$A$ is a $(W,V')$-absorbing set.}\OLD{$A^+, A^-\subseteq E(D)$ are absorbing sets of~$(W,V')$-starting and~$(V',W)$-ending edges for~$D$. Denote~$A\coloneqq A^+\cup A^-$.}
	Suppose~$\cP$ is a partial path decomposition of~$D$ which is consistent with \NEW{$A$}\OLD{$A^+$ and~$A^-$}. Denote~$D'\coloneqq D\setminus \cP$.
	Then, \NEW{$A$ is a $(W,V')$-absorbing set}\OLD{$A^+$ and~$A^-$ are absorbing sets of~$(W,V')$-starting and~$(V',W)$-ending edges} for~$D'$.
\end{fact}

%% file: Constructing_Layouts_Definitions_Statements.tex
	\onlyinsubfile{
		\setcounter{section}{8}
\section{Constructing layouts in general tournaments}}

\NEW{In this section, we discuss how to construct layouts in general tournaments. 
Recall that \cref{lm:approxdecomp} (which constructs an approximate decomposition which respects a given set of layouts) only applies to almost regular robust outexpanders. In general, our tournament $T$ will not be almost regular nor a robust outexpander. In \cref{sec:W}, we discuss how to circumvent this problem.
As discussed in \cref{sec:goodpartialdecomp}, we will need the set of paths obtained with \cref{lm:approxdecomp} to form a good partial path decomposition. In \cref{sec:pathL}, we explain how we can ensure this.
In \cref{sec:cleaninglm}, we discuss the cleaning step.
In \cref{sec:layoutlm}, we state the lemma which guarantees the existence of suitable layouts.
}

\subsection{\texorpdfstring{$W$}{W}-exceptional layouts}\label{sec:W}

\NEW{Let $T$ be a tournament on $n$ vertices with $\exc(T)\leq \varepsilon n^2$ for some small constant $\varepsilon$. Then, there exists a partition of $V(T)$ into $W$ and $V'$ such that $W$ is small and $T[V']$ is almost regular. Our aim is to apply \cref{lm:approxdecomp} to $T[V']$. To do so, we will construct layouts $(L_1, F_1), \dots, (L_\ell, F_\ell)$ so that $E_W(T)\subseteq \bigcup_{i\in [\ell]}F_i$. (Recall from \cref{sec:notation} that $E_W(T)$ denotes the set of edges of $T$ which are incident to $W$.) This will ensure that all the edges in $E_W(T)$ will be contained in the partial path decomposition obtained from the spanning configurations of shapes $(L_1, F_1), \dots, (L_\ell, F_\ell)$.}

\OLD{It turns out to be convenient to first transform a~\emph{$W$-exceptional layout on~$V$} (defined formally below) into a layout on~$V\setminus W$. This allows us to ignore~$W$ when turning layouts into spanning configurations.}

\NEW{Let $V$ be a vertex set and $W\subseteq V$.}\OLD{For~$W\subseteq V$,} We say that a layout~$(L,F)$ is \emph{$W$-exceptional} if~$E_W(L) \subseteq F$%
\COMMENT{i.e. all endpoints (where we need to construct a path) avoid $W$.}.%
\OLD{(Recall that~$E_W(L)$ denotes the set of all those edges of~$L$ which have at least one endpoint in~$W$.)}
Let \OLD{$W\subseteq V$ and}$(L,F)$ be a~$W$-exceptional layout. A multidigraph~$\cH$ on~$V$ is a \emph{$W$-exceptional spanning configuration of shape~$(L,F)$} if~$\cH$ can be decomposed into internally vertex-disjoint paths $\{P_e \mid e \in E(L)\}$ such that each~$P_e$ has shape~$e$;~$P_f = f$ for all~$f \in F$; and \NEW{$\bigcup_{e \in E(L)}V^0(P_e) = V \setminus (V(L)\cup W)$}\OLD{$V^0(\bigcup_{e \in E(L)} \{P_e\} ) = V \setminus (V(L)\cup W)$}.
(Note that the last equality implies that the vertices in~$W\setminus V(L)$ are isolated in~$\cH$.)
Thus, roughly speaking, a~$W$-exceptional spanning configuration of shape~$(L,F)$ is one such that all ``additional" edges (i.e.\ those edges of~$\cH$ that are not in~$F$) are disjoint from~$W$.
\NEW{A path decomposition of $\cH$ is \emph{induced} by $(L,F)$ if it consists of all the paths $P_Q\coloneqq \{P_e\mid e\in E(Q)\}$ where $Q$ is a path in $L$. 
The analogue of \cref{fact:layoutbijection} holds for $W$-exceptional spanning configurations.}

\NEW{\begin{fact}\label{fact:layoutbijectionW}
	Let $V$ be a vertex set and $W\subseteq V$. Let $(L,F)$ be a $W$-exceptional layout on $V$ and $\cH$ be a $W$-exceptional spanning configuration of shape $(L,F)$ on $V$. Let $L'$ denote the set of (non-trivial) paths contained in $L$. Suppose that $\cP$ is a path decomposition of $\cH$ which is induced by $(L,F)$. Then, $V^\pm(\cP)=V^\pm(L')$ and $V^0(\cP)=V^0(L)\cup (V\setminus (W\cup V(L)))$.
\end{fact}}

\NEW{We now show that there is a natural transformation of a $W$-exceptional layout $(L,F)$ into an auxiliary layout $(L^{\upharpoonright W}, F^{\upharpoonright W})$ on $V\setminus W$. Roughly speaking, the auxiliary layout $(L^{\upharpoonright W}, F^{\upharpoonright W})$ is obtained from $(L,F)$ by contracting all the edges in $E_W(L)$ so that $E_W(L^{\upharpoonright W})=\emptyset =E_W(F^{\upharpoonright W})$ and then remove $W$.}

\begin{definition}[Auxiliary layout]\label{def:layoutV'}
	Let~$V$ be a vertex set and~$W\subseteq V$. Suppose~$(L,F)$ is a~$W$-exceptional layout on~$V$. We denote by $(L^{\upharpoonright W}, F^{\upharpoonright W})$ the layout on~$V\setminus W$ obtained from~$(L,F)$ as follows.
	
	Let~$\cP$ be the multiset of maximal paths~$P$ such that~$P\subseteq P'$ for some~$P'\in L$,~$V^0(P)\subseteq W$, and~$V(P)\cap W\neq \emptyset$ (in particular, each isolated vertex~$v\in V(L)\cap W$ is a path in~$\cP$ but no isolated vertex~$v\in V(L)\setminus W$ is a path in~$\cP$).
	Note that, since~$(L,F)$ is~$W$-exceptional, each~$P\in \cP$ satisfies~$E(P)\subseteq F$. 
	Let $P_1, \dots, P_k$ be an enumeration of~$\cP$ and, for each~$i\in [k]$, let~$x_i$ and~$y_i$ denote the starting and ending points of~$P_i$, respectively.
	Then, let~$L^{\upharpoonright W}$ be obtained from~$L$ as follows. For each~$i\in [k]$,
	\begin{itemize}
		\item if both $x_i,y_i\in V\setminus W$, then contract the subpath~$P_i$ into an edge~$x_iy_i$;
		\item otherwise, delete~$E(P_i)$ as well as $V(P_i)\cap W$.
	\end{itemize}
	Note that $V(L^{\upharpoonright W})=V(L)\setminus W\subseteq V\setminus W$.
	Define $F^{\upharpoonright W}\coloneqq \{x_iy_i\mid i\in [k], x_i,y_i\in V\setminus W\}\cup (F\setminus E_W(F))=\{x_iy_i\mid i\in [k], x_i,y_i\in V\setminus W\}\cup (F\setminus E_W(L))$.%
	\COMMENT{Let $P\in L$ be such that $P_i\subseteq P$. If $V(P_i)\subseteq V^0(P)$, then, by maximality of $P_i$, both $x_i,y_i\in V\setminus W$ and so we contract $P_i$ to fixed edge. Suppose $P_i=P$. If both $x_i,y_i\in V\setminus W$, then we contract $P_i$ to a fixed edge; otherwise, we delete $P_i$ and add the vertices in $\{x_i,y_i\}\setminus W$ as isolated vertices. Suppose that $V^+(P)\subseteq V(P_i)\subseteq V^+(P)\cup V^0(P)$. Then, note that, by maximality of $P_i$, $y_i\in V\setminus W$. If $x_i\in V\setminus W$, then we contract $P_i$ to a fixed edge; otherwise, we delete $V(P_i)\cap W$. Suppose that $V^-(P)\subseteq V(P_i)\subseteq V^-(P)\cup V^0(P)$. Then, note that, by maximality of $P_i$, $x_i\in V\setminus W$. If $y_i\in V\setminus W$, then we contract $P_i$ to a fixed edge; otherwise, we delete $V(P_i)\cap W$.}
\end{definition}

\NEW{Note that since $(L,F)$ is $W$-exceptional, each $e\in E(L)\setminus F$ satisfies $V(e)\subseteq V\setminus W$ and so $E(L)\setminus F=E(L^{\upharpoonright W})\setminus F^{\upharpoonright W}$.}

The following \lcnamecref{fact:layoutV'} states that a spanning configuration of shape $(L^{\upharpoonright W}, F^{\upharpoonright W})$ in~$D[V\setminus W]$ can easily be transformed into a~$W$-exceptional spanning configuration of shape~$(L,F)$ in~$D$.
In other words, it allows us to reverse the process described in \cref{def:layoutV'}.

\begin{prop}\label{fact:layoutV'}
	Let~$D$ be a digraph on a vertex set~$V$. Let~$W\subseteq V$ and denote $V'\coloneqq V\setminus W$. Let~$(L,F)$ be a~$W$-exceptional layout on~$V$. Let $(L^{\upharpoonright W}, F^{\upharpoonright W})$ be as in \cref{def:layoutV'}.
	Suppose $\cH^{\upharpoonright W}\subseteq D[V']\cup F^{\upharpoonright W}$ is a spanning configuration of shape $(L^{\upharpoonright W},F^{\upharpoonright W})$.
	Let~$\cH$ be the multidigraph with~$V(\cH)\coloneqq V$ and $E(\cH)\coloneqq (E(\cH^{\upharpoonright W})\setminus F^{\upharpoonright W}) \cup F$. Then, \NEW{$E(\cH)\subseteq E(D[V'])\cup F$}\OLD{$\cH\subseteq D[V']\cup F$} and~$\cH$ is a~$W$-exceptional spanning configuration of shape~$(L,F)$.
\end{prop}

\begin{proof}
	\NEW{Note that $H^{\upharpoonright W}$ can be decomposed into internally vertex-disjoint paths $\{P_e\mid e\in E(L^{\upharpoonright W})\}$ such that $P_e$ has shape $e$ for each $e \in E(L^{\upharpoonright W})$; $P_f=f$ for each $f\in F^{\upharpoonright W}$; and $\bigcup_{e\in E(L^{\upharpoonright W})}V^0(P_e)=V'\setminus V(L)$. Since $E(L)\setminus F=E(L^{\upharpoonright W})\setminus F^{\upharpoonright W}$, $H$ can be decomposed into $\{P_e\mid e\in E(L)\setminus F\}\cup F$. The \lcnamecref{fact:layoutV'} follows by setting $P_f\coloneqq f$ for all $f\in F$.}
\end{proof}

This has the advantage that it suffices to find spanning configurations in an almost regular robust outexpander, which corresponds to the setting of \cref{lm:approxdecomp}. More precisely, if we let~$V'\coloneqq V\setminus W$ be the set of ``non-exceptional vertices" described in \cref{sec:sketch-general} and let~$D'$ be the remainder of the tournament~$T[V']$ after the cleaning step, then~$D'$ is almost complete and almost regular, and hence a robust outexpander. Then, we can split~$D'$ into~$D$ and~$\Gamma$ as required for \cref{lm:approxdecomp}.

\subsection{Path consistent layouts}\label{sec:pathL}

Let~$D$ be a digraph on~$V$. 
When we refer to a spanning configuration of shape~$(L,F)$ in~$D$, we mean that this configuration is contained in the multidigraph $D \cup F$ (as~$F$ may not be in~$D$). 
Let~$\cF$ be a multiset of edges on~$V$ and $(L_1,F_1), \dots, (L_\ell, F_\ell)$ be layouts, where~$F_i\subseteq \cF$ for each~$i\in [\ell]$. 
We would like the union of their spanning configurations to form a good partial path decomposition of~$D \cup \mathcal{F}$.
For this, these layouts will need to satisfy the following properties.
Let $U^*\subseteq U^0(D\cup \cF)$ be such that $|U^*|\leq \texc(D\cup\cF)-\exc(D\cup\cF)$ and define the multiset~$L$ by~$L \coloneqq \bigcup_{i \in [\ell]} L_i$. 
We say $(L_1,F_1), \dots, (L_\ell, F_\ell)$ are \emph{$U^*$-path consistent with respect to~$(D, \mathcal{F})$}, if $\bigcup_{i \in [\ell] } F_i\subseteq\cF$ (counting \NEW{with} multiplicity) and the following hold.
\begin{enumerate}[label=\upshape(P\arabic*$'$)]
	\item For any $v\in V \setminus U^0(D\cup \cF)$,~$v$ is the starting point of at most~$\exc_{D\cup \cF}^+(v)$ non-trivial paths in~$L$ and the ending point of at most~$\exc_{D\cup \cF}^-(v)$ non-trivial paths in~$L$.
	\item For any~$v\in U^*$,~$L$ contains at most one non-trivial path starting at~$v$ and at most one non-trivial path ending at~$v$.
	\item If $v\in U^0(D\cup\cF)\setminus U^*$, then~$v$ is not an endpoint of any non-trivial path in~$L$.
\end{enumerate}
If~$D$ and~$\cF$ are clear from the context, then we omit ``with respect to~$(D, \cF)$''.

\NEW{The following \lcnamecref{fact:layouts} simply states that the union of spanning configurations of $U^*$-path consistent layouts indeed forms a $U^*$-partial path decomposition (as defined in \cref{sec:goodpartialdecomp}). We also track the degrees for later uses.}

\begin{prop}\label{fact:layouts}
	Let~$D$ be a digraph on a vertex set~$V$. Let~$V=W\cup V'$ be a partition of~$V$. Let $U^*\subseteq U^0(D)$ satisfy $|U^*|\leq \texc(D)-\exc(D)$ and $\cF\subseteq E(D)$.
	Let $(L_1,F_1) \dots (L_\ell,F_\ell)$ be~$W$-exceptional layouts.
	For each~$i\in [\ell]$, let~$\cH_i$ be a~$W$-exceptional spanning configuration of shape~$(L_i,F_i)$. Suppose that $\cH_1, \dots, \cH_\ell$ are pairwise edge-disjoint.
	For each~$i\in [\ell]$, denote by~$\cP_i$ a path decomposition of~$\cH_i$ induced by~$(L_i,F_i)$.
	Define the multiset~$L$ by $L\coloneqq \bigcup_{i \in [\ell]}L_i$. Let $F\coloneqq \bigcup_{i \in [\ell]}F_i$, $\cH\coloneqq \bigcup_{i \in [\ell]}\cH_i$, and $\cP\coloneqq \bigcup_{i \in [\ell]}\cP_i$.
	Then, for all~$v \in V$, 
	\[d^{\pm}_\cH(v) = d^{\pm}_L(v)+ |\{i \in [\ell] \mid v \in V'\setminus V(L_i)\}|.\]
	Moreover, if $(L_1,F_1) \dots (L_\ell,F_\ell)$ are~$U^*$-path consistent with respect to~$(D\setminus\cF,\cF)$, then~$\cP$ is a \NEW{$U^*$}-partial path decomposition of~$D$ such that~$|\cP|$ is equal to the number of non-trivial paths in~$L$ \NEW{and} $E_W(\cP)\subseteq F\subseteq \cF$\OLD{and, for each~$v\in U^0(D)$, if~$v$ is an endpoint of a (non-trivial) path in~$\cP$, then~$v\in U^*$}.
\end{prop}

\begin{proof}
		\NEW{By \cref{fact:layoutbijectionW}, each $v\in V$ satisfies
		\begin{align*}
			d^{\pm}_\cH(v)&=\sum_{i\in [\ell]}d_{\cH_i}^\pm(v) 
			=\sum_{i\in [\ell]}(d_{L_i}^\pm(v) +\mathds{1}_{v\notin V(L_i)\cup W})= d^{\pm}_L(v)+ |\{i \in [\ell] \mid v \in V'\setminus V(L_i)\}|,
		\end{align*}
		as desired.
		Moreover, \cref{fact:layoutbijectionW} implies that, for each $v\in V$, the number of paths in $\cP$ which start/end at $v$ is precisely the number of (non-trivial) paths in $L$ which start/end at $v$.		
		By definition of path consistency, this implies that $\cP$ satisfies \cref{def:partialpathdecomp-exc,def:partialpathdecomp-U0}. Moreover, the fact that $|U^*|\leq \texc(D)-\exc(D)$ implies that \cref{def:partialpathdecomp-texc} holds. Thus, $\cP$ is a partial path decomposition of~$D$.
		}
\end{proof}

\subsection{Cleaning}\label{sec:cleaninglm}

\NEW{As discussed in \cref{sec:A}, the leftover from the approximate decomposition will be decomposed using \cref{prop:absorbingedges} and so it needs to have a specific structure: no non-absorbing edge incident to the exceptional set $W$ can be left over, while the non-exceptional vertices in $V'$ must form a digraph which is very close to being regular. One consequence of this is that the degree at $W$ needs to be covered at a much faster rate than the degree at $V'$. Unfortunately, this cannot be achieved via the approximate decomposition. Indeed, \cref{lm:approxdecomp}\cref{lm:approxdecomp-degreev} implies that each $v\in V'$ can only be included as an isolated vertex or covered by a fixed edge in a small proportion of the layouts. Thus, \cref{fact:layouts} implies that the $\ell$ spanning configurations obtained with \cref{lm:approxdecomp} will cover about $\ell$ inedges and $\ell$ outedges at each vertex in $V'$. Therefore, to cover the degree at $W$ at a faster rate than the vertices in $V'$, we would need that each vertex in $W$ belongs (on average) to several paths of each spanning configuration. However, as discussed in \cref{sec:W}, the exceptional vertices will be included via fixed edges and so, in that case, the layouts would be large, while the approximate decomposition only allows small layouts (see \cref{lm:approxdecomp}\cref{lm:approxdecomp-sizeLi})%
	\COMMENT{If $w\in W$ is at the start or end of a path in $L$, then this does not induce any fixed in that corresponding layout on $V'$. However, not all vertices in $W$ will have very large excess.}.}

\NEW{Therefore, we will start with a cleaning procedure which significantly reduces the degree at $W$. To facilitate the construction of layouts, we also cover all the edges inside $W$ in this step. Note that this needs to be done efficiently so that, after the cleaning step, the non-exceptional vertices still form an almost regular oriented graph of very large degree (otherwise, we would not be able to apply \cref{lm:approxdecomp} to obtain an approximate decomposition).}

\OLD{As discussed in the proof overview,
vertices of very high excess and vertices adjacent to absorbing edges will be treated as exceptional vertices throughout the proof of \cref{thm:main}. To be able to cover the edges at these vertices in the approximate decomposition step, we need to start with a cleaning procedure, which will ensure that no two exceptional vertices have an edge between them and that the degree of the exceptional vertices is not too large.
This needs be done efficiently so that, after the cleaning step, the non-exceptional vertices still form an almost regular oriented graph of very large degree.}

\NEW{We now state our cleaning lemma. (The proof is deferred to \cref{sec:cleaning}.) Roughly speaking, \cref{lm:cleaning} says the following. Suppose that $T\notin \cT_{\rm excep}$ (exceptional tournaments have already been decomposed in \cref{sec:annoyingT}). Let $W\subseteq V(T)$ consist of all the vertices of excess at least $\varepsilon n$ and denote $V'\coloneqq V(T)\setminus W$. Let $A^+$ and $A^-$ be small absorbing sets of $(W,V')$-starting and $(V',W)$-ending edges. Let $U^*\subseteq U^0(T)$ satisfy $\texc(T)-\exc(T)$. Then, there exists a good $(U^*,W,A)$-partial path decomposition $\cP$ such that the leftover $D\coloneqq T\setminus \cP$ satisfies the following properties. First, the main objectives of the cleaning step are achieved.
\begin{itemize}
	\item The degree of the exceptional vertices (i.e.\ those in $W$) is significantly lower than the degree of the vertices in $V'$ (compare the bounds in \cref{lm:cleaning}\cref{lm:cleaning-degreeW*,lm:cleaning-degreeV'}). 
	\item All the edges inside the exceptional set are covered (see \cref{lm:cleaning}\cref{lm:cleaning-W}).
\end{itemize}
Moreover, these objectives are achieved very efficiently.
\begin{itemize}[--]
	\item $D$ is still almost complete (see \cref{lm:cleaning}\cref{lm:cleaning-d,lm:cleaning-degreeV',lm:cleaning-degreeW*,lm:cleaning-degreeW0}). Together with \cref{lm:3/8rob}, this will ensure that $D[V']$ is an almost regular robust outexpander, which is needed for the approximate decomposition.
	\item $\texc(D)$ is large compared to number of edges in $D$ (see \cref{lm:cleaning}\cref{lm:cleaning-good,lm:cleaning-exc>d} (by \cref{lm:cleaning}\cref{lm:cleaning-degreeV',lm:cleaning-degreeW*,lm:cleaning-degreeW0}, $d$ is roughly the density of $D$)). If $\texc(D)$ was very small, then $D$ would need to be decomposed with few long paths which would be more difficult and might not even be possible with our strategy%
		\COMMENT{If \cref{lm:cleaning}\cref{lm:cleaning-good} fails, then a decomposition which uses $U^*$ as endpoints will have size $>\texc(D)$.}.
	\item The number of distinct endpoints which can be used to decompose $D$ is roughly the same as for $T$ (see \cref{lm:cleaning}\cref{lm:cleaning-sizeU}). As discussed in \cref{sec:good}, having a large pool of suitable endpoints is very convenient and in fact necessary for the final step of the decomposition if $A^+\cup A^-=\emptyset$ (recall that we aim to apply \cref{prop:absorbingedges} after the approximate decomposition).
\end{itemize}
We now explain and motivate the conditions which are needed for the cleaning strategy to work or to simplify the construction of layouts.
\begin{itemize}
	\item The number of exceptional vertices must be small (see \cref{lm:cleaning}\cref{lm:cleaning-defW}). Otherwise, there would be too many edges to cover within the exceptional set and we would not be able to it efficiently. We are able to assume that $|W|$ is small since otherwise $\exc(T)$ would be large and so \cref{thm:evenlarge} would apply (recall that $W$ consists of vertices of large excess).
	\item The absorbing edges are taken at vertices of as high excess as possible (see \cref{lm:cleaning}\cref{lm:cleaning-A}). This is convenient because it maximises the number of endpoints that are allowed to be used in $\cP$. Indeed, recall that the effect of absorbing edges is to reserve some excess at the vertices of $W$ for the final step of the decomposition, so if the absorbing edges account for all the excess at a vertex $w\in W$, then $w$ cannot be used as an endpoint in the $(U^*,W,A)$-partial path decomposition $\cP$. Taking the absorbing edges at vertices of as high excess as possible ensures that this occurs for as few vertices $w\in W$ as possible.
	\item We distinguish exceptional vertices of very high excess ($W_*$), exceptional vertices of significant but not too high excess ($W_0$), and exceptional vertices which are incident to absorbing edges ($W_A$) (see \cref{lm:cleaning}\cref{lm:cleaning-defW,lm:cleaning-A}). One issue that we have not discussed so far is that the exceptional set for the approximate decomposition and the exceptional set for the final step of the decomposition will have to be different. Indeed, as discussed in \cref{sec:W}, the exceptional set for the approximate decomposition must contain all the vertices with excess at least $\varepsilon n$ to ensure that we apply \cref{lm:approxdecomp} to an almost regular digraph. Thus, $W=W_*\cup W_0$ will be the exceptional set considered during the approximate decomposition. As discussed in \cref{sec:A}, the main role of the exceptional set in the final step of the decomposition is to incorporate the absorbing edges. Thus, $W_A$ will have to be part of the exceptional set when we apply \cref{prop:absorbingedges}. In addition, the vertices of very high excess will also have to be part of this exceptional set because they have almost all of their edges in the same direction and so it would be impossible for them to satisfy the degree conditions of the non-exceptional vertices in \cref{prop:absorbingedges}. Thus, $W_*\cup W_A$ will be the exceptional set used in the final step of the decomposition and $W_0\setminus W_A$ will be incorporated back into the non-exceptional set after the approximate decomposition. (This explains the degree conditions in \cref{lm:cleaning}\cref{lm:cleaning-degreeW0}.) This is necessary because it may not be possible to decrease the degree at $W_0$ significantly during the cleaning step. Indeed, since the excess of the vertices in $W_0$ is not too large, we may have $\texc(T)$ relatively small and $W_0$ relatively large at the same time. In that case, significantly decreasing the degree at $W_0$ during the cleaning step would amount to covering many edges with very few paths, which is not possible.
	\item If $A^+\cup A^-$ is non-empty, then $\texc(T)$ must not be too small (see \cref{lm:cleaning}\cref{lm:cleaning-A}). This will allow us to significantly reduce the degree at $W_A$ during the cleaning step (which is, as discussed above, necessary for applying \cref{prop:absorbingedges}).
\end{itemize}}
\NEW{In addition to our main objectives, we will also achieve the following property.
\begin{itemize}
	\item If $\texc(D)$ is not too large, then we can achieve that all the vertices in $W_*$ have all their edges in the same direction in $D$ (see \cref{lm:cleaning}\cref{lm:cleaning-exc<2d}). This means that, in the decomposition of $D$, no path will need to have a vertex of $W_*$ as an internal vertex. This will be very convenient because if $\texc(D)$ is relatively small but a vertex $w\in W_*$ has large positive excess (say), then almost all of the positive excess of $D$ is concentrated at $w$ and so almost all of the paths in the decomposition of $D$ have to start at $w$. Then, there would be very few paths were $w$ could be incorporated as an internal vertex and so it would be difficult to cover all the inedges at $w$.
\end{itemize}}

\begin{lm}[Cleaning lemma]\label{lm:cleaning}
	Let $0<\frac{1}{n}\ll \varepsilon\ll \eta \ll 1$. Let~$T\notin \cT_{\rm excep}$ be a tournament on a vertex set~$V$ of size~$n$ satisfying the following properties.
	\begin{enumerate}
		\item Let~$W_*\cup W_0\cup V'$ be a partition of~$V$ such that, for each~$w_*\in W_*$, $|\exc_T(w_*)|>(1-20\eta)n$; for each~$w_0\in W_0$, $|\exc_T(w_0)|\leq (1-20\eta)n$; and, for each~$v'\in V'$, $|\exc_T(v')|\leq \varepsilon n$. Let $W\coloneqq W_*\cup W_0$ and suppose $|W|\leq \varepsilon n$.\label{lm:cleaning-defW}
		\item Let $A^+,A^-\subseteq E(T)$ be absorbing sets of~$(W, V')$-starting/$(V',W)$-ending edges for~$T$ of size at most~$\lceil\eta n\rceil$. 
		Denote $A\coloneqq A^+\cup A^-$. 
		Let $W_A^\pm \coloneqq V(A^\pm)\cap W$ and $W_A\coloneqq V(A)\cap W$.
		\NEW{Suppose that the following hold.
		\begin{itemize}
			\item Let $\diamond\in \{+,-\}$. If $|W_A^\diamond|\geq 2$, then $\exc_T^\diamond(v)<\lceil\eta n\rceil$ for each~$v\in V$.
			\item Let $\diamond\in \{+,-\}$. If $|W_A^\diamond|=1$, then $\exc_T^\diamond(v)\leq \exc_T^\diamond(w)$ for each~$v\in V$ and $w\in W_A^\diamond$.
			\item If $W_A\neq \emptyset$, then $\texc(T)\geq \frac{n}{2}+10\eta n$.
		\end{itemize}}\OLD{Suppose that, for each $\diamond\in \{+,-\}$, if $|W_A^\diamond|\geq 2$, then $\exc_T^\diamond(v)<\lceil\eta n\rceil$ for each~$v\in V$ and, if $|W_A^\diamond|=1$, then, for each~$v\in V$ and $w\in W_A^\diamond$, $\exc_T^\diamond(v)\leq \exc_T^\diamond(w)$.
		Moreover, suppose that if $W_A\neq \emptyset$, then $\texc(T)\geq \frac{n}{2}+10\eta n$.}\label{lm:cleaning-A}
		\item Let $U^*\subseteq U^0(T)$ satisfy $|U^*|=\texc(T)-\exc(T)$.\label{lm:cleaning-U*}	 
	\end{enumerate}	
	Then, there exist~$d\in \mathbb{N}$ and a good $(U^*, W, A)$-partial path decomposition~$\cP$ of~$T$ such that the following hold, where $D\coloneqq T\setminus \cP$.
	\begin{enumerate}[label=\rm(\roman*)]
		\item $\left\lceil\frac{n}{2}\right\rceil-10\eta n\leq d\leq \left\lceil\frac{n}{2}\right\rceil-\eta n$. \label{lm:cleaning-d}
		\item Each $v\in U^*\setminus (V^+(\cP)\cup V^-(\cP))$ satisfies $d_D^+(v)=d_D^-(v)\leq \texc(D)-1$\label{lm:cleaning-good}%
		\COMMENT{This will ensure we can decompose~$D$ into~$\texc(D)$ paths. In particular, it implies that~$D$ is not regular.}. 
		\item $E(D[W])=\emptyset$.\label{lm:cleaning-W}
		\item $N^\pm(T)-N^\pm(D)\leq 89\eta n$.\label{lm:cleaning-sizeU}
		\item $\texc(D)\geq d+\lceil\eta n\rceil$.\label{lm:cleaning-exc>d}
		\item If $\texc(D)< 2d+\lceil\eta n\rceil$, then each~$w\in W_*$ satisfies $|\exc_D(w)|=d_D(w)$.\label{lm:cleaning-exc<2d}
		\item For each~$v\in W_*\cup W_A$, $2d-3\sqrt{\eta}n\leq d_D(v)\leq 2d-\lceil\eta n\rceil$.\label{lm:cleaning-degreeW*}
		\item For each~$v\in W_0$, $2d+2\lceil\eta n\rceil-4\sqrt{\eta}n\leq d_D(v)\leq 2d+2\lceil\eta n\rceil$ and $d_D^{\min}(v)\geq \lceil\eta n\rceil$.\label{lm:cleaning-degreeW0}
		\item For each~$v\in V'$, $2d+2\lceil\eta n\rceil-9\sqrt{\varepsilon}n\leq d_D(v)\leq 2d+2\lceil\eta n\rceil$.\label{lm:cleaning-degreeV'}	
	\end{enumerate}
\end{lm}

\OLD{\cref{lm:cleaning} will be proved in \cref{sec:cleaning}.
As mentioned above, the cleaning step has two aims: cleaning the edges inside~$W$ (property \cref{lm:cleaning-W}) and bounding the degree of vertices in~$W_*\cup W_A$ from above (compare the upper bound in \cref{lm:cleaning-degreeW*} to those in \cref{lm:cleaning-degreeV',lm:cleaning-degreeW0}). As discussed in the second paragraph after the statement of \cref{lm:layouts}, property \cref{lm:cleaning-good} is needed to ensure that our final decomposition contains the desired number of paths. 
\cref{lm:cleaning-sizeU} ensures that there are few vertices whose excess has been completely used (this is necessary for the last step of the decomposition). \cref{lm:cleaning-d,lm:cleaning-degreeV',lm:cleaning-degreeW*,lm:cleaning-degreeW0} ensure that the leftover oriented graph~$D$ is still almost complete. \cref{lm:cleaning-exc>d} ensures that the leftover excess is sufficiently large compared to the number of edges left to cover.  \cref{lm:cleaning-exc<2d} is a technical condition that will enable us to cover all edges at~$W_*$ when the total excess is not too large.}

\subsection{Constructing layouts}\label{sec:layoutlm}
\NEW{We now state the lemma which we will use to construct the layouts for the approximate decomposition. (The proof is deferred to \cref{sec:layouts}.)
Roughly speaking, \cref{lm:layouts} says the following. Let $D$ be an oriented graph. Let $W_1\cup W_2\cup V'$ be a partition of $V(D)$ and denote $W\coloneqq W_1\cup W_2$. 
Here, $W$ will be the exceptional set for the approximate decomposition (i.e.\ the same $W$ as in the cleaning lemma), $W_1$ will be the exceptional set for the final step of the decomposition (i.e.\ $W_*\cup W_A$ from the cleaning lemma), and $W_2$ will be the set of exceptional vertices which will be incorporated back into the non-exceptional set after the approximate decomposition (i.e.\ $W_0\setminus W_A$ from the cleaning lemma).
Let $U^*\subseteq U^0(D)$ satisfy $|U^*|=\texc(D)-\exc(D)$. Let $A^+$ and $A^-$ be absorbing sets of $(W_1,V')$-starting and $(V',W_1)$-ending edges.
Suppose that $D$ satisfies the following properties.
\begin{itemize}
	\item The vertices in $V'$ all have small excess (see \cref{lm:layouts}\cref{lm:layouts-degree&excessV'}, this is inherited from \cref{lm:cleaning}\cref{lm:cleaning-defW}). Recall that one of the roles of the layouts is to prescribe the endpoints of the paths we want to construct in the approximate decomposition. Thus, the fact that the vertices in $V'$ have small excess means that each vertex in $V'$ will be an endpoint in only few of the layouts. This is necessary because \cref{lm:approxdecomp} only allows each vertex to be covered by few layouts (see \cref{lm:approxdecomp}\cref{lm:approxdecomp-degreev}).
	\item The vertices in $V'$ all have roughly the same degree (see \cref{lm:layouts}\cref{lm:layouts-degree&excessV'}, this is inherited from \cref{lm:cleaning}\cref{lm:cleaning-degreeV'}). Recall that the role of the isolated vertices in layouts is to specify which vertices need to be avoided in each of the spanning configurations constructed in the approximate decomposition. The fact that the vertices in $V'$ all have roughly the same degree means that they all need to be covered by roughly the same number of spanning configurations. Thus, each vertex in $V'$ will only have to be included as an isolated vertex in few of the layouts. This is necessary because \cref{lm:approxdecomp} only allows each vertex to be included in few layouts (see \cref{lm:approxdecomp}\cref{lm:approxdecomp-degreev}).
	\item The degree at $W_1$ is significantly smaller than the degree at $V'$ (compare \cref{lm:layouts}\cref{lm:layouts-degree&excessV',lm:layouts-degree&excessW1}, this is inherited from \cref{lm:cleaning}\cref{lm:cleaning-degreeV',lm:cleaning-degreeW*}). This will enable us to incorporate all the non-absorbing edges at $W_1$ into the layouts. (Thus, the vertices in $W_1$ will have no non-absorbing edges left over after the approximate decomposition and so we will be able to use $W_1$ as the exceptional set in \cref{prop:absorbingedges}.)
	\item The degree at $W_2$ is comparable or smaller than the degree at $V'$ but every vertex in $W_2$ has a significant number of edges of each direction (see \cref{lm:layouts}\cref{lm:layouts-degree&excessW2}, this is inherited from \cref{lm:cleaning}\cref{lm:cleaning-degreeW0}). 
	This means that the degree at $W_2$ is not too large compared to the number of layouts that we will construct, and the in- and outdegree of each vertex in $W_2$ is larger than the in- and outdegree required to satisfy the non-exceptional degree conditions in \cref{prop:absorbingedges}. Thus, we will be able to incorporate almost all of edges at $W_2$ into the layouts. (The edges left over will be covered using \cref{prop:absorbingedges}.)
	\item $\texc(D)$ is significantly larger than the average degree of $D$ (see \cref{lm:layouts}\cref{lm:layouts-U*,lm:layouts-exc>d}, this is inherited from \cref{lm:cleaning}\cref{lm:cleaning-exc>d,lm:cleaning-good}). Also, if $\texc(D)$ is not very large, then the vertices in $W_1$ satisfy some additional degree conditions (see \cref{lm:layouts}\cref{lm:layouts-degree&excessW1}, this is inherited from \cref{lm:cleaning}\cref{lm:cleaning-exc<2d,lm:cleaning-degreeW0}). Moreover, $D[W]$ is empty (see \cref{lm:layouts,lm:cleaning-W}, this is inherited from \cref{lm:cleaning}\cref{lm:cleaning-W}). As discussed in \cref{sec:cleaninglm}, these conditions will facilitate the construction of the layouts. 
\end{itemize}
If we assume the above conditions, then there exist layouts $(L_1,F_1), \dots, (L_\ell, F_\ell)$ on $V(D)$ which satisfy the following properties, where $L$ denotes the multiset $L\coloneqq \bigcup_{i\in [\ell]}L_i$.
\begin{itemize}[--]
	\item $(L_1,F_1), \dots, (L_\ell, F_\ell)$ are $W$-exceptional. As discussed in \cref{sec:W}, $D$ may not be an almost regular robust outexpander and so we will need to apply \cref{lm:approxdecomp} with $D[V']$ playing the role of $D$. The concept of $W$-exceptional layouts will enable us to incorporate the edges incident to $W$ into the approximate decomposition.
	\item Let $\cF$ consist of all the non-absorbing edges of $D$ which are incident to $W$ (i.e.\ $\cF\coloneqq E_W(D)\setminus (A^+\cup A^-)$) and denote $D'\coloneqq D\setminus \cF$. Then, $(L_1,F_1), \dots, (L_\ell, F_\ell)$ are $U^*$-path consistent with respect to $(D',\cF)$. As discussed in \cref{sec:pathL}, this will ensure that the spanning configurations obtained in the approximate decomposition form a partial path decomposition which does not have any endpoint in $U^0(D)\setminus U^*$ (recall \cref{fact:layouts}). Moreover, the definition of path consistency implies that $F_1,\dots, F_\ell\subseteq \cF\subseteq D\setminus (A^+\cup A^-)$. This will ensure that none of the absorbing edges will be covered during the approximate decomposition (recall from \cref{sec:A} that these edges are reserved for the final step of the decomposition).
	\item The number of layouts is bounded away from the density of $D$ (see \cref{lm:layouts}\cref{lm:layouts-ell}). This is needed for applying \cref{lm:approxdecomp}.
	\item The number of (non-trivial) paths in $L$ is precisely $\texc(D)-r$, where $r$ is the value from \cref{prop:absorbingedges} (see \cref{lm:layouts}\cref{lm:cleaning-good}, in the proof of \cref{thm:main} we will apply \cref{prop:absorbingedges} with $\lceil\eta n\rceil$ playing the role of $r$). By \cref{fact:layouts}, this means that the partial path decomposition obtained in the approximate decomposition step will consist of $\texc(D)-r$ paths. The leftover will then be decomposed into $r$ paths with \cref{prop:absorbingedges} and so, overall, we will obtain a path decomposition of $D$ of size $\texc(D)$, as desired.
	\item Each layout is small (see \cref{lm:layouts}\cref{lm:layouts-size}). This is needed for the approximate decomposition (see \cref{lm:approxdecomp}\cref{lm:approxdecomp-sizeLi}).
	\item Each vertex in $V'$ is included in few of the layouts (see \cref{lm:layouts}\cref{lm:layouts-V'}). This is needed for the approximate decomposition (see \cref{lm:approxdecomp}\cref{lm:approxdecomp-degreev}). 
	\item Each non-absorbing edge incident to $W_1$ is incorporated as a fixed edge into precisely one of the layouts (see \cref{lm:layouts}\cref{lm:layouts-degreeW1}). This implies that after the approximate decomposition, the only remaining edges incident to $W_1$ will be the absorbing edges. This is precisely what we need for applying \cref{prop:absorbingedges} with $W_1$ playing the role of the exceptional set.
	\item All but $r$ inedges and $r$ outedges at each vertex in $W_2$ are included as fixed edges in $L$ (see \cref{lm:layouts}\cref{lm:layouts-degreeW2}). This implies that, after the approximate decomposition, each vertex in $W_2$ will have both its in and outdegree equal $r$, i.e.\ each vertex in $W_2$ will have excess $0$ and total degree $2r$ in the leftover. This means that we will be able to apply \cref{prop:absorbingedges} with $W_2\subseteq X^0$.
	\item  Let $X^+\subseteq (U^+(D)\cup U^*)\cap V'$ and $X^-\subseteq (U^-(D)\cup U^*)\cap V'$. Then, \cref{fact:layouts} and \cref{lm:layouts}\cref{lm:layouts-degreeV'} imply that, for each $v\in V'$, the number of $L_i$ which include $v$ as an isolated vertex and the outdegree of $v$ in $L$ together have precisely the value such that, after the approximate decomposition, the outdegree at $v$ will be $r$ if $v\notin X^-$ and $r-1$ if $v\in X^-$. The analogous statement holds for $X^+$.
	Thus, the vertices in $V'$ will satisfy the degree conditions of \cref{prop:absorbingedges} with $X^+\setminus X^-$, $X^-\setminus X^+$, and $X^+\cap X^-$ playing the roles of $X^+$, $X^-$, and $X^*$, and with $V'\setminus (X^+\cup X^-)\subseteq X^0$.	
\end{itemize}}

\begin{lm}[Layout construction]\label{lm:layouts}
	Let $0<\frac{1}{n}\ll \varepsilon\ll \eta\ll 1$ and $d\in \mathbb{N}$. Let~$D$ be a\NEW{n} oriented graph on a vertex set~$V$ of size~$n$ such that the following hold. 
	\begin{enumerate}
		\item Let $W_1\cup W_2\cup V'$ be a partition of~$V$. Denote $W\coloneqq W_1\cup W_2$. Suppose \NEW{that} $|W|\leq \varepsilon n$ and $E(D[W])=\emptyset$.\label{lm:layouts-W}
		\item Let $U^*\subseteq U^0(D)\setminus W$ be such that $|U^*|=\texc(D)-\exc(D)$. Moreover, each~$v\in U^*$ satisfies $d_D^+(v)=d_D^-(v)\leq \texc(D)-1$.\label{lm:layouts-U*}
		\item Let~$A^+$ and~$A^-$ be absorbing sets of~$(W_1,V')$-starting and~$(V',W_1)$-ending edges for~$D$, respectively, and denote $A\coloneqq A^+\cup A^-$. Suppose $X^\pm \subseteq (U^\pm(D)\cup U^*) \setminus W$ are such that $|A^\pm|+|X^\pm|=\lceil\eta n\rceil$.
		Define $\phi^\pm: V \longrightarrow \{0,1\}$ by
		\begin{equation*}
			\phi^\pm(v)\coloneqq
			\begin{cases}
				1 & \text{if }v\in X^\pm,\\
				0 & \text{otherwise}.\\
			\end{cases}
		\end{equation*} \label{lm:layouts-A}
		\item $d\geq \eta n$.\label{lm:layouts-d}
		\item $\texc(D)\geq d+\lceil\eta n\rceil$.\label{lm:layouts-exc>d}
		\item For all~$v\in W_1$, $10\varepsilon n\leq d_{D\setminus A}(v)\leq 2d -\lceil\eta n\rceil$. 
		Moreover, if $\texc(D)<2d+\lceil\eta n\rceil$, then, for each~$v\in W_1$, one of the following holds:
		\begin{itemize}[--]
			\item $|\exc_D(v)|=d_D(v)$; or
			\item $d_D^{\min}(v)\geq \eta n$ and $|\exc_{D\setminus A}(v)|\leq \lceil\eta n\rceil$; or
			\item $d_D^{\min}(v)\geq \eta n$ and $d_A(v)=\lceil\eta n\rceil$.
		\end{itemize}
		\label{lm:layouts-exc<2d}\label{lm:layouts-degree&excessW1}
		\item For all~$v\in W_2$, $d_D^{\min}(v)\geq \lceil\eta n\rceil$ and $d_D(v)\leq  2d+2\lceil\eta n\rceil$.\label{lm:layouts-degree&excessW2}
		\item For all~$v\in V'$, $2d +2\lceil\eta n\rceil-\varepsilon n \leq d_D(v)\leq 2d +2\lceil\eta n\rceil$ and $|\exc_D(v)|\leq \varepsilon n$.\label{lm:layouts-degree&excessV'}
	\end{enumerate}
	Let $\cF\coloneqq E_W(D)\setminus A$ and $D'\coloneqq D\setminus \cF$. Then, there exist $\ell\in \mathbb{N}$ and~$W$-exceptional layouts $(L_1, F_1), \dots, (L_\ell, F_\ell)$ which are~$U^*$-path consistent with respect to~$(D',\cF)$%
	\COMMENT{Note that we have~$D'$ here rather than~$D$ since otherwise, following our convention from \cref{sec:notation},~$D\cup \cF$ is now a multidigraph which contains two occurrences of each edge in~$\cF$. This not good for the path consistent definition.}
	and satisfy the following, where~$L$ is the multiset defined by \NEW{$L\coloneqq \bigcup_{i\in [\ell]}L_i$}\OLD{$L\coloneqq L_1\cup \dots \cup L_\ell$}.
	\begin{enumerate}[label=\upshape(\roman*)]
		\item $d\leq \ell \leq d+\sqrt{\varepsilon}n$.\label{lm:layouts-ell}
		\item $L$ contains exactly $\texc(D)-\lceil\eta n\rceil$ non-trivial paths.\label{lm:layouts-good}
		\item For all~$v\in W_1$, $d_L^\pm(v)= d^{\pm}_{\mathcal{F}}(v) = d_{D\setminus A}^\pm(v)$.\label{lm:layouts-degreeW1}
		\item For all~$v\in W_2$, $d_L^\pm(v) = d_{\mathcal{F}}^\pm(v)-\lceil\eta n\rceil=d_D^\pm(v)-\lceil\eta n\rceil$.\label{lm:layouts-degreeW2}
		\item For all~$v\in V'$, $d_L^\pm(v) = d_D^\pm(v) - |\{i\in [\ell]\mid v\notin V(L_i)\}|-\lceil\eta n\rceil + \phi^\mp(v)$.\label{lm:layouts-degreeV'}
		\item For all~$i\in [\ell]$, $|V(L_i)|, |E(L_i)|\leq 3\varepsilon^{\frac{1}{3}} n$.\label{lm:layouts-size}
		\item For each~$v\in V'$, $d_L(v)\leq 8\varepsilon n$ and there exist at most~$3\sqrt{\varepsilon}n$ indices~$i\in [\ell]$ such that~$v\in V(L_i)$.\label{lm:layouts-V'} 
	\end{enumerate}
\end{lm}

\NEW{We now motivate the expression appearing in \cref{lm:layouts-degreeV'}. Let $v\in V'$.
	Recall from \cref{fact:layouts} that $d_{\hL}^+(v)+|\{i\in [\ell]\mid v\notin V(\hL_i)\}|$ is precisely the outdegree of $v$ in a set of spanning configurations of shapes $(\hL_1, \hF_1), \dots, (\hL_\ell,\hL_\ell)$. Moreover, as seen in the proof of \cref{thm:main} (see also the explanatory paragraph before the statement of \cref{lm:layouts}), $\lceil\eta n\rceil-\phi^-(v)$ is precisely the leftover outdegree of $v$ that we aim for after the approximate decomposition step.}

\NEW{We now explain in more detail why \cref{lm:cleaning}\cref{lm:cleaning-good} is necessary. This is because a general oriented graph $D$ may not contain sufficiently many zero-excess vertices which satisfy the degree condition of \cref{lm:layouts}\cref{lm:layouts-U*}. For example, if $D$ is regular, then $\texc(D)-\exc(D)>0$ but every $v\in V(D)$ satisfies $d_D^+(v)=d_D^-(v)=\Delta^0(D)=\texc(D)$. This example also illustrates why the degree condition of \cref{lm:layouts}\cref{lm:layouts-U*} is necessary. Indeed, suppose for a contradiction that \cref{lm:layouts} also holds if we omit the ``moreover part" of \cref{lm:layouts}\cref{lm:layouts-U*}. Let $T$ be a regular tournament on $n$ vertices. Let $d\coloneqq \frac{n-1}{2}-\lceil\eta n\rceil$. Let $V'\coloneqq V(T)$ and $W_1\coloneqq W_2\coloneqq A^+\coloneqq A^-\coloneqq \emptyset$. Let $U^*\subseteq V(T)$ satisfy $|U^*|=\frac{n-1}{2}$ and $X^+,X^-\subseteq U^*$ satisfy $|X^+|=\lceil\eta n\rceil=|X^-|$. Then, one can easily verify that \cref{lm:layouts}\cref{lm:layouts-A,lm:layouts-W,lm:layouts-U*,lm:layouts-d,lm:layouts-exc<2d,lm:layouts-exc>d,lm:layouts-degree&excessV',lm:layouts-degree&excessW1,lm:layouts-degree&excessW2} are all fully satisfied except for the ``moreover part" of \cref{lm:layouts}\cref{lm:layouts-U*} and so, by assumption, there exist layouts as in \cref{lm:layouts}. As discussed earlier, this implies that we can use \cref{lm:approxdecomp} to construct a partial path decomposition $\cP$ of $T$ of size $|\cP|=\texc(T)-\lceil\eta n\rceil$ such that $T\setminus \cP$ satisfies the degree conditions of \cref{prop:absorbingedges} with $r\coloneqq \lceil\eta n\rceil$. Thus, \cref{prop:absorbingedges} implies that there exists a path decomposition $\cP'$ of $T\setminus \cP$ of size $\lceil\eta n\rceil$. But, this means that $\cP\cup \cP'$ is a path decomposition of $T$ of size $\texc(T)$. This contradicts \cref{thm:annoyingT}.}

\NEW{The next \lcnamecref{prop:sizeU*good} states that, after the cleaning step, we will be able to find a set $U^*$ of candidates for path endpoints which satisfies \cref{lm:layouts}\cref{lm:layouts-U*}. Let $T\notin \cT_{\rm excep}$. By \cref{prop:sizeU}, there exists $U^*\subseteq U^0(T)$ satisfying $|U^*|=\texc(T)-\exc(T)$. Let $\cP$ be the good partial path decomposition obtained by applying \cref{lm:cleaning}. Denote $D\coloneqq T\setminus \cP$. We now aim to apply \cref{lm:layouts} to $D$ and so we need a new set $U^{**}\subseteq U^0(D)$ which satisfies \cref{lm:layouts}\cref{lm:layouts-U*}. (We cannot use the original $U^*$ since some of the vertices in $U^*$ may have been used as endpoints in $\cP$ and so may have non-zero excess in $D$.) By \cref{prop:sizeU}, there exists $U^{**}\subseteq U^0(T)$ satisfying $|U^{**}|=\texc(D)-\exc(D)$. However, there is no guarantee that the vertices in $U^{**}$ satisfy the desired degree conditions. But by \cref{lm:cleaning}\cref{lm:cleaning-good}, we know that all the vertices in $U^*$ which have not been used as endpoints in $\cP$ have the correct degree conditions for \cref{lm:layouts}\cref{lm:layouts-U*}. Thus, we would like to take $U^{**}\subseteq U^*\setminus (V^+(\cP)\cup V^-(\cP))$ and so we would like $U^*\setminus (V^+(\cP)\cup V^-(\cP))$ to contain at least $\texc(D)-\exc(D)$ vertices of $U^0(D)$. \Cref{prop:sizeU*good} states that this is the case.} 

\begin{prop}\label{prop:sizeU*good}
	Let~$D$ be a digraph and $W,V'\subseteq V(D)$ be disjoint. Suppose $A^+, A^-\subseteq E(D)$ are absorbing sets of~$(W,V')$-starting and~$(V',W)$-ending edges for~$D$. Denote $A\coloneqq A^+\cup A^-$. Let~$U^*\subseteq U^0(D)$ satisfy~$|U^*|= \texc(D)-\exc(D)$.
	Suppose~$\cP$ is a 
	good
	$(U^*,W,A)$-partial path decomposition of~$D$. Let~$U^{**}\coloneqq U^*\setminus (V^+(\cP)\cup V^-(\cP))$. Then, \NEW{$U^{**}\subseteq U^0(D\setminus \cP)$ and} $|U^{**}|=\texc(D\setminus \cP)-\exc(D\setminus \cP)$.
\end{prop}

\NEW{Recall that a $(U^*,W,A)$-partial path decomposition of~$D$ was defined in \cref{def:U*WApartialpathdecomp} and its goodness before \cref{fact:combinefull}.}

\begin{proof}[Proof of \cref{prop:sizeU*good}]
	\NEW{Since the vertices in $U^{**}$ have not been used as endpoints in $\cP$, they still have excess $0$ in $D\setminus \cP$.%
		\COMMENT{More formally, each $v\in U^{**}$ satisfies
			\begin{align*}
				\exc_{D\setminus 	\cP}(v)&\stackrel{\text{\eqmakebox[U][c]{}}}{=} d_{D\setminus\cP}^+(v)-d_{D\setminus\cP}^-(v)= (d_D^+(v)-d_D^-(v))-(d_\cP^+(v)-d_\cP^-(v))\\
				&\stackrel{\text{\eqmakebox[U][c]{$v\in 	U^0(D)$}}}{=}0-(d_\cP^+(v)-d_\cP^-(v))
				\stackrel{v\notin V^+(\cP)\cup V^-(\cP)}{=}0-0=0.
			\end{align*}
			Thus, $U^{**}\subseteq U^0(D\setminus \cP)$, as desired.}}

	\NEW{By definition of a $(U^*, W, A)$-partial path decomposition, no path in $\cP$ has an endpoint in $U^0(D)\setminus U^*$. Therefore,
		\begin{equation}\label{eq:sizeU*good}
			U^0(D)\cap (V^+(\cP)\cup V^-(\cP))= U^*\setminus U^{**}
		\end{equation}
		and so the fact that $\cP$ is a good partial path decomposition implies that
		\begin{align*}\label{eq:sizeU*good-texc}
			\texc(D\setminus \cP)&\stackrel{\text{\eqmakebox[U2][c]{}}}{=}\texc(D)-|\cP|\\
			&\stackrel{\text{\eqmakebox[U2][c]{\cref{prop:excP}}}}{=}(\texc(D)-\exc(D))+\exc(D\setminus \cP)-|U^0(D)\cap (V^+(\cP)\cup V^+(\cP))|\\
			&\stackrel{\text{\eqmakebox[U2][c]{\cref{eq:sizeU*good}}}}{=}|U^*|+\exc(D\setminus \cP)-|U^*\setminus U^{**}|=\exc(D\setminus \cP)+|U^{**}|,
		\end{align*}
		as desired.}
	\OLD{\textbf{Old proof, which uses the auxiliary excess function (not defined at this stage):} We have
		\begin{align*}
			|U^{**}|&\stackrel{\text{\eqmakebox[prop:sizeU*good][c]{\text{\cref{eq:texcA}}}}}{=}\texc_{U^{**}}^+(D\setminus \cP)-\exc(D\setminus \cP)+|A^+|
			\stackrel{\text{\cref{fact:specialdecomp}}}{=}(\texc_{U^*}^+(D)-|\cP|)-\exc(D\setminus \cP)+|A^+|\\
			&\stackrel{\text{\eqmakebox[prop:sizeU*good][c]{\text{\cref{eq:texcA}}}}}{=}(\texc(D)-|A^+|)-|\cP|-\exc(D\setminus \cP)+|A^+|= \texc(D\setminus \cP)-\exc(D\setminus \cP).\\
	\end{align*}}
\end{proof}

%% file: Proof_of_Main_Theorem.tex
	\onlyinsubfile{
		\setcounter{section}{9}
\section{Deriving Theorem \ref{thm:main}}}

\NEW{In this section, we assume that \cref{lm:cleaning,lm:layouts} hold and derive \cref{thm:main}.}
We will proceed as follows. In \cref{step:absorbing}, we select absorbing edges (if they are required).
In \cref{step:cleaning}, we clean up~$T$ by removing a small number of paths using \cref{lm:cleaning}. In \cref{step:approx}, \NEW{we first apply \cref{lm:layouts} to obtain approximate layouts and then apply \cref{lm:approxdecomp} to obtain an approximate decomposition of $T$ based on these layouts}\OLD{we obtain an approximate decomposition of $T$ by \cref{lm:approxdecomp,lm:layouts}}. Finally, in \cref{step:fulldecomp}, we apply \cref{prop:absorbingedges} to decompose the leftover.

\begin{proof}[Proof of \cref{thm:main}] 
	Assume without loss of generality that~$\beta\leq 1$. Fix additional constants such that 
	$0<\frac{1}{n_0}\ll \varepsilon \ll\alpha_1\ll\alpha_2\ll\eta \ll\beta\leq 1$.
	Let~$T\notin\cT_{\rm excep}$ be a tournament on~$n\geq n_0$ vertices satisfying \cref{thm:main-largeexc} or \cref{thm:main-largeU}. By \cref{thm:evenlarge}\cref{thm:even-verylarge}, we may assume that~$\exc(T)\leq \varepsilon^2n^2$. Denote~$V\coloneqq V(T)$.

	\NEW{Recall that $N^\pm(T)=|U^\pm(T)|+\texc(T)-\exc(T)$.}
	If both~$N^\pm(T)\geq \alpha_1 n$, then redefine~$\eta\coloneqq \alpha_1^2$. Suppose not. If both~$N^\pm(T)\leq \alpha_2 n$, then redefine~$\varepsilon\coloneqq \alpha_2$.
	Otherwise, there exists~$\diamond\in \{+,-\}$ such that~$N^\diamond(T)\leq\alpha_1 n$ and~$\circ\in \{+,-\}\setminus\{\diamond\}$ satisfies~$N^\circ(T)\geq \alpha_2n$ and we redefine~$\varepsilon\coloneqq \alpha_1$ and~$\eta\coloneqq \alpha_2^2$.
	Thus, we have defined constants such that \[0<\frac{1}{n}\ll \varepsilon \ll\eta \ll\beta\leq 1,\]
	$\exc(T)\leq \varepsilon^2 n^2$, and, for each~$\diamond\in \{+,-\}$, either~$N^\diamond(T)\geq \sqrt{\eta} n$ or~$N^\diamond(T)\leq \varepsilon n$.
	Define additional constants such that 
	\[0<\frac{1}{n}\ll \varepsilon \ll \nu\ll \tau\ll \gamma \ll \eta\ll\beta\leq 1.\]
	\begin{steps}
		\item \textbf{Choosing absorbing edges.}\label{step:absorbing}
		\NEW{We start by partitioning $V$ into $V'$ and $W$, and selecting a $(W,V')$-absorbing set $A$.}
		Let $r\coloneqq \lceil \eta n\rceil$.
		
		\begin{claim}\label{claim:A}
			\NEW{There exist a partition $W\cup V'$ of $V$ and absorbing sets $A^+,A^-\subseteq E(T)$ of $(W,V')$-starting/$(V',W)$-ending edges for~$T$ such that the following hold, where $A\coloneqq A^+\cup A^-$, $W_A^\pm \coloneqq V(A^\pm)\cap W$, and $W_A\coloneqq V(A)\cap W$.
				\begin{enumerate}[label=\upshape(\roman*)]
					\item $W\cap U^0(T)=\emptyset$.\label{claim:A-U0}
					\item $|W|\leq 4\varepsilon n$.\label{claim:A-W}
					\item For each $v\in V'$, $|\exc_T(v)|\leq \varepsilon n$.\label{claim:A-V'}
					\item Let $\diamond \in \{+,-\}$. If $N^\pm(T)\geq \sqrt{\eta}n$, then $A^\diamond=\emptyset$, otherwise $|A^\diamond|=r$.\label{claim:A-A}
					\item Let $\diamond\in \{+,-\}$. If $|W_A^\diamond|\geq 2$, then $\exc_T^\diamond(v)<r$ for each~$v\in V$.\label{claim:A-A2}
					\item Let $\diamond\in \{+,-\}$. If $|W_A^\diamond|=1$, then $\exc_T^\diamond(v)\leq \exc_T^\diamond(w)$ for each~$v\in V$ and $w\in W_A^\diamond$.\label{claim:A-A1}
					\item If $W_A\neq \emptyset$, then $\texc(T)\geq \frac{n}{2}+10\eta n$.\label{claim:A-A0}
				\end{enumerate}}
		\end{claim}
	
		\begin{proofclaim}
			\NEW{First, we choose suitable sets of endpoints for our absorbing edges.} Let $\diamond\in \{+,-\}$ and define~$W_A^\diamond$ as follows.  \NEW{First, suppose that $N^\diamond(T)\geq \sqrt{\eta}n$. Since we will not need any absorbing edges in that case, we let $W_A^\diamond=\emptyset$. Observe that \cref{claim:A-A1,claim:A-A2} hold for $\diamond$.}\OLD{If~$N^\diamond(T)\geq \sqrt{\eta}n$, then let~$W_A^\diamond=\emptyset$.}
			
			\NEW{Now suppose that $N^\diamond(T)< \sqrt{\eta}n$. By construction, we have $N^\diamond(T)\leq \varepsilon n$. We will need $r$ absorbing edges, so we need to choose a set $W_A^\diamond$ which concentrates a sufficiently large amount of excess but also satisfies \cref{claim:A-A1,claim:A-A2}.}%
			\OLD{Otherwise, by construction,~$N^\diamond(T)\leq \varepsilon n$.}
			By \cref{prop:texcT},
			$\exc(T)\geq \texc(T)-N^\diamond(T)\geq \frac{n}{2}-\varepsilon n\geq r$ \NEW{and so we can} let~$W_A ^\diamond\subseteq U^\diamond(T)$ be a smallest set such that~$\exc_T^\diamond(W_A^\diamond) \ge r$. We further assume that, subject to this,~$\exc_T^\diamond(W_A^\diamond)$ is maximum. Note that \begin{equation}\label{eq:absorbing-size}
				\NEW{|W_A^\diamond|\leq}|U^\diamond(T)|\leq N^\diamond(T)\leq \varepsilon n.
			\end{equation}\OLD{Thus, $|W_A^\diamond|\leq \varepsilon n$.}%
			\NEW{We verify that \cref{claim:A-A0,claim:A-A1,claim:A-A2} are satisfied for $\diamond$. If $|W_A^\diamond|\geq 2$, then the minimality of $|W_A^\diamond|$ implies that each $v\in V$ satisfies $\exc_T^\diamond(v)< r$ and so \cref{claim:A-A2} holds.
			If $|W_A^\diamond|=1$, then the maximality of $\exc_T^\diamond(W_A^\diamond)$ implies that each $v\in V$ and $w\in W_A^\diamond$ satisfy	$\texc_T^\diamond(v)\leq \texc_T^\diamond(w)$, so \cref{claim:A-A1} holds. By assumption, $N^\diamond(T)\leq \varepsilon n< \beta n$. Thus, \cref{thm:main-largeU} does not hold and so \cref{thm:main-largeexc} implies that $\texc(T)\geq \frac{n}{2}+\beta n\geq \frac{n}{2}+10\eta n$, as desired for \cref{claim:A-A0}.}%	
			\OLD{Observe that, if~$|W_A^\diamond|\geq 2$, then, for each~$v\in V$,~$\exc_T^\diamond(v)< r$, and, if~$|W_A^\diamond|=1$, then, for each~$v\in V$ and~$w\in W_A^\diamond$, $\texc_T^\diamond(v)\leq \texc_T^\diamond(w)$, as desired for \cref{lm:cleaning}\cref{lm:cleaning-A}.}
			
			Let $W_A\coloneqq W_A^+\cup W_A^-$.
			\NEW{Define $W\coloneqq W_A\cup \{v\in V\mid |\exc_T(v)|>\varepsilon n\}$ and $V'\coloneqq V\setminus W$. Note that \cref{claim:A-V'} is satisfied. By construction, $W_A^\pm\subseteq U^\pm(T)$ and so \cref{claim:A-U0} holds.}%
			\OLD{Let \[W_*\coloneqq \{v\in V\mid |\exc_T(v)|>(1-20\eta)n\} \quad \text{and} \quad W_0\coloneqq \{v\in V\setminus W_*\mid |\exc_T(v)|>\varepsilon n\}\cup (W_A\setminus W_*).\]
			For each $\diamond\in \{*,0\}$, let $W_\diamond^\pm\coloneqq W_\diamond\cap U^\pm(T)$.
			Let $W^\pm\coloneqq W_*^\pm\cup W_0^\pm$, $W\coloneqq W^+\cup W^-$, and $V'\coloneqq V\setminus W$.}
			\NEW{By \cref{eq:absorbing-size} and since $\exc(T)\leq \varepsilon^2n^2$, we have
				\begin{equation*}
					|W|\leq |W_A^+|+|W_A^-|+|\{v\in V\mid |\exc_T(v)|>\varepsilon n\}|\leq 2\varepsilon n+2\cdot\frac{\exc(T)}{\varepsilon n}\leq 4\varepsilon n.
				\end{equation*} 
				Thus, \cref{claim:A-W} holds.}\OLD{Note that, by construction, $|W|\leq 4\varepsilon n$. (Here we also use that $\exc(T)\leq \varepsilon^2 n^2$.)}
			
			\NEW{We are now ready to choose the absorbing edges.} If $N^+(T)\geq \sqrt{\eta} n$, then let $A^+\coloneqq \emptyset$; otherwise, let $A^+\subseteq E(T)$ be an absorbing set of~$r$~$(W_A^+, V')$-starting edges for~$T$ ($A^+$ exists since, by construction,~$\exc_T^+(W_A^+)\geq r$ and~$d_T^+(v)\geq \frac{n}{2}$ for each~$v\in W_A^+$). Similarly, let~$A^-\subseteq E(T)$ be an absorbing set of~$(V',W_A^-)$-ending edges for~$T$, of size~$0$ if~$N^-(T)\geq \sqrt{\eta}n$ and~$r$ otherwise.
			\NEW{Thus, \cref{claim:A-A} holds.}
			Let $A\coloneqq A^+\cup A^-$.
			\NEW{By minimality of $|W_A^\pm|$,}\OLD{Note that, by construction,} $V(A^\pm)\cap W=W_A^\pm$.%
			\OLD{In particular, each $v\in W_A^\pm$ satisfies
			\begin{equation*}
				d_A^\pm(v)\geq 1.
			\end{equation*}}
			\NEW{This completes the proof.}
		\end{proofclaim}
		\NEW{Let $W_*$ consists of all the vertices $w_*\in W$ for which $|\exc_T(w_*)|>(1-20\eta)n$ and let $W_0$ consist of all the  vertices $w_0\in W$ for which $|\exc_T(w_0)|\leq (1-20\eta)n$.}
		
		\item \textbf{Cleaning.}\label{step:cleaning}
		\NEW{By \cref{claim:A}, \cref{lm:cleaning}\cref{lm:cleaning-A,lm:cleaning-defW} are satisfied with $4\varepsilon$ playing the role of $\varepsilon$.}
		\NEW{By \cref{prop:sizeU}, there exists $U_1^*\subseteq U^0(T)$ which satisfies $|U_1^*|=\texc(T)-\exc(T)$.}\OLD{Let $U_1^*\subseteq U^0(T)$ be such that $|U_1^*|=\texc(T)-\exc(T)$ (this is possible by \cref{prop:sizeU}).}
		\NEW{Then, \cref{lm:cleaning}\cref{lm:cleaning-U*} holds with $U_1^*$ playing the role of $U^*$.
		By \cref{claim:A-U0},
		\begin{equation}\label{eq:cleaning-U1*}
			W\cap U_1^*=\emptyset.
		\end{equation}
		(This will be needed in \cref{step:approx}.)}
		\OLD{Note that $W\cap U_1^*=\emptyset$ since $W\cap U^0(T)=\emptyset$. Moreover, if~$W_A\neq \emptyset$, then, by \cref{step:absorbing}, \cref{thm:main-largeU} does not hold and so $\texc(T)\geq \frac{n}{2}+\beta n\geq \frac{n}{2}+10\eta n$. Thus, the ``moreover part" of \cref{lm:cleaning}\cref{lm:cleaning-A} is satisfied.}
		
		Apply \cref{lm:cleaning} with $U_1^*$ and $4\varepsilon$ playing the roles of $U^*$ and $\varepsilon$ to obtain $d\in \mathbb{N}$ and $\cP_1\subseteq T$ such that the following are satisfied, where $D_1\coloneqq T\setminus \cP_1$.
		\begin{enumerate}[label=\upshape(\greek*)]
			\item $\cP_1$ is a good $(U_1^*, W, A)$-partial path decomposition of~$T$. In particular,~$\cP_1$ is consistent with $A^+$ and $A^-$ and so, by \cref{fact:absorbingedges}, $A^+$ and $A^-$ are absorbing sets of~$(W_A^+,V')$-starting and~$(V',W_A^-)$-ending edges for~$D_1$.\label{thm:cleaning-A}
			\item $\left\lceil\frac{n}{2}\right\rceil-10\eta n\leq d\leq \left\lceil\frac{n}{2}\right\rceil-\eta n$. \label{thm:cleaning-d}
			\item Each $v\in U_1^*\setminus (V^+(\cP_1)\cup V^-(\cP_1))$ satisfies $d_{D_1}^+(v)=d_{D_1}^-(v)\leq \texc(D_1)-1$.\label{thm:cleaning-good}
			\item $E(D_1[W])=\emptyset$.\label{thm:cleaning-W}
			\item $N^\pm(T)-N^\pm(D_1)\leq 89\eta n$.\label{thm:cleaning-sizeU}
			\item $\texc(D_1)\geq d+\lceil\eta n\rceil$.\label{thm:cleaning-exc>d}
			\item If $\texc(D_1)< 2d+\lceil\eta n\rceil$, then each~$w\in W_*$ satisfies $|\exc_{D_1}(w)|=d_{D_1}(w)$.\label{thm:cleaning-exc<2d}
			\item For each~$v\in W_*\cup W_A$, $2d-3\sqrt{\eta}n\leq d_{D_1}(v)\leq 2d-2\lceil\eta n\rceil$.\label{thm:cleaning-degreeW*}
			\item For each~$v\in W_0$, $2d+2\lceil\eta n\rceil-4\sqrt{\eta}n\leq d_{D_1}(v)\leq 2d+2\lceil\eta n\rceil$ and $d_{D_1}^{\min}(v)\geq \lceil\eta n\rceil$.\label{thm:cleaning-degreeW0}
			\item For each~$v\in V'$, $2d+2\lceil\eta n\rceil-18\sqrt{\varepsilon}n\leq d_{D_1}(v)\leq 2d+2\lceil\eta n\rceil$ and $|\exc_{D_1}(v)|\leq \varepsilon n$.\label{thm:cleaning-degreeV'}
		\end{enumerate}
		(The final part of \cref{thm:cleaning-degreeV'} follows from the facts that $\cP_1$ is a \OLD{$(U_1^*,W,A)$-}partial path decomposition of~$T$ and $|\exc_T(v)|\leq \varepsilon n$ for each~$v\in V'$.
		\NEW{Indeed, if $v\in V'\setminus U^0(T)$, then \cref{def:partialpathdecomp-exc} implies that $\cP_1$ contains at most $\exc_T^+(v)$ paths which start at $v$ and at most $\exc_T^-(v)$ paths which end at $v$. Thus, each $v\in V'\setminus U^0(T)$ satisfies $|\exc_{D_1}(v)|\leq |\exc_T(v)|\leq \varepsilon n$. Moreover, \cref{def:partialpathdecomp-U0} implies that each $v\in V'\cap U^0(T)$ is the starting point of at most one path in $\cP_1$ and the ending point of at most one path in $\cP_1$. Thus, each $v\in V'\cap U^0(T)$ satisfies $|\exc_{D_1}(v)|\leq 1\leq \varepsilon n$.})
		
		\item \textbf{Approximate decomposition.}\label{step:approx}
		\NEW{We will approximately decompose $D_1$ as follows. First, we will apply \cref{lm:layouts} to construct $W$-exceptional layouts on $V$. These layouts will then be transformed into auxiliary layouts on $V\setminus W$ via \cref{def:layoutV'}. 
		We will then apply \cref{lm:approxdecomp} to these auxiliary layouts to approximately decompose $D_1[V']$ into auxiliary spanning configurations on $V'$. Finally, we will use \cref{fact:layoutV'} to transform these auxiliary spanning configurations into $W$-exceptional spanning configurations on $V$. By \cref{fact:layouts,lm:layouts}, these will induce a good partial path decomposition of $D_1$ which covers almost all the edges of $D_1$.}
		\OLD{First, we will apply \cref{lm:layouts} to construct layouts for the approximate decomposition.}
		
		\NEW{First, we ensure that all the prerequisites of \cref{lm:layouts} are satisfied.} Let $U_2^*\coloneqq U_1^*\setminus (V^+(\cP_1)\cup V^-(\cP_1))$ and observe that, by \cref{thm:cleaning-A} and \cref{prop:sizeU*good}, \begin{equation}\label{eq:approx-U2*}
			|U_2^*|=\texc(D_1)-\exc(D_1).
		\end{equation}
		
		\begin{claim}\label{claim:Xpm}
			\NEW{There exist $X^\pm\subseteq(U^\pm(D_1)\cup U_2^*)\setminus W$ which satisfy $|X^\pm|=r-|A^\pm|$.}
		\end{claim}
	
		\NEW{\begin{proofclaim}
			Let $\diamond\in \{+,-\}$. First, suppose that $N^\diamond(T)\leq \varepsilon n$. Then, \cref{claim:A-A} implies that $|A^\diamond|=r$ and so we can let $X^\diamond\coloneqq \emptyset$. We may therefore assume that $N^\diamond(T)> \varepsilon n$. By construction, $N^\diamond(T)\geq \sqrt{\eta}n$ and so
			\[|U^\diamond(D_1)\cup U_2^*|=N^\diamond(D_1)\stackrel{\text{\cref{thm:cleaning-sizeU}}}{\geq} \sqrt{\eta}n-89\eta n\geq r + |W|,\]
			as desired.
		\end{proofclaim}}
		
		\OLD{Let $X^\pm\subseteq(U^\pm(D_1)\cup U_2^*)\setminus W$ be such that $|X^\pm|=r-|A^\pm|$. This is possible since, for each $\diamond\in \{+,-\}$, by construction, if $N^\diamond(T)\leq \varepsilon n$, then $|A^\diamond|=r$; otherwise, $N^\diamond(T)\geq \sqrt{\eta}n$ and so, by \cref{thm:cleaning-sizeU}, $|U^\diamond(D_1)\cup U_2^*|=N^\diamond(D_1)\geq \sqrt{\eta}n-89\eta n\geq r + |W|$.}
		Define $\phi^\pm: V\longrightarrow \{0,1\}$ by 
		\[\phi^\pm(v)\coloneqq 
		\begin{cases}
			1 & \text{if }v\in X^\pm,\\
			0 & \text{otherwise.}
		\end{cases}\]
		Denote $\cF\coloneqq E_W(D_1)\setminus A$ and $D_1'\coloneqq D_1\setminus \cF$. Let $W_1\coloneqq W_*\cup W_A$ and $W_2\coloneqq W_0\setminus W_A$.
		
		We now verify that \cref{lm:layouts}\cref{lm:layouts-A,lm:layouts-W,lm:layouts-U*,lm:layouts-d,lm:layouts-exc>d,lm:layouts-degree&excessW2,lm:layouts-degree&excessV',lm:layouts-degree&excessW1} hold with $D_1, U_2^*$, and $18\sqrt{\varepsilon}$ playing the roles of $D, U^*$, and $\varepsilon$. 
		\cref{lm:layouts}\cref{lm:layouts-W} \NEW{follows from \cref{thm:cleaning-W} and \cref{lm:cleaning}\cref{lm:cleaning-defW}}\OLD{holds by \cref{step:absorbing} and \cref{thm:cleaning-W}}. 
		\cref{lm:layouts}\cref{lm:layouts-U*} holds by \cref{thm:cleaning-good}, \NEW{\cref{eq:cleaning-U1*}, and \cref{eq:approx-U2*}}\OLD{and the fact that $U_2^*\subseteq U_1^*\subseteq V'$}.
		\cref{lm:layouts}\cref{lm:layouts-A} holds by \NEW{\cref{claim:Xpm}}\OLD{construction} and \cref{thm:cleaning-A}.
		\cref{lm:layouts}\cref{lm:layouts-d} follows from \cref{thm:cleaning-d} and \cref{lm:layouts}\cref{lm:layouts-exc>d} holds by \cref{thm:cleaning-exc>d}. \cref{lm:layouts}\cref{lm:layouts-degree&excessW2} and \cref{lm:layouts-degree&excessV'} as well as the first part of \cref{lm:layouts}\cref{lm:layouts-degree&excessW1} follow from \cref{thm:cleaning-degreeV',thm:cleaning-degreeW*,thm:cleaning-degreeW0}.
		Finally, we show that the ``moreover part" of \cref{lm:layouts}\cref{lm:layouts-exc<2d} holds. \NEW{Suppose that $\texc(D_1)<2d+\lceil\eta n\rceil$ and $v\in W$. If $v\in W_*$, then \cref{thm:cleaning-exc<2d} implies that $|\exc_{D_1}(v)|=d_{D_1}(v)$. We may therefore assume that $v\in W_1\setminus W_*\subseteq W_A\cap W_0$. Then, \cref{thm:cleaning-degreeW0} implies that $d_{D_1}^{\min}(v)\geq \lceil\eta n\rceil$. Thus, it is enough to show that $|\exc_{D_1\setminus A}(v)|\leq r$ or $d_A(v)=r$. Suppose without loss of generality that $v\in W_A^+$, i.e.\ $v\in U^+(T)$ by \cref{def:A} (similar arguments hold if $v\in W_A^-$). If $|W_A^+|\geq 2$, then,
		\[|\exc_{D_1\setminus A}(v)|\stackrel{\text{\cref{def:A},\cref{thm:cleaning-A}}}{\leq} |\exc_T(v)|\stackrel{\text{\cref{claim:A-A2}}}{\leq}r.\]
		If $W_A^+=\{v\}$, then \cref{claim:A-A1} implies that $d_{A^+}(v)=r$ and \cref{def:A} implies that $d_{A^-}(v)=0$, so $d_A(v)=r$. Therefore, \cref{lm:layouts}\cref{lm:layouts-exc<2d} holds.}\OLD{By \cref{thm:cleaning-exc<2d}, if $\texc(D_1)<2d+\lceil\eta n\rceil$, then each~$v\in W_*$ satisfies $|\exc_{D_1}(v)|=d_{D_1}(v)$. 
		Moreover, by \cref{thm:cleaning-degreeW0}, each~$v\in W_A\setminus W_*\subseteq W_0$ satisfies $d_{D_1}^{\min}(v)\geq \lceil\eta n\rceil$. Finally, by \cref{step:absorbing} and \cref{thm:cleaning-A}, if $|W_A^\pm|\geq 2$, then each~$v\in W_A^\pm$ satisfies $\exc_{D_1\setminus A}^\pm(v)\leq \exc_{D_1}^\pm(v)\leq \exc_T^\pm(v)\leq \lceil \eta n\rceil$ and, if $|W_A^\pm|=1$, then $d_A^\pm(v)=\lceil\eta n\rceil$ for the unique~$v\in W_A^\pm$.}
		
		Apply \cref{lm:layouts} with $D_1, U_2^*$, and $18\sqrt{\varepsilon}$ playing the roles of $D,U^*$, and $\varepsilon$ to obtain $\ell\in \mathbb{N}$ and $W$-exceptional layouts $(L_1, F_1), \dots, (L_\ell, F_\ell)$ which are~$U_2^*$-path consistent with respect to $(D_1', \cF)$ and such that the following hold, \NEW{where $L$ is the multiset defined by $L\coloneqq \bigcup_{i\in [\ell]}L_i$}.\OLD{Let $L\coloneqq \bigcup_{i\in [\ell]}L_i$.}
		\begin{enumerate}[label=\upshape(\Alph*)]
			\item $d\leq \ell \leq d+\varepsilon^{\frac{1}{5}}n$.\label{thm:layouts-ell}
			\item $L$ contains exactly~$\texc(D_1)-r$ non-trivial paths.\label{thm:layouts-good}
			\item For each $v\in W_1$, $d_L^\pm(v)=d_{D_1\setminus A}^\pm(v)$.\label{thm:layouts-degreeW1}
			\item For each~$v\in W_2$, $d_L^\pm(v)=d_{D_1}^\pm(v)-r$.\label{thm:layouts-degreeW2}
			\item For each~$v\in V'$, $d_{D_1}^\pm(v)=d_L^\pm(v)+|\{i\in [\ell]\mid v\notin V(L_i)\}|+r-\phi^\mp(v)$.\label{thm:layouts-degreeV'}
			\newcounter{layoutsV'}
			\setcounter{layoutsV'}{\value{enumi}}
			\item For each~$i\in [\ell]$, $|V(L_i)|, |E(L_i)|\leq \varepsilon^{\frac{1}{7}}n$.\label{thm:layouts-size}
			\item For each~$v\in V'$, $d_L(v)\leq \varepsilon^{\frac{1}{3}}n$ and there exist at most $\varepsilon^{\frac{1}{5}} n$ indices $i\in [\ell]$ such that $v\in V(L_i)$. \label{thm:layouts-V'}
		\end{enumerate}
		
		\NEW{We now transform $(L_1, F_1), \dots, (L_\ell, F_\ell)$ into auxiliary layouts on $V'$.} For each $i\in [\ell]$, let $(L_i^{\upharpoonright W}, F_i^{\upharpoonright W})$ be \NEW{obtained from $(L_i, F_i)$ using the procedure described in}\OLD{as in} \cref{def:layoutV'}.
		Let $L^{\upharpoonright W}\coloneqq \bigcup_{i\in [\ell]}L_i^{\upharpoonright W}$ and $\cF^{\upharpoonright W}\coloneqq \bigcup_{i\in[\ell]}F_i^{\upharpoonright W}$. Then, \NEW{\cref{def:layoutV'} implies that}\OLD{by construction,} $(L_1^{\upharpoonright W}, F_1^{\upharpoonright W}), \dots$, $(L_\ell^{\upharpoonright W}, F_\ell^{\upharpoonright W})$ are layouts on~$V'$ such that \NEW{the following hold.
		\begin{enumerate}[label=\upshape(\Alph*$'$)]
			\setcounter{enumi}{\value{layoutsV'}}
			\item Let $i\in [\ell]$. By \cref{thm:layouts-size}, $|V(L_i^{\upharpoonright W})|\leq |V(L_i)|\leq \varepsilon^{\frac{1}{7}}n\leq 2\varepsilon^{\frac{1}{7}}|V'|$ and, similarly, $|E(L_i^{\upharpoonright W})|\leq 2\varepsilon^{\frac{1}{7}}|V'|$.
			\item Let $v\in V'$. By \cref{thm:layouts-V'}, $d_{L^{\upharpoonright W}}(v)\leq d_L(v)\leq \varepsilon^{\frac{1}{3}}n\leq 2\varepsilon^{\frac{1}{3}}|V'|$ and there exist at most $\varepsilon^{\frac{1}{5}} n\leq 2\varepsilon^{\frac{1}{5}} |V'|$ indices~$i\in [\ell]$ such that~$v\in V(L_i^{\upharpoonright W})$.
		\end{enumerate}
		Thus, \cref{lm:approxdecomp}\cref{lm:approxdecomp-degreev,lm:approxdecomp-sizeLi} are satisfied with $|V'|, \varepsilon^{\frac{1}{29}}$, and $L_1^{\upharpoonright W}, \dots, L_\ell^{\upharpoonright W}$ playing the roles of $n, \varepsilon$, and $L_1, \dots, L_\ell$%
			\COMMENT{$(\varepsilon^{\frac{1}{29}})^2\geq (\varepsilon^{\frac{1}{29}})^4= \varepsilon^{\frac{4}{29}}\geq 2\varepsilon^{\frac{1}{7}}$ (for \cref{lm:approxdecomp}\cref{lm:approxdecomp-sizeLi}), $(\varepsilon^{\frac{1}{29}})^3= \varepsilon^{\frac{3}{29}}\geq 2\varepsilon^{\frac{1}{3}}$ (for the first part of \cref{lm:approxdecomp}\cref{lm:approxdecomp-degreev}), and $(\varepsilon^{\frac{1}{29}})^2=\varepsilon^{\frac{2}{29}}\geq 2\varepsilon^{\frac{1}{5}}$ (for the second part of \cref{lm:approxdecomp}\cref{lm:approxdecomp-degreev})}.}
		\OLD{for each~$i\in [\ell]$, $|V(L_i^{\upharpoonright W})|,|E(L_i^{\upharpoonright W})|\leq \varepsilon^{\frac{1}{7}}n\leq 2\varepsilon^{\frac{1}{7}}|V'|$ and, for each~$v\in V'$, $d_{L^{\upharpoonright W}}(v)\leq \varepsilon^{\frac{1}{3}}n\leq 2\varepsilon^{\frac{1}{3}}|V'|$ and there exist at most $\varepsilon^{\frac{1}{5}} n\leq 2\varepsilon^{\frac{1}{5}} |V'|$ indices~$i\in [\ell]$ such that~$v\in V(L_i^{\upharpoonright W})$.}
		
		\NEW{In order to approximately decompose $D_1[V']$ using \cref{lm:approxdecomp}, we need to partition $D_1[V']$ into a dense almost regular digraph (which will play the role of $D$ in \cref{lm:approxdecomp}) and a sparse almost regular robust expander (which will play the role of $\Gamma$ in \cref{lm:approxdecomp}). We choose $\Gamma$ randomly as follows.}\OLD{We now reserve a random slice~$\Gamma\subseteq D_1[V']$ for \cref{lm:approxdecomp} as follows.}
		Let $\delta\coloneqq \frac{d+r}{n}$.
		Note that, by \cref{thm:cleaning-degreeV'},~$D_1[V']$ is~$(\delta, 10\sqrt{\varepsilon})$-almost regular%
			\COMMENT{By \cref{thm:cleaning-degreeV'} and definition of $\delta$, we have $(2\delta-18\sqrt{\varepsilon})n\leq d_{D_1}(v)\leq 2\delta n$. Therefore, we have $(\delta-9\sqrt{\varepsilon}-\varepsilon)n\leq \frac{d_{D_1}(v)-\exc_{D_1}(v)}{2}\leq d_{D_1}^\pm(v)\leq \frac{d_{D_1}(v)+\exc_{D_1}(v)}{2}\leq (\delta+\varepsilon)n$. Since $|W|\leq 4\varepsilon n$, we have $(\delta-10\sqrt{\varepsilon})|V'|\leq (\delta-9\sqrt{\varepsilon}-\varepsilon)n-4\varepsilon n\leq d_{D_1[V']}^\pm(v)\leq (\delta+\varepsilon)n\leq (\delta+10\sqrt{\varepsilon})|V'|$.}.
		By \cref{lm:3/8rob} and \cref{thm:cleaning-d},~$D_1[V']$ is a robust~$(\nu, \tau)$-outexpander.
		Apply \cref{lm:Gamma} with $D_1[V'], |V'|$, and $10\sqrt{\varepsilon}$ playing the roles of $D, n$, and $\varepsilon$ to obtain $\Gamma\subseteq D_1[V']$ such that~$\Gamma$ is a~$(\gamma, 10\sqrt{\varepsilon})$-almost regular~$(10\sqrt{\varepsilon},|V'|^{-2})$-robust~$(\nu, \tau)$-outexpander and $D_1''\coloneqq D_1[V']\setminus \Gamma$ is~$(\delta-\gamma, 10\sqrt{\varepsilon})$-almost regular.
		
		Observe that, by \cref{thm:layouts-ell}, $\ell\leq (\delta-\frac{\eta}{2})|V'|$%
			\COMMENT{$\ell\leq d+\varepsilon^{\frac{1}{5}}n=\delta n-r+\varepsilon^{\frac{1}{5}}n\leq (\delta-\eta+\varepsilon^{\frac{1}{5}})n\leq (\delta -\frac{\eta}{2})|V'|$.}.
		Apply \cref{lm:approxdecomp} with $D_1'', \cF^{\upharpoonright W}, |V'|, |V'|^{-2}, \frac{\eta}{2}, \varepsilon^{\frac{1}{29}}$,
		and $(L_1^{\upharpoonright W},F_1^{\upharpoonright W}),\dots, (L_\ell^{\upharpoonright W},F_\ell^{\upharpoonright W})$ playing the roles of $D, \cF, n, p, \eta, \varepsilon$, and $(L_1,F_1),\dots, (L_\ell,F_\ell)$ to obtain edge-disjoint $\cH_1^{\upharpoonright W}, \dots, \cH_\ell^{\upharpoonright W}\subseteq D_1''\cup \Gamma\cup \cF^{\upharpoonright W}=D_1[V']\cup \cF^{\upharpoonright W}$ such that, for each~$i\in [\ell]$,~$\cH_i^{\upharpoonright W}$ is a spanning configuration of shape~$(L_i^{\upharpoonright W}, F_i^{\upharpoonright W})$ and the following holds. Let $\cH^{\upharpoonright W}\coloneqq \bigcup_{i\in [\ell]}\cH_i^{\upharpoonright W}$ and $D_2'\coloneqq D_1[V'] \setminus \cH^{\upharpoonright W}$. Then, 
		\begin{enumerate}[label=\upshape(\Roman*)]
			\newcounter{rob}
			\setcounter{rob}{\value{enumi}}
			\item $D_2'$ is a robust~$(\frac{\nu}{2}, \tau)$-outexpander.\label{thm:approxdecomp-rob1}
		\end{enumerate}
			
		\NEW{Next, we transform the auxiliary spanning configurations $\cH_1^{\upharpoonright W}, \dots, \cH_\ell^{\upharpoonright W}$ into edge-disjoint spanning configurations on $V$ of shapes $(L_1, F_1), \dots, (L_\ell, F_\ell)$ as follows.}\OLD{Now construct~$W$-exceptional spanning configurations as follows.}
		For each~$i\in[\ell]$, let~$\cH_i$ be the digraph with $V(\cH_i)\coloneqq V$ and $E(\cH_i)=(E(\cH_i^{\upharpoonright W})\setminus F_i^{\upharpoonright W})\cup F_i$. 
		Denote $\cH\coloneqq \bigcup_{i\in [\ell]}\cH_i$.
		\NEW{Let $i\in[\ell]$. Then, \cref{fact:layoutV'} implies that}\OLD{Then, by \cref{fact:layoutV'}, for each~$i\in[\ell]$,} $\cH_i\subseteq D_1[V']\cup F_i$ and~$\cH_i$ is a~$W$-exceptional spanning configuration of shape~$(L_i,F_i)$.
		Moreover,\OLD{for each~$i\in[\ell]$,} since $F_i\subseteq E_W(D_1)\setminus A$, we have $\cH_i\subseteq D_1\setminus A$ and so \begin{equation}\label{eq:approx-A}
			E(\cH)\cap A=\emptyset.
		\end{equation}
		Furthermore, by definition of~$U_2^*$-path consistency with respect to~$(D_1',\cF)$, the sets $F_1, \dots, F_\ell$ are edge-disjoint. 
		Thus, since $\cH_1^{\upharpoonright W}, \dots, \cH_\ell^{\upharpoonright W}$ are edge-disjoint, $\cH_1, \dots, \cH_\ell$ are edge-disjoint.
		
		\NEW{Finally, we verify that $\cH_1, \dots, \cH_\ell$ induce a good partial path decomposition of $D_1$.} For each~$i\in[\ell]$, let~$\cP_{2,i}$ be a path decomposition of~$\cH_i$ induced by~$(L_i, F_i)$. Let $\cP_2\coloneqq \bigcup_{i\in [\ell]}\cP_{2,i}$ and $D_2\coloneqq D_1\setminus \cP_2$.
		We claim that the following holds.
		
		\begin{claim}\label{claim:P2}
			$\cP_2$ is a $(U_2^*, W, A)$-partial path decomposition of~$D_1$ \NEW{of size $\texc(D_1)-r$}.
		\end{claim}
	
		\begin{proofclaim}
			By \cref{fact:layouts}, \NEW{\cref{thm:layouts-good}, and} since $(L_1, F_1)$, $\dots$, $(L_\ell, F_\ell)$ are~$U_2^*$-path consistent with respect to~$(D_1',\cF)$,~$\cP_2$ is a partial path decomposition of~$D_1$ \NEW{of size $\texc(D_1)-r$} such that 
			\begin{equation}\label{eq:approx-U0}
				U^0(D_1)\cap (V^+(\cP_2)\cup V^-(\cP_2))\subseteq U_2^*.
			\end{equation}
			\NEW{Thus, it only remains to show that $\cP_2$ is consistent with $A^+$ and $A^-$. By \cref{eq:approx-A},}%
			\OLD{Moreover,} $E(\cP_2)\cap A=\emptyset$.
			Thus, it suffices to show that each \NEW{$v\in W$}\OLD{$v\in W_A$} is the starting point of \NEW{at most $\exc_{D_1}^+(v)-d_A^+(v)$ paths in $\cP_2$ and the ending point of at most $\exc_{D_1}^-(v)-d_A^-(v)$.}\OLD{exactly~$\texc_{D_1,U_2^*}^+(v)$ paths in~$\cP_2$ and the ending point of exactly~$\texc_{D_1,U_2^*}^-(v)$ paths in~$\cP_2$.}
			
			\NEW{Let $v\in W$. If $v\in W\setminus W_A$, then \cref{claim:A} implies that $d_A(v)=0$. Moreover, \cref{eq:cleaning-U1*} implies that $W\cap U_2^*=\emptyset$. Thus, the fact that $\cP_2$ is a partial path decomposition satisfying \cref{eq:approx-U0} implies that $\cP_2$ contains at most $\exc_{D_1}^+(v)=\exc_{D_1}^+(v)-d_A^+(v)$ paths which start at $v$ and at most $\exc_{D_1}^-(v)=\exc_{D_1}^-(v)-d_A^-(v)$ paths which end at $v$.
			We may therefore assume that}\OLD{Let} $v\in W_A^+$ (similar argument hold if~$v\in W_A^-$)\OLD{and recall that~$W_A\subseteq W_1$}. \NEW{By \cref{def:A} and \cref{thm:cleaning-A}, we have $\exc_{D_1}^+(v)>0$ and so $v\in U^+(D_1)$. Thus, the fact that $(L_1, F_1)$, $\dots$, $(L_\ell, F_\ell)$ are~$U_2^*$-path consistent with respect to~$(D_1',\cF)$ implies that the number of paths in $\cP_2$ which end at $v$ is $0\leq \exc_{D_1}^-(v)-d_A^-(v)$, as desired.}\OLD{Note that, since $U_2^*\cap W_A=\emptyset$, we have $\texc_{D_1,U_2^*}^-(v)=\exc_{D_1}^-(v)=0$ and $\texc_{D_1,U_2^*}^+(v)=\exc_{D_1}^+(v)-d_A^+(v)$.
			Since $(L_1, F_1)$, $\dots$, $(L_\ell, F_\ell)$ are~$U_2^*$-path consistent with respect to~$(D_1',\cF)$,~$v$ is not the ending point of any path in~$\cP_2$, as desired.}
			In particular,~$v$ is an internal vertex of precisely~$d_{\cP_2}^-(v)$ paths in~$\cP_2$.
			Therefore, the number of paths in~$\cP_2$ \NEW{which start at $v$ is
			\begin{align*}
				d_{\cP_2}^+(v)-d_{\cP_2}^-(v)&
				\stackrel{\text{\eqmakebox[approx][c]{\cref{fact:layouts}}}}{=}d_L^+(v)-d_L^-(v)
				\stackrel{\text{\cref{thm:layouts-degreeW1}}}{=}d_{D_1\setminus A}^+(v)-d_{D_1\setminus A}^-(v)\\
				&\stackrel{\text{\eqmakebox[approx][c]{\cref{def:A}}}}{=}d_{D_1}^+(v)-d_A^+(v)-d_{D_1}^-(v)
				=\exc_{D_1}^+(v)-d_A^+(v),
			\end{align*}}\OLD{whose starting point is~$v$ equals
			\[d_{\cP_2}^+(v)-d_{\cP_2}^-(v)\stackrel{\text{\cref{fact:layouts},\cref{thm:layouts-degreeW1}}}{=}d_{D_1\setminus A}^+(v)-d_{D_1}^-(v)\stackrel{\text{\cref{thm:cleaning-A}}}{=}\exc_{D_1}^+(v)-d_A^+(v)=\texc_{D_1,U_2^*}^+(v),\]}
			as desired.
		\end{proofclaim}
		
		\NEW{It remains to show that $\cP_2$ is good, i.e.\ that $\texc(D_2)=\texc(D_1)-|\cP_2|$.}
		Recall that by \cref{thm:cleaning-A},~$\cP_1$ avoids all the edges in~$A$.
		Thus, by \cref{fact:layouts} and \cref{thm:layouts-degreeW1,thm:layouts-degreeV',thm:layouts-degreeW2}, the following hold.
		\begin{enumerate}[label=\upshape(\Roman*),resume]
			\item For each~$v\in W_1$, $N_{D_2}^\pm(v)=N_A^\pm(v)$.\label{thm:approxdecomp-W1}
			\item For each $v\in V\setminus (W_1\cup X^-)$, $d_{D_2}^+(v)= r$.\label{thm:approxdecomp-d+}
			\item For each $v\in V\setminus (W_1\cup X^+)$, $d_{D_2}^-(v)= r$.\label{thm:approxdecomp-d-}
			\item For each $v\in X^\pm$, $d_{D_2}^\mp(v)=r-1$.\label{thm:approxdecomp-Xpm}
		\end{enumerate}
		\OLD{Note that, since $|X^\pm\cup A^\pm|=r$ and $X^+\cap X^-\subsetneq V\setminus W_1$, \cref{thm:approxdecomp-W1,thm:approxdecomp-Xpm,thm:approxdecomp-d+,thm:approxdecomp-d-} imply that}
		\begin{claim}\label{eq:texcD2}
			\NEW{$\texc(D_2)=r$.}
		\end{claim}
		
		\begin{proofclaim}
			\NEW{First, we check that $\exc(D_2)\leq r$. Let $v\in V$.
			\begin{itemize}
				\item If $v\in X^+\cup X^-$, then \cref{thm:approxdecomp-Xpm} implies that $\exc_{D_2}^+(v)=0$.
				\item If $v\in X^+\setminus X^-$, then \cref{thm:approxdecomp-Xpm,thm:approxdecomp-d+} imply that $\exc_{D_2}^+(v)=1$.
				\item If $v\in X^-\setminus X^+$, then \cref{thm:approxdecomp-Xpm,thm:approxdecomp-d+} imply that $\exc_{D_2}^+(v)=0$.
				\item If $v\in V\setminus (X^+\cup X^-\cup W)$, then \cref{thm:approxdecomp-d-,thm:approxdecomp-d+} imply that $\exc_{D_2}^+(v)=0$.
				\item Suppose that $v\in W_1\setminus W_A^+$. Recall from \cref{step:absorbing} that $A^+$ and $A^-$ are absorbing sets of $(W_A^+,V')$-starting and $(V',W_A^-)$-ending edges for $T$. Thus, $d_A^+(v)=0$ and so \cref{thm:approxdecomp-W1} implies that $\exc_{D_2}^+(v)=0$.
				\item Suppose that $v\in W_1\cup W_A^+$. Then, \cref{def:A} implies that $v\in U^+(T)$ and so $d_A^-(v)=0$. Therefore, \cref{thm:approxdecomp-W1} implies that $\exc_{D_2}^+(v)=d_A^+(v)$.
			\end{itemize}
			Moreover, recall that $W_A^+\subseteq W_1$ and $(X^+\cup X^-)\cap W_1=\emptyset$ (see \cref{claim:Xpm}). Thus,
			\begin{align*}
				\exc(D_2)\stackrel{\text{\cref{eq:exc}}}{=}\sum_{v\in V}\exc_{D_2}^+(v)=|X^+\setminus X^-|+\sum_{v\in W_A^+}d_A^+(v)\stackrel{\text{\cref{def:A}}}{=}|X^+\setminus X^-|+|A^+|\stackrel{\text{\cref{claim:Xpm}}}{\leq} r,
			\end{align*}
			as desired.}
			
			\NEW{Thus, it is enough to show that $\Delta^0(D_2)=r$. By \cref{claim:Xpm}, $|X^-|\leq r$ and, by \cref{claim:A-W}, $|W_1|\leq |W|\leq 4\varepsilon n$. Thus, $V\setminus (W_1\cup X^-)\neq \emptyset$ and so \cref{thm:approxdecomp-d+} implies that $\Delta^0(D_2)\geq r$. Let $v\in V$. We verify that both $d_{D_2}^\pm(v)\leq r$. If $v\notin W_1$, then \cref{thm:approxdecomp-Xpm,thm:approxdecomp-d+,thm:approxdecomp-d-} imply that both $d_{D_2}^\pm(v)\leq r$. We may therefore assume that $v\in W_1\subseteq W$. By \cref{def:A} and \cref{claim:A-A}, $d_A^+(v)=d_{A^+}(v)\leq |A^+|\leq r$ and, similarly, $d_A^-(v)\leq r$. Therefore, $\Delta^0(D_2)\leq r$ and so we are done.} 
		\end{proofclaim}
		\OLD{Indeed, each~$v\in V$ satisfies $d_{D_2}^{\max}(v)\leq r$ and each $v\notin (X^+\cap X^-)\cup W_1$ satisfies $d_{D_2}^{\max}(v)=r$.		
		Thus, $X^+\cap X^-\subsetneq V\setminus W_1$ implies $\Delta^0(D_2)=r$. Moreover, $\exc(D_2)=|(X^\pm\setminus X^\mp)\cup A^\pm|\leq |X^\pm\cup A^\pm|=r$.}
		
		Thus, \NEW{\cref{claim:P2,eq:texcD2} imply}\OLD{\cref{thm:layouts-good} implies} that
		\begin{enumerate}[resume,label=\upshape(\Roman*)]
			\item $\cP_2$ is a good partial path decomposition of~$D_1$.\label{thm:approxdecomp-good}
		\end{enumerate}
	
		\item \textbf{Completing the path decomposition.}\label{step:fulldecomp}
		\NEW{Finally, we will decompose $D_2$ using \cref{prop:absorbingedges}. Recall that to apply \cref{prop:absorbingedges}, all the edges incident to the exceptional set must be absorbing edges. By \cref{thm:approxdecomp-W1,thm:approxdecomp-Xpm,thm:approxdecomp-d+,thm:approxdecomp-d-}, all the absorbing edges are incident to $W_1$ and the vertices in $W_2$ are still incident to some non-absorbing edges in $D_2$. Thus, we will apply \cref{prop:absorbingedges} with $W_1$ and $V'\cup W_2$ playing the roles of $W$ and $V'$.}
		 
		\NEW{First, we check that all the prerequisites of \cref{prop:absorbingedges} are satisfied.}
		Note that~$D_2[V']=D_2'$ \NEW{($D_2'$ was defined just above \cref{thm:approxdecomp-rob1} and $D_2$ was defined just above \cref{claim:P2})}\OLD{and recall that $|W|\leq 4\varepsilon n$}. Thus, by \cref{lm:verticesedgesremovalrob}\NEW{\cref{lm:verticesedgesremovalrob-vertices}}, \NEW{\cref{claim:A-W}}, and \cref{thm:approxdecomp-rob1},
		$D_2-W_1$ is a robust~$(\frac{\nu}{4}, 2\tau)$-outexpander.
		By \cref{fact:absorbingedges}, \NEW{\cref{thm:cleaning-A}}, and \cref{claim:P2}\OLD{and~$\cP_2$ is a $(U_2^*, W, A)$-partial path decomposition of~$D_1$}, 
		$A^+$ and $A^-$ are absorbing sets of~\NEW{$(W_A^+,V')$}\OLD{$(W_1, V')$}-starting and~\NEW{$(V', W_A^-)$}\OLD{$(V', W_1)$}-ending edges for~$D_2$.
		\NEW{In particular, \cref{def:A} and the fact that $W_A\subseteq W_1$ imply that $A^+$ and $A^-$ are absorbing sets of~$(W_1, V'\cup W_2)$-starting and~$(V'\cup W_2, W_1)$-ending edges for~$D_2$.}				
		
		Let $Y^\pm\coloneqq X^\pm\setminus X^\mp$, $Y^*\coloneqq X^+\cap X^-$, and $Y^0\coloneqq V\setminus(Y^+\cup Y^-\cup Y^*\NEW{\cup W_1})$.
		\NEW{We aim to apply to apply \cref{prop:absorbingedges} with $Y^+, Y^-, Y^*$, and $Y^0$ playing the roles of $X^+, X^-, X^*$, and $X^0$.}
		\NEW{First, observe that \cref{thm:approxdecomp-W1,thm:approxdecomp-Xpm,thm:approxdecomp-d+,thm:approxdecomp-d-} imply that the excess and degree conditions of \cref{prop:absorbingedges} hold with $W_1, Y^+, Y^-$, and $Y^*$ playing the roles of $W, X^+, X^-$, and $X^*$. Moreover, \cref{claim:Xpm} implies that $|Y^\pm\cup Y^*|+|A^\pm|=r$. Finally, we claim that $V(A^\pm)\cap (V'\cup W_2)\subseteq Y^\mp\cup Y^0$.
		By \cref{claim:A-A}, $|A^\pm|\in \{0,r\}$. If $|A^\pm|=0$, then $V(A^\pm)\cap (V'\cup W_2)=\emptyset$ and so we are done. If $|A^\pm|=r$, then \cref{claim:Xpm} implies that $X^\pm=\emptyset$ and so $V'\cup W_2=Y^\mp\cup Y^0$. Thus, $V(A^\pm)\cap (V'\cup W_2)\subseteq Y^\mp\cup Y^0$, as desired.}%
		\OLD{Note that, if $A^\pm\neq \emptyset$, then, by \cref{step:absorbing}, $|A^\pm|=r$ and so, by \cref{step:approx}, $X^\pm=\emptyset$. Therefore, $V(A^\pm)\cap (V\setminus W_1)=V(A^\pm)\cap V'\subseteq Y^\mp\cup Y^0$.}
		 
		Apply \cref{prop:absorbingedges} with $D_2, n-|W_1|, V\setminus W_1, W_1, Y^+, Y^-, Y^*, Y^0, \frac{\eta}{2}, \frac{\nu}{4}$, and $2\tau$ playing the roles of $D, n, V', W, X^+, X^-, X^*, X^0, \delta, \nu$, and $\tau$ to obtain a path decomposition~$\cP_3$ of~$D_2$ of size~$r$. Note that, by \cref{eq:texcD2},~$\cP_3$ a perfect path decomposition of~$D_2$.
		Recall that by \cref{thm:cleaning-A} and \cref{thm:approxdecomp-good},~$\cP_1$ and~$\cP_2$ are good.
		Then, by \cref{fact:combinefull}, $\cP\coloneqq\cP_1\cup \cP_2\cup \cP_3$ is a perfect path decomposition of~$T$.
		That is,~$|\cP|= \texc(T)$.
		This completes the proof.\qedhere
	\end{steps}
\end{proof}

%% file: Auxiliary_Excess_Function.tex
	\onlyinsubfile{
				\setcounter{section}{10}
		\section{Auxiliary excess function}}

\NEW{In this section, we introduce some concepts which will be convenient for constructing good $(U^*,W,A)$-partial path decompositions. (Recall these were defined in \cref{def:U*WApartialpathdecomp}.)}

\NEW{Recall that} once we have chosen absorbing edges, we need to ensure that \hypertarget{texc1}{(i)} these edges are not used for other purposes and \hypertarget{texc2}{(ii)} not too many paths have endpoints in~$W$.
Moreover,\OLD{recall from \cref{sec:sketch-general} that} if $\texc(T)>\exc(T)$, then some vertices $v$ will have to be used as starting/ending points of paths more than~$\exc_T^\pm(v)$ times. \NEW{As discussed briefly in \cref{sec:A}}\OLD{For simplicity}, we will choose in advance which vertices will be used as these additional endpoints:\OLD{When choosing these endpoints, we aim to maximise the total number of distinct endpoints available at each step of our decomposition (this is helpful when finding suitable~$X^\pm$ for \cref{prop:absorbingedges} in the case where we are not using absorbing edges). In other words, this means that} we will initially select a set~$U^*$ of~$\texc(T)-\exc(T)$ (distinct) vertices in~$U^0(T)$ as additional endpoints of paths. We then treat each vertex in~$U^*$ as if it has positive and negative excess both equal to one. The following concept of an auxiliary excess function (as defined in \cref{eq:texc}) encapsulates all this -- it also incorporates the constraints given by \hyperlink{texc1}{(i)} and \hyperlink{texc2}{(ii)} above. It will enable us to easily keep track of how many paths remain to be chosen and which vertices can be used as endpoints.

\NEW{Let~$D$ be a digraph and $W,V'\subseteq V(D)$ be disjoint. Suppose that $A\subseteq E(D)$ is a $(W,V')$-absorbing set for~$D$. Let~$U^*\subseteq U^0(D)$ satisfy
$|U^*|= \texc(D)-\exc(D)$.}\OLD{Let~$D$ be a digraph and $U^*\subseteq U^0(D)$ satisf\NEW{y} $|U^*|= \texc(D)-\exc(D)$. Let $A^+,A^-\subseteq E(D)$ be absorbing sets of~$(W,V')$-starting and~$(V',W)$-ending edges.}
Note that, by \cref{def:A}, $(V(A)\cap W)\cap U^*=\emptyset$.
For each~$v\in V(D)$, define
\begin{equation}\label{eq:texc}
\texc_{D,U^*,W,A}^\pm(v)\coloneqq 
\begin{cases}
1 & \text{if }v\in U^*,\\
\exc_D^\pm(v)-d_A^\pm(v)& \text{if }v\in V(A)\cap W,\\
\exc_D^\pm(v)& \text{otherwise}.\\
\end{cases}
\end{equation}
Then, define
\begin{align*}
	\tU_{U^*,W,A}^\pm(D)&\coloneqq \{v\in V(D)\mid
	\texc_{D,U^*,W,A}^\pm(v)>0\};\\
	\tU_{U^*,W,A}^0(D)&\coloneqq V(D)\setminus
	(\tU_{U^*,W,A}^+(D)\cup \tU_{U^*,W,A}^-(D));\\
	\texc_{D,U^*,W,A}^\pm(S)&\coloneqq \sum_{v\in S}\texc_{D,U^*,W,A}^\pm(v) \text{ for each } S\subseteq V(D);\\
	\texc_{U^*,W,A}^\pm(D)&\coloneqq \sum_{v\in V}\texc_{D,U^*,W,A}^\pm(v).
\end{align*}
\NEW{Denote by $A^+$ and $A^-$ the absorbing sets of $(W,V')$-starting and $(V',W)$-ending edges contained in $A$. Then,}\OLD{Note that}
\begin{equation}\label{eq:texcA}
\texc_{U^*,W,A}^\pm(D)=\exc(D)+|U^*|-|A^\pm|=\texc(D)-|A^\pm|.
\end{equation}
\NEW{Note that it is possible that $\texc_{U^*,W,A}^+(D)\neq \texc_{U^*,W,A}^-(D)$.}
For simplicity, when~$A$ and~$W$ are clear from the context, they will be
omitted in the subscripts of the above notation.

Note that the analogue of \cref{fact:partialdecompS} holds for \NEW{the}\OLD{this} auxiliary excess function.

\begin{fact}\label{fact:texcS}
	\NEW{Let~$D$ be a digraph and $W,V'\subseteq V(D)$ be disjoint. Suppose that $A\subseteq E(D)$ is a $(W,V')$-absorbing set for~$D$. Let~$U^*\subseteq U^0(D)$ satisfy	$|U^*|= \texc(D)-\exc(D)$.}\OLD{Let~$D$ be a digraph and~$W,V'\subseteq V(D)$ be disjoint. Suppose $A^+, A^-\subseteq E(D)$ are absorbing sets of~$(W,V')$-starting and~$(V',W)$-ending edges for~$D$. Let $V\coloneqq V(D)$ and let $U^*\subseteq U^0(D)$ satisfy $|U^*|=\texc(D)-\exc(D)$.}
	Then, for any $S\subseteq V\NEW{(D)}$,
	$\texc_{U^*}^\pm(D)=\texc_{D,U^*}^\pm (S)+\texc_{D,U^*}^\pm (V\NEW{(D)}\setminus S)$.
\end{fact}

\OLD{Finally, a set~$\cP$ of edge-disjoint paths of~$D$ is a \emph{$(U^*,W,A)$-partial path decomposition of~$D$} if $\cP\subseteq D\setminus A$ and each~$v\in V(D)$ is the
	starting point of at most~$\texc_{D,U^*,W,A}^+(v)$ paths in~$\cP$ and the ending
	point of at most~$\texc_{D,U^*,W,A}^-(v)$ paths in~$\cP$, i.e.\ if~$\cP$ is a partial path decomposition of~$D$ which is consistent with~$A^+$ and~$A^-$ and such that,
	for each~$v\in U^0(D)$, if~$v$ is an endpoint of a path in~$\cP$, then~$v\in
	U^*$.}

\NEW{Observe that the following holds by definition of the auxiliary excess function.}

\begin{fact}\label{fact:texc}
	\NEW{Let $D, W, V', A$, and $U^*$ satisfy the assumptions of \cref{fact:texcS}.	
	A set~$\cP$ of edge-disjoint paths of~$D$ is a $(U^*,W,A)$-partial path decomposition of~$D$ if and only if $\cP\subseteq D\setminus A$ and each~$v\in V(D)$ is the
	starting point of at most~$\texc_{D,U^*,W,A}^+(v)$ paths in~$\cP$ and the ending
	point of at most~$\texc_{D,U^*,W,A}^-(v)$ paths in~$\cP$.}
\end{fact}

\OLD{Note that, by \cref{eq:texc},}\NEW{Thus,} this auxiliary excess function designates which vertices are still available to use as endpoints and, by \cref{eq:texcA}, it indicates how many paths we are still allowed to \NEW{choose}\OLD{take}. 
For these reasons, fixing~$U^*$ at the beginning will prove very useful in \cref{sec:cleaning}, even though it is not necessary and may look cumbersome at first glance.

\NEW{Let $D, W, V', A$, and $U^*$ satisfy the assumptions of \cref{fact:texcS}. Suppose that $\cP$ is a $(U^*,W,A)$-partial path decomposition of $D$. By \cref{fact:texc}, each path in $\cP$ corresponds to some auxiliary excess and so, by removing $\cP$ (and removing from $U^*$ the vertices which have already been used as endpoints), we reduce the auxiliary positive/negative excess of each $v\in V(D)$ by the number of paths in $\cP$ which start/end at $v$ (\cref{prop:texcv}). 
This implies that the total auxiliary excess of $D$ is decreased by precisely $|\cP|$ when we remove the paths in $\cP$ (\cref{fact:specialdecomp}). The auxiliary excess function will thus be much more convenient to use than the modified excess function introduced in \cref{sec:intro} (compare the bounds in \cref{prop:excpartialdecomp}\cref{prop:excpartialdecomp-general} and \cref{fact:specialdecomp}). Moreover, this implies that good $(U^*,W,A)$-partial path decompositions can be combined to form a larger good $(U^*,W,A)$-partial path decomposition (\cref{fact:combinespecial}).}

\begin{prop}\label{prop:texcv}
	\NEW{Let $D, W, V', A$, and $U^*$ satisfy the assumptions of \cref{fact:texcS}. Suppose that $\cP$ is a $(U^*,W,A)$-partial path decomposition of $D$. Denote $U^{**}\coloneqq U^*\setminus (V^+(\cP)\cup V^-(\cP))$. For each $v\in V(D)$, let $n_\cP^+(v)$ and $n_\cP^-(v)$ denote the number of paths in $\cP$ which start and end at $v$, respectively. Then, each $v\in V(D)$ satisfies $\texc_{D\setminus \cP,U^{**}}^\pm(v)=\texc_{D,U^*}^\pm(v)-n_\cP^\pm(v)$.}
\end{prop}

\begin{proof}
	\NEW{If $v\in V(D)\setminus U^*$, then $\exc_{D\setminus \cP}^\pm(v)=\exc_D^\pm(v)-n_\cP^\pm(v)$ and so the \lcnamecref{prop:texcv} holds.
	We may therefore assume that $v\in U^*$. If $v\in U^{**}$, then $n_\cP^\pm(v)=0$ and so $\texc_{D\setminus \cP,U^{**}}^\pm(v)=1=\texc_{D,U^*}^\pm(v)-n_\cP^\pm(v)$, as desired.
	We may therefore assume that $v\in U^*\setminus U^{**}$. By \cref{def:A}, $v\notin W$. First, suppose that $v\in V^+(\cP)\cap V^-(\cP)$. Note that both $n_\cP^\pm(v)=1$ and $\exc_{D\setminus \cP}(v)=0$.
	Since $v\notin W$, we have $\texc_{D\setminus \cP,U^{**}}^\pm(v)=\exc_{D\setminus \cP}^\pm(v)=0=\texc_{D,U^*}^\pm(v)-n_\cP^\pm(v)$, as desired.
	By symmetry, we may therefore assume that $v\in V^+(\cP)\setminus V^-(\cP)$. Note that $n_\cP^+(v)=1$, $n_\cP^-(v)=0$, and $\exc_{D\setminus \cP}(v)=-1$. Recall that $v\notin W$. Thus, $\texc_{D\setminus \cP}^+(v)=\exc_{D\setminus \cP}^+(v)=0=\texc_{D,U^*}^+(v)-n_\cP^+(v)$ and $\texc_{D\setminus \cP}^-(v)=\exc_{D\setminus \cP}^-(v)=1=\texc_{D,U^*}^-(v)-n_\cP^-(v)$.}
\end{proof}

\begin{cor}\label{fact:specialdecomp}
	\NEW{Let $D, W, V', A$, and $U^*$ satisfy the assumptions of \cref{fact:texcS}.}\OLD{Let~$D$ be a digraph and $W,V'\subseteq V(D)$ be disjoint. Suppose $A^+, A^-\subseteq E(D)$ are absorbing sets of~$(W,V')$-starting and~$(V',W)$-ending edges for~$D$. Denote $A\coloneqq A^+\cup A^-$. Let~$U^*\subseteq U^0(D)$ satisfy
		$|U^*|= \texc(D)-\exc(D)$.}
	Suppose~$\cP$ is a~$(U^*,W,A)$-partial path decomposition of~$D$. Let $U^{**}\coloneqq U^*\setminus (V^+(\cP)\cup V^-(\cP))$. 
	Then, $\texc_{U^{**}}^\pm (D\setminus \cP)=\texc_{U^*}^\pm(D)-|\cP|$.
\end{cor}

\begin{proof}
	\NEW{For each $v\in V(D)$, let $n_\cP^+(v)$ and $n_\cP^-(v)$ denote the number of paths in $\cP$ which start and end at $v$, respectively. Then,
	\begin{align*}
		\texc_{U^{**}}^\pm (D\setminus \cP)=\sum_{v\in V(D)}\texc_{D\setminus \cP,U^{**}}^\pm (v)\stackrel{\text{\cref{prop:texcv}}}{=}\sum_{v\in V(D)}(\texc_{D,U^*}^\pm (v)-n_\cP^\pm(v))=\texc_{U^*}^\pm (D)-|\cP|,
	\end{align*}
	as desired.}
\end{proof}

\OLD{Finally, by \cref{prop:sizeU*good}, the following analogue of \cref{fact:combinepartial} holds.}

\begin{cor}\label{fact:combinespecial}
	Let \NEW{$k\in \mathbb{N}$}\OLD{$k\in \mathbb{N}\setminus \{0\}$}. \NEW{Let $D, W, V', A$, and $U^*$ satisfy the assumptions of \cref{fact:texcS}.}\OLD{Let~$D$ be a digraph and $W,V'\subseteq V(D)$ be disjoint. Suppose $A^+,A^-\subseteq E(D)$ are absorbing sets of~$(W,V')$-starting and~$(V',W)$-ending edges for~$D$. Denote $A\coloneqq A^+\cup A^-$. Let $U^*\subseteq U^0(D)$ satisfy $|U^*|\leq \texc(D)-\exc(D)$.} Denote~$D_0\coloneqq D$ and~$U_0^*\coloneqq U^*$. Suppose that, for each~$i\in [k]$,~$\cP_i$ is a good $(U_{i-1}^*,W,A)$-partial path decomposition of~$D_{i-1}$, $D_i\coloneqq D_{i-1}\setminus \cP_i$, and $U_i^*\coloneqq U_{i-1}^*\setminus (V^+(\cP_i)\cup V^-(\cP_i))$. Let $\cP\coloneqq \bigcup_{i\in [k]}\cP_i$. Then,~$\cP$ is a good $(U^*,W,A)$-partial path decomposition of~$D$ of size $|\cP|=\sum_{i\in [k]}|\cP_i|$.
\end{cor}

\begin{proof}
	\NEW{By induction on $k$, it suffices to prove the case $k=2$. By \cref{fact:texc}, $E(\cP)\subseteq D\setminus A$.
	For each $v\in V(D)$, denote by $n_\cP^+(v)$ and $n_\cP^-(v)$ the number of paths in $\cP$ which start and end at $v$, and define $n_{\cP_1}^\pm(v)$ and $n_{\cP_2}^\pm(v)$ analogously. 
	Then, each $v\in V(D)$ satisfies
	\begin{align*}
		n_\cP^\pm(v)&=n_{\cP_1}^\pm(v)+n_{\cP_2}^\pm(v)\stackrel{\text{\cref{prop:texcv}}}{=}(\texc_{D_0,U_0^*}^\pm(v)-\texc_{D_1,U_1^*}^\pm(v))+(\texc_{D_1,U_1^*}^\pm(v)-\texc_{D_2,U_2^*}^\pm(v))\\
		&\leq \texc_{D_0,U_0^*}^\pm(v).
	\end{align*}
	Thus, each $v\in V(D)$ is the starting point of at most $\texc_{D,U^*}^+(v)$ paths in $\cP$ and the ending point of at most $\texc_{D,U^*}^-(v)$ paths in $\cP$. Therefore, \cref{fact:texc} implies that $\cP$ is a $(U^*,W,A)$-partial path decomposition of $D$. Denote by $A^+$ the absorbing set of starting edges contained in $A$. Then,
	\begin{align*}
		\texc(D\setminus \cP)&\stackrel{\text{\cref{eq:texcA}}}{=}\texc_{U_2^{**}}^+(D\setminus \cP)+|A^+|\stackrel{\text{\cref{fact:specialdecomp}}}{=}\texc_{U^*}^+(D)-|\cP|+|A^+| \stackrel{\text{\cref{eq:texcA}}}{=}\texc(D)-|\cP|.
	\end{align*}
	Therefore, $\cP$ is good.}
\end{proof}

\NEW{Recall from \cref{eq:N} that $N^+(D)$ and $N^-(D)$ denote the maximum number of distinct vertices of $D$ which may be used as starting and ending points in a partial path decomposition of $D$. The next \lcnamecref{cor:Ntexc} states that $N^\pm(D)$ decreases appropriately when a good $(U^*,W,A)$-partial path decomposition is removed from $D$.}

\begin{cor}\label{cor:Ntexc}
	\NEW{Let $D, W, V', A$, and $U^*$ satisfy the assumptions of \cref{fact:texcS}. Suppose that $\cP$ is a good $(U^*,W,A)$-partial path decomposition. Let $X^+$ be the set of vertices $v\in \tU_{U^*}^+(D)$ for which $\cP$ contains precisely $\texc_{D,U^*}^+(v)$ paths starting at $v$. Similarly, let $X^-$ be the set of vertices $v\in \tU_{U^*}^-(D)$ for which $\cP$ contains precisely $\texc_{D,U^*}^-(v)$ paths ending at $v$. Then, both $N^\pm(D)-N^\pm(D\setminus \cP)\leq |X^\pm|$.}
\end{cor}

\begin{proof}
	\NEW{By symmetry, it suffices to show that $N^+(D)-N^+(D\setminus \cP)=|X^+|$. Denote $D'\coloneqq D\setminus \cP$ and $U^{**}\coloneqq U^*\setminus (V^+(\cP)\cup V^-(\cP))$. Let $Y\coloneqq \tU_{U^*}^+(D)\cup \{v\in W\mid \exc_D^+(v)=d_A^+(v)>0\}$ and $Z\coloneqq \tU_{U^{**}}^+(D')\cup \{v\in W\mid \exc_{D'}^+(v)=d_A^+(v)>0\}$. By assumption, $|U^*|=\texc(D)-\exc(D)$ and so
	\[|Y|\stackrel{\text{\cref{eq:texc}}}{=}|U^+(D)|+|U^*|=N^+(D).\]
	Since~$\cP$ is good, \cref{prop:sizeU*good}
	implies that $|U^{**}|=\texc(D')-\exc(D')$ and so, by the same arguments as above, $|Z|=N^+(D')$.
	Thus, it suffices to show that
	\begin{equation}\label{eq:Ntexc}
		Y=X^+\cup Z.
	\end{equation}
	By definition, $X^+\subseteq Y$.
	Next, we show that $Z\subseteq Y$.
	By \cref{prop:texcv}, $\texc_{D',U^{**}}^+(v)\leq \texc_{D,U^*}^+(v)$ for each $v\in V$ and so $U_{U^{**}}^+(D')\subseteq U_{U^*}^+ (D)\subseteq Y$. By definition of absorbing starting edges (\cref{def:A}), each $v\in W$ satisfies $\exc_D^+(v)\geq d_A^+(v)$ and so \cref{eq:texc} implies that $\{v\in W\mid \exc_{D'}^+(v)=d_A^+(v)>0\}\subseteq Y$. Therefore, $X^+\cup Z\subseteq Y$, as desired.}
	
	\NEW{Finally, we prove that $Y\subseteq X^+\cup Z$. Let $v\in Y\setminus X^+$. It is enough to show that $v\in Z$. Suppose first that $v\in \tU_{U^*}^+(D)$. By \cref{fact:texc} and since $v\notin X^+$, $\cP$ contains fewer than $\texc_{D,U^*}^+(v)$ paths which starts at $v$. Then, \cref{prop:texcv} implies that $v\in \tU_{U^{**}}^+(v)\subseteq Z$. We may therefore assume that $v\in \{v'\in W\mid \exc_D^+(v')=d_A^+(v')>0\}$. By definition of absorbing starting edges (\cref{def:A}), $\exc_D^+(v)>0$ and so \cref{eq:texc} implies that $\texc_{D,U^*}^+(v)=0=\texc_{D,U^*}^-(v)$. Therefore, \cref{fact:texc} implies that no path in $\cP$ contains $v$ as an endpoint and so $\exc_{D'}^+(v)=\exc_D(v)=d_A^+(v)>0$. Thus, $v\in Z$ and so $Y\subseteq X^+\cup Z$. Consequently, \cref{eq:Ntexc} holds and we are done.}
\end{proof}

%% file: Cleaning.tex
	\onlyinsubfile{
		\setcounter{section}{11}
\section{The cleaning step: proof of Lemma \ref{lm:cleaning}}}

We now prove \cref{lm:cleaning} using \cref{lm:cleaningallW*,lm:cleaningW0,lm:cleaningdegreeWA} below.
\NEW{Recall from \cref{sec:cleaninglm} that the main goal of the cleaning step is to reduce the degree at the exceptional set $W$ and cover all the edges inside $W$. Moreover, recall that we denote by $W_*$ the set of exceptional vertices of very large excess, by $W_0$ the set of exceptional vertices with not too large excess, and by $W_A$ the set of exceptional vertices which are incident to some absorbing edges.}
 
In \cref{lm:cleaningW0}, we cover all edges of~$T[W_0]$. The \NEW{remaining} edges of~$T[W]$ which are incident to~$W_*$ will be covered in \cref{lm:cleaningallW*}. Since the excess of $T$ is proportional to~$|W_*|$, the edges of~$T[W]$ which are incident to~$W_*$ can be covered one by one with short paths. However, vertices in~$W_0$ may have small excess and so~$T[W_0]$ needs to be covered more efficiently. The idea will be to decompose~$T[W_0]$ into matchings and then tie them into paths.

In \cref{lm:cleaningallW*}, we also decrease the degree of the vertices in~$W_*\cup W_A$ when~$W_*\neq \emptyset$. This is achieved by covering edges at~$W_*$ one by one with short paths until the desired degree is attained. (The endpoints of these paths are chosen via \cref{lm:endpoints}.)
Finally, we will use \cref{lm:cleaningdegreeWA} to decrease the degree at~$W_A$ when~$W_*=\emptyset$. There, we will use long paths to decrease the degree at all vertices in~$W_A$ at the same time. This is necessary because the total excess may be relatively small, so we do not have room to cover the degree at each vertex in~$W_A$ one by one.

\subsection{\NEW{Proof overview}}

\NEW{The proof structure of \cref{lm:cleaningW0,lm:cleaningallW*,lm:cleaningdegreeWA} is similar. In each of these \lcnamecrefs{lm:cleaningW0}, we need to construct a good partial path decomposition. We always proceed inductively to construct the paths either one by one (\cref{lm:cleaningallW*}) or two by two (\cref{lm:cleaningW0,lm:cleaningdegreeWA}). All these paths are constructed using \cref{cor:robpaths}: we use \cref{cor:robpaths}\cref{cor:robpaths-short} when we need short paths (\cref{lm:cleaningallW*,lm:cleaningdegreeWA}) and we use \cref{cor:robpaths}\cref{cor:robpaths-long} when we need long paths (\cref{lm:cleaningW0,lm:cleaningdegreeWA}).
In each application of \cref{cor:robpaths}, we need to specify two main elements.
\begin{enumerate}[label=\upshape(\roman*)]
	\item We need choose which edges we want to cover and which vertices whose degree we want to decrease. (Roughly speaking, these edges will play the roles of $P_2, \dots, P_{k-1}$ in \cref{cor:robpaths}.) E.g.\ in \cref{lm:cleaningW0}, the aim is to cover all the edges inside $W_0$, so for each path we select which of these edges we want to cover. 
	\item We need to choose ``suitable" endpoints for our paths. (Roughly speaking, these will play the roles of $P_1$ and $P_k$ in \cref{cor:robpaths}.) Choosing these endpoints will form the core of the proof as these need to satisfy several requirements. Firstly, they need to be ``compatible" with the edges we want to cover. E.g\ if we want to cover an edge $uv$, we cannot use $v$ as a starting point. Secondly, they need to have an ``appropriate" amount of excess to ensure that the resulting set of paths will form a good partial path decomposition. The auxiliary excess function defined in \cref{sec:texc} (see \cref{eq:texc}) will enable us to keep track of which vertices are allowed to be used as endpoints in each stage.
\end{enumerate}}

\subsection{Covering the edges inside \texorpdfstring{$W_0$}{W0}}
\NEW{The next \lcnamecref{lm:cleaningW0} states that all the edges inside $W_0$ can be covered by a small good partial path decomposition. The idea is to first decompose $T[W_0]$ into matchings using Vizing's theorem. Then, we incorporate each matching into a pair of (almost) spanning paths using \cref{cor:robpaths}\cref{cor:robpaths-long}. We use very long paths so that the maximum semidegree of $T$ is reduced sufficiently quickly to obtain a good partial path decomposition. Moreover, we require two paths to cover each of the matchings because, by definition of a partial path decomposition, we may only construct paths whose starting and ending points belong to $U^+(T)\cup U^*$ and $U^-(T)\cup U^*$, respectively. Indeed, if for example $M$ is a matching such that each of the vertices in $U^+(T)\cup U^*$ is the ending point of an edge in $M$, then we would not be able to construct a path which contains $M$ and starts in $U^+(T)\cup U^*$. Splitting each matching obtained from Vizing's theorem in two ensures that there are always suitable endpoints to cover each of the submatchings.}

\begin{lm}\label{lm:cleaningW0}
	Let $0<\frac{1}{n}\ll \varepsilon\ll \eta \ll 1$. Let~$T\notin\cT_{\rm excep}$ be a tournament on a vertex set~$V$ of size~$n$ satisfying the following properties.
	\begin{enumerate}
		\item Let $W_*\cup W_0\cup V'$ be a partition of~$V$ such that, for each $w_*\in W_*$, $|\exc_T(w_*)|>(1-20\eta)n$; for each $w_0\in W_0$, $|\exc_T(w_0)|\leq (1-20\eta)n$; and, for each $v'\in V'$, $|\exc_T(v')|\leq \varepsilon n$. Let $W\coloneqq W_*\cup W_0$ and suppose $|W|\leq \varepsilon n$.\label{lm:cleaningW0-W}
		\item Let $A^+,A^-\subseteq E(T)$ be absorbing sets of~$(W, V')$-starting/$(V',W)$-ending edges for~$T$ of size at most~$\lceil\eta n\rceil$. Denote $A\coloneqq A^+\cup A^-$.\label{lm:cleaningW0-A}
		\item Let $U^*\subseteq U^0(T)$ satisfy $|U^*|=\texc(T)-\exc(T)$.
	\end{enumerate} 
	Then, there exists a good $(U^*, W,A)$-partial path decomposition~$\cP$ of~$T$ such that the following hold, where~$D\coloneqq T\setminus \cP$.
	\begin{enumerate}[label=\rm(\roman*)]
		\item $|\cP|\leq 2\varepsilon n$.\label{lm:cleaningW0-size}
		\item $E(D[W_0])=\emptyset$.\label{lm:cleaningW0-cleaningW0}
		\item If $|U^+(D)|=|U^-(D)|=1$, then $e(U^-(D),U^+(D))=0$ or $\texc(D)-\exc(D)\geq 2$.\label{lm:cleaningW0-annoyingD}
		\item Each $v\in U^*\setminus (V^+(\cP)\cup V^-(\cP))$ satisfies $d_D^+(v)=d_D^-(v)\leq \texc(D)-1$.\label{lm:cleaningW0-U*}
	\end{enumerate}
\end{lm}

\NEW{Property \cref{lm:cleaningW0-U*} will ensure that \cref{lm:cleaning}\cref{lm:cleaning-good} is satisfied at the end of the cleaning.}
One can use the \cref{lm:cleaningW0-annoyingD} to ensure that the leftover oriented graph~$D$ does not have all its positive and negative excess concentrated on two vertices~$v^+$ and~$v^-$, respectively, with an edge~$v^-v^+$ between them. Otherwise, we would encounter a similar problem as with the tournaments in the class~$\cT_{\rm apex}$ (recall \cref{prop:annoyingT,fact:annoyingT}). 

\begin{proof}[Proof of \cref{lm:cleaningW0}]
	If $W_0=\emptyset$, then we can set~$\cP\coloneqq \emptyset$ and, by \cref{fact:annoyingT,prop:texcT}, we are done. Thus, we may assume that~$W_0\neq \emptyset$.
	Fix additional constants such that $\varepsilon\ll \nu\ll \tau\ll \eta$.
	Let $W^\pm\coloneqq W\cap U^\pm(T)$ and, for each $\diamond\in \{*,0\}$, denote $W_\diamond^\pm \coloneqq W_\diamond\cap U^\pm(T)$.
	
	Fix a\OLD{n optimal} matching decomposition $M_1, \dots, M_m$ of~$T[W_0]$\OLD{(where~$m\leq \varepsilon n$ by Vizing's theorem)}. \NEW{By Vizing's theorem, we may assume that $m\leq |W_0|\leq \varepsilon n$.}
	Assume inductively that for some $0\leq k\leq m$, we have constructed edge-disjoint paths $P_{1,1},P_{1,2}, P_{2,1} \dots,P_{k,2}\subseteq T$ such that $\cP_k\coloneqq \{P_{i,j}\mid i\in [k], j\in [2]\}$ is a $(U^*,W,A)$-partial path decomposition of~$T$ and the following hold. 
	\begin{enumerate}[label=\upshape(\greek*)]
		\item For each~$i\in [k]$, $E(P_{i,1}\cup P_{i,2})\cap E(T[W])=M_i$.\label{lm:cleaningW0-matching}
		\item For each $i\in [k]$ and $j\in[2]$, $V\setminus W_*\subseteq V(P_{i,j})$\label{lm:cleaningW0-spanning}%
			\COMMENT{Not strictly necessary but an easy way to ensure $\cP$ is a good partial decomposition. An alternative would be covering all vertices that ``need'' to be covered at each stage but that would require more technical details so, since we take very few paths, we may as well take spanning paths all the time. (Taking short paths as we do e.g.\ in \cref{lm:cleaningallW*} would not be good enough if $\texc(T)$ is very small.)}.
		\item For each $\diamond\in \{+,-\}$, if $W_*^\diamond\neq \emptyset$, then $V^\diamond(\cP_k)\subseteq W_*^\diamond$.\label{lm:cleaningW0-W*}
	\end{enumerate}	
	
	Denote $D_k\coloneqq T\setminus \cP_k$. Then, following holds.
	
	\begin{claim}\label{eq:good1}
		We have $\texc(D_k)=\texc(T)-2k\geq 2\eta n$. In particular, $\cP_k$ is a good $(U^*, W, A)$-partial path decomposition of $T$.
	\end{claim}

	\begin{proofclaim}		
		First, note that $\texc(T)-2k\geq 2\eta n$ holds by \cref{prop:texcT} and since $k\leq \varepsilon n$.
		\NEW{By \cref{prop:excpartialdecomp}\cref{prop:excpartialdecomp-general},}\OLD{To show $\texc(D_k)=\texc(T)-2k$, note that, by \cref{prop:excpartialdecomp} and since~$\cP_k$ is a partial path decomposition of~$T$, we have $\texc(D_k)\geq \texc(T)-2k$ and $\exc(D_k)\leq \texc(T)-2k$. Thus,} 
		it is enough to show that $\Delta^0(D_k)\leq \texc(T)-2k$.
		
		If there exists $\diamond\in \{+,-\}$ such that $|W_*^\diamond|\geq 2$, then $\exc(D_k)\geq 2(1-20\eta)n-2k\geq n$ and so $\Delta^0(D_k)\leq \exc(D_k)\leq\texc(T)-2k$, as desired.
		\NEW{We may therefore assume that}\OLD{Suppose} both $|W_*^\pm|\leq 1$.
		\NEW{By \cref{lm:cleaningW0-spanning,lm:cleaningW0-W*}, $P_k$ consists of Hamilton paths. Since $\cP_k$ is a partial path decomposition, each $v\in U^+(T)$ satisfies
		\begin{align*}
			d_{D_k}^+(v)=d_T^+(v)-2k\leq \Delta^0(T)-2k\leq \texc(T)-2k
		\end{align*}
		and
		\begin{align*}
			d_{D_k}^-(v)\leq d_T^-(v)-(2k-\exc_T^+(v))\stackrel{\text{\cref{fact:exc}\cref{fact:exc-excv}}}{=}d_T^+(v)-2k\leq \texc(T)-2k.
		\end{align*}
		Similarly, each $v\in U^-(T)$ satisfies both $d_{D_k}^\pm(v)\leq \texc(T)-2k$.
		Hence, we may assume that $U^0(T)\neq \emptyset$ and so $n$ is odd by definition of $\exc_T(v)$ or \cref{prop:even}. Let $v\in U^0(T)$.}%
		\OLD{If~$v\in U^\pm(T)$, then, by \cref{lm:cleaningW0-spanning}, \cref{lm:cleaningW0-W*}, and since~$\cP_k$ is a $(U^*,W,A)$-partial path decomposition of~$T$, $v\in V^\pm(P)\cup V^0(P)$ for each $P\in \cP_k$ and so, $d_{D_k}^{\max}(v)=d_T^{\max}(v)-2k\leq \texc(T)-2k$.
		Similarly, if $v\in U^0(T)\setminus U^*$, then $v\in V^0(P)$ for each $P\in \cP_k$ and so, $d_{D_k}^{\max}(v)=d_T^{\max}(v)-2k\leq \texc(T)-2k$.
		Thus, if $U^*=\emptyset$, then $\Delta^0(D_k)\leq \texc(T)-2k$, as desired. 
		Suppose there exists $v\in U^*\subseteq U^0(T)$. It is enough to show that both $d_{D_k}^\pm(v)\leq \texc(T)-2k$. Note that, by definition of $U^*$, $\texc(T)>\exc(T)$ and, thus, by \cref{prop:even}, $n$ is odd.}
		By \cref{lm:cleaningW0-spanning}, \cref{prop:texcT}, and \NEW{\cref{def:partialpathdecomp-U0}, $v\in V^\pm(P)\cup V^0(P)$}\OLD{since $\cP_k$ is a $(U^*,W,A)$-partial path decomposition of $T$, $v\in V^+(P)\cup V^0(P)$} for all but at most one path $P\in \cP_k$ and so\OLD{$d_{D_k}^+(v)\leq d_T^+(v)-(2k-1)=\frac{n+1}{2}-2k\leq \texc(T)-2k$. Similarly, $d_{D_k}^-(v)\leq \texc(T)-2k$.}
		\NEW{\begin{align*}
			d_{D_k}^\pm(v)\leq d_T^\pm(v)-(2k-1)=\frac{n+1}{2}-2k\stackrel{\text{\cref{prop:texcT}}}{\leq} \Delta^0(T)-2k\leq \texc(T)-2k.
		\end{align*}}
		Thus, $\Delta^0(D_k)\leq \texc(T)-2k$, as desired.
	\end{proofclaim}
	
	If $k=m$, then let $\cP\coloneqq \cP_m$. \NEW{We verify that all the assertions of \cref{lm:cleaningW0} hold. By construction and \cref{eq:good1}, $\cP$ is a good $(U^*,W,A)$-partial path decomposition of $T$.}\OLD{Then,~$\cP$ is good by \cref{eq:good1}.} Moreover, by construction, $|\cP|=2m\leq 2\varepsilon n$ and~$D\coloneqq T\setminus \cP$ satisfies $E(D[W_0])=\emptyset$ by \cref{lm:cleaningW0-matching}. Thus, \cref{lm:cleaningW0-size,lm:cleaningW0-cleaningW0} hold.	
	For \cref{lm:cleaningW0-annoyingD}, suppose both $|U^\pm(D)|=1$, say $U^\pm(D)=\{u_\pm\}$, and assume that $u_-u_+\in E(D)$. \NEW{We need to show that $\texc(D)-\exc(D)\geq 2$.}
	If there exists $\diamond\in \{+,-\}$ such that~$u_\diamond\notin W_*$, then, \OLD{as $|\cP|\leq 2\varepsilon n$, $d_D^{\min}(u_\diamond)\geq 9\eta n$ and so, by \cref{fact:exc}\cref{fact:exc-excv}, $\texc(D)-\exc(D)\geq 2$.}
	\NEW{\begin{align*}
		\texc(D)-\exc(D)&\stackrel{\text{\eqmakebox[W0][c]{\cref{eq:exc}}}}{=}\texc(D)-\exc_D^\diamond(u_\diamond)\stackrel{\text{\cref{fact:exc}\cref{fact:exc-excv}}}{\geq}d_D^{\min}(u_\diamond)\geq d_T^{\min}(u_\diamond)-|\cP|\\
		&\stackrel{\text{\eqmakebox[W0][c]{\cref{lm:cleaningW0-W},\cref{fact:exc}\cref{fact:exc-dmin}}}}{\geq} \frac{20\eta n-1}{2}-2\varepsilon n\geq 2.
	\end{align*}}
	We may therefore assume that both~$u_\pm\in W_*$. Then, by \cref{lm:cleaningW0-W*}, all paths in~$\cP$ start in~$W_*^+\subseteq U^+(T)$ and end in~$W_*^-\subseteq U^-(T)$. Thus, \NEW{\cref{prop:excP} implies that}\OLD{as~$\cP$ is a partial path decomposition of~$T$, we have} $\exc(D)=\exc(T)-|\cP|$%
		\COMMENT{Each path from $U^+(T)$ to $U^-(T)$ decreases the excess by one. This calculation was done in \cref{prop:excpartialdecomp} (case $k=|S|=0$). Could add a new fact if you want.}
	and, since $k\leq \varepsilon n$ and each~$v\in W_*$ satisfies $|\exc_T(v)|\geq (1-20\eta)n$, we have $U^\pm(D)=U^\pm(T)$. Therefore, by \cref{eq:good1} and since~$T\notin \cT_{\rm excep}$, \cref{fact:annoyingT} implies that $\texc(D)-\exc(D)=\texc(T)-\exc(T)\geq 2$ and so \cref{lm:cleaningW0-annoyingD} holds. 
	Finally, for \cref{lm:cleaningW0-U*}, suppose that $v\in U^*\setminus (V^+(\cP)\cup V^-(\cP))$.
	Then, note that $v\notin W_*$ and, by \cref{prop:even}, $n$ is odd. Thus,\OLD{by \cref{prop:texcT}, \cref{lm:cleaningW0-spanning}, and \cref{eq:good1}, $d_D^+(v)=d_D^-(v)=d_T^-(v)-|\cP|=\frac{n-1}{2}-|\cP|\leq \texc(T)-1-|\cP|\leq \texc(D)-1$}  \[\NEW{d_D^+(v)=d_D^-(v)\stackrel{\text{\cref{lm:cleaningW0-spanning}}}{=}d_T^-(v)-|\cP|=\frac{n-1}{2}-|\cP|\stackrel{\text{\cref{prop:texcT}}}{\leq} \texc(T)-1-|\cP|\stackrel{\text{\cref{eq:good1}}}{=} \texc(D)-1}\]
	and so \cref{lm:cleaningW0-U*} holds.
	
	If $k<m$, then construct~$P_{k+1,1}$ and~$P_{k+1,2}$ as follows.
	Denote $M_{k+1}\coloneqq \{x_1^+x_1^-, \dots, x_\ell^+ x_\ell^-\}$. For each $\diamond\in \{+,-\}$, let $X^\diamond\coloneqq \{x_i^\diamond\mid i\in[\ell]\}$.
	Let $U_k^*\coloneqq U^*\setminus (V^+(\cP_k)\cup V^-(\cP_k))$. Note that, $U_k^*\subseteq U^0(D_k)$.
	Moreover, since $|U^*| = \texc(T) - \exc(T)$, \cref{eq:good1,prop:sizeU*good} imply that $|U_k^*| = \texc(D_k) - \exc(D_k)$.
	
	We claim that there exist suitable endpoints $v_1^\pm$ and $v_2^\pm$ for~$P_{k+1,1}$ and~$P_{k+1,2}$. More precisely, we want to find $v_1^+,v_2^+,v_1^-,v_2^-\in V$ such that the following hold.
	\begin{enumerate}[label=\upshape(\Alph*)]
		\item For each $i\in[2]$, $v_i^+\neq v_i^-$ and $v_i^\pm\in \tU_{U_k^*}^\pm(D_k)$. Moreover, for each~$\diamond\in \{+,-\}$, if~$v_1^\diamond=v_2^\diamond$, then $\texc_{D_k, U_k^*}^\diamond(v_1^\diamond)\geq2$.\label{lm:cleaningW0-endpoints-partialdecomp}
		\item For each $\diamond\in \{+,-\}$, if $W_*^\diamond\neq \emptyset$, then, for each~$i\in [2]$,~$v_i^\diamond\in W_*^\diamond$.\label{lm:cleaningW0-endpoints-W*}
		\item There exists a partition $M_{k+1}=M_{k+1,1}\cup M_{k+1,2}$ such that, for each~$i\in [2]$ and~$x^+x^-\in M_{k+1,i}$, \NEW{we have} $x^\mp\neq v_i^\pm$ and~$x^+x^-\neq v_i^+v_i^-$.
		\label{lm:cleaningW0-endpoints-M}
	\end{enumerate}
	\NEW{Property} \cref{lm:cleaningW0-endpoints-partialdecomp} will ensure that~$\cP_{k+1}$ is a $(U^*, W,A)$-partial path decomposition of~$T$ and \cref{lm:cleaningW0-endpoints-W*} will ensure that \cref{lm:cleaningW0-W*} holds. Finally, \cref{lm:cleaningW0-endpoints-M} will ensure that all edges of~$M_{k+1}$ can be covered by~$P_{k+1,1}\cup P_{k+1,2}$.%
		\COMMENT{Indeed, a $(v^+,v^-)$-path cannot cover an edge ending at $v^+$ nor an edge starting at $v^-$ and, if it covers the edge $v^+v^-$, then it has length one and cannot cover any other edge. Thus, given a partition $M_{k+1}=M_{k+1,1}\cup M_{k+1,2}$ as in \cref{lm:cleaningW0-endpoints-M}, we will be able to construct, for each $i\in [2]$, a $(v_i^+,v_i^-)$-path covering $M_{k+1,i}$, as desired for \cref{lm:cleaningW0-matching}.}
	\NEW{(We will cover $M_{k+1,1}$ with a $(v_1^+,v_1^-)$-path $P_{k+1,1}$ and cover $M_{k+1,2}$ with a $(v_2^+,v_2^-)$-path $P_{k+1,2}$.)}
	
	To find $v_1^+,v_2^+,v_1^-,v_2^-\in V$ satisfying \cref{lm:cleaningW0-endpoints-M,lm:cleaningW0-endpoints-W*,lm:cleaningW0-endpoints-partialdecomp}, we will need the following claim.

	\begin{claim}\label{claim:cleaningW0}
		There exist $v_1^+,v_2^+,v_1^-,v_2^-\in V$ such that, for each~$\diamond\in \{+,-\}$ and~$i,j\in[2]$, the following hold.
		\begin{enumerate}[label=\rm(\Roman*)]
			\item If $W_*^\diamond\neq \emptyset$, then $v_1^\diamond=v_2^\diamond\in W_*^\diamond$ and $\texc_{D_k,U_k^*}^\diamond(v_1^\diamond)\geq2$; otherwise, $v_1^\diamond,v_2^\diamond\in \tU_{U_k^*}^\diamond(D_k)$ are distinct. \label{claim:cleaningW0-W*}
			\item Both $v_1^-v_2^+, v_2^-v_1^+\notin M_{k+1}$.\label{claim:cleaningW0-nocrossedge}
			\item $v_i^+\neq v_j ^-$.\label{claim:cleaningW0-notbothstartend}
		\end{enumerate}
	\end{claim}
	
	\begin{proofclaim}
		If $W_*^+\neq \emptyset$, then pick~$v_1^+\in W_*^+$ and let~$v_2^+\coloneqq v_1^+$ and note that, since~$k\leq \varepsilon n$, $\texc_{D_k, U_k^*}^+ (v_1^+)\geq\exc_T^+(v_1^+)-d_{A^+}(v_1^+)-2k\geq (1-20\eta)n-\lceil\eta n\rceil-2\varepsilon n\geq 2$, as desired.
		Assume that~$W_*^+=\emptyset$. We claim that $|\tU_{U_k^*}^+(D_k)|\geq 2$. Assume not. Then, since $|A^+|\leq \lceil\eta n\rceil$, by \cref{eq:good1,eq:texcA}, $\texc_{U_k^*}^+(D_k)=\texc(D_k)-|A^+|\geq 2\eta n-\lceil\eta n\rceil>1$ and so there exists $v\in V$ such that $\tU_{U_k^*}^+(D_k)=\{v\}$ and $\texc_{D_k,U_k^*}^+(v)\geq2$. Then, note that $v\notin U^0(D_k)$.
		As $v\notin W_*^+$ and $k\leq \varepsilon n$,\OLD{$d_{D_k}^-(v)\geq 9\eta n$.
		Thus, by \cref{fact:exc}\cref{fact:exc-excv}, $\texc(D_k)-\exc_{D_k}^+(v)\geq d_{D_k}^-(v)\geq 9\eta n$}
		\[\NEW{\texc(D_k)-\exc_{D_k}^+(v)\stackrel{\text{\cref{fact:exc}\cref{fact:exc-excv}}}{\geq}d_{D_k}^{\min}(v)\geq d_T^{\min}(v)-|\cP|\stackrel{\text{\cref{lm:cleaningW0-W},\cref{fact:exc}\cref{fact:exc-dmin}}}{\geq} \frac{20\eta n-1}{2}-2\varepsilon n\geq 9\eta n}\]
		 and so, since $v\notin U^0(D_k)$, we have
		\NEW{\begin{align*}
				0&\stackrel{\text{\eqmakebox[claim:cleaningW0][c]{}}}{=}\texc_{D_k,U_k^*}^+(V\setminus \{v\})= \texc_{U_k^*}^+(D_k)-\texc_{D_k, U_k^*}^+(v)\\
				&\stackrel{\text{\eqmakebox[claim:cleaningW0][c]{\text{\cref{eq:texc},\cref{eq:texcA}}}}}{\geq} \texc(D_k)-|A^+|-\exc_{D_k}^+(v)\\
				&\stackrel{\text{\eqmakebox[claim:cleaningW0][c]{}}}{\geq} 9\eta n-\lceil\eta n\rceil\geq 7\eta n,
		\end{align*}}\OLD{\begin{align*}
			0&\stackrel{\text{\eqmakebox[claim:cleaningW0][c]{}}}{=}\texc_{D_k,U_k^*}^+(V\setminus \{v\})= \texc_{U_k^*}^+(D_k)-\texc_{D_k, U_k^*}^+(v)\\
			&\stackrel{\text{\eqmakebox[claim:cleaningW0][c]{\text{\cref{eq:texc},\cref{eq:texcA}}}}}{=} \texc(D_k)-|A^+|-(\exc_{D_k}^+(v)-d_{A^+}(v))\\
			&\stackrel{\text{\eqmakebox[claim:cleaningW0][c]{}}}{\geq} (\texc(D_k)-\exc_{D_k}^+(v)) - |A^+|\geq 9\eta n-\lceil\eta n\rceil\geq 7\eta n,
		\end{align*}}%
		a contradiction. Thus, $|\tU_{U_k^*}^+(D_k)|\geq 2$ and so we can pick distinct $v_1^+,v_2^+\in \tU_{U_k^*}^+(D_k)$.
		
		Then, proceed similarly as above to pick $v_1^-,v_2^-\in \tU_{U_k^*}^-(D_k)\setminus \{v_1^+,v_2^+\}$%
.
		Note that this is possible since, for each~$i\in [2]$, $\texc_{D_k,U_k^*}^-(v_i^+)\leq 1$.
		\NEW{Then, \cref{claim:cleaningW0-W*,claim:cleaningW0-notbothstartend} are satisfied.}
		By relabelling~$v_1^-$ and~$v_2^-$ if necessary, we can ensure \cref{claim:cleaningW0-nocrossedge} holds.
		\NEW{Indeed, suppose that $v_1^-v_2^+\in M_{k+1}$ (the case $v_2^-v_1^+\in M_{k+1}$ is similar). It suffices to show that $v_2^-v_2^+, v_1^-v_1^+\notin M_{k+1}$. 
		Note that, since $V(M_{k+1})\subseteq W_0$, we have $v_1^-,v_2^+\in W_0$. Thus, by \cref{claim:cleaningW0-W*}, $v_1^-\neq v_2^-$ and so, as $M_{k+1}$ is a matching $v_2^-v_2^+\notin M_{k+1}$. Similarly, $v_1^+\neq v_2^+$ and so $v_1^-v_1^+\notin M_{k+1}$.}
		This completes the proof.	  
	\end{proofclaim}
	
	Fix $v_1^+,v_2^+,v_1^-,v_2^-\in V$ satisfying properties \cref{claim:cleaningW0-W*,claim:cleaningW0-notbothstartend,claim:cleaningW0-nocrossedge} of \cref{claim:cleaningW0}.
	We claim that \cref{lm:cleaningW0-endpoints-W*,lm:cleaningW0-endpoints-partialdecomp,lm:cleaningW0-endpoints-M} hold. Indeed, \cref{lm:cleaningW0-endpoints-W*,lm:cleaningW0-endpoints-partialdecomp} follow immediately from \cref{claim:cleaningW0-W*}.
	Recall the notation $M_{k+1}=\{x_i^+x_i^-\mid i\in [\ell]\}$.
	For \cref{lm:cleaningW0-endpoints-M}, let $M_{k+1,2}\coloneqq \{x_i^+x_i^-\in M_{k+1} \mid v_1^+=x_i^- \text{ or } v_1^-=x_i^+\text{ or } v_1^+v_1^-=x_i^+x_i^-\}$ and $M_{k+1,1}\coloneqq M_{k+1}\setminus M_{k+1,2}$. We claim that the partition $M_{k+1}=M_{k+1,1}\cup M_{k+1,2}$ witnesses that \cref{lm:cleaningW0-endpoints-M} holds. By definition,~$M_{k+1,1}$ clearly satisfies the desired properties, so it is enough to show that $M_{k+1,2} \subseteq M_{k+1}\setminus \{x_i^+x_i^-\in M_{k+1} \mid v_2^+=x_i^- \text{ or } v_2^-=x_i^+ \text{ or } v_2^+v_2^-=x_i^+x_i^-\}$.
	If $v_1^+v_1^-=x_i^+x_i^-$ for some~$i\in [\ell]$, then, by \cref{claim:cleaningW0-W*}, \cref{claim:cleaningW0-notbothstartend}, and the fact that $V(M_{k+1})\subseteq W_0$, we have $v_2^+,v_2^-\notin \{x_i^+,x_i^-\}$. Moreover, if $v_1^+=x_i^-$ for some~$i\in [\ell]$, then, by \cref{claim:cleaningW0-W*},~$v_2^+\neq x_i^-$, by \cref{claim:cleaningW0-notbothstartend},~$v_2^-\neq v_1^+$ and so $v_2^+v_2^-\neq x_i^+x_i^-$, and, by \cref{claim:cleaningW0-nocrossedge}, $v_2^-\neq x_i^+$. Similarly, if~$v_1^-=x_i^+$ for some~$i\in [\ell]$, then $v_2^-\neq x_i^+$, $v_2^+\neq x_i^-$, and $v_2^+v_2^-\neq x_i^+x_i^-$.
	Therefore, \cref{lm:cleaningW0-endpoints-M} holds, as desired.
	
	We will now construct, for each~$i\in[2]$, a~$(v_i^+,v_i^-)$-path~$P_{k+1,i}$ covering~$M_{k+1,i}$. The idea is to join together the edges in~$M_{k+1,i}$ via~$V'$. In order to satisfy \cref{lm:cleaningW0-spanning}, we also incorporate the vertices in~$W_0\setminus V(M_{k+1,i})$ in a similar fashion. This will be done using \cref{cor:robpaths}\NEW{\cref{cor:robpaths-long}} as follows. 
	
	Denote \[k'\coloneqq 2+|M_{k+1,1}\setminus\{e\in M_{k+1,1}\mid V^+(e)=\{v_1^+\} \text{ or } V^-(e)=\{v_1^-\}\}|+|W_0\setminus V(M_{k+1,1})|\] ($k'$ will play the role of~$k$ in \cref{cor:robpaths}\NEW{\cref{cor:robpaths-long}}). 
	\NEW{Since $M_{k+1,1}$ is a matching on $W_0$, we have
	\begin{align}\label{eq:cleaningW0-k'}
		k'\leq 2+|W_0|\stackrel{\text{\cref{lm:cleaningW0-W}}}{\leq}2\varepsilon n.
	\end{align}}
	Now construct the~$k'$ paths for \cref{cor:robpaths}\NEW{\cref{cor:robpaths-long}} as follows.
	If $v_1^+\notin V(M_{k+1,1})$, let $Q_1\coloneqq v_1^+$; otherwise, let~$Q_1$ be the (unique) edge~$e\in M_{k+1,1}$ such that $V^+(e)=\{v_1^+\}$. Similarly, if $v_1^-\notin V(M_{k+1,1})$, let $Q_{k'}\coloneqq v_1^-$; otherwise, let~$Q_{k'}$ be the (unique) edge~$e\in M_{k+1,1}$ such that~$V^-(e)=\{v_1^-\}$. Let $Q_2, \dots, Q_{k'-1}$ be an enumeration of $(M_{k+1,1}\setminus \{Q_1, Q_{k'}\})\cup (W_0\setminus V(M_{k+1,1}))$.
	Recall that $V(M_{k+1,1})\subseteq W_0$. Thus, $V'\cap \left(\bigcup_{i\in [k']}V(Q_i)\right)\subseteq \{v_1^+, v_1^-\}$ and so, since $k\leq \varepsilon n$, \cref{lm:3/8rob} implies that $D_k[V'\setminus\bigcup_{i\in [k']}V(Q_i)]$ is a robust~$(\nu, \tau)$-outexpander%
		\COMMENT{\label{com:Dkrob}For each $v\in V'$, $|\exc_T(v)|\leq \varepsilon n$ so $d_{D_k[V'\setminus\bigcup_{j\in [k']}V(Q_j)]}^\pm(v)\geq \frac{n-1-\varepsilon n}{2}-2k-|W|-2\geq \frac{(1-\sqrt{\varepsilon})n}{2}\geq \frac{3|V'|}{7}$.}.
	\NEW{In order to apply \cref{cor:robpaths}, we first need to check that the endpoints of the paths $Q_1, \dots, Q_{k'}$ have sufficiently many neighbours.}
		
	\begin{claim}\label{claim:cleaningW0-endpoints}
		\NEW{For all $i\in [k'-1]$, the ending point $v$ of $Q_i$ satisfies $|N_{D_k\setminus A}^+(v)\cap (V'\setminus \bigcup_{j\in [k']}V(Q_j))|\geq 2k'$ and the starting point $v'$ of $Q_{i+1}$ satisfies $|N_{D_k\setminus A}^-(v')\cap (V'\setminus \bigcup_{j\in [k']}V(Q_j))|\geq 2k'$.}
	\end{claim}
	
	\begin{proofclaim}
		\NEW{Let $i\in [k'-1]$. By symmetry, it is enough to show that the ending point $v$ of $Q_i$ satisfies $N\coloneqq |N_{D_k\setminus A}^+(v)\cap (V'\setminus \bigcup_{j\in [k']}V(Q_j))|\geq 2k'$.}
		
		\NEW{First, observe that $v\in V\setminus W_*^-$. Indeed, $v\in \{v_1^+\}\cup V(M_{k+1})\cup W_0$. By construction,
		$V(M_{k+1})\subseteq W_0$. Moreover,
		\begin{align*}
			v_1^+\stackrel{\text{\cref{lm:cleaningW0-endpoints-partialdecomp}}}{\in} \tU_{U_k^*}^+(D_k)
			\stackrel{\text{\cref{prop:texcv}}}{\subseteq}\tU_{U^*}^+(T)
			\stackrel{\text{\cref{eq:texc}}}{\subseteq} U^+(D)\cup U^*
			\subseteq V\setminus W_*^-.
		\end{align*}
		Thus, $v\in V\setminus W_*^-$.}
		
		\NEW{Since $V(M_{k+1})\subseteq W_0$, we have
		\begin{align*}
			N&\geq d_T^+(v)-k-|A|-|W|-\left|\bigcup_{j\in [k']}V(Q_j)\right|
			\stackrel{\text{\cref{lm:cleaningW0-W},\cref{lm:cleaningW0-A}}}{\geq} d_T^+(v)-\varepsilon n-2\lceil\eta n\rceil-\varepsilon n-2\\
			&\geq d_T^+(v)-3\eta n
		\end{align*}
		and so \cref{eq:cleaningW0-k'} implies that it is enough to show that $d_T^+(v)\geq 4\eta n$.
		If $v\in V\setminus U^-(T)$, then $d_T^+(v)\geq \frac{n-1}{2}$ and so we are done. Suppose that $v\in U^-(T)$. Recall that we have shown that $v\notin W_*^-$ and so
		\begin{align*}
			d_T^+(v)\stackrel{\text{\cref{fact:exc}\cref{fact:exc-dmin}}}{=}\frac{d_T(v)-|\exc_T(v)|}{2}\stackrel{\text{\cref{lm:cleaningW0-W}}}{\geq}\frac{20\eta n-1}{2}\geq 4\eta n.
		\end{align*}
		This completes the proof.}
	\end{proofclaim}
		 	
	\NEW{Thus, all the conditions of \cref{cor:robpaths} are satisfied.} Apply \cref{cor:robpaths}\cref{cor:robpaths-long} with $D_k\setminus A, V\setminus\bigcup_{i\in [k']}V(Q_i), k', \frac{3}{8}, W_*\setminus \{v_1^+,v_1^-\}$, and $Q_1, \dots, Q_{k'}$ playing the roles of $D, V', k, \delta, S$, and $P_1, \dots, P_k$ to obtain a~$(v_1^+,v_1^-)$-path~$P_{k+1,1}$ covering the edges in~$M_{k+1, 1}$ and all vertices in~$V'\cup W_0$.%
	\OLD{Note that all the conditions of \cref{cor:robpaths}\NEW{\cref{cor:robpaths-long}} are satisfied since, by construction, the ending points of $Q_1, \dots, Q_{k'-1}$ lie in $V\setminus W_*$ and so, by  \cref{lm:cleaningW0-W}, \cref{lm:cleaningW0-A}, and \cref{fact:exc}\cref{fact:exc-dmin} they have, in $D_k\setminus A$, at least $\frac{20\eta n\NEW{-1}}{2}-2k-\lceil\eta n\rceil-|W|-2k'\geq 8\eta n\geq 2k'$ outneighbours in $V'\setminus \bigcup_{i\in [k']} V(Q_i)$. Similarly, the starting points of $Q_2,\dots, Q_{k'}$ have, in $D_k\setminus A$, at least~$2k'$ inneighbours in $V'\setminus \bigcup_{i\in [k']} V(Q_i)$, as desired.}
	Construct $P_{k+1,2}$ similarly, but deleting the edges in $P_{k+1,1}$ before applying \cref{cor:robpaths}\NEW{\cref{cor:robpaths-long}} (this will ensure that $P_{k+1,1}$ and $P_{k+1,2}$ are edge-disjoint).
	Thus, \cref{lm:cleaningW0-W*,lm:cleaningW0-matching,lm:cleaningW0-spanning} hold with $k$ replaced by $k+1$. This completes the proof.	
\end{proof}

\subsection{Covering the remaining edges inside \texorpdfstring{$W$}{W} and decreasing the degree at \texorpdfstring{$W_*$}{W*}} 

Since the vertices in~$W_*$ have almost all their edges in the same direction, \NEW{it is not possible to cover}\OLD{one could not have covered} the remaining edges in~$W$ with a similar approach as in \cref{lm:cleaningW0}. \NEW{(In order to incorporate an edge $uv$ into a longer path using \cref{cor:robpaths}, we need $u$ to have many inneighbours and $v$ to have many outneighbours, but this may not be the case if $u\in W_*^+$ or $v\in W_*^-$.)}
However, since the vertices in~$W_*$ have very large excess,~$W_*\neq \emptyset$ implies that the excess of the tournament is very large and so we have room to cover each remaining edge in~$W$ one by one. \NEW{Moreover}, one can also decrease the degree at~$W_*\cup W_A$ with a similar approach. This is achieved in the next \lcnamecref{lm:cleaningallW*}.

\NEW{Note that in \cref{lm:cleaningallW*}\cref{lm:cleaningallW*-defW}, the definition of $W_*$ is adjusted so that the vertex partition $W_*\cup W_0\cup V'$ can be chosen to be the same as in \cref{lm:cleaningW0}.}

\begin{lm}\label{lm:cleaningallW*}
	Let $0<\frac{1}{n}\ll \varepsilon\ll \eta \ll 1$. Let~$D$ be an oriented graph on a vertex set~$V$ of size~$n$ such that $\delta(D)\geq (1-\varepsilon)n$ 
	and the following properties are satisfied.
	\begin{enumerate}
		\item Let $W_*\cup W_0\cup V'$ be a partition of~$V$ such that, for each $w_*\in W_*$, $|\exc_D(w_*)|>(1-21\eta)n$; for each $w_0\in W_0$, $|\exc_D(w_0)|\leq (1-20\eta)n$; and, for each $v'\in V'$, $|\exc_D(v')|\leq \varepsilon n$.
		Let $W\coloneqq W_*\cup W_0$ and suppose $|W|\leq \varepsilon n$ and $W_*\neq \emptyset$.\label{lm:cleaningallW*-defW}
		\item Let $A^+,A^-\subseteq E(D)$ be absorbing sets of~$(W, V')$-starting/$(V',W)$-ending edges for~$D$ of size at most~$\lceil\eta n\rceil$. Denote $A\coloneqq A^+\cup A^-$.
		Let $W_A^\pm \coloneqq V(A^\pm)\cap W$ and $W_A\coloneqq V(A)\cap W$.
		\NEW{Suppose that the following hold.
			\begin{itemize}
				\item Let $\diamond\in \{+,-\}$. If $|W_A^\diamond|\geq 2$, then $\exc_D^\diamond(v)<\lceil\eta n\rceil$ for each~$v\in V$.
				\item Let $\diamond\in \{+,-\}$. If $|W_A^\diamond|=1$, then $\exc_D^\diamond(v)\leq \exc_D^\diamond(w)+\varepsilon n$ for each~$v\in V$ and $w\in W_A^\diamond$.
		\end{itemize}}%
		\OLD{Suppose that, for each $\diamond\in \{+,-\}$, if $|W_A^\diamond|\geq 2$, then, for each~$v\in V$, $\exc_D^\diamond(v)<\eta n$ and, if $|W_A^\diamond|=1$, then, for each~$v\in V$ and $w_A\in W_A^\diamond$, $\exc_D^\diamond(v)\leq \exc_D^\diamond(w_A)+\varepsilon n$.}\label{lm:cleaningallW*-A}
		\item  Let $U^*\subseteq U^0(D)$ satisfy $|U^*|=\texc(D)-\exc(D)$.\label{lm:cleaningallW*-defU*}
		\item Suppose $E(D[W_0])=\emptyset$ and, if $|U^+(D)|=|U^-(D)|=1$, then $e(U^-(D),U^+(D))=0$ or $\texc(D)-\exc(D)\geq 2$.\label{lm:cleaningallW*-annoyingD}
	\end{enumerate}
	Then, there exists a good $(U^*, W, A)$-partial path decomposition~$\cP$ of~$D$ such that $D'\coloneqq D\setminus \cP$ satisfies the following.
	\begin{enumerate}[label=\rm(\roman*)]
		\item $E(D'[W])=\emptyset$.\label{lm:cleaningallW*-cleaningW*}
		\item $N^\pm(D)-N^\pm(D')\leq 88\eta n$.\label{lm:cleaningallW*-exc0}
		\item For each $v\in W_*\cup W_A$, $(1-3\sqrt{\eta})n\leq d_{D'}(v)\leq (1-4\eta)n$.\label{lm:cleaningallW*-degreeW*}
		\item For each $v\in W_0$, $d_{D'}(v)\geq (1-3\sqrt{\eta})n$ and $d_{D'}^{\min}(v)\geq 5\eta n$.\label{lm:cleaningallW*-degreeW0}
		\item For each $v\in V'$, $d_{D'}(v)\geq (1-3\sqrt{\varepsilon})n$.\label{lm:cleaningallW*-degreeV'}
		\item If $|W_*^+|,|W_*^-|\leq 1$, then each~$v\in W_*$ satisfies $|\exc_{D'}(v)|=d_{D'}(v)$.\label{lm:cleaningallW*-exc<2d}
		\item Each $v\in U^*\setminus (V^+(\cP)\cup V^-(\cP))$ satisfies $d_{D'}^+(v)=d_{D'}^-(v)\leq \texc(D')-1$.\label{lm:cleaningallW*-U*}
	\end{enumerate}
\end{lm}

\NEW{Property} \cref{lm:cleaningallW*-exc<2d} will enable us to satisfy \cref{lm:cleaning}\cref{lm:cleaning-exc<2d}.
As mentioned above, the strategy in the proof of \cref{lm:cleaningallW*} is to cover the remaining edges of~$D[W]$ one by one. To decrease the degree at~$W_*$, we further fix some additional edges that will be covered with short paths in the same way. The degree at~$W_A\setminus W_*$ will be decreased by incorporating these vertices in some of these paths.

Similarly as in the proof of \cref{lm:cleaningW0}, given an edge that needs to be covered, we need to find suitable endpoints, that is, endpoints of the correct excess and which are ``compatible'' with the edge \NEW{$e=uv$} that needs to be covered (\NEW{the starting point cannot be $v$ and the ending point cannot be $u$}).
The next \lcnamecref{lm:endpoints} enables us to find, given a set~$H$ of edges to be covered, pairs of suitable endpoints to cover each of these edges with a path.

\begin{lm}\label{lm:endpoints}
	Let $0<\frac{1}{n}\ll \varepsilon\ll \eta \ll 1$. Let~$D$ be an oriented graph on a vertex set~$V$ of size~$n$
	such that \cref{lm:cleaningallW*}\cref{lm:cleaningallW*-defU*,lm:cleaningallW*-defW,lm:cleaningallW*-A,lm:cleaningallW*-annoyingD} are satisfied. 
	Let~$H\subseteq D$ satisfy $\Delta (H) \leq 11\eta n$ and $k\coloneqq |E(H)|\leq 11\eta n |W_*|$.
	Let $w_1^+w_1^-, \dots, w_k^+w_k^-$ be an enumeration of~$E(H)$.
	Then, there exist pairs of vertices $(v^+_1, v^-_1),  \dots, (v^+_k, v^-_k)$ such that the following hold.	
	\begin{enumerate}[label=\rm(\roman*)]
		\item For each $v \in V$ and $\diamond\in \{+,-\}$, there exist at most $\min \{2\sqrt{\eta}n,  \texc_{D,U^*}^\diamond (v)\}$ indices~$i \in [k]$ such that~$v = v_i^\diamond$. \label{lm:endpoints-partialdecomp}
		\item For all $i \in [k]$, if $w_i^{\pm} \in W_*^{\pm}$, then $v ^{\pm}_i = w^{\pm}_i$.\label{lm:endpoints-W*}
		\item For all $i \in [k]$, if there exists $\diamond\in \{+,-\}$ such that $w_i^\diamond \in V'$, then $(v_i^+,v_i^-)\neq (w_i^+,w_i^-)$\label{lm:endpoints-coverW0}%
			\COMMENT{This ensures that we are able to cover vertices of $W_0$ which have high degree in \cref{lm:cleaningallW*}.}.
		\item For each $i \in [k]$, $\{v_i^+, w_i^+\} \cap \{ v_i^-, w_i^-\} = \emptyset$.  \label{lm:endpoints-consistent}
		\item For each $\diamond\in \{+,-\}$, there exist at most~$88\eta n$ vertices \NEW{$v\in \tU_{U^*}^\diamond(D)$}\OLD{$v\in V$} such that there exist exactly $\texc_{D,U^*}^\diamond(v)$ indices~$i\in [k]$ such that~$v_i^\diamond= v$.\label{lm:endpoints-excess0}
		\item \NEW{Denote $V^\pm\coloneqq \{v\in V\mid d_D^\pm(v) \geq \texc(D) - 22\eta n\}$. Then, both $V^\pm\subseteq \{w_i^+, w_i^-, v_i^\pm\}\setminus \{v_i^\mp\}$ for all~$i \in [k]$.}%
		\OLD{If $v \in V$ satisfies $d_D^+(v) \geq \texc(D) - 22\eta n$, then $v \in \{w_i^+, w_i^-, v_i^+\}\setminus \{v_i^-\}$ for all~$i \in [k]$. Similarly, if~$v \in V$ satisfies $d_D^-(v) \geq \texc(D) - 22\eta n$, then $v \in \{w_i^+, w_i^-, v_i^-\}\setminus \{v_i^+\}$ for all $i \in [k]$.}\label{lm:endpoints-good}
	\end{enumerate}
\end{lm}

\NEW{Property \cref{lm:endpoints-partialdecomp} will ensure that a vertex is not used as an endpoint too many times.}
Since vertices in~$W_*$ have most of their edges in the same direction, \cref{lm:endpoints-W*} will ensure that we will be able to tie up edges to the designated endpoints of the path. \NEW{Property} \cref{lm:endpoints-coverW0} will ensure that some of the paths will have length more than one, which will enable us to cover a significant number of edges at~$W_A\setminus W_*$.
\NEW{Property} \cref{lm:endpoints-consistent} implies that the chosen endpoints are distinct, the chosen starting point for the path is not the ending point of the edge we want to cover and, similarly, that the chosen ending point is not the starting point of the edge we want to cover.
\NEW{Moreover, \cref{lm:endpoints-excess0} will ensure that \cref{lm:cleaningallW*}\cref{lm:cleaningallW*-exc0} is satisfied.}
Together with \cref{prop:Delta0}, property \cref{lm:endpoints-good} will ensure that the partial path decomposition constructed with this set of endpoints will be good.
\NEW{Note that the main difficulties in the proof of \cref{lm:endpoints} arise from the cases where $D$ is ``close" to being a tournament from $\cT_{\rm apex}$ (defined in \cref{sec:intro}).}

\NEW{First, we suppose that \cref{lm:endpoints} holds and derive \cref{lm:cleaningallW*}.}

\begin{proof}[Proof of \cref{lm:cleaningallW*}]
	Recall that, by assumption,~$W_*\neq\emptyset$.
	Fix additional constants such that $\varepsilon\ll \nu\ll \tau\ll \eta$.
	Let $W^\pm\coloneqq W\cap U^\pm(D)$ and, for each $\diamond\in \{*,0\}$, denote $W_\diamond^\pm \coloneqq W_\diamond\cap U^\pm(D)$.
	
	We now define a subdigraph $H \subseteq D$, whose edges will be covered by $\cP$. If $\max\{|W_*^+|,|W_*^-|\}\geq 2$, then let $H\subseteq D\setminus A$ be obtained from~$D[W]$ by adding, for each~$v\in W_*$,~$\lceil4\eta n\rceil$ edges of~$D\setminus A$ between~$v$ and~$V'$ (of either direction). Otherwise, let~$H\subseteq D\setminus A$ be obtained from~$D[W]$ by adding, for each~$v\in W_*^\pm$, $\max\{d_D^\mp(v),\lceil4\eta n\rceil\}$ edges of~$D\setminus A$ between~$v$ and~$V'$ (of either direction) such that $d_H^\mp(v)=d_{D\setminus A}^\mp(v)=d_D^\mp(v)$.
	Note that each $v\in V\setminus W_*$ satisfies \OLD{$d_H(v)\leq |W_*|\leq \varepsilon n$.}
	\NEW{\begin{equation}\label{eq:cleaningallW*-dHv}
			d_H(v)\leq |W|\stackrel{\text{\cref{lm:cleaningallW*-defW}}}{\leq} \varepsilon n.
	\end{equation}}
	\NEW{Moreover, each $v\in W_*$ satisfies}
	\OLD{$\Delta(H)\leq |W|+\max\left\{\max_{v\in W^*}d_D^{\min}(v),\lceil4\eta n\rceil\right\}\leq \varepsilon n+\frac{21 \eta n}{2} \leq 11\eta n$.}
	\NEW{\begin{align}
		\lceil4\eta n\rceil \leq d_H(v)&\stackrel{\text{\eqmakebox[DeltaH][c]{}}}{\leq} |W|+\max\left\{\max_{v\in W^*}d_D^{\min}(v),\lceil4\eta n\rceil\right\}\label{eq:cleaningallW*-dHw}\\
		&\stackrel{\text{\eqmakebox[DeltaH][c]{\text{\cref{fact:exc}\cref{fact:exc-dmin}}}}}{\leq} |W|+\max\left\{\max_{v\in W^*}\frac{d_D(v)-|\exc_D(v)|}{2},\lceil4\eta n\rceil\right\}\nonumber\\ 
		&\stackrel{\text{\eqmakebox[DeltaH][c]{\text{\cref{lm:cleaningallW*-defW}}}}}{\leq} \varepsilon n+\frac{21 \eta n}{2}
		\leq 11\eta n.\label{eq:cleaningallW*-DeltaH}
	\end{align}}
	\NEW{Observe that}\OLD{and}
	\begin{equation}\label{eq:k}
		\lceil4\eta n\rceil\leq k\coloneqq |E(H)|\stackrel{\text{\cref{eq:cleaningallW*-DeltaH}}}{\leq} 11\eta n|W_*|.
	\end{equation}
	
	Let $w_1^+w_1^-, \dots, w_k^+w_k^-$ be an enumeration of~$E(H)$. Recall that $W_*\neq \emptyset$ and, for each~$w\in W^*$, $|N_H(w)\cap V'|\geq \lceil 4\eta n\rceil$. Thus, we may assume without loss generality that, for each $i\in [\lceil4\eta n\rceil]$, $w_i^+w_i^-\notin E(D[W])$. 
	Apply \cref{lm:endpoints} to obtain pairs of vertices $(v_1^+,v_1^-), \dots, (v_k^+,v_k^-)$ such that the following hold.
	\begin{enumerate}[label=\upshape(\greek*)]
		\item For each $v \in V$ and $\diamond\in \{+,-\}$, there exist at most $\min \{2\sqrt{\eta}n,  \texc_{D,U^*}^\diamond (v)\}$ indices~$i \in [k]$ such that~$v = v_i^\diamond$. \label{lm:cleaningallW*-partialdecomp}
		\item For all $i \in [k]$, if $w_i^{\pm} \in W_*^{\pm}$, then $v ^{\pm}_i = w^{\pm}_i$.\label{lm:cleaningallW*-W*}
		\newcounter{coverW0}
		\setcounter{coverW0}{\value{enumi}}
		\item For all $i \in [k]$, if there exists $\diamond\in \{+,-\}$ such that $w_i^\diamond \in V'$, then $(v_i^+,v_i^-)\neq (w_i^+,w_i^-)$.\label{lm:cleaningallW*-coverW0}
		\item For each $i \in [k]$, $\{v_i^+, w_i^+\} \cap \{ v_i^-, w_i^-\} = \emptyset$.  \label{lm:cleaningallW*-consistent}
		\item For each $\diamond\in \{+,-\}$, there exist at most~$88\eta n$ vertices \NEW{$v\in \tU_{U^*}^\diamond(D)$}\OLD{$v\in V$} such that there exist exactly $\texc_{D,U^*}^\diamond(v)$ indices~$i\in [k]$ such that~$v_i^\diamond= v$.\label{lm:cleaningallW*-excess0}
		\item \NEW{Denote $V^\pm\coloneqq \{v\in V\mid d_D^\pm(v) \geq \texc(D) - 22\eta n\}$. Then, both $V^\pm\subseteq \{w_i^+, w_i^-, v_i^\pm\}\setminus \{v_i^\mp\}$ for all~$i \in [k]$.}%
		\OLD{If $v \in V$ satisfies $d_D^+(v) \geq \texc(D) - 22\eta n$, then $v \in \{w_i^+, w_i^-, v_i^+\}\setminus \{v_i^-\}$ for all~$i \in [k]$. Similarly, if~$v \in V$ satisfies $d_D^-(v) \geq \texc(D) - 22\eta n$, then $v \in \{w_i^+, w_i^-, v_i^-\}\setminus \{v_i^+\}$ for all $i \in [k]$.}
		\label{lm:cleaningallW*-good}
	\end{enumerate}
	By assumption on our ordering of~$E(H)$, \cref{lm:cleaningallW*-coverW0} implies that the following holds.
	\begin{enumerate}[label=\upshape(\greek*$'$)]
		\setcounter{enumi}{\value{coverW0}}
		\item For all $i \in [\lceil4\eta n\rceil]$, $(v_i^+,v_i^-)\neq (w_i^+,w_i^-)$.\label{lm:cleaningallW*-coverW0'}
	\end{enumerate}
	
	We \NEW{will} now cover each edge~$w_i^+w_i^-$ with a short~$(v_i^+,v_i^-)$-path inductively. \NEW{In the first few paths, we also cover the vertices in $W_A\setminus W_*$ whose degree is too high. More precisely, we proceed} as follows.
	Suppose that for some $0\leq \ell\leq k$ we have constructed edge-disjoint paths $P_1, \dots, P_\ell\subseteq D\setminus A$. For each $0\leq i\leq \ell$, let $D_i\coloneqq D\setminus\bigcup_{j\in[i]}P_i$ and~$S_i$ be the set of vertices $w\in W_A\setminus W_*$ such that $d_{D_i}(w)>(1-4\eta)n$. (Note that~$S_\ell$ corresponds to the set of vertices in~$W_A\setminus W_*$ whose degree is currently too high.) Suppose furthermore that the following hold for each~$i\in [\ell]$.
	\begin{enumerate}[label=\upshape(\Roman*)]
		\item $P_i$ is a $(v_i^+,v_i^-)$-path.\label{lm:cleaningallW*-endpoints}
		\item $w_i^+w_i^-\in E(P_i)$.\label{lm:cleaningallW*-edge}
		\item $S_{i-1}\subseteq V(P_i)$.\label{lm:cleaningallW*-Si}
		\item For each $v\in V'$, there exist at most~$\sqrt{\varepsilon} n$ indices~$j\in [\ell]$ such that $v\in V^0(P_j)\setminus \{w_j^+,w_j^-\}=V(P_j)\setminus\{v_j^+,w_j^+,w_j^-,v_j^-\}$.\label{lm:cleaningallW*-V'}
		\item For each $v\in V(P_i)\cap W$, $v\in \{v_i^+,v_i^-,w_i^+,w_i^-\}\cup S_{i-1}$.\label{lm:cleaningallW*-W}
		\item $e(P_i)\leq 7\nu^{-1}(|S_{i-1}|+1)$.\label{lm:cleaningallW*-short}
	\end{enumerate}
	
	\NEW{First, suppose that $\ell=k$.}\OLD{If $\ell=k$, then} Let $\cP\coloneqq \bigcup_{i\in[\ell]}P_i$ and $D'\coloneqq D_\ell$. 
	
	\begin{claim}
		\NEW{$\cP$ is a good $(U^*,W,A)$-partial path decomposition of~$D$. Moreover, \cref{lm:cleaningallW*-cleaningW*,lm:cleaningallW*-exc<2d,lm:cleaningallW*-degreeW*,lm:cleaningallW*-degreeW0,lm:cleaningallW*-degreeV',lm:cleaningallW*-U*,lm:cleaningallW*-exc0} are satisfied.}
	\end{claim}
	
	\begin{proofclaim}
		\NEW{By assumption, $\cP\subseteq D\setminus A$. Moreover, \cref{lm:cleaningallW*-partialdecomp} and \cref{lm:cleaningallW*-endpoints} imply that each $v\in V$ is the starting point of at most $\texc_{D, U^*}^+(v)$ paths in $\cP$ and the ending point of at most $\texc_{D, U^*}^-(v)$ paths in $\cP$. Thus, \cref{fact:texc} implies that $\cP$ is a $(U^*,W,A)$-partial path decomposition.}
		
		\NEW{We now verify \cref{lm:cleaningallW*-exc0}. By \cref{lm:cleaningallW*-excess0} and \cref{lm:cleaningallW*-endpoints}, there are at most $88\eta n$ vertices $v\in \tU_{U^*}^+(D)$ for which $\cP$ contains precisely $\texc_{D,U^*}^+(v)$ paths starting at $v$ and at most $88\eta n$ vertices $v'\in \tU_{U^*}^-(D)$ for which $\cP$ contains precisely $\texc_{D,U^*}^-(v')$ paths ending at $v'$. Thus, \cref{lm:cleaningallW*-exc0} follows from \cref{cor:Ntexc}.}
		
		\NEW{Next, we show that $\cP$ is good. By \cref{lm:cleaningallW*-good}, \cref{lm:cleaningallW*-endpoints}, and \cref{lm:cleaningallW*-edge}, both $V^\pm\subseteq V^\pm(P_i)\cup V^0(P_i)$ for each $i\in [k]$. If $k\leq 22\eta n$, then \cref{prop:Delta0}\cref{prop:Delta0-good} implies that $\cP$ is good. We may therefore assume that $k>22\eta n$. Then, \cref{eq:k} implies that $\max\{|W_*^+|,|W_*^-|\}\geq 2$ and so, by \cref{lm:cleaningallW*-defW} and \cref{eq:exc}, we have $\exc(D)\geq 2(1-21\eta)n\geq n\geq \Delta^0(D)$. Therefore, $\texc(D)=\exc(D)$ and so \cref{prop:excpartialdecomp}\cref{prop:excpartialdecomp-exc} implies that 
		\[\exc(D')=\exc(D)-|\cP|\stackrel{\text{\NEW{\cref{lm:cleaningallW*-defW},\cref{eq:exc}},\cref{eq:k}}}{\geq} \max\{|W_*^+|,|W_*^-|\}(1-21\eta)n-11\eta n|W_*|\geq n\geq \Delta^0(D').\]
		Thus, $\texc(D')=\exc(D')=\exc(D)-|\cP|=\texc(D)-|\cP|$ and so~$\cP$ is good.}
	
		\NEW{By \cref{lm:cleaningallW*-edge}, $\cP$ covers all the edges of $H$. By construction, $H$ contains all the edges of $D$ which lie inside $W$ and so \cref{lm:cleaningallW*-cleaningW*} holds. Moreover, if both $|W_*^\pm|\leq 1$, then by construction $H$ contains all the edges which end in $W_*^+$ as well as all the edges which start in $W_*^-$, so \cref{lm:cleaningallW*-exc<2d} is satisfied.}
		
		\NEW{We now verify \cref{lm:cleaningallW*-degreeW*}. Let $v\in W_*\cup W_A$. If $v\in W_*$, then \cref{eq:cleaningallW*-dHw} implies that $d_H(v)\geq \lceil4\eta n\rceil$. Thus, \cref{lm:cleaningallW*-edge} implies that the upper bound in \cref{lm:cleaningallW*-degreeW*} holds if $v\in W_*$. If $v\in W_A$, then \cref{lm:cleaningallW*-Si} implies that $v\notin S_\ell$ for each $\ell\geq 4\eta n$. In particular, \cref{eq:k} implies that $v\notin S_k$ and so the upper bound in \cref{lm:cleaningallW*-degreeW*} also holds.
		Moreover, \cref{lm:cleaningallW*-endpoints}, \cref{lm:cleaningallW*-W}, \cref{eq:cleaningallW*-dHv}, and \cref{eq:cleaningallW*-DeltaH} imply that there are at most $d_H(v)+4\eta n\leq 15\eta n$ paths in $\cP$ which contain $v$ as an internal vertex.
		By \cref{lm:cleaningallW*-partialdecomp} and \cref{lm:cleaningallW*-endpoints}, there are at most $2\sqrt{\eta}n$ paths in $\cP$ which have $v$ as an endpoint. 
		Thus,
		\[d_{D'}(v)\geq d_D(v)-2\sqrt{\eta}n-30\eta n\geq (1-\varepsilon)n-2\sqrt{\eta}n-30\eta n\geq (1-3\sqrt{\eta})n,\]
		and so the lower bound in \cref{lm:cleaningallW*-degreeW*} holds.
		Therefore, \cref{lm:cleaningallW*-degreeW*} is satisfied.}
		
		\NEW{Next, we verify \cref{lm:cleaningallW*-degreeW0}. Let $v\in W_0$. 
		By the same arguments as for \cref{lm:cleaningallW*-degreeW*}, there are at most $2\sqrt{\eta}n$ paths in $\cP$ which contain $v$ as an endpoint and at most \[d_H(v)+4\eta n\stackrel{\text{\cref{eq:cleaningallW*-dHv}}}{\leq} (4\eta +\varepsilon)n\] paths which contain $v$ as an internal vertex.
		Thus,
		\[d_{D'}(v)\geq d_D(v)-2\sqrt{\eta}n-2(4\eta +\varepsilon)n\geq (1-\varepsilon)n-(2\sqrt{\eta}+8\eta +2\varepsilon)n\geq (1-3\sqrt{\eta})n,\]
		as desired.
		It remains to show that $d_{D'}^{\rm min}(v)\geq 5\eta n$. Suppose without loss of generality that $d_D^+(v)\geq d_D^-(v)$, i.e.\ that $d_D^{\rm min}(v)=d_D^-(v)$. Then, $v\in U^+(D)\cup U^0(D)$, so \cref{def:partialpathdecomp-exc,def:partialpathdecomp-U0} imply that $\cP$ contains at most $\max\{\exc_D^+(v),1\}$ paths which start at $v$ and at most one path which ends at $v$. Therefore, $d_{D'}^+(v)\geq d_{D'}^-(v)-1$ and so it is enough to show that $d_{D'}^-(v)>5\eta n$. Since there is at most one path in $\cP$ which ends at $v$ and at most $(4\eta +\varepsilon)n$ paths in $\cP$ which contain $v$ as an internal vertex, we have
		\begin{align*}
			d_{D'}^-(v)&\stackrel{\text{\eqmakebox[dD'-][c]{}}}{\geq} d_D^-(v)-1-(4\eta +\varepsilon)n\geq d_D^{\rm min}(v)-(4\eta +2\varepsilon)n\\
			&\stackrel{\text{\eqmakebox[dD'-][c]{\text{\cref{fact:exc}\cref{fact:exc-dmin}}}}}{=}\frac{d_D(v)-|\exc_D(v)|}{2}-(4\eta +2\varepsilon)n
			\stackrel{\text{\cref{lm:cleaningallW*-defW}}}{\geq}\frac{(1-\varepsilon)n-(1-20\eta)n}{2}-(4\eta +2\varepsilon)n\\
			&\stackrel{\text{\eqmakebox[dD'-][c]{}}}{>}5\eta n.
		\end{align*}
		Thus, \cref{lm:cleaningallW*-degreeW0} holds.}
	
		\NEW{For \cref{lm:cleaningallW*-degreeV'}, let $v\in V'$. By \cref{lm:cleaningallW*-defW}, $|\exc_D(v)|\leq \varepsilon n$ and so, as $\cP$ is a partial path decomposition, $v$ is an endpoint of at most $\varepsilon n$ paths in $\cP$. Moreover, \cref{lm:cleaningallW*-V'} implies that there are at most \[\sqrt{\varepsilon}n+d_H(v)\stackrel{\text{\cref{eq:cleaningallW*-dHv}}}{\leq} (\sqrt{\varepsilon}+\varepsilon)n\]
		paths in $\cP$ which contain $v$ as an internal vertex.
		Thus,
		\[d_{D'}(v)\geq d_D(v)-\varepsilon n- 2(\sqrt{\varepsilon}+\varepsilon)n\geq (1-\varepsilon)n-(2\sqrt{\varepsilon}+3\varepsilon)n\geq (1-3\sqrt{\varepsilon})n\]
		and so \cref{lm:cleaningallW*-degreeV'} holds.}
	
		\NEW{Finally, we verify \cref{lm:cleaningallW*-U*}. Let $v\in U^*\setminus (V^+(\cP)\cup V^-(\cP))$. By \cref{lm:cleaningallW*-defU*}, $v\in U^0(D)$ and so $d_{D'}^+(v)=d_{D'}^-(v)=\frac{d_{D'}(v)}{2}\leq \frac{n-1}{2}$. Thus, it is enough to show that $\texc(D')\geq \frac{n+1}{2}$. By assumption, there exists $w\in W_*$ and so
		\begin{align*}
			\texc(D')\geq \exc(D')\stackrel{\text{\cref{eq:exc}}}{\geq} |\exc_{D'}(w)|\stackrel{\text{\cref{lm:cleaningallW*-partialdecomp},\cref{lm:cleaningallW*-endpoints}}}{\geq} |\exc_D(w)|-2\sqrt{\eta}n \stackrel{\text{\cref{lm:cleaningallW*-defW}}}{\geq}(1-21\eta)n-2\sqrt{\eta}n>\frac{n+1}{2}.
		\end{align*}
		Thus, \cref{lm:cleaningallW*-U*} holds.}
	\end{proofclaim}

	\NEW{We may therefore assume that}\OLD{Suppose} $\ell< k$. Note that, by \cref{lm:cleaningallW*-Si}, if $\lceil4\eta n\rceil <i\leq \ell$, then~$S_i=\emptyset$.
	We construct $P_{\ell+1}$ as follows. If $(v_{\ell+1}^+,v_{\ell+1}^-)=(w_{\ell+1}^+,w_{\ell+1}^-)$, then let $P_{\ell+1}\coloneqq w_{\ell+1}^+w_{\ell+1}^-$. 
	Note that, in this case, by \cref{lm:cleaningallW*-coverW0'},~$S_\ell=\emptyset$. Thus, \cref{lm:cleaningallW*-Si,lm:cleaningallW*-V',lm:cleaningallW*-W*,lm:cleaningallW*-edge,lm:cleaningallW*-endpoints,lm:cleaningallW*-short,lm:cleaningallW*-W} hold with~$\ell+1$ playing the role of~$\ell$ and we are done.
	We may therefore assume that $(v_{\ell+1}^+,v_{\ell+1}^-)\neq(w_{\ell+1}^+,w_{\ell+1}^-)$. We construct~$P_{\ell+1}$ using \cref{cor:robpaths} as follows.
	Let~$X$ be the set of vertices $v\in V'\setminus\{v_{\ell+1}^+,v_{\ell+1}^-,w_{\ell+1}^+,w_{\ell+1}^-\}$ such that there exist $\lfloor\sqrt{\varepsilon}n\rfloor$ indices~$i\in [\ell]$ such that $v\in V^0(P_i)\setminus \{w_i^+,w_i^-\}$. 
	Note that, by \cref{lm:cleaningallW*-short}, 
	\begin{equation}\label{eq:cleaningallW*-X}
		|X|\leq \frac{7\nu^{-1}(\sum_{i\in [\ell]}|S_i|+2\ell)}{\lfloor\sqrt{\varepsilon}n\rfloor}\leq\frac{8\nu^{-1}(\lceil4\eta n\rceil\cdot\varepsilon n+22\eta n\cdot\varepsilon n)}{\sqrt{\varepsilon}n}\leq \varepsilon^{\frac{1}{3}}n.
	\end{equation}
	Recall that each~$v\in V'$ satisfies $|\exc_D(v)|\leq \varepsilon n$. Thus, by \cref{lm:3/8rob}, \cref{lm:cleaningallW*-partialdecomp}, \cref{lm:cleaningallW*-endpoints}, and \cref{lm:cleaningallW*-V'}, $D_\ell[V'\setminus (X\cup \{v_{\ell+1}^+, w_{\ell+1}^+,w_{\ell+1}^-,v_{\ell+1}^-\})]$ is a robust~$(\nu, \tau)$-outexpander.
	\NEW{The idea is to use \cref{cor:robpaths}\cref{cor:robpaths-short} to tie together the edge $w_{\ell+1}^+w_{\ell+1}^1$ and the vertices in $S_\ell$ into a short $(v_{\ell+1}^+,v_{\ell+1}^-)$-path via the vertices in $V'\setminus X$.}
	
	Let $u_1,\dots, u_s$ be an enumeration of $S_\ell\setminus \{v_{\ell+1}^+, w_{\ell+1}^+,w_{\ell+1}^-,v_{\ell+1}^-\}$.
	If both $v_{\ell+1}^\pm\neq w_{\ell+1}^\pm$, then let $m\coloneqq s+3$, $Q_1\coloneqq v_{\ell+1}^+$, $Q_2\coloneqq w_{\ell+1}^+w_{\ell+1}^-$, $Q_m\coloneqq v_{\ell+1}^-$, and, for each~$i\in [s]$, let $Q_{i+2}\coloneqq u_i$ \NEW{($m, Q_1, \dots, Q_m$ will play the roles of $k, P_1, \dots, P_k$ in \cref{cor:robpaths}\cref{cor:robpaths-short})}.
	If $v_{\ell+1}^+\neq w_{\ell+1}^+$ and $v_{\ell+1}^-= w_{\ell+1}^-$, then let $m\coloneqq s+2$, $Q_1\coloneqq v_{\ell+1}^+$, $Q_m\coloneqq w_{\ell+1}^+w_{\ell+1}^-$, and, for each~$i\in [s]$, let $Q_{i+1}\coloneqq u_i$.
	Similarly, if $v_{\ell+1}^+= w_{\ell+1}^+$ and $v_{\ell+1}^-\neq w_{\ell+1}^-$, then let $m\coloneqq s+2$, $Q_1\coloneqq w_{\ell+1}^+w_{\ell+1}^-$, $Q_m\coloneqq v_{\ell+1}^-$, and, for each~$i\in [s]$, let $Q_{i+1}\coloneqq u_i$.
	Note that, by \cref{lm:cleaningallW*-consistent}, this covers all possible cases.
	Moreover, \NEW{\cref{lm:cleaningallW*-defW} implies that} we always have
	\begin{equation}\label{eq:cleaningallW*-m}
		m\leq |S_\ell|+3\NEW{\leq |W|+3\leq 2\varepsilon n}.
	\end{equation}
	\NEW{In order to apply \cref{cor:robpaths}, we first need to check that that endpoints of the paths $Q_1, \dots, Q_m$ have sufficiently many neighbours. The proof is similar to that of \cref{claim:cleaningW0-endpoints} in the proof of \cref{lm:cleaningW0}.}
	
	\begin{claim}\label{claim:cleaningallW*-endpoints}
		\NEW{For each $i\in [m-1]$, the ending point $v$ of $Q_i$ satisfies $|N_{D_\ell\setminus A}^+(v)\cap (V'\setminus (\bigcup_{j\in [m]}V(Q_j)\cup X))|\geq 2m$ and the starting point $v'$ of $Q_{i+1}$ satisfies $|N_{D_\ell\setminus A}^-(v')\cap (V'\setminus (\bigcup_{j\in [m]}V(Q_j)\cup X))|\geq 2m$.}
	\end{claim}
	
	\begin{proofclaim}
		\NEW{Let $i\in [m-1]$. By symmetry, it is enough to show that the ending point $v$ of $Q_i$ satisfies $N\coloneqq|N_{D_\ell\setminus A}^+(v)\cap (V'\setminus (\bigcup_{j\in [m]}V(Q_j)\cup X))|\geq 2m$.}
		
		\NEW{First, we show that $v\in V\setminus W_*^-$. 
		By construction, $S_\ell\subseteq W\setminus W_*$. We may therefore assume that $v\in \{v_i^+,w_i^-\}$. 
		By \cref{lm:cleaningallW*-partialdecomp}, $v_i^+\in \tU_{U^*}^+(v)\subseteq V\setminus U^-(D)\subseteq V\setminus W_*^-$. We may therefore assume that $v=w_i^-$. Since $i<m$, we have $w_i^-\neq v_i^-$ and so 
		\cref{lm:cleaningallW*-W*} implies that $v\notin W_*^-$.
		Thus, $v\in V\setminus W_*^-$, as desired.}
		
		\NEW{Next, observe that
		\begin{align*}
			N&\stackrel{\text{\eqmakebox[N+v1][c]{}}}{\geq} d_D^+(v)-d_{\cP_\ell}^+(v)-|A|-|W|-\left|V'\cap\bigcup_{j\in [m]}V(Q_j)\right|-|X|\\
			&\stackrel{\text{\eqmakebox[N+v1][c]{\text{\cref{lm:cleaningallW*-defW},\cref{lm:cleaningallW*-A},\cref{eq:cleaningallW*-X}}}}}{\geq} d_D^+(v)-d_{\cP_\ell}^+(v)-2\lceil\eta n\rceil-\varepsilon n-2-\varepsilon^{\frac{1}{3}}n \geq d_D^+(v)-d_{\cP_\ell}^+(v)-3\eta n
		\end{align*}
		and so \cref{eq:cleaningallW*-m} implies that it is enough to show that $d_D^+(v)-d_{\cP_\ell}^+(v)\geq \frac{7\eta n}{2}$.}

		\NEW{Suppose first that $v\in W\setminus U^-(D)$. Then, note that $d_D^+(v)=d_D^{\max}(v)$. By \cref{lm:cleaningallW*-partialdecomp}, there are at most $2\sqrt{\eta}n$ indices $j\in [\ell]$ for which $v\in \{v_j^+,v_j^-\}$. Moreover, \cref{lm:cleaningallW*-Si} implies that there are at most $4\eta n$ indices $j\in [\ell]$ for which $v\in S_{j-1}$. Therefore,
		\begin{align*}
			d_D^+(v)-d_{\cP_\ell}^+(v)&\stackrel{\text{\eqmakebox[N+v2][c]{\text{\cref{lm:cleaningallW*-W}}}}}{\geq} \frac{d_D(v)}{2}-(2\sqrt{\eta}n+d_H(v)+4\eta n)\\
			& \stackrel{\text{\eqmakebox[N+v2][c]{\text{\cref{eq:cleaningallW*-DeltaH}}}}}{\geq} \frac{(1-\varepsilon)n}{2}-(2\sqrt{\eta}n+11\eta n+4\eta n)
			\geq \frac{7\eta n}{2}.
		\end{align*}
		Next, suppose that $v\in W\cap U^-(D)$. Note that $d_D^+(v)=d_D^{\rm \min}(v)$. By \cref{lm:cleaningallW*-partialdecomp}, there is no index $j\in [\ell]$ for which $v=v_j^+$. Moreover, \cref{lm:cleaningallW*-Si} implies that there are at most $4\eta n$ indices $j\in [\ell]$ for which $v\in S_{j-1}$. 
		Recall from \cref{lm:cleaningallW*-endpoints} that $v_j^-$ is the ending point of $P_j$ for each $j\in [\ell]$. Moreover, we have shown that $v\notin W_*^-$ and so we have $v\in W_0$.
		Thus,
		\begin{align*}
			d_D^+(v)-d_{\cP_\ell}^+(v)&\stackrel{\text{\eqmakebox[N+v3][c]{\text{\cref{fact:exc}\cref{fact:exc-dmin},\cref{lm:cleaningallW*-W}}}}}{\geq} \frac{d_D(v)-|\exc_D(v)|}{2}-(d_H(v)+4\eta n)\\
			& \stackrel{\text{\eqmakebox[N+v3][c]{\text{\cref{lm:cleaningallW*-defW}},\cref{eq:cleaningallW*-dHv}}}}{\geq} \frac{(20\eta-\varepsilon)n}{2}-(\varepsilon n+4\eta n)
			\geq \frac{7\eta n}{2}.
		\end{align*}
		We may therefore assume that $v\in V'$. Then, \cref{lm:cleaningallW*-defW}, \cref{lm:cleaningallW*-partialdecomp}, and \cref{lm:cleaningallW*-endpoints} imply that there are at most $\varepsilon n$ indices $j\in [\ell]$ such that $v$ is the starting point of $P_j$. Moreover, \cref{lm:cleaningallW*-V'} implies that there are at most $\sqrt{\varepsilon}n$ indices $j\in [\ell]$ for which $v$ is an internal vertex of $P_j$. Thus,
		\begin{align*}
			d_D^+(v)-d_{\cP_\ell}^+(v)&\stackrel{\text{\cref{fact:exc}\cref{fact:exc-dmin}}}{\geq} \frac{d_D(v)-|\exc_D(v)|}{2}-(\varepsilon n+\sqrt{\varepsilon}n)
			\stackrel{\text{\cref{lm:cleaningallW*-defW}}}{\geq} \frac{(1-2\varepsilon)n}{2}-2\sqrt{\varepsilon} n
			\geq \frac{7\eta n}{2}.
		\end{align*}
		This completes the proof of \cref{claim:cleaningallW*-endpoints}.}
	\end{proofclaim}

	\NEW{Thus, all the conditions of \cref{cor:robpaths} are satisfied. Let $S_\ell'\coloneqq S_\ell\cup \{v_{\ell+1}^+, w_{\ell+1}^+,w_{\ell+1}^-,v_{\ell+1}^-\}$.} Apply \cref{cor:robpaths}\cref{cor:robpaths-short} with \NEW{$D_\ell\setminus A, V\setminus S_\ell', m, \frac{3}{8}, X\cup (W\setminus S_\ell')$}\OLD{$D_\ell\setminus A, V\setminus (S_\ell\cup \{v_{\ell+1}^+, w_{\ell+1}^+,w_{\ell+1}^-,v_{\ell+1}^-\}), m, \frac{3}{8}, X\cup (W\setminus(S_\ell\cup \{v_{\ell+1}^+, w_{\ell+1}^+,w_{\ell+1}^-,v_{\ell+1}^-\}))$}, and $Q_1, \dots, Q_m$ playing the roles of $D,V',k, \delta, S$, and $P_1, \dots, P_k$ to obtain a~$(v_{\ell+1}^+,v_{\ell+1}^-)$-path~$P_{\ell+1}$ of length at most $2\nu^{-1}m+1\leq 2\nu^{-1}(|S_\ell|+3)+1\leq 7\nu^{-1}(|S_\ell|+1)$ which covers~$w_{\ell+1}^+w_{\ell+1}^-$ and the vertices in~$S_\ell$ and avoids the vertices in \NEW{$X\cup (W\setminus S_\ell')$}\OLD{$X\cup (W\setminus(S_\ell\cup \{v_{\ell+1}^+, w_{\ell+1}^+,w_{\ell+1}^-,v_{\ell+1}^-\}))$}. 

	One can easily verify that \cref{lm:cleaningallW*-Si,lm:cleaningallW*-V',lm:cleaningallW*-W*,lm:cleaningallW*-edge,lm:cleaningallW*-endpoints,lm:cleaningallW*-short,lm:cleaningallW*-W} hold with~$\ell+1$ playing the role of~$\ell$.
\end{proof}

\NEW{We now prove \cref{lm:endpoints}.}

\begin{proof}[Proof of \cref{lm:endpoints}]
	Let $W^\pm\coloneqq W\cap U^\pm(D)$ and, for each $\diamond\in \{*,0\}$, denote $W_\diamond^\pm \coloneqq W_\diamond\cap U^\pm(D)$.
	Observe that the following holds.
	
	\begin{claim}\label{property:notannoyingD}
		There are no distinct $v_+,v_-,v_0\in V$ such that $v_+v_-\in E(D)$ and both $\tU_{U^*}^\pm(D)=\{v_\mp,v_0\}$.
	\end{claim}

	\begin{proofclaim}
		Suppose for a contradiction that $v_+,v_-,v_0\in V$ are distinct and such that $v_+v_-\in E(D)$ and both $\tU_{U^*}^\pm(D)=\{v_\mp,v_0\}$. We now show that $U^\pm(D)=\{v_\mp\}$, which implies that $|U^\pm(D)|=1$, $e(U^-(D),U^+(D))\neq 0$ and $\texc(D)-\exc(D)=|U^*|=|\tU_{U^*}^+(D)\cap \tU_{U^*}^-(D)|=1<2$, a contradiction to \cref{lm:cleaningallW*-annoyingD}.
		\NEW{By \cref{lm:cleaningallW*-defW},}\OLD{Note that, since $W_*\neq \emptyset$,} $\exc(D)\geq (1-21\eta)n$. Moreover, since $v_\mp\in \tU_{U^*}^\pm(D)\setminus \tU_{U^*}^\mp(D)$, we have $v_\mp\in U^\pm(D)$.
		Thus,\OLD{by \cref{fact:texcS} and \cref{eq:texcA},} 
		\begin{align*}
			\exc_D^\pm(v_\mp)&\stackrel{\text{\NEW{\cref{eq:texc}}}}{\geq} \texc_{D, U^*}^\pm(v_\mp)\stackrel{\text{\NEW{\cref{fact:texcS}}}}{=}\texc_{U^*}^\pm(D)-\texc_{D,U^*}^\pm(v_0)\\
			&\stackrel{\text{\NEW{\cref{eq:texcA}}}}{\geq} \NEW{(\exc(D)-|A^\pm|)-1
			\stackrel{\text{\cref{lm:cleaningallW*-A}}}{\geq}}(\exc(D)-\lceil\eta n\rceil)-1\geq 2\eta n.
		\end{align*} 
		Thus,\OLD{by} \cref{lm:cleaningallW*-A} \NEW{implies that} $|W_A^\pm|\leq 1$.% 
		\OLD{Moreover, $W_A^\pm\subseteq \{v_\mp\}$ since, if $v\in W_A^\pm \setminus \{v_\mp\}$, then $v\notin \tU_{U^*}^\pm(D)$ and so, by \cref{lm:cleaningallW*-A}, $\exc_D^\pm(v_\mp)\leq \exc_D^\pm(v)+\varepsilon n\leq \lceil\eta n\rceil+\varepsilon n< 2\eta n$, a contradiction. Thus, $U^\pm(D)=\tU_{U^*}^\pm(D)\setminus U^*=\{v_\mp\}$, as desired.} 
		\NEW{If both $W_A^\pm\subseteq \{v_\mp\}$, then \cref{eq:exc} implies that both $U^\pm(D)=\tU_{U^*}^\pm(D)\setminus U^*=\{v_\mp\}$ and so we are done. We may therefore assume without loss of generality that there exists $v\in W_A^+\setminus \{v_-\}$. We find a contradiction. By \cref{def:A}, we have $v\notin U^*$ and so $v\notin \tU_{U^*}(D)$. Thus, $\texc_{D,U^*}^+(v)=0$ and so \cref{lm:cleaningallW*-A} implies that
		\begin{align*}
			\exc_D^\pm(v_\mp)\leq \exc_D^\pm(v)+\varepsilon n\stackrel{\text{\cref{eq:texc}}}{\leq} |A^+|+\varepsilon n\leq \lceil\eta n\rceil+\varepsilon n< 2\eta n,
		\end{align*}
		a contradiction. This completes the proof of \cref{property:notannoyingD}.}
	\end{proofclaim}
	
	Let $\tW_*^\pm\coloneqq \{v\in V\mid \exc_D^\pm(v)\geq (1-86\eta)n\}$. For technical reasons, we will ensure that, for any $i\in [k]$ and $\diamond\in \{+,-\}$, if $w_i^\diamond\in \tW_*^\diamond$, then $v_i^\diamond=w_i^\diamond$. Note that this will imply \cref{lm:endpoints-W*}, as $W_*^\pm\subseteq \tW_*^\pm$.
	Without loss of generality, we may assume that \NEW{$E(H)$ is ordered so that}, if $|\tW_*^+|=|\tW_*^-|=1$ and there exists $i\in [k]$ such that $w_i^+\in \tW_*^-$ and $w_i^-\in \tW_*^+$, then~$i=1$.
	
	Suppose that, for some $0\leq \ell\leq k$, we have already constructed pairs $(v^+_1, v^-_1),  \dots, (v^+_\ell, v^-_\ell)$ such that the following hold.
	For each~$v\in V$, define $\hexc_\ell^\pm(v)\coloneqq \texc_{D, U^*}^\pm(v)-|\{i\in [\ell]\mid v_i^\pm=v\}|$. Denote $\hU_\ell^\pm(D)\coloneqq \{v\in V\mid \hexc_\ell^\pm(v)>0\}$.	
	\begin{enumerate}[label=\upshape(\greek*)]
		\item For each $v \in V$, $\hexc_\ell^\pm(v)\geq 0$.\label{lm:endpoints-IH-partialdecomp}
		\item For each $v \in V$ and $\diamond\in\{+,-\}$, there exist at most $\sqrt{\eta}n$ indices $i \in [\ell]$ such that $v_i^\diamond =v\neq w_i^\diamond$. \label{lm:endpoints-IH-forbidden}
		\item For all $i \in [\ell]$ and $\diamond\in \{+,-\}$, if $w_i^\diamond \in \tW_*^\diamond$, then $v_i^\diamond = w_i^\diamond$. \label{lm:endpoints-IH-W*}
		\item For all $i \in [\ell]$, if there exists $\diamond\in \{+,-\}$ such that $w_i^\diamond \in V'$, then $(v_i^+,v_i^-)\neq (w_i^+,w_i^-)$.\label{lm:endpoints-IH-coverW0}
		\item For each $i \in [\ell]$, $\{v_i^+, w_i^+\} \cap \{v_i^-, w_i^-\} = \emptyset$.  \label{lm:endpoints-IH-consistent}
		\item For each $\diamond\in \{+,-\}$, if $\tU_{U^*}^\diamond(D)\setminus \hU_\ell^\diamond(D)\neq \emptyset$, then \NEW{both $|W_*^\pm|\leq 4$}\OLD{$|W_*^+|\leq 4$ and $|W_*^-|\leq 4$}.\label{lm:endpoints-IH-excess0}
		\item \NEW{Recall that $V^\pm= \{v\in V\mid d_D^\pm(v) \geq \texc(D) - 22\eta n\}$. Then, both $V^\pm\subseteq \{w_i^+, w_i^-, v_i^\pm\}\setminus \{v_i^\mp\}$ for all~$i \in [\ell]$.}\OLD{If $v \in V$ satisfies $d_D^+(v) \geq \texc(D) - 22\eta n$, then $v \in \{w_i^+, w_i^-, v_i^+\}\setminus \{v_i^-\}$ for all $i \in [\ell]$. Similarly, if $v \in V$ satisfies $d_D^-(v) \geq \texc(D) - 22\eta n$, then $v \in \{w_i^+, w_i^-, v_i^-\}\setminus \{v_i^+\}$ for all $i \in [\ell]$.}
		\label{lm:endpoints-IH-good}
	\end{enumerate}
	
	Assume $\ell=k$. \NEW{Since}\OLD{Then, since, by assumption,} $\Delta(H)\leq 11\eta n$, \cref{lm:endpoints-partialdecomp} follows from \cref{lm:endpoints-IH-forbidden,lm:endpoints-IH-partialdecomp}. \NEW{Moreover,} \cref{lm:endpoints-W*,lm:endpoints-coverW0,lm:endpoints-consistent,lm:endpoints-good} hold by \cref{lm:endpoints-IH-consistent,lm:endpoints-IH-coverW0,lm:endpoints-IH-W*,lm:endpoints-IH-good}, respectively.
	\NEW{It remains to verify \cref{lm:endpoints-excess0}. 
	By definition, we need to show that both $|\tU_{U^*}^\pm(D)\setminus \hU_k^\pm|\leq 88\eta n$.	
	If $\max\{|W_*^+|,|W_*^-|\}\geq 5$, then \cref{lm:endpoints-IH-excess0} implies that both $|\tU_{U^*}^\pm(D)\setminus \hU_k^\pm|=0$. If both $|W_*^\pm|\leq 4$, then $|\tU_{U^*}^\pm(D)\setminus \hU_\ell^\pm(D)|\leq k\leq 11\eta n|W_*|\leq 88\eta n$
	and so \cref{lm:endpoints-excess0} holds.}
	\OLD{By \cref{lm:endpoints-IH-excess0}, \cref{lm:endpoints-excess0} holds if $\max\{|W_*^+|,|W_*^-|\}\geq 5$. Otherwise, $|\tU_{U^*}^\pm(D)\setminus \hU_\ell^\pm(D)|\leq k\leq 88\eta n$ and \cref{lm:endpoints-excess0} is also satisfied.}
	
	Suppose $\ell<k$.
	First, observe that, by definition of $\hexc_\ell^\pm(v)$, the following hold.
	\begin{align}
	\hexc_\ell^\pm(D)&\stackrel{\text{\eqmakebox[eq:hexc/texc][c]{}}}{\coloneqq} \sum_{v\in V}\hexc_\ell^\pm(v)
	=\texc_{U^*}^\pm(D)-\ell\nonumber\\
	&\stackrel{\text{\eqmakebox[eq:hexc/texc][c]{\text{\cref{eq:texcA}}}}}{=}\texc(D)-|A^\pm|-\ell\label{eq:hexc/texc}\\
	&\stackrel{\text{\eqmakebox[eq:hexc/texc][c]{}}}{\geq} \exc(D)-|A^\pm|-\ell
	\geq \max\{|W_*^+|,|W_*^-|\}(1-21\eta)n-\lceil\eta n\rceil-11\eta n|W_*|\nonumber\\
	&\stackrel{\text{\eqmakebox[eq:hexc/texc][c]{}}}{\geq} \max\{|W_*^+|,|W_*^-|\}(1-46\eta)n.\label{eq:hexc}
	\end{align}	
	
	Let $X^\pm$ be the set of vertices~$v \in V\setminus\{w_{\ell+1}^\pm\}$ such that $\texc_{D,U^*}^\pm(v)-\hexc_\ell^\pm(v)=\lfloor\sqrt{\eta} n\rfloor$. \NEW{Note that each $v\in \tW_*^\pm\subseteq V\setminus U^*$ satisfies
	\begin{align}
		\hexc_\ell^\pm(v)&\stackrel{\text{\cref{lm:endpoints-IH-forbidden}}}{\geq} \texc_{D, U^*}^\pm(v)- \sqrt{\eta}n-d_H^\pm(v)
		\stackrel{\text{\cref{eq:texc}}}{\geq}
		\exc_D^\pm(v)-|A^\pm|-\sqrt{\eta}n-\Delta(H)\nonumber\\
		&\stackrel{\text{\cref{lm:cleaningallW*-defW}}}{\geq} (1-86\eta)n-\lceil\eta n\rceil-\sqrt{\eta}n-11\eta n\nonumber\\
		&\geq 2.\label{property:extW*}
	\end{align}}	
	\OLD{Also observe that,
	by \cref{lm:endpoints-IH-forbidden} and the fact that $\Delta(H)\leq 11\eta n$, the following hold.
	\begin{property}{\ddagger}
		Each $v\in \tW_*^\pm$ satisfies $\hexc_\ell^\pm(v)\geq \exc_D^\pm(v)-|A^\pm|-\sqrt{\eta}n-\Delta(H)\geq 2$.
	\end{property}}

	\begin{claim}
		\OLD{One can easily verify that}It is enough to find distinct $v_{\ell+1}^+,v_{\ell+1}^-\in V$ such that the following hold.%
		\COMMENT{\cref{lm:endpoints-IH-partialdecomp,lm:endpoints-IH-forbidden} follow from \cref{lm:endpoints-IS-partialdecomp}.\\
			\cref{lm:endpoints-IH-W*} follows from \cref{lm:endpoints-IS-W*}.\\
			\cref{lm:endpoints-IH-coverW0} follows from \cref{lm:endpoints-IH-partialdecomp} and \cref{lm:endpoints-IS-coverW0}.\\
			\cref{lm:endpoints-IH-consistent} follows from \cref{lm:endpoints-IS-partialdecomp} and the fact that $v_{\ell+1}^+\neq v_{\ell+1}^-$.\\
			\cref{lm:endpoints-IH-excess0} follows from \cref{lm:endpoints-IS-excess2}.\\
			By \cref{lm:endpoints-IS-good}, each $v\notin \{w_{\ell+1}^+, w_{\ell+1}^-\}$ satisfies \cref{lm:endpoints-IH-good}. The vertices $w_{\ell+1}^\pm$ satisfy \cref{lm:endpoints-IH-good} by \cref{lm:endpoints-IS-partialdecomp} and the fact that $d^\mp(w_{\ell+1}^\pm)\geq \texc(D)-22\eta n$ implies $\hexc_\ell^\pm(w_{\ell+1}^\pm)=0$.}
		\begin{enumerate}[label=\rm(\Roman*)]
			\item $v_{\ell+1}^\pm\in \hU_\ell^\pm(D)\setminus (X^\pm\cup \{w_{\ell+1}^\mp\})$.\label{lm:endpoints-IS-partialdecomp}
			\item If $w_{\ell+1}^\pm\in \tW_*^\pm$, then $v_{\ell+1}^\pm=w_{\ell+1}^\pm$.\label{lm:endpoints-IS-W*}
			\item If $\max\{|W_*^+|,|W_*^-|\}\geq 5$, then both $\hexc_\ell^\pm(v_{\ell+1}^\pm)\geq 2$.\label{lm:endpoints-IS-excess2}
			\item If $w_{\ell+1}^+\in \hU_\ell^+(D)$, $w_{\ell+1}^-\in \hU_\ell^-(D)$, then for each $\diamond \in \{+,-\}$ such that $w_{\ell+1}^\diamond\in V'$, we have $v_{\ell+1}^\diamond\neq w_{\ell+1}^\diamond$.\label{lm:endpoints-IS-coverW0}
			\item \OLD{For each $\diamond\in \{+,-\}$,}If \NEW{$v\in V^\pm\setminus \{w_{\ell+1}^+, w_{\ell+1}^-\}$}\OLD{$v\in V\setminus\{w_{\ell+1}^+,w_{\ell+1}^-\}$ satisfies $d_D^\diamond(v)\geq \texc(D)-22\eta n$}, then \NEW{$v_{\ell+1}^\pm=v$}\OLD{$v_{\ell+1}^\diamond=v$}.\label{lm:endpoints-IS-good}
		\end{enumerate}
	\end{claim}

	\begin{proofclaim}
		\NEW{Suppose that $v_{\ell+1}^+,v_{\ell+1}^-\in V$ are distinct and satisfy \cref{lm:endpoints-IS-W*,lm:endpoints-IS-coverW0,lm:endpoints-IS-excess2,lm:endpoints-IS-good,lm:endpoints-IS-partialdecomp}. We show that \cref{lm:endpoints-IH-W*,lm:endpoints-IH-consistent,lm:endpoints-IH-coverW0,lm:endpoints-IH-excess0,lm:endpoints-IH-forbidden,lm:endpoints-IH-good,lm:endpoints-IH-partialdecomp} hold with $\ell+1$ playing the role of $\ell$.}
		
		\NEW{First, \cref{lm:endpoints-IH-partialdecomp,lm:endpoints-IH-forbidden} follow from \cref{lm:endpoints-IS-partialdecomp}, while \cref{lm:endpoints-IH-W*} follows from \cref{lm:endpoints-IS-W*}. Moreover, \cref{lm:endpoints-IH-excess0} follows from \cref{lm:endpoints-IS-excess2}, while \cref{lm:endpoints-IH-consistent} follows from \cref{lm:endpoints-IS-partialdecomp} and the fact that $v_{\ell+1}^+\neq v_{\ell+1}^-$.}
		
		\OLD{Indeed,}In order to verify \cref{lm:endpoints-IH-coverW0}, note that, if both $w_{\ell+1}^+,w_{\ell+1}^-\in W$, then \cref{lm:endpoints-IH-coverW0} holds vacuously with~$\ell+1$ playing the role of~$\ell$ and, if there exists $\diamond\in \{+,-\}$ such that $w_{\ell+1}^\diamond\notin \hU_\ell^\diamond(D)$, then \cref{lm:endpoints-IH-coverW0} follows from \cref{lm:endpoints-IH-partialdecomp}.
		In the remaining cases, \cref{lm:endpoints-IH-coverW0} holds by \cref{lm:endpoints-IS-coverW0}.
		
		In order to verify \cref{lm:endpoints-IH-good}, first note that each $v\notin \{w_{\ell+1}^+, w_{\ell+1}^-\}$ satisfies \cref{lm:endpoints-IH-good} by \cref{lm:endpoints-IS-good}. To check that the vertices~$w_{\ell+1}^\pm$ satisfy \cref{lm:endpoints-IH-good}, first note that, by \cref{prop:Delta0}\cref{prop:Delta0-excv}, if \NEW{$w_{\ell+1}^\pm\in V^\mp$}\OLD{$d^\mp(w_{\ell+1}^\pm)\geq \texc(D)-22\eta n$}, then $\exc_D^\mp(w_{\ell+1}^\pm)\geq 2$ and so \NEW{\cref{eq:texc} implies that} $\hexc_\ell^\pm(w_{\ell+1}^\pm)\leq \texc_{D, U^*}^\pm(w_{\ell+1}^\pm)=0$.
		Thus, \cref{lm:endpoints-IH-good} for the vertices~$w_{\ell+1}^\pm$ follows from \cref{lm:endpoints-IS-partialdecomp}.
	\end{proofclaim}

	The following observation\OLD{s} will enable us to ensure that \cref{lm:endpoints-IS-W*,lm:endpoints-IS-good,lm:endpoints-IS-partialdecomp,lm:endpoints-IS-excess2,lm:endpoints-IS-coverW0} are satisfied simultaneously.
	
	\begin{claim}\label{property:endpoints}
		Let $v^\pm\in V\NEW{^\pm}$\OLD{satisfy $d_D^\pm(v^\pm)\geq \texc(D)-22\eta n$}. Then $\hexc_\ell^\pm(v^\pm)\geq 2$ and, if $v^\pm\neq w_{\ell+1}^\mp$, then $v_{\ell+1}^\pm\coloneqq v^\pm$ satisfies (the ``$\pm$ part" of) \cref{lm:endpoints-IS-W*,lm:endpoints-IS-good,lm:endpoints-IS-excess2,lm:endpoints-IS-partialdecomp,lm:endpoints-IS-coverW0}.
		In particular, if~$w_{\ell+1}^\pm\in \tW_*^\pm$ and~$v^\pm\neq w_{\ell+1}^\mp$, then~$v^\pm=w_{\ell+1}^\pm$.	
	\end{claim}
	
	\begin{proofclaim}
		Let $\diamond\in \{+,-\}$ and suppose~$v\in V\NEW{^\diamond}$\OLD{satisfies $d_D^\diamond(v)\geq \texc(D)-22\eta n$}. By \cref{prop:Delta0}\cref{prop:Delta0-excD}, $\texc(D)\leq (1+22\eta)n$ and so $|\tW_*^\pm|\leq 1$. In particular, $|W_*^\pm|\leq1$ and so,~$k\leq 22\eta n$. Thus, both~$X^\pm=\emptyset$.
		By \cref{prop:Delta0}\cref{prop:Delta0-excv},~$v\in \tW_*^\diamond$. Thus,~$\tW_*^\diamond=\{v\}$ and, by \cref{property:extW*},~$\hexc_\ell^\diamond(v)\geq 2$. 
		By \cref{prop:Delta0}\cref{prop:Delta0-excv}, each $u\in V\setminus\{v\}$ satisfies
		\NEW{$v\notin V^\diamond$}\OLD{$d_D^\diamond(u)< \texc(D)-22\eta n$} (otherwise $\exc_D^\diamond(u)\geq (1-86\eta)n$ and thus $|\tW_*^\diamond|>1$, a contradiction). 
	\end{proofclaim}

	To find $v_{\ell+1}^\pm$ when each \NEW{$v^\pm\in V\setminus (V^\pm\cup \{w_{\ell+1}^\mp\})$}\OLD{$v^\pm\in V\setminus\{w_{\ell+1}^\mp\}$ satisfies $d_D^\pm(v^\pm)< \texc(D)-22\eta n$}, we will use the following claim. For each~$S\subseteq V$, denote $\hexc_\ell^\pm(S)\coloneqq \sum_{v\in S}\hexc_\ell^\pm(v)$.
	
	\begin{claim}\label{claim:endpoints}
		The following hold.
		\begin{enumerate}[label=\rm(\Alph*)]
			\item If $\max\{|W_*^+|,|W_*^-|\}\geq 5$, then there exists $v\in \hU_\ell^+(D) \setminus (X^+ \cup \{w_{\ell+1}^+,w_{\ell+1}^-\})$ satisfying $\hexc_\ell^+(v)\geq 2$.\label{claim:hexc>2}
			\item Suppose $|W_*^+|,|W_*^-|\leq 4$. Then, $X^+=\emptyset$ and $\hU_\ell^+(D)\setminus\{w_{\ell+1}^-\}\neq \emptyset$. Moreover, if $w_{\ell+1}^-\in \hU_\ell^-(D)$ and $w_{\ell+1}^+\in V'$, then $\hexc_\ell^+(V\setminus \{w_{\ell+1}^+,w_{\ell+1}^-\})\geq 2$ (and thus, in particular, $\hU_\ell^+(D) \setminus (X^+ \cup \{w_{\ell+1}^+,w_{\ell+1}^-\})\neq \emptyset$).\label{claim:hexc1}
		\end{enumerate}
		Both statements also hold if~$+$ and~$-$ are swapped. 
	\end{claim}
	
	\begin{proofclaim}	
		For \cref{claim:hexc>2}, suppose that $\max\{|W_*^+|,|W_*^-|\}\geq 5$.		
		Assume for a contradiction that each $v\in \hU_\ell^+(D)\setminus (X^+\cup \{w_{\ell+1}^+,w_{\ell+1}^-\})$ satisfies $\hexc_\ell^+(v)=1$.
		Note that $\lfloor\sqrt{\eta} n\rfloor|X^+|\leq \ell < k\leq 11\eta n|W_*|$ and so $|X^+|\leq 23\sqrt{\eta} \max\{|W_*^+|,|W_*^-|\}$.
		Thus,
		\begin{align*}
		\hexc_\ell^+(D)&\leq |X^+\cup \{w_{\ell+1}^+,w_{\ell+1}^-\}|n+|\hU_\ell^+(D)|
		\leq \left(23\sqrt{\eta}\max\{|W_*^+|,|W_*^-|\}+2\right)n+n\\
		&\leq \max\{|W_*^+|,|W_*^-|\}23\sqrt{\eta}n+3n.
		\end{align*}
		But, by \cref{eq:hexc}, $\hexc_\ell^+(D)\geq \max\{|W_*^+|,|W_*^-|\}(1-46\eta)n$, so $\max\{|W_*^+|,|W_*^-|\}\leq \frac{3}{1-24\sqrt{\eta}}\leq 4$, a contradiction.
		
		For \cref{claim:hexc1}, assume that $|W_*^+|,|W_*^-|\leq 4$. Then, $\ell\leq k\leq \NEW{11\eta n|W_*|\leq }88\eta n$ and so $X^+=\emptyset$.
		If $w_{\ell+1}^-\in \hU_\ell^-(D)\subseteq \tU_{U^*}^-(D)$ and~$w_{\ell+1}^+\in V'$, then \begin{align}\label{eq:hexc+w-}
			\hexc_\ell^+(w_{\ell+1}^-)\leq \texc_{D, U^*}^+(w_{\ell+1}^-)\stackrel{\text{\NEW{\cref{eq:texc}}}}{\leq} 1
		\end{align} 
		and 
		\begin{align}\label{eq:hexc+w+}
			\hexc_\ell^+(w_{\ell+1}^+)\leq \texc_{D,U^*}^+(w_{\ell+1}^+)\stackrel{\text{\NEW{\cref{eq:texc}}}}{\leq}\max\{\exc_D^+(w_{\ell+1}^+),1\}\stackrel{\text{\NEW{\cref{lm:cleaningallW*-defW}}}}{\leq} \varepsilon n.
		\end{align} 
		Hence,\OLD{by \cref{eq:hexc}, }
		\begin{align*}
			\hexc_\ell^+(V\setminus\{w_{\ell+1}^+,w_{\ell+1}^-\})&\stackrel{\text{\eqmakebox[hexc+][c]{}}}{=}\hexc_\ell^+(D)-\hexc_\ell^+(w_{\ell+1}^-)-\hexc_\ell^+(w_{\ell+1}^+)\\
			&\stackrel{\text{\eqmakebox[hexc+][c]{\text{\NEW{\cref{eq:hexc},\cref{eq:hexc+w-},\cref{eq:hexc+w+}}}}}}{\geq} (1-46\eta)n-1-\varepsilon n\geq 2,
		\end{align*}
		as desired.
		It only remains to show that $\hU_\ell^+(D)\setminus \{w_{\ell+1}^-\}\neq \emptyset$.
		By \cref{eq:hexc}, $\hexc_\ell^+(D)>0$ and so $\hU_\ell^+(D)\neq \emptyset$. Suppose for a contradiction that $\hU_\ell^+(D)=\{w_{\ell+1}^-\}$.
		Note that\OLD{by \cref{eq:hexc},} 
		\begin{align*}
			\texc_{D,U^*}^+(w_{\ell+1}^-)\geq \hexc_\ell^+(w_{\ell+1}^-)=\hexc_\ell^+(D)\stackrel{\text{\NEW{\cref{eq:hexc}}}}{\geq} (1-46\eta)n.
		\end{align*} 
		Thus, $w_{\ell+1}^-\notin U^0(D)$ and so 
		\begin{align}\label{eq:exc+w-}
			\exc_D^+(w_{\ell+1}^-)\stackrel{\text{\NEW{\cref{eq:texc}}}}{\geq} \texc_{D,U^*}^+(w_{\ell+1}^-)\geq (1-46\eta)n.
		\end{align}
		Thus, $w_{\ell+1}^-\in \tW_*^+$ and, by \cref{lm:cleaningallW*-A}, $|W_A^+|\leq 1$. Moreover, each $v\in V\setminus \{w_{\ell+1}^-\}$ satisfies 
		\begin{align*}
			\exc_D^+(v)&\stackrel{\text{\eqmakebox[exc+v][c]{\text{\NEW{\cref{eq:texc}}}}}}{\leq} \texc_{D, U^*}^+(v)+d_{A^+}^+(v)\leq (\hexc_\ell^+(v)+\ell)+|A^+|\\
			&\stackrel{\text{\eqmakebox[exc+v][c]{\text{\NEW{\cref{lm:cleaningallW*-defW}}}}}}{\leq} 0+88\eta n+\lceil\eta n\rceil \stackrel{\text{\NEW{\cref{eq:exc+w-}}}}{<} \exc_D^+(w_{\ell+1}^-)-\varepsilon n.
		\end{align*} 
		Thus, by \cref{lm:cleaningallW*-A}, we have $W_A^+\subseteq \{w_{\ell+1}^-\}$. Therefore, $d_{A^+}^+(w_{\ell+1}^-)=|A^+|$ and so \NEW{the fact that $w_{\ell+1}^-\notin U^0(D)$ implies that}
		\begin{equation}\label{eq:texcw-}
			\texc_{D,U^*}^+(w_{\ell+1}^-)
			\stackrel{\text{\NEW{\cref{eq:texc}}}}{=}\exc_D^+(w_{\ell+1}^-)-d_{A^+}^+(w_{\ell+1}^-)=\exc_D^+(w_{\ell+1}^-)-|A^+|.
		\end{equation}
		Then,
		\begin{align}
			\hexc_\ell^+(D)&=\hexc_\ell^+(w_{\ell+1}^-)
			\leq \texc_{D,U^*}^+(w_{\ell+1}^-)
			\stackrel{\text{\cref{eq:texcw-}}}{=}\exc_D^+(w_{\ell+1}^-)-|A^+|
			\leq n,\label{eq:hexcD}
		\end{align}
		and so
		\begin{equation}\label{eq:hexw-}
			\texc(D)\stackrel{\text{\cref{eq:hexc/texc}}}{=}\hexc_\ell^+(D)+|A^+|+\ell\stackrel{\text{\cref{eq:hexcD}}}{\leq} \exc_D^+(w_{\ell+1}^-)+\ell.
		\end{equation}
		Suppose first that \NEW{$w_{\ell+1}^-\notin V^+$, that is,} $d_D^+(w_{\ell+1}^-)< \texc(D)-22\eta n$.
		By \cref{eq:hexc,eq:hexcD}, both~$|W_*^\pm|\leq1$. Thus,~$\ell\NEW{\leq 11\eta n|W_*|}\leq 22\eta n$ and so\OLD{by \cref{eq:hexw-},} 
		\begin{align*}
			d_D^+(w_{\ell+1}^-)<\texc(D)-22\eta n\leq \texc(D)-\ell\stackrel{\text{\NEW{\cref{eq:hexw-}}}}{\leq} \exc_D^+(w_{\ell+1}^-),
		\end{align*} 
		a contradiction.
		Therefore, \NEW{$w_{\ell+1}^-\in V^+$}\OLD{$d_D^+(w_{\ell+1}^-)\geq \texc(D)-22\eta n$}. 
		\NEW{Observe that
		\begin{equation}\label{eq:endpoints-i}
			\{i\in [\ell]\mid v_i^+\neq w_{\ell+1}^-\}=\{i\in [\ell]\mid w_i^-=w_{\ell+1}^-\}.
		\end{equation}
		Indeed, \cref{lm:endpoints-IH-good} implies that $w_{\ell+1}^-\in \{w_i^+, w_i^-, v_i^+\}$ for each $i\in [\ell]$ and \cref{lm:endpoints-IH-W*} implies that, for each $i\in [\ell]$, if $w_{\ell+1}^-=w_i^+$, then $w_{\ell+1}^-=v_i^+$. Thus, $w_{\ell+1}^-\in \{w_i^-, v_i^+\}$ for each $i\in [\ell]$. But \cref{lm:endpoints-IH-consistent} implies that $v_i^+\neq w_i^-$ for each $i\in [\ell]$. Therefore, \cref{eq:endpoints-i} holds and so}
		\OLD{Then, by \cref{lm:endpoints-IH-W*}, \cref{lm:endpoints-IH-consistent}, \cref{lm:endpoints-IH-good}, and since~$w_{\ell+1}^-\in \tW_*^+$, we have $\{i\in [\ell]\mid v_i^+\neq w_{\ell+1}^-\}=\{i\in [\ell]\mid w_i^-=w_{\ell+1}^-\}$%
			%\COMMENT{By \cref{lm:endpoints-IH-good}, for each $i\in [\ell]$, $w_{\ell+1}^-\in \{w_i^+, w_i^-, v_i^+\}$. By \cref{lm:endpoints-IH-W*}, for each $i\in [\ell]$, if $w_{\ell+1}^-=w_i^+$, then $w_{\ell+1}^-=v_i^+$. Thus, for each $i\in [\ell]$, $w_{\ell+1}^-\in \{w_i^-, v_i^+\}$. By \cref{lm:endpoints-IH-consistent}, $v_i^+\neq w_i^-$.}.
		Thus,}%
			\COMMENT{Strict inequality since $w_{\ell+1}^-$ is the ending point of $w_{\ell+1}^+w_{\ell+1}^-$. Indeed, $d_H^-(w_{\ell+1}^-)$ counts the number of indices $i\in [k]$ such that $w_i^-=w_{\ell+1}^-$, while $|\{i\in [\ell]\mid w_i^-=w_{\ell+1}^-\}|$ only counts such indices $i\leq \ell$.}
		\begin{align}	
			\texc_{U^*}^+(D)&\stackrel{\text{\eqmakebox[texc+D][c]{}}}{=} \hexc_\ell^+(D)+\ell=\hexc_\ell^+(w_{\ell+1}^-)+\ell\nonumber\\
			&\stackrel{\text{\eqmakebox[texc+D][c]{}}}{=}\texc_{D,U^*}^+(w_{\ell+1}^-)-|\{i\in [\ell]\mid v_i^+= w_{\ell+1}^-\}|+\ell\nonumber\\			
			&\stackrel{\text{\eqmakebox[texc+D][c]{\text{\NEW{\cref{eq:endpoints-i}}}}}}{=}\texc_{D,U^*}^+(w_{\ell+1}^-)+|\{i\in [\ell]\mid w_i^-=w_{\ell+1}^-\}|<\texc_{D,U^*}^+(w_{\ell+1}^-)+d_H^-(w_{\ell+1}^-)\nonumber\\
			&\stackrel{\text{\eqmakebox[texc+D][c]{}}}{\leq} \texc_{D,U^*}^+(w_{\ell+1}^-)+d_D^-(w_{\ell+1}^-)\label{eq:texcD}.
		\end{align}
		Therefore, we have
		\begin{align*}
			\NEW{\Delta^0(D)}&\stackrel{\text{\eqmakebox[claim:endpoints][c]{}}}{\leq}\texc(D)\stackrel{\text{\cref{eq:texcA}}}{=}\texc_{U^*}^+(D)+|A^+|
			\stackrel{\text{\cref{eq:texcD}}}{<}(\texc_{D,U^*}^+(w_{\ell+1}^-)+d_D^-(w_{\ell+1}^-))+|A^+|\\
			&\stackrel{\text{\eqmakebox[claim:endpoints][c]{\text{\cref{eq:texcw-}}}}}{=}\exc_D^+(w_{\ell+1}^-)+d_D^-(w_{\ell+1}^-)=d_D^+(w_{\ell+1}^-)\leq \Delta^0(D),
		\end{align*}
		 a contradiction.
		
		The same arguments hold with~$+$ and~$-$ swapped. This concludes the proof of \cref{claim:endpoints}.
	\end{proofclaim}
	
	We are now ready to choose distinct~$v_{\ell+1}^\pm\in V$ such that \cref{lm:endpoints-IS-W*,lm:endpoints-IS-good,lm:endpoints-IS-excess2,lm:endpoints-IS-partialdecomp,lm:endpoints-IS-coverW0} are satisfied.
	Without loss of generality, suppose that $|\hU_\ell^+(D)\setminus(X^+\cup\{w_{\ell+1}^-\})|\leq |\hU_\ell^-(D)\setminus(X^-\cup\{w_{\ell+1}^+\})|$.
	We start by picking~$v_{\ell+1}^+$ as follows \NEW{(where we assume in each case that the previous ones do not apply)}.
	\begin{case}
		\item \NEW{\textbf{$w_{\ell+1}^+\in \tW_*^+$ or $V^+\setminus\{w_{\ell+1}^-\}\neq\emptyset$.}} If $w_{\ell+1}^+\in \tW_*^+$, then let $v_{\ell+1}^+\coloneqq w_{\ell+1}^+$ and if there exists $v\in V\NEW{^+}\setminus\{w_{\ell+1}^-\}$\OLD{such that $d_D^+(v)\geq \texc(D)-22\eta n$}, then let~$v_{\ell+1}^+\coloneqq v$. (Note that~$v_{\ell+1}^+$ is well defined by the ``in particular" part of \cref{property:endpoints}.) Then, by \cref{property:extW*} and \cref{property:endpoints}, \NEW{(the ``$+$ part" of)}\OLD{properties} \cref{lm:endpoints-IS-W*,lm:endpoints-IS-good,lm:endpoints-IS-excess2,lm:endpoints-IS-partialdecomp,lm:endpoints-IS-coverW0} hold for~$v_{\ell+1}^+$.
		\item \NEW{\textbf{$\max\{|W_*^+|,|W_*^-|\}\geq 5$.} Let $v_{\ell+1}^+\in \hU_\ell^+(D)\setminus (X^+\cup \{w_{\ell+1}^+,w_{\ell+1}^-\})$ satisfy $\hexc_\ell^+(v_{\ell+1}^+)\geq 2$ ($v_{\ell+1}^+$ exists by \cref{claim:hexc>2}). Then, (the ``$+$ part" of) \cref{lm:endpoints-IS-W*,lm:endpoints-IS-good,lm:endpoints-IS-excess2,lm:endpoints-IS-partialdecomp,lm:endpoints-IS-coverW0} are clearly satisfied for~$v_{\ell+1}^+$.}
		\item \NEW{\textbf{$w_{\ell+1}^-\notin \hU_\ell^-(D)$ or $w_{\ell+1}^+\notin V'$.} Let $v_{\ell+1}^+\in \hU_\ell^+(D)\setminus (X^+\cup \{w_{\ell+1}^-\})$ ($v_{\ell+1}^+$ exists by \cref{claim:hexc1}). Then, (the ``$+$ part" of) \cref{lm:endpoints-IS-W*,lm:endpoints-IS-good,lm:endpoints-IS-excess2,lm:endpoints-IS-partialdecomp,lm:endpoints-IS-coverW0} are clearly satisfied for~$v_{\ell+1}^+$.}
		\item \NEW{\textbf{$w_{\ell+1}^-\in \hU_\ell^-(D)$ and $w_{\ell+1}^+\in V'$.} Let $v_{\ell+1}^+\in \hU_\ell^+(D)\setminus (X^+\cup \{w_{\ell+1}^+,w_{\ell+1}^-\})$ ($v_{\ell+1}^+$ exists by the ``moreover part" of \cref{claim:hexc1}). Then, (the ``$+$ part" of) \cref{lm:endpoints-IS-W*,lm:endpoints-IS-good,lm:endpoints-IS-excess2,lm:endpoints-IS-partialdecomp,lm:endpoints-IS-coverW0} are clearly satisfied for~$v_{\ell+1}^+$.}
	\end{case}%
	\OLD{Otherwise, by using \cref{claim:endpoints} and distinguishing the cases when $\max\{|W_*^+|,|W_*^-|\}\geq 5$ and when $|W_*^+|, |W_*^-|\leq 4$, it is easy to check that we can choose $v_{\ell+1}^+\in \hU_\ell^+(D)\setminus (X^+\cup \{w_{\ell+1}^-\})$ satisfying \cref{lm:endpoints-IS-W*,lm:endpoints-IS-good,lm:endpoints-IS-excess2,lm:endpoints-IS-partialdecomp,lm:endpoints-IS-coverW0}.}%
	Note that, since~$v_{\ell+1}^+$ satisfies \cref{lm:endpoints-IS-partialdecomp}, \NEW{we have $v_{\ell+1}^+\neq w_{\ell+1}^-$. Moreover, \cref{lm:endpoints-IS-partialdecomp} and \cref{eq:texc} imply that} $\hexc_\ell^-(v_{\ell+1}^+)\leq 1$. 
	Therefore, by \cref{property:extW*},~$v_{\ell+1}^+\notin \tW_*^-$ and, by \cref{property:endpoints}, \NEW{$v_{\ell+1}^+\notin V^-$}\OLD{$d_D^-(v_{\ell+1}^+)<\texc(D)-22\eta n$} (otherwise $\hexc_\ell^-(v_{\ell+1}^+)\geq 2$, a contradiction).
	Thus, if $\hU_\ell^-(D)\setminus (X^-\cup \{w_{\ell+1}^+,v_{\ell+1}^+\})\neq \emptyset$, then, we can proceed similarly as for~$v_{\ell+1}^+$ to obtain~$v_{\ell+1}^-\NEW{\neq v_{\ell+1}^+}$ satisfying \cref{lm:endpoints-IS-W*,lm:endpoints-IS-good,lm:endpoints-IS-excess2,lm:endpoints-IS-partialdecomp,lm:endpoints-IS-coverW0}.
	\NEW{More precisely, we proceed as follows.}
	\begin{case}
		\item \NEW{\textbf{$w_{\ell+1}^-\in \tW_*^-$ or $V^-\setminus\{w_{\ell+1}^+\}\neq\emptyset$.} If $w_{\ell+1}^-\in \tW_*^-$, then let $v_{\ell+1}^-\coloneqq w_{\ell+1}^-$ and if there exists $v\in V^-\setminus\{w_{\ell+1}^+\}$, then let~$v_{\ell+1}^+\coloneqq v$. (Note that~$v_{\ell+1}^-$ is well defined by the ``in particular" part of \cref{property:endpoints}.)  
		Then, by \cref{property:extW*} and \cref{property:endpoints}, (the ``$-$ part" of) \cref{lm:endpoints-IS-W*,lm:endpoints-IS-good,lm:endpoints-IS-excess2,lm:endpoints-IS-partialdecomp,lm:endpoints-IS-coverW0} hold for~$v_{\ell+1}^-$.
		Moreover, we have shown above that $v_{\ell+1}^+\notin \{w_{\ell+1}^-\}\cup V^-$, thus $v_{\ell+1}^-\neq v_{\ell+1}^+$ and so we are done.}
		\item \NEW{\textbf{$\max\{|W_*^+|,|W_*^-|\}\geq 5$.} Let $v_{\ell+1}^-\in \hU_\ell^-(D)\setminus (X^-\cup \{w_{\ell+1}^+,w_{\ell+1}^-\})$ satisfy $\hexc_\ell^-(v_{\ell+1}^-)\geq 2$ ($v_{\ell+1}^-$ exists by (the ``$-$ analogue" of) \cref{claim:hexc>2}). Then, (the ``$-$ part" of) \cref{lm:endpoints-IS-W*,lm:endpoints-IS-good,lm:endpoints-IS-excess2,lm:endpoints-IS-partialdecomp,lm:endpoints-IS-coverW0} are clearly satisfied for~$v_{\ell+1}^-$. Moreover, we have shown above that $\hexc_\ell^-(v_{\ell+1}^+)\leq 1$, thus $v_{\ell+1}^-\neq v_{\ell+1}^+$ and so we are done.}
		\item \NEW{\textbf{$w_{\ell+1}^+\notin \hU_\ell^+(D)$ or $w_{\ell+1}^-\notin V'$.} Let $v_{\ell+1}^-\in \hU_\ell^-(D)\setminus (X^-\cup \{w_{\ell+1}^+,v_{\ell+1}^+\})$ ($v_{\ell+1}^-$ exists by assumption). Then, (the ``$-$ part" of) \cref{lm:endpoints-IS-W*,lm:endpoints-IS-good,lm:endpoints-IS-excess2,lm:endpoints-IS-partialdecomp,lm:endpoints-IS-coverW0} are clearly satisfied for~$v_{\ell+1}^-$ and so we are done.}
		\item \NEW{\textbf{$w_{\ell+1}^+\in \hU_\ell^+(D)$ and $w_{\ell+1}^-\in V'$.} Recall that $\hexc_\ell^-(v_{\ell+1}^+)\leq 1$, so (the ``$-$ analogue" of) \cref{claim:hexc1} implies that $X^-=\emptyset$ and $\hexc_\ell^-(V\setminus \{w_{\ell+1}^+,w_{\ell+1}^-,v_{\ell+1}^+\})=\hexc_\ell^-(V\setminus \{w_{\ell+1}^+,w_{\ell+1}^-\})-\hexc_\ell^-(v_{\ell+1}^+)\geq 2-1>0$.
		We can thus let $v_{\ell+1}^-\in \hU_\ell^-(D)\setminus (X^-\cup \{w_{\ell+1}^+,w_{\ell+1}^-,v_{\ell+1}^+\})$. Then, (the ``$-$ part" of) \cref{lm:endpoints-IS-W*,lm:endpoints-IS-good,lm:endpoints-IS-excess2,lm:endpoints-IS-partialdecomp,lm:endpoints-IS-coverW0} are clearly satisfied for~$v_{\ell+1}^-$ and so we are done.}
	\end{case}%
	\OLD{(To see that this can be done in the case when $|W_*^+|,|W_*^-|\leq 4$, note that, by \cref{claim:endpoints}\cref{claim:hexc1}, if $|W_*^+|,|W_*^-|\leq 4$, $w_{\ell+1}^+\in \hU_\ell^+(D)$, and $w_{\ell+1}^-\in V'$, then $\hexc_\ell^-(V\setminus\{w_{\ell+1}^+,w_{\ell+1}^-,v_{\ell+1}^+\})\geq 2-\hexc_\ell^-(v_{\ell+1}^+)\geq 1$ and so $\hU_\ell^-(D) \setminus \{w_{\ell+1}^+,w_{\ell+1}^-,v_{\ell+1}^+\}\neq \emptyset$.)}%
	We may therefore assume that $\hU_\ell^-(D)\setminus (X^-\cup \{w_{\ell+1}^+,v_{\ell+1}^+\})= \emptyset$.
	But, by \cref{claim:endpoints}, $\hU_\ell^-(D)\setminus (X^-\cup \{w_{\ell+1}^+\})\neq \emptyset$ and so $\hU_\ell^-(D)\setminus (X^-\cup \{w_{\ell+1}^+\})= \{v_{\ell+1}^+\}$. Thus, by assumption, \[|\hU_\ell^+(D)\setminus (X^+\cup \{w_{\ell+1}^-\})|\leq |\hU_\ell^-(D)\setminus (X^-\cup \{w_{\ell+1}^+\})|=|\{v_{\ell+1}^+\}|=1.\]
	Then, since~$v_{\ell+1}^+$ satisfies \cref{lm:endpoints-IS-partialdecomp}, $\hU_\ell^+(D)\setminus (X^+\cup \{w_{\ell+1}^-\})=\{v_{\ell+1}^+\}$. We will find a contradiction.
	
	Note that, $v_{\ell+1}^+\in \hU_\ell^+(D)\cap \hU_\ell^-(D)\subseteq U^0(D)$ and so 
	\begin{align}\label{eq:hexc1}
		\hexc_\ell^\pm(v_{\ell+1}^+)\stackrel{\text{\NEW{\cref{eq:texc}}}}{=}1.
	\end{align} 
	\NEW{Since $v_{\ell+1}^+$ satisfies (the ``$+$ part" of) \cref{lm:endpoints-IS-excess2}, this implies that $|W_*^+|,|W_*^-|\leq 4$. By \cref{claim:hexc1} (and its ``$-$ analogue"), we have $X^\pm=\emptyset$ and so $\{v_{\ell+1}^+\}\subseteq \hU_\ell^\pm(D)\subseteq \{v_{\ell+1}^+,w_{\ell+1}^\mp\}$. Hence,
	\begin{align*}
		\exc_{D,U^*}^\pm(w_{\ell+1}^\mp)&\stackrel{\text{\eqmakebox[excw][c]{\text{\cref{eq:texc}}}}}{\geq}\texc_{D, U^*}^\pm(w_{\ell+1}^\mp)-1\geq \hexc_\ell^\pm(w_{\ell+1}^\mp)-1\\
		&\stackrel{\text{\eqmakebox[excw][c]{}}}{=} \hexc_\ell^\pm(D)-\hexc_\ell^\pm(v_{\ell+1}^+)-1
		\stackrel{\text{{\text{\cref{eq:hexc},\cref{eq:hexc1}}}}}{\geq}(1-46\eta)n-1-1\geq (1-47\eta)n.
	\end{align*}
	Together with \cref{property:extW*} and \cref{eq:hexc1}, this implies that $\tW_*^\pm=\{w_{\ell+1}^\mp\}$ and so}%
	\OLD{Then, since $v_{\ell+1}^+$ satisfies \cref{lm:endpoints-IS-excess2}, we have $|W_*^+|,|W_*^-|\leq 4$. Thus, by \cref{claim:endpoints}\cref{claim:hexc1} (and its analogue with~$+$ and~$-$ swapped),~$X^\pm=\emptyset$. Therefore, we have $\hU_\ell^\pm(D)\setminus \{w_{\ell+1}^\mp\}=\{v_{\ell+1}^+\}$.
	Moreover, by \cref{eq:hexc}, $\hexc_\ell^\pm(D)\geq (1-46\eta)n$ and so, as $\hexc_\ell^\pm(v_{\ell+1}^+)= 1$, we must have $\hU_\ell^\pm(D)=\{w_{\ell+1}^\mp, v_{\ell+1}^+\}$.
	Then, $\texc_{D, U^*}^\pm(w_{\ell+1}^\mp)\geq \hexc_\ell^\pm(w_{\ell+1}^\mp)= \hexc_\ell^\pm(D)-1\geq (1-47\eta)n$. Thus, $w_{\ell+1}^\mp\notin U^*$ and so $\exc_D^\pm(w_{\ell+1}^\mp)\geq \texc_{D, U^*}^\pm(w_{\ell+1}^\mp)\geq (1-47\eta)n$.
	Therefore, $w_{\ell+1}^\mp\in \tW_*^\pm$. But, by \cref{property:extW*}, $\tW_*^\pm\subseteq \hU_\ell^\pm(D)$ and so $\tW_*^\pm=\{w_{\ell+1}^\mp\}$. Thus}, by assumption on our ordering of~$E(H)$ \NEW{made after \cref{property:notannoyingD}}, it follows that~$\ell=0$. Therefore, $\tU_{U^*}^\pm(D)=\hU_\ell^\pm(D)=\{w_{\ell+1}^\mp, v_{\ell+1}^+\}$, contradicting \cref{property:notannoyingD}.
\end{proof}

\subsection{Decreasing the degree at \texorpdfstring{$W_A$}{WA}} 

Note that, if~$W_*=\emptyset$, then the excess of our tournament may be relatively small and so we do not have room to proceed similarly as in \cref{lm:cleaningallW*} to decrease the degree of the vertices in~$W_A$. The strategy is to find a partial path decomposition~$\cP$ such that each vertex in~$W_A$ is covered by each of the paths in~$\cP$ and such that each vertex in~$V'$ is covered by half of the paths in~$\cP$. In this way, the degree at~$W_A$ is decreased faster than the degree at~$V'$. Decreasing the degree at~$V'$ will ensure that the leftover excess is not too small compared to degree of the leftover oriented graph (recall \cref{lm:cleaning}\cref{lm:cleaning-exc>d}).

\begin{lm}\label{lm:cleaningdegreeWA}
	Let $0<\frac{1}{n}\ll \varepsilon\ll \eta \ll 1$. Let~$D$ be an oriented graph on a vertex set~$V$ of size~$n$ satisfying $\delta(D)\geq (1-\varepsilon)n$, $\texc(D)\geq \frac{n}{2} +9\eta n$, and the following properties.
	\begin{enumerate}
		\item Let $W\cup V'$ be a partition of~$V$ such that, for each $v\in V'$, $|\exc_D(v)|\leq \varepsilon n$ and, for each $v\in W$, $|\exc_D(v)|\leq (1-20\eta)n$.\label{lm:cleaningdegreeWA-W}	
		Suppose $E(D[W])=\emptyset$ and $|W|\leq \varepsilon n$.
		\item Let $A^+,A^-\subseteq E(T)$ be absorbing sets of~$(W, V')$-starting/$(V',W)$-ending edges for~$D$ of size at most~$\lceil\eta n\rceil$.
		Let $A\coloneqq A^+\cup A^-$, $W_A^\pm \coloneqq V(A^\pm)\cap W$, and $W_A\coloneqq V(A)\cap W$. Assume $A\neq \emptyset$, i.e.\ $W_A\neq \emptyset$.
		\item Let $U^*\subseteq U^0(D)$ satisfy $|U^*|=\texc(D)-\exc(D)$.
	\end{enumerate} 
	Then, there exists a good $(U^*,W,A)$-partial path decomposition~$\cP$ of~$D$ such that $|\cP|=8\lceil\eta n\rceil$ and $D'\coloneqq D\setminus \cP$ satisfies the following.
	\begin{enumerate}[label=\rm(\roman*)]
		\item For each $v\in W$, $d_{D'}(v)\leq d_D(v)-12\lceil\eta n\rceil$.\label{lm:cleaningdegreeWA-degreeW}
		\item For each $v\in V'$, $d_D(v)-8\lceil\eta n\rceil\leq d_{D'}(v)\leq d_D(v)-8\lceil\eta n\rceil+1$.\label{lm:cleaningdegreeWA-degreeV'}
		\item Each $v\in U^*\setminus (V^+(\cP)\cup V^-(\cP))$ satisfies $d_{D'}^+(v)=d_{D'}^-(v)\leq \texc(D')-1$. \label{lm:cleaningdegreeWA-U*}
	\end{enumerate}
\end{lm}

\begin{proof}
	Fix additional constants such that $\varepsilon \ll\nu\ll\tau\ll\eta$. 	
	Let $k\coloneqq 4\lceil\eta n\rceil$. Assume inductively that, for some $0\leq\ell\leq k$, we have constructed edge-disjoint paths $P_{1,1},P_{1,2},P_{2,1}, \dots, P_{\ell,2}\subseteq D$ such that $\cP_\ell\coloneqq \{P_{i,j}\mid i\in [\ell], j\in [2]\}$ is a $(U^*, W, A)$-partial path decomposition of~$D$ such that the following hold, where $D_\ell\coloneqq D\setminus \cP_\ell$.
	\begin{enumerate}[label=\upshape(\greek*)]
		\item For each $i\in [\ell]$ and~$v\in W$, $v\in V(P_{i,1})\cap V(P_{i,2})$.\label{lm:cleaningdegreeWA-spanningW}
		\item For each $i\in [\ell]$ and~$v\in W$,~$v$ is \NEW{an}\OLD{the} endpoint of at most one of~$P_{i,1}$ and~$P_{i,2}$.\label{lm:cleaningdegreeWA-endpointW}
		\item For each $i\in [\ell]$ and $v\in V'$, either $v\in V(P_{i,1})\triangle V(P_{i,2})$ or~$v$ is an endpoint of both~$P_{i,1}$ and~$P_{i,2}$.\label{lm:cleaningdegreeWA-spanningV'}
		\item For each $v\in V'$, there is at most one~$i\in [\ell]$ such that~$v$ is an endpoint of exactly one of~$P_{i,1}$ and~$P_{i,2}$. Moreover, for each~$v\in V'$, if there exists~$i\in [\ell]$ such that~$v$ is an endpoint of exactly one of~$P_{i,1}$ and~$P_{i,2}$, then $\exc_{D_\ell}(v)=0$. \label{lm:cleaningdegreeWA-endpointV'}
	\end{enumerate}
	If $\ell=k$, then let $\cP\coloneqq \cP_k$ and $D'\coloneqq D\setminus \cP$.
	
	\begin{claim}\label{claim:good}
		$\cP$ is a good partial path decomposition of~$D$, \NEW{i.e.\ $\texc(D')=\texc(D)-2k=\texc(D)-8\lceil\eta n\rceil$}.
	\end{claim}

	Note that if \NEW{$\ell=k$ and} \cref{claim:good} holds, then we are done. 
	\NEW{Indeed, \cref{lm:cleaningdegreeWA-spanningW,lm:cleaningdegreeWA-endpointW} imply that each $w\in W$ satisfies $d_\cP(w)\geq 3k=12\lceil\eta n\rceil$, while \cref{lm:cleaningdegreeWA-endpointV',lm:cleaningdegreeWA-spanningV'} imply that each $v\in V'$ satisfies $8\lceil\eta n\rceil-1=2k-1\leq d_\cP(v)\leq 2k=8\lceil\eta n \rceil$. Thus, \cref{lm:cleaningdegreeWA-degreeW,lm:cleaningdegreeWA-degreeV'} hold.}%
	\OLD{Indeed, \cref{lm:cleaningdegreeWA-degreeW} holds by \cref{lm:cleaningdegreeWA-spanningW,lm:cleaningdegreeWA-endpointW} while \cref{lm:cleaningdegreeWA-degreeV'} holds by \cref{lm:cleaningdegreeWA-endpointV',lm:cleaningdegreeWA-spanningV'}.}
	Finally, \cref{lm:cleaningdegreeWA-U*} follows from \cref{claim:good}. Indeed, for each $v\in U^*\setminus(V^+(\cP)\cup V^-(\cP))$, we have $d_{D'}^+(v)=d_{D'}^-(v)\leq \frac{n-1}{2}<\frac{n}{2}+9\eta n-2k \leq \texc(D)-2k=\texc(D')$, as desired.\OLD{Thus, it suffices to prove \cref{claim:good}.}
	
	\begin{proofclaim}[Proof of \cref{claim:good}]
		\NEW{By \cref{prop:excpartialdecomp}\cref{prop:excpartialdecomp-general},}\OLD{Observe that, as~$\cP$ is a partial path decomposition of~$D$, by \cref{prop:excpartialdecomp}, $\texc(D')\geq \texc(D)-|\cP|$ and $\exc(D')\leq \texc(D)-|\cP|$. 
		Thus,} it is enough to show that $\Delta^0(D')\leq \texc(D)-|\cP|$. 
		
		Let $v\in V$. By \cref{lm:cleaningdegreeWA-spanningW},  \cref{lm:cleaningdegreeWA-spanningV'}, and since~$\cP$ is a partial path decomposition of~$D$, the following hold.
		\begin{itemize}[--]
			\item If $v\in U^\pm(D)\cap W$, then $v\in V^\pm(P)\cup V^0(P)$ for each~$P\in \cP$.
			\item If $v\in U^0(D)\cap W$, then for each $\diamond\in\{+,-\}$, $v\in V^\diamond(P)\cup V^0(P)$ for all but at most one~$P\in \cP$.
			\item If $v\in U^\pm(D)\cap V'$, then $v\in V^\pm(P)\cup V^0(P)$ for at least $k$ paths $P\in \cP$.
			\item If $v\in U^0(D)\cap V'$, then, for each $\diamond\in\{+,-\}$, $v\in V^\diamond(P)\cup V^0(P)$ for at least $k-1\NEW{=|\cP|-k-1}$ paths $P\in \cP$.
		\end{itemize}
		Thus, since $\cP$ is a partial path decomposition of $D$, we have 
		\begin{equation}\label{eq:good2}
		d_{D'}^{\max}(v)\leq 
		\begin{cases}
		d_D^{\max}(v) -|\cP| & \text{if }v\in W\setminus U^0(D),\\
		d_D^{\max}(v) -|\cP|+k+1 & \text{if }v\in V'\cup U^0(D).\\
		\end{cases}
		\end{equation}
		\OLD{By \cref{fact:exc}\cref{fact:exc-dmax},}For each $v\in V'\cup U^0(D)$, we have 
		\begin{align*}
			d_D^{\max}(v)\stackrel{\text{\NEW{\cref{fact:exc}\cref{fact:exc-dmax}}}}{=} \frac{d_D(v)+|\exc_{D}(v)|}{2}\stackrel{\text{\NEW{\cref{lm:cleaningdegreeWA-W}}}}{\leq} \frac{n-1+\varepsilon n}{2}\leq \frac{n}{2}+9\eta n-4\lceil\eta n\rceil -1\leq \texc(D)-k-1
		\end{align*}
		and so, by \cref{eq:good2}, $\Delta^0(D')\leq \texc(D)-|\cP|$. Thus,~$\cP$ is a good partial path decomposition of~$D$, as desired.
	\end{proofclaim} 
	
	If $\ell<k$, then let $D_\ell\coloneqq D\setminus \cP_\ell$ and $U_\ell^*\coloneqq U^*\setminus (V^+(\cP_\ell)\cup V^-(\cP_\ell))$.
	We claim that there exist suitable endpoints $v_1^+,v_1^-,v_2^+,v_2^-\in V$ for~$P_{\ell+1,1}$ and~$P_{\ell+1,2}$.
	
	\begin{claim}\label{claim:endpoints2}
		There exist $v_1^+,v_1^-,v_2^+,v_2^-\in V$ such that the following hold.
		\begin{enumerate}[label=\rm(\Roman*)]
			\item For each $i\in [2]$, $v_i^+\neq v_i^-$ and $v_i^\pm\in \tU_{U_\ell^*}^\pm(D_\ell)$. Moreover, for each $\diamond\in\{+,-\}$, if $v_1^\diamond=v_2^\diamond$, then $\texc_{D_\ell, U_\ell^*}^\diamond(v_1^\diamond)\geq 2$.\label{lm:cleaningdegreeWA-partialdecomp}
			\item For each $v\in W$, there exists at most one pair $(i,\diamond)\in [2]\times\{+,-\}$ such that $v_i^\diamond=v$.\label{lm:cleaningdegreeWA-Wendpoint}
			\item For each $v\in V'$, if there exists exactly one pair $(i,\diamond)\in [2]\times\{+,-\}$ such that $v_i^\diamond=v$, then $\exc_{D_\ell}^\diamond(v)=1$.\label{lm:cleaningdegreeWA-V'endpoint}
		\end{enumerate}
	\end{claim}

	Before proving \cref{claim:endpoints2}, let us first apply it to construct~$P_{\ell+1,1}$ and~$P_{\ell+1,2}$. Let $v_1^+,v_1^-,v_2^+,v_2^-\in V$ be as in \cref{claim:endpoints2}. We construct a~$(v_1^+,v_1^-)$-path~$P_{\ell+1,1}$ and a~$(v_2^+,v_2^-)$-path~$P_{\ell+1,2}$ using \cref{cor:robpaths} as follows.
	Observe that, by \cref{lm:3/8rob}, $D_\ell[V'\setminus \{v_1^+,v_1^-,v_2^+,v_2^-\}]$ is a robust~$(\nu, \tau)$-outexpander.
	\NEW{Moreover, each $v\in V$ satisfies
	\begin{align*}
		d_{D_\ell}^{\min}(v)&
		\geq d_D^{\min}(v)-|\cP_\ell|
		\stackrel{\text{\cref{fact:exc}\cref{fact:exc-dmin}}}{\geq}\frac{d_D(v)-|\exc_D(v)|}{2}-2\ell
		\stackrel{\text{\cref{lm:cleaningdegreeWA-W}}}{\geq}\frac{(20\eta -\varepsilon)n}{2}-8\lceil\eta n\rceil\\
		&\geq \eta n \stackrel{\text{\cref{lm:cleaningdegreeWA-W}}}{\geq} 2(|W|+2)+2.
	\end{align*}}%
	\NEW{Let $\delta\coloneqq \frac{3}{8}$ and $S\coloneqq \{v_2^+,v_2^-\}\setminus (W\cup \{v_1^+,v_1^-\})$. For each $i\in [2]$, let $V_i'\coloneqq V'\setminus \{v_i^+,v_i^-\}$ and $k_i\coloneqq |W\cup \{v_i^+,v_i^-\}|$.} Apply \cref{cor:robpaths}\cref{cor:robpaths-short} with 
	\NEW{\begin{center}
	\begin{tabular}{r|c|c|c|c|c|c}
		  & $D_\ell\setminus A$ & $V_1'$ & $k_1$ & $v_1^+$ & $W\setminus \{v_1^+,v_1^-\}$ &$v_1^-$\\
		 \hline
		 playing the role of & $D$& $V'$ & $k$ & $P_1$ & $\{P_2, \dots, P_{k-1}\}$ & $P_k$
	\end{tabular}
	\end{center}}%
	\noindent\OLD{$D_\ell\setminus A, V'\setminus \{v_1^+,v_1^-\}, |W\cup \{v_1^+,v_1^-\}|, \frac{3}{8}, \{v_2^+,v_2^-\}\setminus (W\cup \{v_1^+,v_1^-\}),v_1^+, W\setminus \{v_1^+,v_1^-\}$, and $v_1^-$ playing the roles of $D, V', k, \delta, S, P_1, \{P_2, \dots, P_{k-1}\}$, and $P_k$}to obtain a~$(v_1^+,v_1^-)$-path~$P_{\ell+1,1}$ of length at most~$\sqrt{\varepsilon}n$ which covers~$W$ and avoids $\{v_2^+, v_2^-\}\setminus (W\cup \{v_1^+, v_1^-\})$.
	Let $D_\ell'\coloneqq D_\ell \setminus P_{\ell+1,1}$ and observe that, by \cref{lm:3/8rob}, $D_\ell'[V'\setminus (V(P_{\ell+1,1})\cup \{v_2^+,v_2^-\})]$ is still a robust~$(\nu,\tau)$-outexpander.
	Then, \NEW{let $S'\coloneqq V(P_{\ell+1,1})\setminus(W\cup \{v_2^+,v_2^-\})$ and} apply \cref{cor:robpaths}\cref{cor:robpaths-long} with 
	\NEW{\begin{center}
			\begin{tabular}{r|c|c|c|c|c|c|c}
				& $D_\ell\setminus A$ & $V_2'$ & $k_2$ & $S'$ & $v_2^+$ & $W\setminus \{v_2^+,v_2^-\}$ &$v_2^-$\\
				\hline
				playing the role of & $D$& $V'$ & $k$ & $S$ & $P_1$ & $\{P_2, \dots, P_{k-1}\}$ & $P_k$
			\end{tabular}
	\end{center}}%
	\noindent\OLD{$D_\ell', V'\setminus \{v_2^+,v_2^-\}, |W\cup \{v_2^+,v_2^-\}|, \frac{3}{8}, V(P_{\ell+1,1})\setminus(W\cup \{v_2^+,v_2^-\}),v_2^+, W\setminus \{v_2^+,v_2^-\}$, and $v_2^-$ playing the roles of $D, V', k, \delta, S, P_1$, $\{P_2, \dots, P_{k-1}\}$, and $P_k$}to obtain a~$(v_2^+,v_2^-)$-path~$P_{\ell+1,2}$ satisfying $V\setminus V(P_{\ell+1,2})=V(P_{\ell+1,1})\setminus (W\cup \{v_2^+,v_2^-\})$. Then,
	note that, by \cref{lm:cleaningdegreeWA-partialdecomp},~$\cP_{\ell+1}$ is a $(U^*,W,A)$-partial path decomposition of~$D$ and, by \cref{lm:cleaningdegreeWA-Wendpoint,lm:cleaningdegreeWA-V'endpoint}, \cref{lm:cleaningdegreeWA-endpointV',lm:cleaningdegreeWA-endpointW} are satisfied with~$\ell+1$ playing the role of~$\ell$, respectively. Finally, by construction of~$P_{\ell+1,1}$ and~$P_{\ell+1,2}$,  \cref{lm:cleaningdegreeWA-spanningW,lm:cleaningdegreeWA-spanningV'} are satisfied.
	
	It remains to prove \cref{claim:endpoints2}.
	
	\begin{proof}[Proof of \cref{claim:endpoints2}]\renewcommand{\qedsymbol}{\rotatebox[origin=c]{45}{$\square$}}
		Since~$\cP_\ell$ is a $(U^*,W,A)$-partial path decomposition of~$D$, $|A^\pm|\leq\lceil\eta n\rceil$, and $2\ell\leq 2k\leq 8\lceil\eta n\rceil$,
		\begin{align}
		\texc_{U_\ell^*}^\pm(D_\ell)&\stackrel{\text{\eqmakebox[eq:cleaningdegreeWA-texc][c]{\text{\cref{fact:specialdecomp}}}}}{=}\texc_{U^*}^\pm(D)-2\ell\nonumber\\
		&\stackrel{\text{\eqmakebox[eq:cleaningdegreeWA-texc][c]{\text{\cref{eq:texcA}}}}}{=}\texc(D)-|A^\pm|-2\ell\label{eq:cleaningdegreeWA-texc1}\\
		&\stackrel{\text{\eqmakebox[eq:cleaningdegreeWA-texc][c]{}}}{\geq} \frac{n}{2}-\eta n\label{eq:cleaningdegreeWA-texc2}.
		\end{align}	
		Thus, we can choose endpoints $v_1^+, v_1^-,v_2^+,v_2^-\in V$ satisfying \cref{lm:cleaningdegreeWA-V'endpoint,lm:cleaningdegreeWA-Wendpoint,lm:cleaningdegreeWA-partialdecomp} as follows.
		
		If $|U_\ell^*\cap V'|\geq 2$, then pick distinct $u_1,u_2\in U_\ell^*\cap V'$ and let $v_1^+\coloneqq u_1, v_2^+\coloneqq u_2, v_1^-\coloneqq u_2$, and $v_2^-\coloneqq u_1$. Then, \cref{lm:cleaningdegreeWA-V'endpoint,lm:cleaningdegreeWA-Wendpoint,lm:cleaningdegreeWA-partialdecomp} are satisfied, as desired.
		
		We may therefore assume that $|U_\ell^*\cap V'|\leq 1$. 
		We first pick $v_1^+, v_2^+\in \tU_{U_\ell^*}^+(D_\ell) \setminus (U_\ell^*\cap V')$ as follows.		
		If $\texc_{D_\ell, U_\ell^*}^+(V'\setminus U_\ell^*)\geq 2$, then pick $v_1^+\in \tU_{U_\ell^*}^+(D_\ell)\cap (V'\setminus U_\ell^*)$ such that $\texc_{D_\ell, U_\ell^*}^+(v_1^+)$ is maximum. If $\texc_{D_\ell, U_\ell^*}^+(v_1^+)\geq 2$, then let $v_2^+\coloneqq v_1^+$; otherwise, let $v_2^+\in (\tU_{U_\ell^*}^+(D_\ell)\cap (V'\setminus U_\ell^*))\setminus \{v_1^+\}$%
			\COMMENT{Note that in the latter case, since $v_1^+,v_2^+\notin U_\ell^*$, we have $\exc_{D_\ell}^+(v_i^+)=\texc_{D_\ell, U_\ell^*}^+(v_i)=1$ for each $i\in [2]$, as desired for \cref{lm:cleaningdegreeWA-V'endpoint}.}.
		If $\texc_{D_\ell, U_\ell^*}^+(V'\setminus U_\ell^*)=1$, then, note that by \cref{fact:texcS} and \cref{eq:cleaningdegreeWA-texc2}, $\texc_{D_\ell, U_\ell^*}^+(W)\geq \texc_{U_\ell^*}^+(D_\ell)-\texc_{D_\ell, U_\ell^*}^+(V'\setminus U_\ell^*)-\texc_{D_\ell, U_\ell^*}^+(U_\ell^*\cap V')\geq \frac{n}{2}-\eta n-1-1\geq 1$.
		Thus, we can let $v_1^+\in \tU_{U_\ell^*}^+(D_\ell)\cap (V'\setminus U_\ell^*)$ and $v_2^+\in \tU_{U_\ell^*}^+(D_\ell)\cap W$%
			\COMMENT{Again, since $v_1^+\notin U_\ell^*$, we have $\exc_{D_\ell}^+(v_1^+)=\texc_{D_\ell, U_\ell^*}^+(v_1)=1$, as required for \cref{lm:cleaningdegreeWA-V'endpoint}.}.
		If $\texc_{D_\ell, U_\ell^*}^+(V'\setminus U_\ell^*)=0$, then it is enough to show that $|\tU_{U_\ell^*}^+(D_\ell)\cap W|\geq 2$ (so that we can take distinct $v_1^+,v_2^+\in \tU_{U_\ell^*}^+(D_\ell)\cap W$, as desired for \cref{lm:cleaningdegreeWA-Wendpoint}). 
		
		Note that, by \cref{fact:texcS} and \cref{eq:cleaningdegreeWA-texc2}, $\texc_{D_\ell, U_\ell^*}^+(W)\geq \texc_{U_\ell^*}^+(D_\ell)-\texc_{D_\ell, U_\ell^*}^+(V'\setminus U_\ell^*)-\texc_{D_\ell, U_\ell^*}^+(U_\ell^*\cap V')\geq \frac{n}{2}-\eta n-0-1\geq 2$ and so, in particular, $\tU_{U_\ell^*}^+(D_\ell)\cap W\neq \emptyset$. Assume for a contradiction that $\tU_{U_\ell^*}^+(D_\ell)\cap W=\{v\}$ for some~$v\in W$. Note that since $\texc_{D_\ell, U_\ell^*}^+(v)=\texc_{D_\ell, U_\ell^*}^+(W)\geq 2$, $v\notin U^0(D_{\ell})$ and so, $U_\ell^*\subseteq V'$ and $\exc_D^+(v)\geq \exc_{D_\ell}^+(v)\geq \texc_{D_\ell, U_\ell^*}^+(v)$.
		Thus, since $|U_\ell^*|=|U_\ell^*\cap V'|\leq 1$, $|A^+|\leq\lceil\eta n\rceil$, and $2\ell\leq 2k\leq 8\lceil \eta n\rceil$,%
		\begin{align*}
			d_D^-(v)+\exc_D^+(v)&\stackrel{\text{\eqmakebox[d+exc][c]{\text{\cref{fact:exc}\cref{fact:exc-excv}}}}}{\leq} \texc(D)
			\stackrel{\text{\cref{eq:cleaningdegreeWA-texc1}}}{=}\texc_{U_\ell^*}^+(D_\ell)+2\ell+|A^+|\\
			&\stackrel{\text{\eqmakebox[d+exc][c]{\text{\cref{fact:texcS}}}}}{=}\texc_{D_\ell,U_\ell^*}^+(v)+|U_\ell^*|+2\ell+|A^+|
			\leq \exc_D^+(v)+1+9\lceil\eta n\rceil.
		\end{align*}
		But, by \cref{lm:cleaningdegreeWA-W} and \cref{fact:exc}\cref{fact:exc-dmin}, $d_D^-(v)=\frac{d_D(v)-\exc_D^+(v)}{2}\geq \frac{(20\eta-\varepsilon)n}{2}> 9\lceil\eta n\rceil+1$, a contradiction. Thus, $|\tU_{U_\ell^*}^+(D_\ell)\cap W|\geq 2$ and we can let $v_1^+,v_2^+\in \tU_{U_\ell^*}^+(D_\ell)\cap W$ be distinct.
		
		Now proceed analogously to pick $v_1^-,v_2^-\in \tU_{U_\ell^*}^-(D_\ell)\setminus((U_\ell^*\cap V')\cup \{v_1^+,v_2^+\})$ (this is possible since, for each~$i\in[2]$, $\texc_{D, U_\ell^*}^-(v_i^+)\leq 1$)%
.
		One can easily verify that \cref{lm:cleaningdegreeWA-V'endpoint,lm:cleaningdegreeWA-Wendpoint,lm:cleaningdegreeWA-partialdecomp} are satisfied.
	\end{proof}
	
	This completes the proof.
\end{proof}

\subsection{Deriving Lemma \ref{lm:cleaning}}

\begin{proof}[Proof of \cref{lm:cleaning}]
	Successively apply \cref{lm:cleaningW0,lm:cleaningallW*,lm:cleaningdegreeWA} as follows.
	\begin{steps}
		\item \textbf{Covering~$T[W_0]$.} 
		First, apply \cref{lm:cleaningW0} to obtain a good $(U^*, W, A)$-partial path decomposition~$\cP_1$ of~$T$ such that the following hold. Denote $D_1\coloneqq T\setminus \cP_1$ and $U_1^*\coloneqq U^*\setminus (V^+(\cP_1)\cup V^-(\cP_1))$.
		\begin{enumerate}[label=\upshape(\greek*)]
			\item $\texc(D_1)=\texc(T)-|\cP_1|$.\label{lm:cleaning-W0-good}
			\item $|\cP_1|\leq 2\varepsilon n$. \label{lm:cleaning-W0-size}
			\item $E(D_1[W_0])=\emptyset$. \label{lm:cleaning-W0-cleaningW0}
			\item If $|U^+(D_1)|=|U^-(D_1)|=1$, then $e(U^-(D_1),U^+(D_1))=0$ or $\texc(D_1)-\exc(D_1)\geq 2$.
			\item Each $v\in U_1^*$ satisfies $d_{D_1}^+(v)=d_{D_1}^-(v)\leq \texc(D_1)-1$.\label{lm:cleaning-W0-U*}
		\end{enumerate}
		In particular, observe that, by \cref{prop:texcT}, \cref{lm:cleaning}\cref{lm:cleaning-A}, and \cref{lm:cleaning-W0-size}, the following hold.
		\begin{enumerate}[resume, label=\upshape(\greek*)]
			\item $\texc(D_1)\geq \frac{n}{2}-2\varepsilon n$ and, if $W_A\neq \emptyset$, then $\texc(D_1)\geq \frac{n}{2}+9\eta n$.\label{lm:cleaning-W0-excess}
			\item For each $v\in V$, $d_{D_1}(v)\geq (1-5\varepsilon)n$. \label{lm:cleaning-W0-degree}
			\item For each $v\in W_*$, $|\exc_{D_1}(v)|>(1-21\eta)n$.\label{lm:cleaning-W0-excessW}
			\item For each $\diamond\in \{+,-\}$, if $|W_A^\diamond|\geq 2$, then $\exc_{D_1}^\diamond(v)< \eta n$ for each $v\in V$ and, if $|W_A^\diamond|=1$, then, for each $v\in V$ and $w\in W_A^\diamond$, $\exc_{D_1}^\diamond(v)\leq \exc_{D_1}^\diamond(w)+5\varepsilon n$. \label{lm:cleaning-W0-WA}
		\end{enumerate}

		\item \textbf{Covering the remaining edges of~$T[W]$ and decreasing the degree of the vertices in $W_*\cup W_A$ when $W_*\neq \emptyset$.}
		If $W_*=\emptyset$, then let $\cP_2\coloneqq \emptyset$. Otherwise, note that by \cref{prop:sizeU*good}, $|U_1^*|=\texc(D_1)-\exc(D_1)$ and let~$\cP_2$ be the good $(U_1^*, W,A)$-partial path decomposition of~$D_1$ obtained by applying \cref{lm:cleaningallW*} with $D_1$, $U_1^*$, and $5\varepsilon$ playing the roles of $D$, $U^*$, and $\varepsilon$.
		Denote $D_2\coloneqq D_1\setminus \cP_2$ and $U_2^*\coloneqq U_1^*\setminus (V^+(\cP_2)\cup V^-(\cP_2))$.
		Then, note that, if $W_*\neq 0$, then the following hold.
		\begin{enumerate}[label=\upshape(\Roman*)]
			\item $\texc(D_2)=\texc(D_1)-|\cP_2|$. \label{lm:cleaning-allW*-good}
			\item $E(D_2[W])=\emptyset$. \label{lm:cleaning-allW*-cleaningW*}
			\item $N^\pm(D_1)-N^\pm(D_2)\leq 88\eta n$.\label{lm:cleaning-allW*-size0}
			\item For each $v\in W_*\cup W_A$, $(1-3\sqrt{\eta})n\leq d_{D_2}(v)\leq (1-4\eta)n$.\label{lm:cleaning-allW*-degreeW*}
			\item For each $v\in W_0$, $d_{D_2}(v)\geq (1-3\sqrt{\eta})n$ and $d_{D_2}^{\min}(v)\geq 5\eta n$.\label{lm:cleaning-allW*-degreeW0}
			\item For each $v\in V'$, $d_{D_2}(v)\geq (1-8\sqrt{\varepsilon})n$.\label{lm:cleaning-allW*-degreeV'}
			\item If $|W_*^+|,|W_*^-|\leq 1$, then each $v\in W_*$ satisfies $|\exc_{D_2}(v)|=d_{D_2}(v)$.\label{lm:cleaning-allW*-exc=d}
			\item Each $v\in U_2^*$ satisfies $d_{D_2}^+(v)=d_{D_2}^-(v)\leq \texc(D_2)-1$.\label{lm:cleaning-allW*-U*}
		\end{enumerate}
		Note that, by \cref{lm:cleaning-allW*-degreeW*}, the following holds.
		\begin{enumerate}[resume,label=\upshape(\Roman*)]
			\item Each $v\in W_*$ satisfies $|\exc_{D_2}(v)|\geq |\exc_T(v)|-3\sqrt{\eta}n\geq (1-4\sqrt{\eta})n$.\label{lm:cleaning-allW*-excessW*}
		\end{enumerate}
		Thus, \cref{lm:cleaning-allW*-exc=d} implies the following.
		\begin{enumerate}[resume,label=\upshape(\Roman*)]
			\item If $\texc(D_2)\leq 2\left\lceil\frac{n}{2}\right\rceil-\lceil\eta n\rceil$, then $|\exc_{D_2}(v)|=d_{D_2}(v)$ for each $v\in W_*$.\label{lm:cleaning-allW*-exc<2d}
		\end{enumerate}
		
		\item \textbf{Decreasing the degree of the vertices in~$W_A$ when $W_*=\emptyset$.}
		If $W_*\neq\emptyset$ or $W_A=\emptyset$, then let $\cP_3\coloneqq \emptyset$. Assume $W_*=\emptyset$ and $W_A\neq \emptyset$. Recall that, by construction, $D_2=D_1$ and $U_2^*=U_1^*$. In particular, \cref{lm:cleaning-W0-cleaningW0,lm:cleaning-W0-degree,lm:cleaning-W0-excess} are satisfied and $|U_2^*|=\texc(D_2)-\exc(D_2)$. Let~$\cP_3$ be the good $(U_2^*, W, A)$-partial path decomposition of~$D_2$ obtained by applying \cref{lm:cleaningdegreeWA} with $D_2$, $U_2^*$, and $5\varepsilon$ playing the roles of $D$, $U^*$, and $\varepsilon$.
		Denote $D_3\coloneqq D_2\setminus \cP_3$ and note that, if $W_*=\emptyset$ and $W_A\neq \emptyset$, then the following hold.
		\begin{enumerate}[label=\upshape(\Alph*)]
			\item $\texc(D_3)=\texc(D_2)-|\cP_3|$.\label{lm:cleaning-degreeWA-good}
			\item Each $v\in U_2^*\setminus (V^+(\cP_3)\cup V^-(\cP_3))$ satisfies $d_{D_3}^+(v)=d_{D_3}^-(v)\leq \texc(D_3)-1$.\label{lm:cleaning-degreeWA-U*}
			\item $|\cP_3|=8\lceil\eta n\rceil$.\label{lm:cleaning-degreeWA-size}
			\item For each $v\in W$, $d_{D_3}(v)\leq d_{D_2}(v)-12\lceil\eta n\rceil$.\label{lm:cleaning-degreeWA-degreeW}
			\item For each $v\in V'$, $d_{D_2}(v)-8\lceil\eta n\rceil\leq d_{D_3}(v)\leq d_{D_2}(v)-8\lceil\eta n\rceil+1$.\label{lm:cleaning-degreeWA-degreeV'}
		\end{enumerate}
		
		\item \textbf{Checking the assertions of \cref{lm:cleaning}.}
	Let $\cP\coloneqq \bigcup_{i\in [3]}\cP_i$ and $D\coloneqq T\setminus \cP=D_3$. If $W_*\neq \emptyset$ or $W_A=\emptyset$, then let $d\coloneqq \min\{\left\lceil\frac{n}{2}\right\rceil-\lceil\eta n\rceil,\texc(D)-\lceil\eta n\rceil\}$.
	If $W_*=\emptyset$ and $W_A\neq \emptyset$, then let \NEW{$d\coloneqq \left\lceil\frac{n}{2}\right\rceil -5\lceil\eta n\rceil$}\OLD{$d\coloneqq \min \{\left\lceil\frac{n}{2}\right\rceil -5\lceil\eta n\rceil, \texc(D)-\lceil\eta n\rceil\}$}. In both cases,~$\cP$ is a good $(U^*,W,A)$-partial path decomposition of~$T$ by \cref{fact:combinespecial}. \NEW{Note that \cref{lm:cleaning-good} follows immediately from \cref{lm:cleaning-W0-U*}, \cref{lm:cleaning-allW*-U*}, and \cref{lm:cleaning-degreeWA-U*}, while \cref{lm:cleaning-W} follows immediately from \cref{lm:cleaning-W0-cleaningW0} and \cref{lm:cleaning-allW*-cleaningW*}. If $W_*=\emptyset$, then \cref{lm:cleaning-exc<2d} holds vacuously, otherwise \cref{lm:cleaning-exc<2d} follows from \cref{lm:cleaning-allW*-exc<2d} and the definition of $d$.}
	
	\NEW{We now verify \cref{lm:cleaning-exc>d}. If $W_*\neq\emptyset$ or $W_A=\emptyset$, then \cref{lm:cleaning-exc>d} holds immediately by definition of $d$. Suppose that $W_*=\emptyset$ and $W_A\neq \emptyset$. Then, $D_2=D_1$, so \cref{lm:cleaning-W0-excess}, \cref{lm:cleaning-degreeWA-good}, and \cref{lm:cleaning-degreeWA-size} imply that $\texc(D)\geq \frac{n}{2}+9\eta n-|\cP_3| \geq \lceil\frac{n}{2}\rceil-4\lceil\eta n\rceil=d+\lceil\eta n\rceil$. Therefore, \cref{lm:cleaning-exc>d} holds.}
	
	\OLD{Moreover, one can easily verify that \cref{lm:cleaning-d,lm:cleaning-good,lm:cleaning-W,lm:cleaning-sizeU,lm:cleaning-exc>d,lm:cleaning-exc<2d} are satisfied and, if $d\neq \texc(D)-\lceil\eta n\rceil$, \cref{lm:cleaning-degreeV',lm:cleaning-degreeW*,lm:cleaning-degreeW0} hold too.}%
		\COMMENT{\textbf{Case $W_*\neq \emptyset$ or $W_A=\emptyset$.}\\
		For \cref{lm:cleaning-d}, note that, if both $W_*,W_A=\emptyset$, then, by \cref{lm:cleaning-W0-excess}, $\texc(D_3)\geq \left\lceil\frac{n}{2}\right\rceil-2\varepsilon n$ and, if $W_*\neq \emptyset$, by \cref{lm:cleaning-allW*-degreeW*}, $\texc(D_3)\geq (1-3\sqrt{\eta})n$.\\ 
		\cref{lm:cleaning-good} holds by \cref{lm:cleaning-W0-U*,lm:cleaning-allW*-U*}.\\
		\cref{lm:cleaning-W} holds by \cref{lm:cleaning-W0-cleaningW0} and \cref{lm:cleaning-allW*-cleaningW*}.\\ 
		\cref{lm:cleaning-sizeU} holds by \cref{lm:cleaning-W0-size} and \cref{lm:cleaning-allW*-size0}.\\
		\cref{lm:cleaning-exc>d} holds by our choice of $d$.\\
		\cref{lm:cleaning-exc<2d} holds by \cref{lm:cleaning-allW*-exc<2d}.\\
		If $d=\left\lceil\frac{n}{2}\right\rceil-\lceil\eta n\rceil$, then the upper bounds in \cref{lm:cleaning-degreeV',lm:cleaning-degreeW0} are clearly satisfied, the lower bound in \cref{lm:cleaning-degreeV'} follow from \cref{lm:cleaning-W0-degree} and \cref{lm:cleaning-allW*-degreeV'}, the lower bounds in \cref{lm:cleaning-degreeW0} follow from \cref{lm:cleaning-W0-degree} and \cref{lm:cleaning-allW*-degreeW0}, and \cref{lm:cleaning-degreeW*} holds by \cref{lm:cleaning-allW*-degreeW*}.\\
		}%

	\NEW{Next, we check \cref{lm:cleaning-d}. If both $W_*,W_A=\emptyset$, then \cref{lm:cleaning-W0-excess} implies that $\texc(D)\geq \left\lceil\frac{n}{2}\right\rceil-2\varepsilon n$ and so
	$\left\lceil\frac{n}{2}\right\rceil-\lceil\eta n\rceil -2\varepsilon n\leq d\leq \left\lceil\frac{n}{2}\right\rceil-\lceil\eta n\rceil$.
	If $W_*\neq \emptyset$, then $D=D_2$, so \cref{lm:cleaning-allW*-excessW*} implies that $\texc(D)\geq (1-4\sqrt{\eta})n$ and thus $d=\left\lceil\frac{n}{2}\right\rceil-\lceil\eta n\rceil$. If $W_*=\emptyset$ and $W_A\neq \emptyset$, then $d=\left\lceil\frac{n}{2}\right\rceil-5\lceil\eta n\rceil$ and so \cref{lm:cleaning-d} holds.}

	\NEW{For \cref{lm:cleaning-sizeU}, note that by \cref{cor:Ntexc} and \cref{lm:cleaning-W0-size}, $N^\pm(T)-N^\pm(D_1)\leq |\cP_1|\leq 2\varepsilon n\leq \eta n$. Thus, if $W_*\neq \emptyset$ or $W_A=\emptyset$, then \cref{lm:cleaning-sizeU} follows from \cref{lm:cleaning-allW*-size0}. We may therefore assume that $W_*=\emptyset$ and $W_A\neq \emptyset$. Then, $D_2=D_1$ and so \cref{cor:Ntexc} and \cref{lm:cleaning-degreeWA-size} imply that $N^\pm(D_1)-N^\pm(D)\leq |\cP_3|\leq 8 \lceil\eta n\rceil\leq 88\eta n$. Therefore, \cref{lm:cleaning-sizeU} holds.}
	
	\NEW{It remains to check \cref{lm:cleaning-degreeV',lm:cleaning-degreeW*,lm:cleaning-degreeW0}. First, suppose that $W_*\neq \emptyset$ or $W_A=\emptyset$, and $d=\left\lceil\frac{n}{2}\right\rceil-\lceil\eta n\rceil$. 
	Then, the upper bounds in \cref{lm:cleaning-degreeV',lm:cleaning-degreeW0} are clearly satisfied. 
	The lower bound in \cref{lm:cleaning-degreeV'} follows from \cref{lm:cleaning-W0-degree} and \cref{lm:cleaning-allW*-degreeV'}, while \cref{lm:cleaning-degreeW*} holds by \cref{lm:cleaning-allW*-degreeW*}. By \cref{lm:cleaning-defW} and \cref{fact:exc}\cref{fact:exc-dmin}, each $v\in W_0$ satisfies $d_T^{\min}(v)\geq \frac{20\eta n-1}{2}\geq 9\eta n$ and so the lower bounds in \cref{lm:cleaning-degreeW0} follow from \cref{lm:cleaning-W0-degree} and \cref{lm:cleaning-allW*-degreeW0}. Thus, \cref{lm:cleaning-degreeV',lm:cleaning-degreeW*,lm:cleaning-degreeW0} hold if $W_*\neq \emptyset$ or $W_A=\emptyset$, and $d=\left\lceil\frac{n}{2}\right\rceil-\lceil\eta n\rceil$.}
	
	\NEW{Next,} assume that \NEW{$W_*\neq \emptyset$ or $W_A=\emptyset$, and} $d=\texc(D)-\lceil\eta n\rceil$. \NEW{Then}\OLD{By construction}, $2d+2\lceil\eta n\rceil=2\texc(D)\NEW{\geq 2\Delta^0(D)}\geq d_{D}(v)$ for each $v\in V$ and so the upper bounds of \cref{lm:cleaning-degreeV',lm:cleaning-degreeW0} hold.
	Moreover, $\texc(D)\leq \left\lceil\frac{n}{2}\right\rceil$ and so, by \cref{lm:cleaning-allW*-excessW*}, $W_*=\emptyset$.
	\NEW{By assumption, this implies that}\OLD{If} $W_A=\emptyset$ \NEW{and so}\OLD{then} \cref{lm:cleaning-degreeW*} holds vacuously. \NEW{Moreover,}\OLD{and by \cref{lm:cleaning-W0-degree},} each $v\in V$ satisfies
	\begin{align*}
		d_D^{\min}(v)\stackrel{\text{\NEW{\cref{lm:cleaning-W0-degree}}}}{\geq} d_T^{\min}(v)-5\varepsilon n
		\stackrel{\text{\NEW{\cref{fact:exc}\cref{fact:exc-dmin}}}}{=}\frac{n-1-|\exc_T(v)|}{2}-5\varepsilon n
		\stackrel{\text{\NEW{\cref{lm:cleaning-defW}}}}{\geq}\frac{20\eta n-1}{2}-5\varepsilon n
		\geq \lceil\eta n\rceil+1
	\end{align*} 
	and
	\begin{align*}
		d_D(v)\stackrel{\text{\NEW{\cref{lm:cleaning-W0-degree}}}}{\geq} (1-5\varepsilon)n\geq 2\left(\left\lceil \frac{n}{2}\right\rceil-\lceil\eta n\rceil\right)+2\lceil\eta n\rceil-6\varepsilon n\geq 2d+2\lceil\eta n\rceil-6\varepsilon n,
	\end{align*} 
	so the lower bounds in \cref{lm:cleaning-degreeV',lm:cleaning-degreeW0} hold. \NEW{Therefore, \cref{lm:cleaning-degreeV',lm:cleaning-degreeW*,lm:cleaning-degreeW0} hold if $W_*\neq \emptyset$ or $W_A=\emptyset$, and $d=\texc(D)-\lceil\eta n\rceil$.}
	
	\NEW{We may therefore assume that $W_*=\emptyset$ and $W_A\neq \emptyset$. Recall that, in this case, $d=\left\lceil\frac{n}{2}\right\rceil -5\lceil\eta n\rceil$ and $\cP_2=\emptyset$.
	First, observe that each $v\in W$ satisfies 
	\begin{align*}
		d_D(v)\geq n-1- 2|\cP|\stackrel{\text{\cref{lm:cleaning-W0-size},\cref{lm:cleaning-degreeWA-size}}}{\geq} n-17\lceil\eta n\rceil
	\end{align*}
	and if $v\in W_0$, then 
	\begin{align*}
		d_D^{\min}(v)&\geq d_T^{\min}(v)-|\cP|\stackrel{\text{\cref{fact:exc}\cref{fact:exc-dmin}}}{\geq} \frac{n-1-|\exc_T(v)|}{2}-|\cP|\stackrel{\text{\cref{lm:cleaning-defW},\cref{lm:cleaning-W0-size},\cref{lm:cleaning-degreeWA-size}}}{\geq}
		\frac{20\eta n-1}{2}-2\varepsilon n-8\lceil\eta n\rceil\\
		&\geq \lceil\eta n\rceil.
	\end{align*}
	Thus, the lower bounds in \cref{lm:cleaning-degreeW*,lm:cleaning-degreeW0} are satisfied. Moreover, each $v\in W$ satisfies
	\begin{align*}
		d_D(v)\stackrel{\text{\cref{lm:cleaning-degreeWA-degreeW}}}{\leq} d_{D_2}(v)-12\lceil\eta n\rceil\leq n-1-12\lceil\eta n\rceil\leq 2d-1-2\lceil\eta n\rceil
	\end{align*}
	and so the upper bounds in \cref{lm:cleaning-degreeW*,lm:cleaning-degreeW0} hold. 
	Finally, note that each $v\in V'$ satisfies
	\begin{align*}
		d_D(v)\stackrel{\text{\cref{lm:cleaning-degreeWA-degreeV'}}}{\leq} d_{D_2}(v)-8\lceil\eta n\rceil+1
		\leq n-8\lceil\eta n\rceil\leq 2d+2\lceil\eta n\rceil
	\end{align*}
	and
	\begin{align*}
		d_D(v)&\stackrel{\text{\eqmakebox[cleaning][c]{\text{\cref{lm:cleaning-degreeWA-degreeV'}}}}}{\geq} d_{D_2}(v)-8\lceil\eta n\rceil\geq d_T(v)-2|\cP_1| -8\lceil\eta n\rceil\\
		&\stackrel{\text{\eqmakebox[cleaning][c]{\text{\cref{lm:cleaning-W0-size}}}}}{\geq} n-1-4\varepsilon n -8\lceil\eta n\rceil
		\geq 2d-5\varepsilon n+2\lceil\eta n\rceil,
	\end{align*}
	so \cref{lm:cleaning-degreeV'} holds. 
	Therefore, \cref{lm:cleaning-degreeV',lm:cleaning-degreeW*,lm:cleaning-degreeW0} hold if $W_*= \emptyset$ and $W_A\neq\emptyset$. This completes the proof of \cref{lm:cleaning}.}%
\qedhere
	\end{steps}
\end{proof}

%% file: Constructing_Layouts_Proof.tex
	\onlyinsubfile{
		\setcounter{section}{12}
\section{Constructing layouts: proof of Lemma \ref{lm:layouts}}}

We will prove \cref{lm:layouts} as follows. In \cref{step:endpoints}, we choose a set~$\tE$ of auxiliary edges which ``neutralise" the excess of the vertices in~$D$. In \cref{step:bulidlayouts}, we then subdivide these edges into paths which form a layout~$(\tL, \tF)$. In \cref{step:coverW}, we subdivide the paths in~$\tL$ further to obtain layouts $(\hL_1, \hF_1), \dots, (\hL_\ell, \hF_\ell)$ which cover the edges of~$D\setminus A$ at~$W$ in such a way that \cref{lm:layouts-degreeW1,lm:layouts-degreeW2} are satisfied. Finally, in \cref{step:degreeV'}, we adjust the degrees of the vertices in~$V'$ \NEW{so that they satisfy \cref{lm:layouts-degreeV'}}\OLD{as follows}. \NEW{To achieve this, we proceed as follows.} For those vertices~$v\in V'$ where the current layouts would result in a degree which is too small after the approximate decomposition, we add~$v$ as an isolated vertex to some of the layouts. For vertices~$v\in V'$ whose degree would be too large, we subdivide two edges from a suitable layout and include~$v$ into both of the resulting paths. \NEW{Recall that the relevant definitions involving layouts were introduced in \cref{sec:approxdecomp,sec:W,sec:pathL}.}

\begin{proof}[Proof of \cref{lm:layouts}]
	Let $W^\pm\coloneqq W\cap U^\pm(D)$.
	Denote $\phi\coloneqq \phi^++\phi^-$.
	\NEW{Note that since $A$ is a $(W_1,V')$-absorbing set, \cref{def:A} implies that $A$ does not contain any edge incident to $W_2$ and so $A$ is also a $(W,V')$-absorbing set. Thus, for simplicity, we can let $W$ (rather than $W_1$) play the role of $W$ in the auxiliary excess notation.} 
	For each~$v\in V$, define 
	\begin{equation}\label{eq:layouts-hexc}
		\hexc^\pm(v)\coloneqq \texc_{D,U^*,\NEW{W,A}}^\pm(v)-\phi^\pm(v)
	\end{equation} \NEW{to be the excess at $v$ that we want to cover with the layouts (i.e.\ the number of paths which we want to start/end at $v$)}.
	Let $\hU^\pm\coloneqq \{v\in V\mid \hexc^\pm(v)>0\}$. 
	Note that \cref{eq:texcA} \NEW{and \cref{lm:layouts-A} imply}\OLD{implies} that \[\hexc(D)\coloneqq\sum_{v\in V}\hexc^+(v)=\sum_{v\in V}\hexc^-(v)=\texc(D)-\lceil\eta n\rceil.\]
	If \NEW{$\hexc(D)-d<\sqrt{\varepsilon}n$}\OLD{$(\texc(D)-\lceil\eta n\rceil)-d<\sqrt{\varepsilon}n$}, \NEW{then} let \NEW{$\ell\coloneqq \hexc(D)$}\OLD{$\ell\coloneqq \texc(D)-\lceil\eta n\rceil$}; otherwise, let $\ell\coloneqq d$. Note that, by \cref{lm:layouts-exc>d}, $\ell\geq d$. Thus, \cref{lm:layouts-ell} holds, as desired.
	Observe that either
	\begin{equation}\label{eq:layouts-ell}
		\hexc(D)-\ell\geq \sqrt{\varepsilon}n \quad  \text{or} \quad \hexc(D)-\ell=0.
	\end{equation}
	\OLD{\begin{equation*}
		(\texc(D)-\lceil\eta n\rceil)-\ell\geq \sqrt{\varepsilon}n \quad  \text{or} \quad (\texc(D)-\lceil\eta n\rceil)-\ell=0.
	\end{equation*}}
	
	We claim that each~$v\in V'$ satisfies
	\begin{equation}\label{eq:U*}
		d_D^\pm(v)\leq \texc(D)-\texc_{D,U^*}^\mp(v).
	\end{equation}
	Indeed, if~$v\in U^*$, then \NEW{\cref{eq:texc} implies that $\texc_{D,U^*}^\pm(v)=1$}\OLD{$\texc_{D,U^*}^\pm(v)\leq 1$} and so \cref{eq:U*} holds by \cref{lm:layouts-U*}. We may therefore assume that~$v\notin U^*$. Suppose without loss of generality that $d_D^+(v)\geq d_D^-(v)$. Then, \NEW{\cref{eq:texc} implies that} $\texc_{D, U^*}^-(v)\NEW{=\exc_D^-(v)}=0$ and so $\texc(D)-\texc_{D,U^*}^-(v)\geq \Delta^0(D)\geq d_D^+(v)$. Finally, \NEW{\cref{eq:texc} implies that $\texc_{D, U^*}^+(v)= \exc_D^+(v)$}\OLD{$\texc_{D, U^*}^+(v)= \exc_D^+(v)$} and so $\texc(D)-\texc_{D, U^*}^+(v)\geq d_D^+(v)-\exc_D^+(v)=d_D^-(v)$, as desired.
	
	Note that throughout this proof, given a multiset~$L'$ of paths, the corresponding edge set~$F'$ in the layout~$(L',F')$ we construct will always satisfy $F'=E(L') \cap \mathcal{F} = E(L') \cap (E_W(D) \setminus A)$. 
	
	\begin{steps}
		\item \textbf{Choosing suitable endpoints.}\label{step:endpoints}
		Let $s\coloneqq \hexc(D)=\texc(D)-\lceil\eta n\rceil$.
		In this step, we will select suitable endpoints for the $s$ (non-trivial) paths in~$L$. \NEW{Each $v\in V'$ will be used precisely $\hexc^+(v)$ times as a starting point and $\hexc^-(v)$ times as an ending point. We now fix the endpoints of each path by defining a multidigraph $\tE$ on $V$ such that each $e\in \tE$ corresponds to a path of shape $e$. Hence, $|E|=s$. Note that $\hexc^+(v)\neq 0\neq \hexc^-(v)$ if and only if $\hexc^+(v)=1=\hexc^-(v)$. We now formalise this in the following paragraph.}
		
		For each $\diamond\in \{+,-\}$, let $v_1^\diamond, \dots, v_s^\diamond\in \hU^\diamond$ be such that, for each~$v\in \hU^\diamond$, there exist exactly~$\hexc^\diamond(v)$ indices~$i\in [s]$ for which~$v=v_i^\diamond$. Since $s\geq d >1$ \NEW(by \cref{lm:layouts-d,lm:layouts-exc>d}) and each $v\in \hU^+\cap \hU^-$ satisfies $\hexc^+(v)=\hexc^-(v)=1$ \NEW{(by \cref{eq:texc})}, we may assume without loss of generality that, for each~$i\in [s]$, $v_i^+\neq v_i^-$.
		Let $\tE\coloneqq \{v_j^+v_j^-\mid j \in [s]\}$.
		Note that,
		\begin{equation}\label{eq:layouts-endpoints-size}
			|\tE|=\hexc(D)=\texc(D)-\lceil\eta n\rceil
		\end{equation} 
		and, for each~$v \in V$, 
		\begin{align}\label{eq:layouts-endpoints-excess}
		d_{\tE}^\pm(v)&=\hexc^\pm(v).
		\end{align}
		\NEW{Let $v\in W$. By \cref{def:A}, we have $d_A^\pm(v)\leq \exc_D^\pm(v)$ and so $\exc_{D\setminus A}^\pm(v)=\exc_D^\pm(v)-d_A^\pm(v)$. Moreover, \cref{lm:layouts-U*} implies that $v\notin U^*$ and \cref{lm:layouts-A} implies that $\phi(v)=0$. Therefore,
		\begin{align}\label{eq:layouts-endpoints-excessA}
			\exc_{D\setminus A}^\pm (v)=\exc_D^\pm(v)-d_A^\pm(v)\stackrel{\text{\cref{eq:texc}}}{=}\texc_{D,U^*,W,A}^\pm(v)=\hexc^\pm(v)\stackrel{\text{\cref{eq:layouts-endpoints-excess}}}{=} d_{\tE}^\pm(v).
		\end{align}}
		
		\item \textbf{Constructing layouts.}\label{step:bulidlayouts}
		\NEW{In \cref{step:bulidlayouts,step:coverW}, we will transform $\tE$ into a $W$-exceptional layout $(\hL, \hF)$. Initially, we set $(\hL,\hF)=(\tE,\emptyset)$ where each edge in $\tE$ is considered as a path. To be a $W$-exceptional layout, each path in $\hL$ requires an edge entirely in $V'$ and $\tF$ must contain $E_W(L)$. For this, we proceed roughly as follows.}
		
		\NEW{Suppose that the path $v^+v^-$ does not lie entirely in $V'$, say $v^+\in W$ and $v^-\in V'$. We pick $u\in N_{D\setminus A}(v^+)$ and replace the path $v^+v^-$ with the subdivided path $v^+uv^-$ and add $v^+u$ into the set $\hF$ of fixed edges. (Note that $uv^-$ lies entirely in $V'$.)}
		
		\NEW{More precisely,} recall that $\cF= E_W(D)\setminus A$ and $D'= D\setminus \cF$.
		In this step, we will use~$\tE$ to construct a layout~$(\tL,\tF)$ such that the following hold. 
		\begin{enumerate}[label=\upshape(\greek*)]
			\item $(\tL,\tF)$ is a $W$-exceptional layout such that $\tF\subseteq \cF$ and~$\tL$ contains no isolated vertex.\label{lm:layouts-buildlayouts-layouts}
			\item $\tL$ has shape~$\tE$ (and thus, by \cref{eq:layouts-endpoints-excess},~$(\tL,\tF)$ is~$U^*$-path consistent with respect to~$(D',\cF)$).\label{lm:layouts-buildlayouts-shape}
			\item For each~$P\in \tL$ and~$v\in V(P)\cap V'$,~$v$ is an endpoint of~$P$ or has (in~$P$) a neighbour in~$W$.\label{lm:layouts-buildlayouts-V'}
			\item For each~$P\in \tL$,~$V^0(P)\subseteq V'$.\label{lm:layouts-buildlayouts-internalvertices}
			\item Each~$P\in \tL$ has at most~$4$ vertices and contains an edge which lies entirely in~$V'$.\label{lm:layouts-buildlayouts-pathlength}%
		\end{enumerate}
		\NEW{Properties \cref{lm:layouts-buildlayouts-V',lm:layouts-buildlayouts-internalvertices,lm:layouts-buildlayouts-pathlength} mean that the paths in $\tL$ are only obtained by subdividing, with vertices of $V'$, the edges in $\tE$ which are incident to $W$. Property \cref{lm:layouts-buildlayouts-V'} will ensure that no vertex of $V'$ belongs to too many layouts (as desired for \cref{lm:layouts-V'}) and \cref{lm:layouts-buildlayouts-pathlength} will ensure that the layouts are not too large (as desired for \cref{lm:layouts-size}). Property \cref{lm:layouts-buildlayouts-internalvertices} means that we have not yet incorporated the exceptional vertices as internal vertices. This will give us more flexibility in \cref{step:coverW} when we cover the remaining edges incident to $W$.} 
		
		Initially, let $\tL^0\coloneqq \tE$ and $\tF^0\coloneqq \emptyset$.
		Let $w_1,\dots, w_k$ be an enumeration of $W\cap V(\tL^0)$. 
		\NEW{We will consider each $w_i$ in turn and, at each stage $i$, subdivide all the edges incident to $w_i$. Let $i\in [k]$. By \cref{eq:layouts-endpoints-excessA}, $\exc_{D\setminus A}(w_i)\neq 0$ and so recall from \cref{eq:Nmax} that $N_{D\setminus A}^{\rm max}(w_i)$ denotes the outneighbourhood of $w_i$ in $D\setminus A$ if $\exc_{D\setminus A}(w_i)>0$ and the inneigbourhood of $w_i$ in $D\setminus A$ otherwise}\OLD{Note that, by \cref{eq:layouts-endpoints-excess}, for each~$i\in [k]$, $\exc_{D\setminus A}(w_i)\neq 0$}.
		
		Assume inductively that for some $0\leq m\leq k$, we have constructed, for each~$i\in [m]$, a multiset of paths~$\tL^i$ and a set of edges~$\tF^i$ such that the following are satisfied.
		\begin{enumerate}[label=\upshape(\Roman*)]
			\item Let $i\in [m]$. Let $S_i\coloneqq \{e\in E(\tL^{i-1}) \mid w_i\in V(e)\}$. Then,			
			$\tL^i$ is the multiset of paths obtained from~$\tL^{i-1}$ by subdividing each edge~$e\in S_i$ with some vertex $z_e\in N_{D\setminus A}^{\max}(w_i)\cap V'$, where the vertices~$z_e$ are distinct for different edges~$e\in S_i$. \NEW{(I.e.\ $\tL^i$ is obtained by subdividing, with a neighbour of $w_i$, all the edges of $\tL^{i-1}$ which are incident to $w_i$.)} \label{lm:layouts-buildlayouts-subdivisionL}
			\item For each~$i \in [m]$, $\tF^i = \tF^{i-1}\cup E_{\{w_i\}}(\tL^i)$. \label{lm:layouts-buildlayouts-subdivisionF}
		\end{enumerate}
		\NEW{Note that \cref{lm:layouts-buildlayouts-subdivisionL} and \cref{eq:layouts-endpoints-excess} imply that the following holds.
		\begin{enumerate}[resume,label=\upshape(\Roman*)]
			\item For all $i\in [m]$ and $v\in W$, $d_{\tL^i}^\pm(v)=d_{\tE}^\pm(v)=\hexc^\pm(v)$.\label{lm:layouts-buildlayouts-degreeW}
		\end{enumerate}
		Moreover,} \cref{lm:layouts-buildlayouts-subdivisionL,lm:layouts-buildlayouts-subdivisionF} imply that, for each~$i\in[m]$,~$\tF^i$ is a set of edges obtained from~$\tF^{i-1}$ by adding all the edges of the form~$w_iz_e$ or~$z_ew_i$ from \cref{lm:layouts-buildlayouts-subdivisionL}. In particular, $\tF^i \subseteq E_{\{w_j\mid j\in [i]\}}(D) \setminus A\subseteq \cF$ is satisfied.
		
		If~$m=k$, then let~$\tL\coloneqq \tL^k$ and~$\tF\coloneqq \tF^k$. 
		Observe that \cref{lm:layouts-buildlayouts-shape,lm:layouts-buildlayouts-layouts,lm:layouts-buildlayouts-internalvertices,lm:layouts-buildlayouts-V',lm:layouts-buildlayouts-pathlength} hold. \NEW{Indeed, \cref{lm:layouts-buildlayouts-subdivisionL} implies that $\tL$ is obtained by subdividing $\tL^0=\tE$, so \cref{lm:layouts-buildlayouts-shape} holds. By \cref{lm:layouts-buildlayouts-subdivisionL}, all these subdivisions are done with vertices of $V'$. Thus, \cref{lm:layouts-buildlayouts-V',lm:layouts-buildlayouts-internalvertices,lm:layouts-buildlayouts-pathlength} are satisfied. (For \cref{lm:layouts-buildlayouts-pathlength}, note that at each stage the paths all contain at most two vertices of $W$, so each edge in $\tL^0=\tE$ is subdivided at most twice.) Moreover, this implies that $\tL$ is a set of non-trivial paths which all contain at least one edge whose endpoints are both in $V'$. By \cref{lm:layouts-buildlayouts-subdivisionF}, $\tF$ consists of the edges of $\tL$ which are incident to $W$. Altogether, this implies that $(\tL,\tF)$ is a $W$-exceptional layout. As mentioned above, \cref{lm:layouts-buildlayouts-subdivisionL,lm:layouts-buildlayouts-subdivisionF} imply that $\tF\subseteq \cF$ and so \cref{lm:layouts-buildlayouts-layouts} is satisfied, as desired.}
		
		We may therefore assume that~$m<k$.
		\NEW{By assumption and \cref{lm:layouts-buildlayouts-degreeW}, $S_{m+1}\coloneqq \{e\in E(\tL^m)\mid w_{m+1}\in V(e)\}\neq \emptyset$ and so $w_{m+1}\notin U^0(D)$.}%
		\OLD{Note that, by \cref{eq:layouts-endpoints-excess},~$w_{m+1}\notin U^0(D)$.}
		We may therefore assume without loss of generality that~$w_{m+1}\in W^+$%
			\COMMENT{$w_{m+1}\in V(\tL^0)=V(\tE)$ implies $w_{m+1}\notin U^0(D)$.}.
		This implies that $\exc_D(w_{m+1})=|\exc_D(w_{m+1})|=\exc_D^+(w_{m+1})$ and $d_D^{\min}(w_{m+1})=d_D^-(w_{m+1})$.
		\OLD{Moreover, by \cref{def:A}, \cref{eq:texc}, and since, by \cref{lm:layouts-A}, $w_{m+1}\notin X^+$, we have		
		\begin{equation*}
			\hexc^+(w_{m+1})=\texc_{D, U^*}^+(w_{m+1})=\exc_D^+(w_{m+1})-d_A^+(w_{m+1})= \exc_{D\setminus A}^+(w_{m+1}).
		\end{equation*}
		Construct $(\tL^{m+1},\tF^{m+1})$ as follows.		
		Define the multiset~$X$ by $X\coloneqq \{v \mid w_{m+1}v \in E(\tL^m)\}$.}\NEW{We now subdivide all the edges in $S_{m+1}$ using Hall's theorem as follows. Let $Y\coloneqq N_{D\setminus A}^+(w_{m+1})$ and observe that, by \cref{lm:layouts-W}, $Y\subseteq V'$.}
		Construct an auxiliary bipartite graph~$G$ on vertex classes \NEW{$S_{m+1}$}\OLD{$X$} and \NEW{$Y$}\OLD{$Y\coloneqq N_{D\setminus A}^+(w_{m+1})\subseteq V'$} by joining \NEW{$e\in S_{m+1}$}\OLD{$v \in X$} and~$u\in Y$ if and only if \NEW{$u\notin V(e)$}\OLD{$u \neq v$}.
		\NEW{Note that 
		\[|S_{m+1}|\stackrel{\text{\cref{lm:layouts-buildlayouts-degreeW}}}{=}\hexc^+(w_{m+1})\stackrel{\text{\cref{eq:layouts-endpoints-excessA}}}{=}\exc_{D\setminus A}^+(w_{m+1})\leq |N_{D\setminus A}^+(w_{m+1})|=|Y|\]}%
		\OLD{Observe that, by \cref{eq:layouts-endpoints-excess} and \cref{lm:layouts-buildlayouts-subdivisionL}, 
		\[|X|=\hexc^+(w_{m+1})\leq |N_{D\setminus A}^+(w_{m+1})|=|Y|\]}%
		and\OLD{by \cref{lm:layouts-degree&excessW1,lm:layouts-degree&excessW2},}
		\[|Y|\NEW{=d_{D\setminus A}^{\max}(w_{m+1})\geq \frac{d_{D\setminus A}(w_{m+1})}{2}}\stackrel{\text{\NEW{\cref{lm:layouts-degree&excessW1},\cref{lm:layouts-degree&excessW2}}}}{\geq} 5\varepsilon n.\]
		\NEW{Since $Y$ is a set (rather than a multiset), each $e\in S_{m+1}$ satisfies $d_G(e)\geq |Y|-1$.}\OLD{Clearly, for each~$v\in X$, $d_G(v)\geq |Y|-1$.}
		Note that \cref{lm:layouts-buildlayouts-subdivisionL} implies that if~$v \in V'$ is contained in a path~$P$ in~$\tL^m$, then~$v$ is an endpoint of~$P$ or has (in~$P$) a neighbour in~$W$.
		Hence,\OLD{together with~\cref{eq:layouts-endpoints-excess} and \cref{lm:layouts-degree&excessV'}, we have, for} each $v \in Y \subseteq V'$ \NEW{satisfies
		\begin{equation*}
			d_G(v)\stackrel{\text{\cref{eq:layouts-endpoints-excess},\cref{lm:layouts-buildlayouts-subdivisionL}}}{\geq} |S_{m+1}|-\hexc^-(v)- |W|\stackrel{\text{\cref{eq:layouts-hexc}}}{\geq} |S_{m+1}|-\texc_{D,U^*,W,A}^-(v)-|W|
			\stackrel{\text{\cref{lm:layouts-W},\cref{lm:layouts-degree&excessV'}}}{\geq} |S_{m+1}|-2\varepsilon n.
		\end{equation*}}% 
		\OLD{\begin{equation*}
		d_G(v)\geq |X|-\hexc^-(v)- |W|
		\geq |X|-2\varepsilon n.
		\end{equation*}}%
		Thus, applying \cref{cor:Hall} with \NEW{$S_{m+1}$}\OLD{$X$} and~$Y$ playing the roles of~$A$ and~$B$ gives a matching~$M$ of~$G$ covering \NEW{$S_{m+1}$}\OLD{$X$}. 
		
		Let~$\tL^{m+1}$ be obtained from~$\tL^m$ by subdividing, for each~$vu\in M$ (with~$v \in X$ and~$u \in Y$), the edge $w_{m+1}v\in E(\tL^m)$ into the path~$w_{m+1}uv$.
		Note that this is a valid subdivision since \cref{lm:layouts-buildlayouts-subdivisionL} implies that the path~$P \in\tL^m$ containing~$w_{m+1}v$ satisfies $V' \cap V(P) \subseteq \{v\}$. 
		Let $\tF^{m+1}\coloneqq \tF^m \cup E_{w_{m+1}}(\tL^{m+1})$. 
		Clearly, \cref{lm:layouts-buildlayouts-subdivisionL,lm:layouts-buildlayouts-subdivisionF} are satisfied with~$m+1$ playing the role of~$m$, as desired. This completes \cref{step:bulidlayouts}.
		
		\item \textbf{Covering additional edges incident to~$W$.}\label{step:coverW}
		\NEW{Now that we have constructed suitable paths, we need to partition them into $\ell$ small layouts. Moreover, recall that we need these layouts to cover all the non-absorbing edges incident to $W_1$ (see \cref{lm:layouts-degreeW1}), as well as a prescribed number of edges incident to $W_2$ (see \cref{lm:layouts-degreeW2}). These are the goals of \cref{step:coverW}. First, we will ensure that \cref{lm:layouts-degreeW1,lm:layouts-degreeW2} are satisfied as follows. For each $w\in W$ in turn, we will subdivide some of the paths constructed in \cref{step:bulidlayouts} to incorporate $w$ as an internal vertex. Since the layout needs to remain $W$-exceptional, we will once again need to prescribe the new edges incident to $W$. More precisely, suppose that we want to incorporate $w\in W$ as an internal vertex in a path $P\in \tL$. First, we will choose an unfixed edge $uv\in E(P)$ (which exists by \cref{lm:layouts-buildlayouts-pathlength}) and edges $u'w, wv'\in D\setminus A$ which are not already covered by $\tL$ and such that $u',v'\notin V(P)$. Then, we will subdivide the edge $uv$ in $P$ into the path $uu'wv'v$ and consider the edges $u'w$ and $wv'$ as fixed edges. We will repeat this procedure until $\tL$ covers the desired amount of edges incident to $W$ (as prescribed by \cref{lm:layouts-degreeW1,lm:layouts-degreeW2}). Once this achieved, we will split $(\tL,\tF)$ into $\ell$ layouts.}
		\OLD{We now proceed similarly as above to ensure \cref{lm:layouts-degreeW1,lm:layouts-degreeW2} are satisfied.}
		
		More precisely, we construct $(\hL_1,\hF_1), \dots, (\hL_\ell,\hF_\ell)$ such that the following hold, where \NEW{$\tL$ is the multiset defined by} $\hL\coloneqq \bigcup_{i\in [\ell]}\hL_i$ and $\hF\coloneqq \bigcup_{i\in [\ell]}\hF_i$.
		\begin{enumerate}[label=\upshape(\greek*$'$)]
			\item $(\hL_1,\hF_1), \dots, (\hL_\ell,\hF_\ell)$ are $W$-exceptional layouts such that~$\hF\subseteq \cF$ and~$\hL$ contains no isolated vertex.\label{lm:layouts-coverW-layouts}
			\item $\hL$ is a subdivision of~$\tL$ (and thus, by \cref{lm:layouts-buildlayouts-shape},~$\hL$ has shape~$\tE$ and $(\hL_1,\hF_1), \dots, (\hL_\ell,\hF_\ell)$ are~$U^*$-path consistent with respect to~$(D', \cF)$).\label{lm:layouts-coverW-shape}
			\item For each~$v\in V'$, $d_{\hL}(v)\leq |\exc_D(v)|-\phi(v)+2+2|W|$.\label{lm:layouts-coverW-V'}
			\item For each~$i\in[\ell]$, $|V(\hL_i)|\leq 5\sqrt{\varepsilon}n$ and $|E(\hL_i)|\leq 4\sqrt{\varepsilon}n$. 
			Moreover, each path~$P \in \hL$ contains an edge which lies entirely in~$V'$.\label{lm:layouts-coverW-size}
			\item Either $\ell=\texc(D)-\lceil\eta n\rceil$ or there exist at least~$\sqrt{\varepsilon}n$ indices~$i\in [\ell]$ such that~$\hL_i$ contains at least~$2$ paths.
			Moreover, if $\ell=\texc(D)-\lceil\eta n\rceil$, then for all~$i \in [\ell]$, $|\hL_i|=1$.
			\label{lm:layouts-coverW-size2}
			\item For each~$v\in W_1$, $d_{\hL}^\pm(v) = d_{D\setminus A}^\pm (v)$.\label{lm:layouts-coverW-degreeW1}
			\item For each~$v\in W_2$, $d_{\hL}^\pm(v)= d_D^\pm(v)-\lceil\eta n\rceil$.\label{lm:layouts-coverW-degreeW2}
		\end{enumerate}
		\NEW{Properties \cref{lm:layouts-coverW-layouts,lm:layouts-coverW-shape,lm:layouts-coverW-degreeW1,lm:layouts-coverW-degreeW2,lm:layouts-coverW-size} mean that our main objectives for \cref{step:coverW} are achieved: $(\hL_1,\hF_1), \dots, (\hL_\ell,\hF_\ell)$ are $U^*$-consistent $W$-exceptional layouts which are small and cover the desired amount of edges incident to $W$ (and so satisfy \cref{lm:layouts-size,lm:layouts-degreeW1,lm:layouts-degreeW2}). Moreover, \cref{lm:layouts-coverW-V'} will ensure that each vertex of $V'$ is covered by only few of the layouts (as desired for \cref{lm:layouts-V'}). Finally, \cref{lm:layouts-coverW-size2} is a technical property which will enable us, in \cref{step:degreeV'}, to adjust the layouts in order for \cref{lm:layouts-degreeV'} to be satisfied.}
		
		Let $\hL^0\coloneqq \tL$ and $\hF^0\coloneqq \tF$.
		Let $w_1, \dots, w_k$ be an enumeration of~$W$. \NEW{We will consider each $w_i$ in turn and, at each stage $i$, subdivide the required the number of paths with $w_i$.}
		Let~$\hQ^0$ be a set of paths in~$\hL^0$ of size $\min \{2 \ell, |\hL^0|\}$.
		\NEW{We will restrict ourselves to only subdivide the paths in $\hQ^0$. This will ensure that only few of the final paths are long, which will enable us to form small layouts.}
		Assume inductively that, for some $0\leq m\leq k$, we have constructed, for each~$i\in [m]$, two multisets of paths~$\hL^i$ and~$\hQ^i$, and a set of edges~$\hF^i$ such that the following hold.
		\begin{enumerate}[label=\upshape(\Roman*$'$)]
			\item Let $i\in [m]$. Then, for each~$P\in \hQ^i$, either~$P \in \hQ^{i-1}$ or there exist~$P'\in \hQ^{i-1}$, an edge $e=u_ev_e\in E(P') \setminus \hF^{i-1}$ with $u_e,v_e\in V'$, and distinct~$u_e',v_e'\in V'\setminus V(P')$ such that~$P$ is obtained from~$P'$ by subdividing the edge $e= u_ev_e$ into the path $u_eu_e'w_iv_e'v_e$, where $u_e'w_i, w_iv_e'\in E(D\setminus A)\setminus \tF$ and $\{u_e',v_e'\}\cap \{u_{e'}',v_{e'}'\}=\emptyset$ whenever $e,e'\in E(\hQ^{i-1})$ are distinct edges to be subdivided in order to form~$\hQ^i$. Moreover, $\hL^i = (\hL^{i-1}\setminus \hQ^{i-1} )\cup \hQ^i$\label{lm:layouts-coverW-subdivisionL}. \NEW{(I.e.\ $\hL^i$ is obtained from $\hL^{i-1}$ by incorporating $w$ as an internal vertex in some of the paths in $\hQ^{i-1}\subseteq \hL^i$.)}%
				\COMMENT{I.e.\ the paths in $\hL^0\setminus \hQ^0$ are left untouched in this step.}
			\item For each~$i \in [m]$, $\hF^i= \hF^{i-1} \cup E_{\{w_i\}}(\hL^i)$.\label{lm:layouts-coverW-subdivisionF}
			\item Let~$i\in [m]$. If~$w_i\in W_1$, then $N_{\hL^m}^\pm(w_i)= N_{D\setminus A}^\pm(w_i)$ and, if~$w_i\in W_2$, then $N_{\hL^m}^\pm(w_i)\subseteq N_D^\pm(w_i)$ and $d_{\hL^m}^\pm(w_i)= d_D^\pm(w_i)-\lceil\eta n\rceil$.\label{lm:layouts-coverW-cover}
		\end{enumerate}
		Note that \cref{lm:layouts-buildlayouts-pathlength} and \cref{lm:layouts-coverW-subdivisionL} imply the following.
		\begin{enumerate}[resume, label=\upshape(\Roman*$'$)]
			\item For each~$i \in [m]$,  each~$P \in \hL^i$ contains an edge which lies entirely in~$V'$.
			\label{lm:layouts-coverW-edgeinV'}
		\end{enumerate}
		Also note that, by \cref{lm:layouts-coverW-subdivisionL,lm:layouts-coverW-subdivisionF}, for each~$i\in [m]$,~$\hF^i$ is a set of edges (rather than a multiset) and is obtained from~$\hF^{i-1}$ by adding all the edges of the form~$u_e'w_i$ and~$w_iv_e'$ in \cref{lm:layouts-coverW-subdivisionL}. In particular, \begin{equation}\label{eq:hF}
		\hF^{i-1}\subseteq \hF^i =E_W(\hL^i)\subseteq E_W(D)\setminus A=\cF.
		\end{equation} 
		
		\NEW{First, suppose that $m=k$. Then, \cref{lm:layouts-coverW-cover} implies that we have finished incorporating all the desired edges incident to $W$. By \cref{lm:layouts-coverW-subdivisionF,lm:layouts-coverW-subdivisionL}, $(\hL^k, \hF^k)$ is obtained by subdividing some of the paths in $\tL$ and adding all the new edges incident to $W$ to $\hF^k$. Thus, $(\hL^k, \hF^k)$ is still a $U^*$-path consistent $W$-exceptional layout and there only remains to partition $(\hL^k, \hF^k)$ into $\ell$ small layouts. We will do so by splitting the paths in $\hL^k$ as evenly as possible across the $\ell$ layouts and, subject that, also distribute the paths in $\hQ^k$ as evenly as possible. This will ensure that not all of the long paths belong to the same layout (recall that we want our layouts to be small).}%
		\OLD{If~$m=k$, then note that, by \cref{eq:layouts-endpoints-size}, we have $|\hL^k| = |\tL| = |\tE| = \texc(D)-\lceil\eta n\rceil$.}
		
		\NEW{More precisely,} we partition~$\hL^k$ into $\hL_1, \dots,\hL_\ell$ such that, for each $i,j\in [\ell]$,
		\begin{equation}\label{eq:layouts-hQ}
			||\hL_i|-|\hL_j||\leq 1\quad \text{and} \quad |\hL_i\cap \hQ^k|\leq 2
		\end{equation}
		\NEW{(this is possible since $|\hQ^k|=\min\{2\ell, |\hL^0|\}=\min\{2\ell, |\hL^k|\}$)}.
		Note that \cref{eq:layouts-endpoints-size}, \NEW{\cref{lm:layouts-buildlayouts-shape}, and \cref{lm:layouts-coverW-subdivisionL}} imply that 
		\begin{equation}\label{eq:layouts-hL}
			|\hL^k| = |\tL| = |\tE| \stackrel{\text{\cref{eq:layouts-endpoints-size}}}{=} \texc(D)-\lceil\eta n\rceil
		\end{equation}
		and so, for each~$i \in [\ell]$, 
		\begin{align}
			|\hL_i|&\stackrel{\text{\eqmakebox[hLi][c]{ }}}{\leq}\left\lceil\frac{\texc(D)-\lceil\eta n\rceil}{\ell}\right\rceil
			\stackrel{\text{\cref{lm:layouts-d},\cref{lm:layouts-ell}}}{\leq} \frac{\max\{\Delta^0(D),\exc(D)\}-\lceil\eta n\rceil}{\eta n}+1\nonumber\\
			&\stackrel{\text{\eqmakebox[hLi][c]{\text{\cref{lm:layouts-degree&excessV'}}}}}{\leq} \frac{\max\{n,n|W|+ \varepsilon n|V'|\}-\lceil\eta n\rceil}{\eta n}+1
			\stackrel{|W|\leq \varepsilon n}{\leq} \frac{\max\{n,2 \varepsilon n^2\}-\lceil\eta n\rceil}{\eta n}+1\nonumber\\
			&\stackrel{\text{\eqmakebox[hLi][c]{ }}}{\leq} \sqrt{\varepsilon} n.\label{eq:layouts-sizehLi}
		\end{align}
		For each $i\in [\ell]$, define $\hF_i\coloneqq E(\hL_i)\cap \cF=E(\hL_i)\cap \hF^k$. \NEW{Let $\hL$ be the multiset defined by $\hL\coloneqq \bigcup_{i\in[\ell]}\hL_i=\hL^k$ and denote $\hF\coloneqq \bigcup_{i\in[\ell]}\hF_i=\hF^k$ (as mentioned after \cref{lm:layouts-coverW-edgeinV'}, $\hF$ is subset of $\tF$ rather than a multiset).}
		
		\begin{claim}
			\NEW{Properties \cref{lm:layouts-coverW-V',lm:layouts-coverW-degreeW1,lm:layouts-coverW-degreeW2,lm:layouts-coverW-layouts,lm:layouts-coverW-shape,lm:layouts-coverW-size,lm:layouts-coverW-size2} are satisfied.}
		\end{claim}

		\begin{proofclaim}
			\NEW{First, observe that \cref{lm:layouts-coverW-degreeW1,lm:layouts-coverW-degreeW2} follow immediately from \cref{lm:layouts-coverW-cover}. For \cref{lm:layouts-coverW-size2}, suppose that $\ell\neq \texc(D)-\lceil\eta n\rceil$. By \cref{eq:layouts-hL} and \cref{eq:layouts-ell}, $|\hL|=\texc(D)-\lceil\eta n\rceil\geq \ell +\sqrt{\varepsilon}n$. Thus, \cref{eq:layouts-hQ} implies that there exist at least~$\sqrt{\varepsilon}n$ indices~$i\in[\ell]$ such that~$|\hL_i|\geq 2$. Thus, \cref{lm:layouts-coverW-size2} holds.} 
				
			\NEW{By \cref{lm:layouts-coverW-subdivisionL}, $\hL$ is obtained by subdividing some of the paths in $\tL$ and so \cref{lm:layouts-coverW-shape} holds. To check \cref{lm:layouts-coverW-layouts}, we need to verify \cref{def:layout-F,def:layout-L,def:layout-unfixededge} as defined in \cref{sec:approxdecomp}. First, \cref{lm:layouts-buildlayouts-layouts} implies that $\hL_1, \dots, \hL_k$ are multisets of non-trivial paths on $V$ and so they satisfy \cref{def:layout-L}. By \cref{lm:layouts-buildlayouts-layouts} and \cref{lm:layouts-coverW-subdivisionF}, $\hF_i$ consists, for each $i\in [k]$, of all the edges in $\hL_i$ which are incident to $W$. Therefore, $(\hL_1, \hF_1),\dots, (\hL_\ell, \hF_\ell)$ satisfy \cref{def:layout-F}, and \cref{lm:layouts-coverW-edgeinV'} implies that \cref{def:layout-unfixededge} holds. Therefore, $(\hL_1, \hF_1),\dots, (\hL_\ell, \hF_\ell)$ are $W$-exceptional layouts. As observed above, $\hF\subseteq \cF$ and $\hL$ only consists of non-trivial paths. Thus, \cref{lm:layouts-coverW-layouts} holds.}
		
			\NEW{Next, we verify \cref{lm:layouts-coverW-V'}.
			Let $v\in V'$. 
			By \cref{lm:layouts-buildlayouts-shape} and \cref{lm:layouts-coverW-shape}, $\hL$ is a subdivision of $\tE$ and so \cref{eq:layouts-endpoints-excess} implies that $v$ is the starting point of $\hexc^+(v)$ paths in $\hL$ and the ending point of $\hexc^-(v)$ paths in $\hL$. 
			Suppose that $v\in V^0(P)$ for some $P\in \hL$. Then, \cref{lm:layouts-buildlayouts-V'} and \cref{lm:layouts-coverW-subdivisionL} imply that $P$ contains an edge $e\in E(D)$ between $v$ and a vertex of $W$. By \cref{lm:layouts-coverW-layouts}, $(\hL_1, \hF_1), \dots, (\hL_\ell,\hF_\ell)$ are $W$-exceptional and so $e\in \hF$. Moreover, since $\hF$ is a set rather than a multiset, we have $e\notin E(P')$ for all $P'\in \hL\setminus \{P\}$. Therefore,
			\begin{align*}
				d_{\hL}(v)&\stackrel{\text{\eqmakebox[coverW][c]{}}}{\leq} \hexc^+(v)+\hexc^-(v)+2|N_D(v)\cap W|\\
				&\stackrel{\text{\eqmakebox[coverW][c]{}}}{\leq}(\texc_{D,U^*}^+(v)-\phi^+(v))+(\texc_{D,U^*}^-(v)-\phi^-(v))+2|W|\\
				&\stackrel{\text{\eqmakebox[coverW][c]{\cref{eq:texc}}}}{\leq} (\exc_D^+(v)+1)+(\exc_D^-(v)+1)-\phi(v)+2|W|
				\leq |\exc_D(v)|-\phi(v)+2+2|W|.
			\end{align*}
			Thus, \cref{lm:layouts-coverW-V'} holds.}
			
			\NEW{Finally, we verify \cref{lm:layouts-coverW-size}. The ``moreover part" holds by~\cref{lm:layouts-buildlayouts-pathlength,lm:layouts-coverW-subdivisionL}.
			Moreover, \cref{lm:layouts-coverW-subdivisionL} implies that whenever we incorporate a vertex of $W$ as an internal vertex to a path in $\hQ^0$, we add three vertices to that path. Therefore, \cref{lm:layouts-buildlayouts-pathlength} implies that each $Q \in \hQ^k$ satisfies $|V(Q)|\leq 4+3|W|$. Thus, \cref{lm:layouts-W}, \cref{eq:layouts-hQ}, and \cref{eq:layouts-sizehLi} imply that, for each $i\in [\ell]$, we have $|V(\hL_i)|\leq 4|\hL_i|+6|W|\leq 4\sqrt{\varepsilon}n+6\varepsilon n\leq 5\sqrt{\varepsilon}n$. Similarly, each $Q \in \hQ^k$ satisfies $|E(Q)|\leq 3+3|W|$ and so each $i\in [\ell]$ satisfies $|E(\hL_i)|\leq 3|\hL_i|+6|W|\leq 3\sqrt{\varepsilon}n+6\varepsilon n\leq 4\sqrt{\varepsilon}n$. Thus, \cref{lm:layouts-coverW-size} holds.}
		\end{proofclaim}
		
		\smallskip
				
		\NEW{We may therefore assume that $m<k$. Suppose}\OLD{If $m<k$, then assume} without loss of generality that $w_{m+1}\notin W^-$. Thus, $w_{m+1}\in U^+(D)\cup U^0(D)$ and so $\exc_D(w_{m+1})=|\exc_D(w_{m+1})|=\exc_D^+(w_{m+1})$ and $d_D^{\min}(w_{m+1})=d_D^-(w_{m+1})$. Moreover, by \cref{def:A}, $d_{D\setminus A}^-(w_{m+1})=d_D^-(w_{m+1})$.
		Finally, by assumptions \cref{lm:layouts-U*,lm:layouts-A}, $\phi^\pm(w_{m+1})=0$ and $U^*\cap W=\emptyset$. Thus, \begin{equation}\label{eq:hexcwm+1}
			\hexc^+(w_{m+1})=\texc_{D, U^*}^+(w_{m+1})=\exc_D^+(w_{m+1})-d_A^+(w_{m+1})\stackrel{\text{\NEW{\cref{def:A}}}}{=}\exc_{D\setminus A}^+(w_{m+1}).
		\end{equation}
		
		\NEW{We will construct $\hQ^{m+1}$ as follows. First, we will pair each inneighbour of $w_{m+1}$ in $D\setminus (A\cup \hF^m)$ to an outneighbour of $w_{m+1}$ in $D\setminus (A\cup \hF^m)$. Let $X$ denote the set of these pairs. We will use these pairs to incorporate $w_{m+1}$ as an internal vertex in some of the paths in $\hQ^m$ as follows. 
		Let $Y$ be the set of paths in $\cQ^m$ which do not already contain $w_{m+1}$. Form an auxiliary bipartite graph by joining each $(u',v')\in X$ and $P\in Y$ if and only if both $u',v'\notin V(P)$ (if $u'\in V(P)$ or $v'\in V(P)$, then we cannot use $(u',v')$ to incorporate $w_{m+1}$ as an internal vertex in $P$). Then, we will use Hall's theorem to find a large matching $M$ in this auxiliary graph. For each $(u',v')P\in M$, we will subdivide an unfixed edge $uv\in E(P)$ into the path $uu'w_{m+1}v'v$.}
		
		\OLD{First, note that}By \cref{eq:layouts-endpoints-excess}, \ref{lm:layouts-buildlayouts-shape}, \cref{lm:layouts-buildlayouts-internalvertices}, and \cref{lm:layouts-coverW-subdivisionL}, 
		\NEW{none of the paths in $\hL^m$ have $w_{m+1}$ as internal vertex or ending point and so}
		we have $d_{\hL^m}^-(w_{m+1})=0$.
		Fix a bijection $\sigma: N_{D\setminus A}^-(w_{m+1}) \longrightarrow N_{D\setminus (A \cup \hF^m)}^+(w_{m+1})$. Note that this is possible since
		\begin{align*}
			d_{\hF^m}^+(w_{m+1})&\stackrel{\text{\eqmakebox[eq1][c]{\text{\cref{eq:hF}}}}}{=}d_{\hL^m}^+(w_{m+1})\stackrel{\text{\cref{lm:layouts-coverW-subdivisionL}}}{=}d_{\tL}^+(w_{m+1})\stackrel{\text{\cref{lm:layouts-buildlayouts-shape},\cref{lm:layouts-buildlayouts-internalvertices}}}{=}d_{\tE}^+(w_{m+1})\\
			&\stackrel{\text{\eqmakebox[eq1][c]{\text{\cref{eq:layouts-endpoints-excess}}}}}{=}\hexc^+(w_{m+1})\stackrel{\text{\cref{eq:hexcwm+1}}}{=}\exc_{D\setminus A}^+(w_{m+1}).
		\end{align*}
		Thus,
		\begin{align*} 
			d_{D\setminus (A\cup\hF^m)}^+(w_{m+1})&\stackrel{\text{\eqmakebox[eq2][c]{\text{\cref{eq:hF}}}}}{=}d_{D\setminus A}^+(w_{m+1})-d_{\hF^m}^+(w_{m+1})=d_{D\setminus A}^+(w_{m+1})- \exc_{D\setminus A}^+(w_{m+1})\\
			&\stackrel{\text{\eqmakebox[eq2][c]{}}}{=}d_{D\setminus A}^-(w_{m+1}),
		\end{align*}
		as desired.
		
		Let $X\coloneqq \{(u,\sigma(u))\mid u\in N_{D\setminus A}^-(w_{m+1})\}$. 
		Let $Y \subseteq \hQ^m$ be obtained from~$\hQ^m$ by deleting all the paths that contain~$w_{m+1}$.
		Define an auxiliary bipartite graph~$G$ with vertex classes~$X$ and~$Y$ by joining~$(u,v)\in X$ and~$P\in Y$ if and only if both~$u,v\notin V(P)$.
		
		\begin{claim}\label{claim:coverW}
			If~$w_{m+1}\in W_1$, then~$G$ contains a matching~$M$ covering~$X$. If~$w_{m+1}\in W_2$, then~$G$ contains a matching~$M$ of size $|X|-\lceil\eta n\rceil$.
		\end{claim}
	
		Let~$M$ be as in \cref{claim:coverW}. We obtain~$\hQ^{m+1}$ from~$\hQ^{m}$ by subdividing, for each~$(u',v')P \in M$, an edge~$uv\in P$ that lies entirely in~$V'$ (which exists by~\cref{lm:layouts-coverW-edgeinV'}) into the path $uu'w_{m+1}v'v$.
		Let $\hL^{m+1} \coloneqq (\hL^m\setminus \hQ^m )\cup \hQ^{m+1}$ and $\hF^{m+1} \coloneqq \hF^m\cup E_{\{w_{m+1}\}}(\hL^{m+1})$.
		One can easily verify that \cref{lm:layouts-coverW-subdivisionL,lm:layouts-coverW-subdivisionF,lm:layouts-coverW-cover} are satisfied with~$m+1$ playing the role of~$m$. There only remains to show \cref{claim:coverW}.
		\NEW{We will need the following observation.}
		
		\begin{claim}\label{claim:XY}
			\NEW{If $X\neq \emptyset$, then
			\begin{align*}
				\max\{|X|,|Y|\}\geq \eta n\geq 10\varepsilon n \quad \text{and} \quad \min\{|X|,|Y|\}\geq
				\begin{cases}
					|X| & \text{if }w_{m+1}\in W_1,\\
					|X|-\lceil\eta n\rceil & \text{if }w_{m+1}\in W_2.
				\end{cases}
			\end{align*}}
		\end{claim}
	
		\NEW{We first assume that \cref{claim:XY} holds and derive \cref{claim:coverW}.}
		
		\begin{proofclaim}[Proof of \cref{claim:coverW}]
			Clearly, we may assume that~$X \ne \emptyset$.
			\NEW{The goal is to use \cref{cor:Hall}. We start by checking that the degree of each vertex in $G$ is large.
			First, observe that, by \cref{lm:layouts-W}, we have~$u,v\in V'$ for each~$(u,v)\in X$.
			By \cref{lm:layouts-coverW-subdivisionL}, $\hL^m$ is obtained from $\tL$ by repeated subdivisions and, in each subdivision, we incorporate a vertex of $W$ using only two new vertices of $V'$. Thus, each $P\in Y$ satisfies}
			\OLD{By~\cref{lm:layouts-buildlayouts-pathlength} and \cref{lm:layouts-coverW-subdivisionL},
			we have, for each~$P\in Y$,} 
			\begin{equation*}
				d_G(P)\geq |X|-|V(P)\cap V'|\stackrel{\text{\cref{lm:layouts-buildlayouts-pathlength},\NEW{\cref{lm:layouts-coverW-subdivisionL}}}}{\geq} |X|- (4+2|W|) \stackrel{\text{\cref{lm:layouts-W}}}{\geq} |X|-3\varepsilon n\NEW{\stackrel{\text{\cref{claim:XY}}}{\geq}|X|-\frac{\max\{|X|,|Y|\}}{2}.}
			\end{equation*}	
			\NEW{Let $(u,v)\in X$. We count the number of paths in $\hL^m$ which contain $u$. By \cref{lm:layouts-coverW-subdivisionL}, $\hL^m$ is a subdivision of $\tL$ and so \cref{lm:layouts-buildlayouts-shape} and \cref{eq:layouts-endpoints-excess} imply that $u$ is an endpoint of precisely $\hexc^+(u)+\hexc^-(u)$ paths in $\hL^m$. Suppose that $P\in \hL^m$ contains $u$ as an internal vertex. By \cref{lm:layouts-buildlayouts-V'} and \cref{lm:layouts-coverW-subdivisionL}, $P$ contains an edge $e$ between $u$ and a vertex of $W$. By \cref{lm:layouts-buildlayouts-layouts} and \cref{lm:layouts-coverW-subdivisionL}, $e\in \hF^m\subseteq E(D\setminus A)$. In particular, the fact that $\hF^m$ is a set (rather than a multiset) implies that $e\notin E(\hL^m\setminus \{P\})$. Thus, there are at most $|N_{D\setminus A}(u)\cap W|$ paths in $\hL^m$ which contain $u$ as an internal vertex. Similarly, there are at most $\hexc^+(v)+\hexc^-(v)+|N_{D\setminus A}(v)\cap W|$ paths in $\hL^m$ which contain $v$ (as an endpoint or internal vertex). Thus,
			\begin{align*}
				d_G((u,v))&\stackrel{\text{\eqmakebox[Hall][c]{}}}{\geq} |Y|-(\hexc^+(u)+\hexc^-(u)+|N_{D\setminus A}(u)\cap W|)\\
				&\qquad\qquad -(\hexc^+(v)+\hexc^-(v)+|N_{D\setminus A}(v)\cap W|)\\
				&\stackrel{\text{\eqmakebox[Hall][c]{}}}{\geq} |Y|-(\texc_{D,U^*}^+(u)+\texc_{D,U^*}^-(u)-\phi(u))-(\texc_{D,U^*}^+(v)+\texc_{D,U^*}^-(v)-\phi(v))\\
				&\qquad\qquad-2|W|\\
				&\stackrel{\text{\eqmakebox[Hall][c]{\cref{eq:texc}}}}{\geq} |Y|-(|\exc_D(u)|+2-\phi(u))-(|\exc_D(v)|+2-\phi(v))-2|W|\\
				&\stackrel{\text{\eqmakebox[Hall][c]{}}}{\geq} |Y|-|\exc_D(u)|-|\exc_D(v)|-4-2|W|
				\stackrel{\text{\cref{lm:layouts-W},\cref{lm:layouts-degree&excessV'}}}{\geq} |Y|-5\varepsilon n\\
				&\stackrel{\text{\eqmakebox[Hall][c]{\cref{claim:XY}}}}{\geq}|Y|-\frac{\max\{|X|,|Y|\}}{2}.
			\end{align*}}%
			\OLD{Note that \cref{lm:layouts-buildlayouts-V'} and \cref{lm:layouts-coverW-subdivisionL} and imply that if~$v \in V'$ is contained in a path~$P\in \tL^m$, then~$v$ is an endpoint of~$P$ or has (in~$P$) a neighbour in~$W$. Since 
			Moreover, by \cref{lm:layouts-W}, we have~$u,v\in V'$ for each~$(u,v)\in X$.
			Hence, together with~\cref{eq:layouts-endpoints-excess} and \cref{lm:layouts-degree&excessV'}, we have, for each~$(u,v)\in X$,
			\begin{align*}
			d_G((u,v))&\geq |Y|-(\hexc^+(u)+\hexc^-(u))-2|N_{D\setminus A}(u)\cap W|\\
			&\qquad\qquad -(\hexc^+(v)+\hexc^-(v))-2|N_{D\setminus A}(v)\cap W|\\
			&\geq |Y|-(\texc_{D,U^*}^+(u)+\texc_{D,U^*}^-(u)-\phi(u))-(\texc_{D,U^*}^+(v)+\texc_{D,U^*}^-(v)-\phi(v))-4|W|\\
			&\stackrel{\text{\cref{eq:texc}}}{\geq} |Y|-(|\exc_D(u)|+2-\phi(u))-(|\exc_D(v)|+2-\phi(v))-4|W|\\
			&\geq |Y|-|\exc_D(u)|-|\exc_D(v)|-4-4|W|
			\stackrel{\text{\cref{lm:layouts-W},\cref{lm:layouts-degree&excessV'}}}{\geq} |Y|-7\varepsilon n.
			\end{align*}}%
			\NEW{Thus, $G$ satisfies the degree conditions of \cref{cor:Hall}, applied with~$\{X,Y\}$ playing the roles of~$\{A,B\}$ (with~$|A| \le |B|$). Therefore, $G$ contains a matching $M$ of size $\min\{|X|,|Y|\}$. By \cref{claim:XY}, we may assume that $|M|=|X|$ if $w_{m+1}\in W_1$ and $|M|=|X|-\lceil\eta n\rceil$ if $w_{m+1}\in W_2$.}%
			\OLD{Suppose that 
				\begin{align}
					|Y| \ge \begin{cases}
						\max\{|X|,\eta n\} & \text{if $w_{m+1} \in W_1$,}\\
						|X| - \lceil \eta n  \rceil & \text{if $w_{m+1} \in W_2$.}
					\end{cases}
					\label{eqn:|Y|'}
				\end{align}
				Then, we are done by \cref{cor:Hall}, applied with~$\{X,Y\}$ playing the roles of~$\{A,B\}$ with~$|A| \le |B|$%
%				\COMMENT{To be precise, if~$w \in W_2$, then let 
%					\begin{equation*}
%						X'\coloneqq
%						\begin{cases}
%							X& \text{if }|X|\leq |Y|,\\
%							Y& \text{otherwise};
%						\end{cases}
%						\quad
%						\text{and}
%						\quad
%						Y'\coloneqq
%						\begin{cases}
%							Y& \text{if }|X|\leq |Y|,\\
%							X& \text{otherwise}.
%						\end{cases}
%					\end{equation*}
%					Thus, by \cref{cor:Hall}, applied with $X'$ and $Y'$ playing the roles of $A$ and $B$, there exists a matching covering $X'$. Since $|X'|\geq |X|-\lceil\eta n\rceil$, the desired matching exists.}.
				To complete the proof of the claim, it suffices to prove~\cref{eqn:|Y|'}.}
		\end{proofclaim}
		
		\NEW{Finally, it remains to prove \cref{claim:XY}.}
		
		\begin{proofclaim}[Proof of \cref{claim:XY}]
			By \cref{lm:layouts-degree&excessW2}, if $w_{m+1}\in W_2$, then $|X|=d_{D\setminus A}^-(w_{m+1})=d_D^{\min}(w_{m+1})\geq \lceil \eta n\rceil$.
			\NEW{Thus, it is enough to show that 
			\begin{align}
				|Y| \ge \begin{cases}
					\max\{|X|,\eta n\} & \text{if $w_{m+1} \in W_1$,}\\
					|X| - \lceil \eta n  \rceil & \text{if $w_{m+1} \in W_2$.}
				\end{cases}
				\label{eqn:|Y|}
			\end{align}}%
			Note that, by \cref{lm:layouts-buildlayouts-internalvertices} and \cref{lm:layouts-coverW-subdivisionL},~$w_{m+1}$ is not an internal vertex of any path in~$\hQ^m$. \NEW{Moreover, \cref{lm:layouts-coverW-subdivisionL} implies that $\hL^m$ is a subdivision of $\tL$ and so \cref{lm:layouts-buildlayouts-shape} and \cref{eq:layouts-endpoints-excess} imply that}\OLD{Thus, by \cref{lm:layouts-buildlayouts-shape} and \cref{eq:layouts-endpoints-excess},}
			\begin{equation}\label{eq:Y}
				|Y|\geq |\hQ^m|-\hexc^+(w_{m+1})\stackrel{\text{\cref{eq:hexcwm+1}}}{=}|\hQ^m|-\exc_{D\setminus A}^+(w_{m+1}).
			\end{equation}		
			\NEW{If}\OLD{Assume first that} $|\hQ^m| \ge 2 d$ and~$w_{m+1} \in W_1$, then
			\begin{align*}
				|Y|& \stackrel{\text{\eqmakebox[Y1][c]{\text{\cref{eq:Y}}}}}{\geq} |\hQ^m|-\exc_{D\setminus A}^+(w_{m+1})
				\geq 2d -\exc_{D\setminus A}^+(w_{m+1})\\
				&\stackrel{\text{\eqmakebox[Y1][c]{\text{\cref{lm:layouts-degree&excessW1}}}}}{\geq} d_{D\setminus A}(w_{m+1})+\eta n-\exc_{D\setminus A}^+(w_{m+1})
				=2d_{D\setminus A}^-(w_{m+1})+\eta n
				\geq |X|+\eta n.
			\end{align*}
			Similarly, if $|\hQ^m| \geq 2 d$ and~$w_{m+1} \in W_2$, then
			\begin{align*}
				|Y|& \stackrel{\text{\eqmakebox[Y2][c]{\text{\cref{eq:Y}}}}}{\geq} |\hQ^m|-\exc_{D\setminus A}^+(w_{m+1})
				\stackrel{w_{m+1}\notin V(A)}{=} |\hQ^m|-\exc_D^+(w_{m+1})\\
				&\stackrel{\text{\eqmakebox[Y2][c]{}}}{\geq} 2d-\exc_D^+(w_{m+1})
				\stackrel{\text{\cref{lm:layouts-degree&excessW2}}}{\geq} d_D(w_{m+1})-2\lceil\eta n\rceil-\exc_D^+(w_{m+1})\\
				&\stackrel{\text{\eqmakebox[Y2][c]{}}}{=}2d_D^-(w_{m+1})-2\lceil\eta n\rceil
				=|X|-\lceil\eta n\rceil+d_D^-(w_{m+1})-\lceil\eta n\rceil
				\stackrel{\text{\cref{lm:layouts-degree&excessW2}}}{\geq} |X|-\lceil\eta n\rceil.
			\end{align*}%
			\OLD{as desired.}%
			\NEW{We may therefore}\OLD{Next,} assume that $|\hQ^m| < 2d$.
			Since, by~\cref{lm:layouts-ell},~$d \le \ell$, we have $|\hL^m|=|\hL^0|=|\hQ^0|=|\hQ^m|$ and so 
			\begin{equation}
				\texc(D)-\lceil\eta n\rceil  \stackrel{\text{\cref{eq:layouts-endpoints-size}}}{=} |\tE|\stackrel{\text{\cref{lm:layouts-buildlayouts-shape},\cref{lm:layouts-coverW-subdivisionL}}}{=}|\hL^m| = |\hQ^m| < 2d. \label{eqn:W1}
			\end{equation}
			Thus,
			\begin{align}
				|Y|& \stackrel{\text{\eqmakebox[Y3][c]{\text{\cref{eq:Y}}}}}{\geq}|\hQ^m|-\exc_{D\setminus A}^+(w_{m+1})\nonumber\\
				&\stackrel{\text{\eqmakebox[Y3][c]{\text{\cref{eqn:W1}}}}}{=}\texc(D)-\lceil\eta n\rceil-\exc_{D\setminus A}^+(w_{m+1}) \label{eqn:Y1} \\
				&\stackrel{\text{\eqmakebox[Y3][c]{}}}{\geq} d_D^+(w_{m+1}) -\lceil\eta n\rceil-\exc_D^+(w_{m+1})
				= d_D^-(w_{m+1})-\lceil\eta n\rceil
				\stackrel{w_{m+1}\notin V(A^-)}{=}|X|-\lceil\eta n\rceil.\nonumber
			\end{align}
			We may therefore assume that~$w_{m+1}\in W_1$ and~$|\hQ^m|<2d$. We need to show that $|Y|\geq \max\{|X|, \eta n\}$.
			Recall that $d_D^-(w_{m+1})=d_{D\setminus A}^-(w_{m+1}) = |X| >0$. Then, \NEW{\cref{fact:exc}\cref{fact:exc-dmax}} implies that $d_D^+(w_{m+1})>\exc_D^+(w_{m+1})$.
			Thus, by~\cref{lm:layouts-exc<2d} and~\cref{eqn:W1}, we have $|X|\geq \eta n$ and one of the following holds: $\exc_{D\setminus A}^+(w_{m+1})\leq \lceil\eta n\rceil$ or $d_A^+(w_{m+1})=d_A(w_{m+1})=\lceil\eta n\rceil$%
				\COMMENT{This is just restating the last two bullet points of \cref{lm:layouts-exc<2d}.}.
			Thus, it suffices to show that $|Y| \ge |X|$. 
			If $\exc_{D\setminus A}^+(w_{m+1})\leq \lceil\eta n\rceil$, then, by \cref{eqn:Y1} and \cref{lm:layouts-exc>d},
			\begin{align*}
				|Y|& 
				\stackrel{\text{\eqmakebox[Y4][c]{}}}{\geq}
				d- \exc_{D\setminus A}^+(w_{m+1})
				\stackrel{\text{\cref{lm:layouts-degree&excessW1}}}{\geq}
				\frac{d_{D\setminus A}(w_{m+1})+\lceil\eta n\rceil}{2}-\exc_{D\setminus A}^+(w_{m+1})\\
				&\stackrel{\text{\eqmakebox[Y4][c]{\text{\cref{fact:exc}\cref{fact:exc-dmin}}}}}{=}
				d_{D\setminus A}^-(w_{m+1})+\frac{\lceil\eta n\rceil-\exc_{D\setminus A}^+(w_{m+1})}{2}
				\geq |X|,
			\end{align*} 
			as desired.
			If $d_A^+(w_{m+1})=\lceil\eta n\rceil$, then \cref{eqn:Y1} implies that 
			\begin{align*}		
				|Y|\geq d_D^+(w_{m+1})-\lceil\eta n\rceil-\exc_{D\setminus A}^+(w_{m+1})=d_{D\setminus A}^+(w_{m+1})-\exc_{D\setminus A}^+(w_{m+1})=d_{D\setminus A}^-(w_{m+1})=|X|,
			\end{align*}
			as desired.
		\end{proofclaim}
		
		\item \textbf{Adjusting the degree of the vertices in~$V'$.}\label{step:degreeV'}
		\NEW{Recall that, in \cref{step:coverW}, we constructed $\ell$ $W$-exceptional layouts which are $U^*$-path consistent and satisfy \cref{lm:layouts-coverW-degreeW1,lm:layouts-coverW-degreeW2}, and thus satisfy \cref{lm:layouts-degreeW1,lm:layouts-degreeW2}. We will now adjust these layouts to ensure that \cref{lm:layouts-degreeV'} is satisfied.}
		\OLD{Finally, we add isolated vertices and subdivide paths to ensure \cref{lm:layouts-degreeV'} is satisfied. }
		
		Let $v_1, \dots, v_k$ be an enumeration of~$V'$ and, for each~$i\in[k]$, define 
		\begin{equation}\label{eq:defni}
			n_i\coloneqq d_{\hL}^+(v_i)+|\{j\in [\ell]\mid v_i\notin V(\hL_j)\}|+\lceil\eta n\rceil-\phi^-(v_i)-d_D^+(v_i).
		\end{equation}
		\NEW{Note that together with \cref{eq:nipm} below, \cref{lm:layouts-degreeV'} holds if $n_i=0$ for all $i\in [k]$.}
		
		\begin{claim}\label{eq:nipm}
			For each~$i\in [k]$, 
			\begin{equation*}
			d_D^\pm(v_i)=d_{\hL}^\pm(v_i)+ |\{j\in [\ell]\mid v_i\notin V(\hL_j)\}|-n_i+ \lceil\eta n\rceil -\phi^\mp(v).
			\end{equation*} 
		\end{claim}
		
		\begin{proofclaim}
			Let~$i\in[k]$. The equality for~$+$ holds immediately by definition of~$n_i$. One can easily verify that, in order to show that the equality for~$-$ holds, it is enough to prove that \begin{equation}\label{eq:ni-}
			d_{\hL}^+(v_i)-\phi^-(v_i)-d_D^+(v_i)=d_{\hL}^-(v_i)-\phi^+(v_i)-d_D^-(v_i). 
			\end{equation}
			We now show that \cref{eq:ni-} is satisfied. First, note that, by \cref{eq:layouts-endpoints-excess} and \cref{lm:layouts-coverW-shape}, \NEW{$v_i$ is the starting point of precisely $\hexc^+(v_i)$ paths in $\hL$ and the ending point of precisely $\hexc^-(v_i)$ paths in $\hL$, so}
			\begin{equation}\label{eq:dL}
			d_{\hL}^+(v_i)=(d_{\hL}^-(v_i)-\hexc^-(v_i))+\hexc^+(v_i).
			\end{equation}
			Assume without loss of generality that $d_D^+(v_i)\geq d_D^-(v_i)$. 
			Suppose first that $v_i\notin U^*$. Then, $\hexc^\pm(v_i)=\exc_D^\pm(v_i)-\phi^\pm(v_i)$.		
			Moreover, $\exc_D^-(v_i)=0$ and so $\hexc^-(v_i)=-\phi^-(v_i)$%
			\COMMENT{(in fact, $\phi^-(v_i)=0$)}. 
			Thus, by \cref{eq:dL}, 
			\begin{align*}
			d_{\hL}^+(v_i)&=(d_{\hL}^-(v_i)+\phi^-(v_i))+(\exc_D^+(v_i)-\phi^+(v_i))\\
			&=d_{\hL}^-(v_i)+\phi^-(v_i)+(d_D^+(v_i)-d_D^-(v_i))-\phi^+(v_i),
			\end{align*}
			so \cref{eq:ni-} holds, as desired.
			Now suppose that~$v_i\in U^*$. Then, $\hexc^\pm(v_i)=1-\phi^\pm(v_i)$ and $d_D^+(v_i)=d_D^-(v_i)$. Thus, by \cref{eq:dL},
			\begin{align*}
			d_{\hL}^+(v_i)&=(d_{\hL}^-(v_i)-1+\phi^-(v_i))+(1-\phi^+(v_i))
			=d_{\hL}^-(v_i)+\phi^-(v_i)-\phi^+(v_i)+(d_D^+(v_i)-d_D^-(v_i)),
			\end{align*}
			so \cref{eq:ni-} holds, as desired.
		\end{proofclaim}
		
		If~$n_i>0$,  then, in order to satisfy \cref{lm:layouts-degreeV'}, it is enough to add~$v_i$ as an isolated vertex to exactly~$n_i$ of the sets of paths $\hL_1, \dots, \hL_\ell$ that do not contain~$v_i$ \NEW{(this will decrease $|\{j\in [\ell]\mid v_i\notin V(\hL_j)\}|$ by $n_i$ and so we will be done by \cref{eq:nipm})}. 
		If~$n_i<0$, then it is enough to find~$-n_i$ indices~$j\in [\ell]$ such that~$v_i\notin V(\hL_j)$ and~$|\hL_j|\geq 2$, and add~$v_i$ as an internal vertex in exactly two paths in~$\hL_j$ \NEW{(this will decrease $|\{j\in [\ell]\mid v_i\notin V(\hL_j)\}|$ by $-n_i$ but increase both $d_{\hL}^\pm(v_i)$ by $-2n_i$ and so we will be done by \cref{eq:nipm}).}%
		\OLD{We do so inductively as follows.}
		
		\NEW{We now bound $n_i$ with the following \lcnamecref{eq:ni}.}
		
		\begin{claim}\label{eq:ni}
		For each~$i\in [k]$,
			\begin{equation*}
			-2\varepsilon n \leq n_i\leq 2\sqrt{\varepsilon}n.
			\end{equation*}
		\end{claim}
	
		\begin{proofclaim}
			Let~$i\in [k]$. We have 
			\[2n_i\stackrel{\text{\cref{eq:nipm}}}{=}d_{\hL}(v_i)+2|\{j\in [\ell]\mid v_i\notin V(\hL_j)\}|+2\lceil\eta n\rceil-\phi(v_i)-d_D(v_i).\]
			Thus,
			\begin{align*}
				2n_i&\stackrel{\text{\eqmakebox[nilow][c]{}}}{\geq} d_{\hL}(v_i)+ 2(\ell-d_{\hL}(v_i))+2\lceil\eta n\rceil-\phi(v_i)-d_D(v_i)\\
				&\stackrel{\text{\eqmakebox[nilow][c]{}}}{=}2\ell-d_{\hL}(v_i)+2\lceil\eta n\rceil-\phi(v_i)-d_D(v_i)\\
				&\stackrel{\text{\eqmakebox[nilow][c]{\text{\cref{lm:layouts-coverW-V'}}}}}{\geq} 2\ell-(|\exc_D(v_i)|-\phi(v_i)+2+2|W|)+2\lceil\eta n\rceil-\phi(v_i)-(2d+2\lceil\eta n\rceil)\\
				&\stackrel{\text{\eqmakebox[nilow][c]{}}}{=} 2(\ell-d)-|\exc_D(v_i)|-2-2|W|
				\stackrel{\text{\cref{lm:layouts-W},\cref{lm:layouts-degree&excessV'},\cref{lm:layouts-ell}}}{\geq} -4\varepsilon n.
			\end{align*}
			Similarly,
			\begin{align*}
			2n_i&\stackrel{\text{\eqmakebox[niup][c]{\text{\cref{lm:layouts-coverW-V'}}}}}{\leq} (|\exc_D(v_i)|-\phi(v_i)+2+2|W|)+2\ell+2\lceil\eta n\rceil-\phi(v_i)-(2d+2\lceil\eta n\rceil-\varepsilon n)\\
			&\stackrel{\text{\eqmakebox[nilow][c]{}}}{\leq} 2(\ell-d)+|\exc_D(v_i)|+2+2|W|+\varepsilon n
			\stackrel{\text{\cref{lm:layouts-W},\cref{lm:layouts-degree&excessV'},\cref{lm:layouts-ell}}}{\leq} 4\sqrt{\varepsilon} n,
			\end{align*}
			which proves the claim.
		\end{proofclaim}

		Assume without loss of generality \NEW{that} $(n_i)_{i \in [k]}$ is an increasing sequence and so, for any~$i,j\in [k]$, if~$n_i<0$ but~$n_j\geq 0$, then~$i<j$. 
		For each~$i\in [\ell]$, let $L_i^0\coloneqq \hL_i$. 
		Assume inductively that, for some $0\leq m\leq k$, we have constructed, for each~$i\in [\ell]$ and~$j\in [m]$, a multiset~$L_i^j$ of paths and isolated vertices such that the following are satisfied, where \NEW{$L^j$ is the multiset defined by} $L^j\coloneqq \bigcup_{i\in[\ell]}L_i^j$ for each~$j\in [m]$.
		\begin{enumerate}[label=\upshape(\Roman*$''$)]
			\item For each~$j\in [m]$, if~$n_j<0$, then there exists~$N_j\subseteq [\ell]$ such that~$|N_j|=-n_j$ and the following hold.
			For each~$i\in N_j$,~$v_j\notin V(L_i^{j-1})$ and there exist two paths $P_1, P_2\in L_i^{j-1}$ such that~$L_i^j$ is obtained from~$L_i^{j-1}$ by subdividing, for each~$s\in [2]$, an edge $uw\in E(P_s)\setminus E_W(P_s)$%
				\COMMENT{$uw$ exists by \cref{lm:layouts-coverW-edgeinV'}.}
			into the path~$uv_jw$. For each~$i\in [\ell]\setminus N_j$,~$L_i^j= L_i^{j-1}$.\label{lm:layouts-adjustdegrees-internalvertex}
			\item For each~$j\in [m]$, if~$n_j\geq 0$, then there exists~$N_j\subseteq [\ell]$ such that~$|N_j|=n_j$ and the following hold. For each~$i\in N_j$,~$v_j\notin V(L_i^{j-1})$ and~$L_i^j$ is obtained from~$L_i^{j-1}$ by adding~$v_j$ as an isolated vertex. For each~$i\in [\ell]\setminus N_j$,~$L_i^j=L_i^{j-1}$.\label{lm:layouts-adjustdegrees-isolatedvertex}
			\item For each~$i\in [\ell]$ and~$j\in [m]$, $|V(L_i^j)\setminus V(\hL_i)|\leq \varepsilon^{\frac{1}{3}} n$.\label{lm:layouts-adjustdegrees-sizeV}
		\end{enumerate}
		\NEW{By \cref{lm:layouts-adjustdegrees-internalvertex,lm:layouts-adjustdegrees-isolatedvertex}, each $L^j$ is obtained from $L^{j-1}$ either by subdividing $2$ edges with $v_j$ in $|n_j|$ layouts which did not already cover $v_j$, or by adding $v_j$ as an isolated vertex in $|n_j|$ layouts which did not already cover $v_j$. Thus, the following holds.}\OLD{Note that \cref{lm:layouts-adjustdegrees-internalvertex,lm:layouts-adjustdegrees-isolatedvertex} imply that the following hold.}
		\begin{enumerate}[resume,label=\upshape(\Roman*$''$)]
			\item For each~$i\in [\ell]$ and~$j\in [m]$, $|E(L_i^j)\setminus E(\hL_i)| \le 2|V(L_i^j)\setminus V(\hL_i)|$.\label{lm:layouts-adjustdegrees-sizeE}
			\item For each~$j \in [m]$, $\sum_{i \in [\ell]} |V(L_i^j)| = \sum_{i \in [\ell]} |V(\hL_i)|+ \sum_{j' \in [j]} |n_{j'}|$.\label{lm:layouts-adjustdegrees-sizeV2}
		\end{enumerate}
		
		\NEW{First, assume that $m=k$. For each $i\in [\ell]$, let $L_i\coloneqq L_i^k$ and $F_i\coloneqq \hF_i$. Denote by $L$ the multiset $L\coloneqq \bigcup_{i\in[\ell]}L_i$ and let $F\coloneqq \bigcup_{i\in[\ell]} F_i=\hF$.}
		\NEW{Recall that $D'=D\setminus \cF$.}
	
		\begin{claim}
			\NEW{$(L_1, F_1), \dots, (L_\ell,F_\ell)$ are~$W$-exceptional~$U^*$-path consistent layouts with respect to $(D',\cF)$.
			Moreover, \cref{lm:layouts-ell,lm:layouts-good,lm:layouts-degreeV',lm:layouts-degreeW1,lm:layouts-degreeW2,lm:layouts-size,lm:layouts-V'} hold.}
		\end{claim}
	
		\begin{proofclaim}
			\NEW{By \cref{lm:layouts-adjustdegrees-internalvertex} and \cref{lm:layouts-adjustdegrees-isolatedvertex}, $(L_1, F_1), \dots, (L_\ell,F_\ell)$ are obtained from $(\hL_1, \hF_1), \dots, (\hL_\ell,\hF_\ell)$ by adding isolated vertices and subdividing, with vertices of $V'$, edges whose endpoints are both in $V'$. In particular, \cref{lm:layouts-coverW-layouts,lm:layouts-coverW-shape} imply that $(L_1, F_1), \dots, (L_\ell,F_\ell)$ are still $W$-exceptional~$U^*$-path consistent layouts with respect to~$(D',\cF)$. Moreover, the number of non-trivial paths in $L$ is precisely $|\hL|=|\tE|$ and so \cref{lm:layouts-good} follows from \cref{eq:layouts-endpoints-size}. Furthermore, each $v\in W$ satisfies both $d_L^\pm(v)=d_{\hL}^\pm(v)$ and so \cref{lm:layouts-degreeW1,lm:layouts-degreeW2} follows from \cref{lm:layouts-coverW-degreeW1,lm:layouts-coverW-degreeW2}.}
			
			\NEW{We have already shown before \cref{step:endpoints} that \cref{lm:layouts-ell} holds. For each $i\in [\ell]$,
			\begin{align*}
				|V(L_i)|=|V(\hL_i)|+|V(L_i)\setminus V(\hL_i)|\stackrel{\text{\cref{lm:layouts-coverW-size},\cref{lm:layouts-adjustdegrees-sizeV}}}{\leq}5\sqrt{\varepsilon}n+\varepsilon^{\frac{1}{3}}n\leq 3\varepsilon^{\frac{1}{3}}n
			\end{align*}
			and 
			\begin{align*}
				|E(L_i)|&\stackrel{\text{\eqmakebox[ELi][c]{}}}{=}|E(\hL_i)|+|E(L_i)\setminus E(\hL_i)|\stackrel{\text{\cref{lm:layouts-coverW-size},\cref{lm:layouts-adjustdegrees-sizeE}}}{\leq}4\sqrt{\varepsilon}n+2|V(L_i)\setminus V(\hL_i)|\\
				&\stackrel{\text{\eqmakebox[ELi][c]{\cref{lm:layouts-adjustdegrees-sizeV}}}}{\leq }4\sqrt{\varepsilon}n+2\varepsilon^{\frac{1}{3}}n\leq 3\varepsilon^{\frac{1}{3}}n.
			\end{align*}
			Thus, \cref{lm:layouts-size} holds.}
				
			\NEW{We now verify \cref{lm:layouts-degreeV'}. Recall that $v_1, \dots, v_k$ is an enumeration of $V'$. Let $i\in [k]$. First, suppose that $n_i\geq 0$. Then, \cref{lm:layouts-adjustdegrees-isolatedvertex,lm:layouts-adjustdegrees-internalvertex} imply that $|\{j\in [\ell]\mid v_i\notin V(L_j)\}|=|\{j\in [\ell]\mid v_i\notin V(\hL_j)\}|-n_i$ and both $d_L^\pm(v_i)=d_{\hL}^\pm(v_i)$. Thus, \cref{lm:layouts-degreeV'} follows from \cref{eq:nipm}. Next, suppose that $n_i<0$. Then, \cref{lm:layouts-adjustdegrees-internalvertex,lm:layouts-adjustdegrees-isolatedvertex} imply that $\{j\in [\ell]\mid v_i\notin V(L_j)\}=\{j\in [\ell]\mid v_i\notin V(\hL_j)\}+n_i$ and both $d_L^\pm(v_i)=d_{\hL}^\pm(v_i)-2n_i$. Therefore, \cref{lm:layouts-degreeV'} follows from \cref{eq:nipm}.}
			
			\NEW{Finally, we check \cref{lm:layouts-V'}. Let $i\in [k]$. If $n_i\geq 0$, then
			\begin{align*}
				d_L(v_i)\stackrel{\text{\cref{lm:layouts-adjustdegrees-internalvertex},\cref{lm:layouts-adjustdegrees-isolatedvertex}}}{=}d_{\hL}(v)\stackrel{\text{\cref{lm:layouts-coverW-V'}}}{\leq}|\exc_D(v)|-\phi(v)+2+2|W|\stackrel{\text{\cref{lm:layouts-W},\cref{lm:layouts-degree&excessV'}}}{\leq}8\varepsilon n,
			\end{align*}
			as desired. If $n_i<0$, then
			\begin{align*}
				d_L(v_i)\stackrel{\text{\cref{lm:layouts-adjustdegrees-internalvertex},\cref{lm:layouts-adjustdegrees-isolatedvertex}}}{=}d_{\hL}(v_i)+2|n_i|\stackrel{\text{\cref{lm:layouts-coverW-V'},\cref{eq:ni}}}{\leq}(|\exc_D(v)|-\phi(v)+2+2|W|)+4\varepsilon n \stackrel{\text{\cref{lm:layouts-W},\cref{lm:layouts-degree&excessV'}}}{\leq}8\varepsilon n,
			\end{align*}
			as desired. Moreover, \cref{lm:layouts-adjustdegrees-internalvertex,lm:layouts-adjustdegrees-isolatedvertex} imply that the number of indices $j\in [\ell]$ for which $v_i\in V(L_j)$ is at most
			\begin{align*}
				d_L(v_i)+|n_i|\stackrel{\text{\cref{eq:ni}}}{\leq}8\varepsilon n+2\sqrt{\varepsilon}n\leq 3\sqrt{\varepsilon}n.
			\end{align*}
			Therefore, \cref{lm:layouts-V'} holds.}
		\end{proofclaim}
		
		\NEW{We may therefore assume that}\OLD{Assume} $m<k$. By \cref{eq:nipm}, $n_i=d_{\hL}^\pm(v_i)+|\{j\in [\ell]\mid v_i\notin V(\hL_j)\}|+\lceil\eta n\rceil-\phi^\mp(v_i)-d_D^\pm(v_i)$ and so we may suppose without loss of generality that~$v_{m+1}\notin U^-(D)$. Let~$X$ be the set of indices~$i\in [\ell]$ such that $|V(L_i^m)\setminus V(\hL_i)|=\lfloor\varepsilon^{\frac{1}{3}} n\rfloor$. \NEW{(Thus, $X$ is the set of indices $i\in [\ell]$ for which cannot modify $L_i^m$ anymore (otherwise \cref{lm:layouts-adjustdegrees-sizeV} would not be satisfied with $m+1$ playing the role of $m$).)}		
		Let~$Z$ be the set of indices~$i\in [\ell]$ such that~$v_{m+1}\in V(L_i^m)$. \NEW{(Thus, $Z$ is the set of indices $i\in [m]$ where $v_{m+1}$ cannot be added (because it is already present).)}	
		\NEW{Observe that
		\begin{align}\label{eq:layouts-Z}
			|Z|\leq d_{L^m}(v_{m+1})\stackrel{\text{\cref{lm:layouts-adjustdegrees-internalvertex},\cref{lm:layouts-adjustdegrees-isolatedvertex}}}{=} d_{\hL}(v_{m+1})\stackrel{\text{\cref{lm:layouts-coverW-V'}}}{\leq} |\exc_D(v_{m+1})|-\phi(v)+2+2|W| \stackrel{\text{\cref{lm:layouts-W},\cref{lm:layouts-degree&excessV'}}}{\leq} 4 \varepsilon n.
		\end{align}}\OLD{By \cref{lm:layouts-W}, \cref{lm:layouts-degree&excessV'}, \cref{lm:layouts-coverW-layouts}, \cref{lm:layouts-coverW-V'}, \cref{lm:layouts-adjustdegrees-internalvertex}, and \cref{lm:layouts-adjustdegrees-isolatedvertex}, $|Z|\leq d_{L^m}(v_{m+1})= d_{\hL}(v_{m+1})\leq |\exc_D(v_{m+1})|-\phi(v)+2+2|W| \leq 4 \varepsilon n$.}
		
		If~$n_{m+1}<0$, then proceed as follows. 		
		Let~$Y$ be the set of indices $i\in [\ell]\setminus (X \cup Z)$ such that $|L_i^m|\geq 2$. \NEW{(Thus, $Y$ is precisely the set of indices $i\in [\ell]$ for which we could incorporate $v_{m+1}$ as an internal vertices in two of the paths in $L_i^m$.)}
		We claim that $|Y| \ge -n_{m+1}$. 
		\NEW{By our choice of ordering $v_1,\dots, v_k$ of~$V'$, we have $n_i<0$ for each $i\in [m]$. Thus,
		\begin{align}\label{eq:layouts-X}
			|X|\stackrel{\text{\cref{lm:layouts-adjustdegrees-internalvertex},\cref{lm:layouts-adjustdegrees-isolatedvertex}}}{\leq} \frac{\sum_{i \in [\ell]}|V(L_i^m)\setminus V(\hL_i)|}{\lfloor\varepsilon^{\frac{1}{3}}n\rfloor} \stackrel{\text{\cref{lm:layouts-adjustdegrees-sizeV2}}}{=}\frac{\sum_{i \in [m]}|n_i|}{\lfloor\varepsilon^{\frac{1}{3}}n\rfloor}\stackrel{\text{\cref{eq:ni}}}{\leq}\frac{2\varepsilon n^2}{\lfloor\varepsilon^{\frac{1}{3}}n\rfloor} \leq 3\varepsilon^{\frac{2}{3}}n.
		\end{align}}%
		\OLD{By our choice of ordering $v_1,\dots, v_k$ of~$V'$,
			%\COMMENT{$n_i < 0$ for all $i \in [m]$ so $ |n_i| \le 2\eps n $}
		\cref{eq:ni}, and \cref{lm:layouts-adjustdegrees-sizeV2},
		we have $|X|\leq\frac{n\cdot 2\varepsilon n}{\lfloor \varepsilon^{\frac{1}{3}}n\rfloor}\leq 3\varepsilon^{\frac{2}{3}}n$%
			%\COMMENT{$|X|\leq\frac{n\cdot 2\varepsilon n}{\lfloor \varepsilon^{\frac{1}{3}}n\rfloor}\leq \frac{2\varepsilon n^2}{\varepsilon^{\frac{1}{3}}n-1}\leq \frac{3\varepsilon n}{\varepsilon^{\frac{1}{3}}}$.}
		.}%
		\NEW{If $\ell\neq\texc(D)-\lceil\eta n\rceil$, then
		\begin{align*}
			|Y|\stackrel{\text{\cref{lm:layouts-coverW-size2}}}{\geq} \sqrt{\varepsilon}n-|X|-|Z|\stackrel{\text{\cref{eq:layouts-Z},\cref{eq:layouts-X}}}{\geq} \sqrt{\varepsilon} n-3\varepsilon^{\frac{2}{3}}n-4\varepsilon n\stackrel{\text{\cref{eq:ni}}}{\geq} -n_{m+1}
		\end{align*}
		and so we are done. It is therefore enough to show that $\ell\neq\texc(D)-\lceil\eta n\rceil$. Suppose not.}%
		\OLD{Assume for a contradiction that $\ell=\texc(D)-\lceil\eta n\rceil$.} Then, by \cref{lm:layouts-coverW-size2}, $|\hL_i|=1$ for each~$i\in [\ell]$. 
		\NEW{Thus, $d_{\hL}^+(v_{m+1})$ is precisely the number of indices $i\in [\ell]$ for which $v_{m+1}$ is the starting point or an internal vertex of the unique path contained in $\hL_i$. Similarly, \cref{eq:layouts-endpoints-excess} and \cref{lm:layouts-coverW-shape} imply that there are precisely $\hexc^-(v_{m+1})$ indices $i\in [\ell]$ for which $v_{m+1}$ is the ending point of the unique path contained in $\hL_i$. Altogether, this implies that $d_{\hL}^+(v_{m+1})+|\{i\in [\ell]\mid v_{m+1}\notin V(\hL_i)\}|=\ell-\hexc^-(v_{m+1})$.}%
		\OLD{Thus, by \cref{eq:layouts-endpoints-excess} and \cref{lm:layouts-coverW-shape}, $d_{\hL}^+(v_{m+1})+|\{i\in [\ell]\mid v_{m+1}\notin V(\hL_i)\}|=\ell-\hexc^-(v_{m+1})$%
			%\COMMENT{Note that, since $|\hL_i|=1$ for each $i\in [\ell]$, $d_{\hL}^+(v_{m+1})$ counts the number of $\hL_i$ for which $v_{m+1}$ is either the starting point or a middle vertex. $\hexc^-(v_{m+1})$ counts the number of $\hL_i$ for which $v_{m+1}$ is the ending point.}
		.}
		Therefore,
		\begin{align*}
			n_{m+1}&\stackrel{\text{\eqmakebox[nm+1][c]{\text{\cref{eq:defni}}}}}{=} (\ell+\lceil\eta n\rceil)-(\hexc^-(v_{m+1})+\phi^-(v_{m+1}))-d_D^+(v_{m+1})\\
			&\stackrel{\text{\eqmakebox[nm+1][c]{\text{\cref{eq:layouts-hexc}}}}}{=}\texc(D)-\texc_{D, U^*}^-(v_{m+1})-d_D^+(v_{m+1})
			\stackrel{\text{\cref{eq:U*}}}{\geq} 0,
		\end{align*}
		a contradiction. \NEW{Consequently, $|Y|\geq -n_{m+1}$, as desired.}%
		\OLD{Therefore, by \cref{lm:layouts-coverW-size2} and \cref{eq:ni}, $|Y|\geq \sqrt{\varepsilon} n-3\varepsilon^{\frac{2}{3}}n-4\varepsilon n\geq -n_{m+1}$, as desired.}
		
		Let~$N_{m+1}\subseteq Y$ be such that~$|N_{m+1}|=-n_{m+1}$ and, for each~$i\in N_{m+1}$, fix two paths $P_{i,1}, P_{i,2}\in \hL_i^m$. For each~$i\in N_{m+1}$ and~$j\in [2]$, let $u_{i,j}w_{i,j}\in E(P_{i,j})\setminus E_W(P_{i,j})$, which exists by~\cref{lm:layouts-coverW-size} and \cref{lm:layouts-adjustdegrees-internalvertex}. 
		For each $i\in [\ell]\setminus N_{m+1}$, let $L_i^{m+1}\coloneqq L_i^m$. For each~$i\in N_{m+1}$, let~$L_i^{m+1}$ be obtained from~$L_i^m$ by subdividing, for each~$j\in [2]$, the edge~$u_{i,j}w_{i,j}$ in~$P_{i,j}$ into the path~$u_{i,j}v_{m+1}w_{i,j}$. Then, \cref{lm:layouts-adjustdegrees-internalvertex,lm:layouts-adjustdegrees-sizeV,lm:layouts-adjustdegrees-isolatedvertex} are satisfied with~$m+1$ playing the role of~$m$.
		
		If~$n_{m+1}\geq 0$, then proceed as follows. 
		Let $Y\coloneqq [\ell] \setminus (X \cup Z)$.%
		\OLD{We claim that~$|Y|\geq n_{m+1}$.}
		\NEW{Note that
		\begin{align}\label{eq:layouts-X2}
			|X|\stackrel{\text{\cref{lm:layouts-adjustdegrees-internalvertex},\cref{lm:layouts-adjustdegrees-isolatedvertex}}}{\leq} \frac{\sum_{i \in [\ell]}|V(L_i^m)\setminus V(\hL_i)|}{\lfloor\varepsilon^{\frac{1}{3}}n\rfloor} \stackrel{\text{\cref{lm:layouts-adjustdegrees-sizeV2}}}{=}\frac{\sum_{i \in [m]}|n_i|}{\lfloor\varepsilon^{\frac{1}{3}}n\rfloor}\stackrel{\text{\cref{eq:ni}}}{\leq}\frac{2\sqrt{\varepsilon} n^2}{\lfloor\varepsilon^{\frac{1}{3}}n\rfloor} \leq 3\varepsilon^{\frac{1}{6}}n
		\end{align}
		and so
		\begin{align*}
			|Y|\stackrel{\text{\cref{lm:layouts-d},\cref{lm:layouts-ell}}}{\geq}\eta n-|X|-|Z|\stackrel{\text{\cref{eq:layouts-Z},\cref{eq:layouts-X2}}}{\geq}\eta n-3\varepsilon^{\frac{1}{6}}n-4\varepsilon n\stackrel{\text{\cref{eq:ni}}}{\geq}n_{m+1}.
		\end{align*}}%
		\OLD{By \cref{lm:layouts-adjustdegrees-sizeV2} and \cref{eq:ni}, $|X|\leq \frac{2\sqrt{\varepsilon}n\cdot n}{\lfloor\varepsilon^{\frac{1}{3}}n\rfloor}\leq 3\varepsilon^{\frac{1}{6}}n$.
		Recall that $|Z| \le 4 \eps n$.
		Thus, by \cref{lm:layouts-d}, \cref{lm:layouts-ell}, and \cref{eq:ni}, $|Y|\geq \eta n-3\varepsilon^{\frac{1}{6}}n-4\varepsilon n\geq n_{m+1}$.}%
		Let $N_{m+1}\subseteq Y$ satisfy $|N_{m+1}|= n_{m+1}$. For each $i\in [\ell]\setminus N_{m+1}$, let $L_i^{m+1}\coloneqq L_i^m$ and, for each~$i\in N_{m+1}$, let~$L_i^{m+1}$ be obtained from~$L_i^m$ by adding~$v_{m+1}$ as an isolated vertex.
		Clearly, \cref{lm:layouts-adjustdegrees-isolatedvertex,lm:layouts-adjustdegrees-internalvertex,lm:layouts-adjustdegrees-sizeV} hold with~$m+1$ playing the role of~$m$, as desired.		
		\qedhere
	\end{steps}
\end{proof}

%% file: Concluding_Remarks.tex
	\onlyinsubfile{
		\setcounter{section}{13}
\section{Concluding remarks}}

\subsection{Approximate Hamilton decompositions of robust outexpanders}
In~\cite{osthus2013approximate}, Osthus and Staden showed that any regular robust outexpander of linear semidegree can be approximately decomposed into Hamilton cycles. This was used as a tool in~\cite{kuhn2013hamilton} to prove that such graphs actually have a Hamilton decomposition.

\begin{thm}[{\cite{osthus2013approximate}}]\label{thm:approxHamdecomp}
	Let $0<\frac{1}{n}\ll \tau \ll \alpha \leq 1$ and $0\leq \frac{1}{n}\ll\varepsilon\ll \nu, \eta\leq 1$. If~$D$ is an~$(\alpha,\varepsilon)$-almost regular robust~$(\nu,\tau)$-outexpander on~$n$ vertices, then~$D$ contains at least~$(\alpha-\eta)n$ edge-disjoint Hamilton cycles.
\end{thm}

\cref{lm:approxdecomp} can also be used to construct approximate Hamilton decompositions of (almost) regular robust outexpanders. In fact, our tools also enable us to assign some specific edges to each element of our approximate decomposition and so can be used to find approximate decompositions with prescribed edges.

\begin{thm}\label{cor:approxHamdecomp}
	Let $0<\frac{1}{n}\ll \tau\ll \alpha \leq 1$ and  $0<\frac{1}{n}\ll\varepsilon\ll\eta,\nu\leq 1$. Let $\ell \leq (\alpha-\eta)n$.	
	Suppose~$D$ is an~$(\alpha,\varepsilon)$-almost regular~$(\varepsilon,n^{-2})$-robust~$(\nu,\tau)$-outexpander on~$n$ vertices. Suppose that, for each~$i\in[\ell]$,~$F_i$ is a \NEW{linear forest}%
		\OLD{matching}
	on~$V(D)$%
		\COMMENT{fictive edges}
	\NEW{satisfying $e(F_i)\leq \varepsilon n$ and}%
		\OLD{of size at most~$\varepsilon n$}
	such that, for each~$v\in V(D)$, there exist at most~$\varepsilon n$ indices~$i\in[\ell]$ such that~$v\in V(F_i)$. Define a multiset~$\cF$ by $\cF\coloneqq \bigcup_{i\in [\ell]}F_i$.
	Then, there exist edge-disjoint Hamilton cycles $C_1,\dots, C_\ell \subseteq D\cup \cF$ such that, for each~$i\in [\ell]$,~$F_i\subseteq C_i$.
\end{thm}

\begin{proof}
	By \cref{lm:robparameters}, we may assume without loss of generality that
	\[0<\frac{1}{n}\ll\varepsilon\ll \nu\ll \tau\ll \eta\ll \alpha \leq 1.\]
	Define an additional constant~$\gamma$ such that $\tau\ll\gamma\ll \eta$. For each~$i\in [\ell]$, let $v_{i,1},v_{i,2}\in V(D)\setminus V(F_i)$ be distinct and such that, for any~$v\in V$, there exists at most two $(i,j)\in [\ell]\times [2]$ such that~$v=v_{i,j}$%
		\COMMENT{Construct an auxiliary bipartite graph $G$ on vertex classes $V(D)$ and $[n]$ with an edge $vi\in V(D)\times [n]$ if and only if $v\notin V(F_i)$, where $F_i\coloneqq \emptyset$ if $i>\ell$. Then, for each $i\in [n]$, $d_G(i)=n-|V(F_i)|\geq \frac{n}{2}$ and, for each $v\in V(D)$, $d_G(v)=n-d_\cF(v)\geq \frac{n}{2}$. By \cref{cor:Hall}, $G$ contains a perfect matching $M$. For each $i\in [n]$, denote by $v_{i,1}$ the neighbour of $i$ in $M$.\\
		Define the $v_{i,2}$ similarly. Let $G'$ be the bipartite graph on vertex classes $V(D)$ and $[n]$ with an edge $vi\in V(D)\times [n]$ if and only if $v\notin V(F_i)\cup \{v_{i,1}\}$, where $F_i\coloneqq \emptyset$ if $i>\ell$. Then, for each $i\in [n]$, $d_{G'}(i)=n-|V(F_i)\cup \{v_{i,1}\}|\geq \frac{n}{2}$ and, for each $v\in V(D)$, $d_{G'}(v)=n-d_\cF(v)-1\geq \frac{n}{2}$. By \cref{cor:Hall}, $G'$ contains a perfect matching $M'$. For each $i\in [n]$, denote by $v_{i,2}$ the neighbour of $i$ in $M$.}.
	\NEW{For each $i\in [\ell]$, denote by $P_{i,1}, \dots, P_{i,f_i}$ the (non-trivial) components of $F_i$ and, for $j\in [f_i]$, denote by $u_{i,j}$ and $w_{i,j}$ the starting and ending points of $P_{i,j}$.
	For each $i\in [\ell]$, let $L_i\coloneqq \{v_{i,1}u_{i,1}P_{i,1}w_{i,1}u_{i,2}P_{i,2}w_{i,2}u_{i,3}\dots w_{i,f_i}v_{i,2},v_{i,2}v_{i,1}\}$.}\OLD{For each $i\in [\ell]$, denote $F_i\coloneqq \{u_{i,1}w_{i,1}, \dots, u_{i,f_i}w_{i,f_i}\}$ and let $L_i\coloneqq \{v_{i,1}u_{i,1}w_{i,1}u_{i,2}\dots w_{i,f_i}v_{i,2},v_{i,2}v_{i,1}\}$.} Denote $L\coloneqq \bigcup_{i\in [\ell]}L_i$.
	Note that $(L_1,F_1), \dots, (L_\ell,F_\ell)$ are layouts such that, for each~$i\in [\ell]$, $V(L_i)\subseteq V$, $|V(L_i)|\leq 3\varepsilon n$ and $|E(L_i)|\leq 3\varepsilon n$. Moreover, for each~$v\in V(D)$, $d_L(v)\leq 3\varepsilon n$ and there exist at most~$2\varepsilon n$ indices~$i\in [\ell]$ such that~$v\in V(L_i)$.
		
	By similar arguments as in \cref{lm:Gamma}%
		\COMMENT{Same proof with $\frac{\gamma\nu}{2}$ playing the role of $\nu$ and last paragraph omitted.},
	there exists a spanning subdigraph~$\Gamma\subseteq D$ such that~$\Gamma$ is a~$(\gamma,\varepsilon)$-almost regular~$(\varepsilon,n^{-2})$-robust~$(\frac{\nu\gamma}{2}, \tau)$-outexpander and~$D'\coloneqq D\setminus \Gamma$ is~$(\alpha-\gamma, \varepsilon)$-almost regular.
	
	Apply \cref{lm:approxdecomp} with $D', \alpha-\gamma, \frac{\nu \gamma}{2}, \varepsilon^{\frac{1}{5}}$, and $\frac{\eta}{2}$ playing the roles of $D, \delta, \nu, \varepsilon$, and $\eta$ to obtain edge-disjoint $C_1,\cdots, C_\ell\subseteq D\cup \cF$ such that, for each~$i\in [\ell]$,~$C_i$ is a spanning configuration of shape~$(L_i, F_i)$. Then, by construction, for each~$i\in [\ell]$,~$C_i$ is a Hamilton cycle of~$D\cup \cF$ such that~$F_i\subseteq E(C_i)$.
\end{proof}

Recall that, by \cref{lm:approxdecomp}\cref{lm:approxdecomp-all}, the leftover from \cref{cor:approxHamdecomp} is actually still a robust~$(\frac{\nu\gamma}{4}, \tau)$-outexpander of linear minimum semidegree at least $\frac{\eta n}{2}$.
Thus, if~$D \cup \cF$ is regular, we can actually obtain a Hamilton decomposition of~$D \cup \cF$ so that for all
$i \in [\ell]$, the edges of
$F_i$ are contained in~$C_i$ (indeed, it suffices to apply to \cref{thm:robHamdecomp} to the leftover from \cref{cor:approxHamdecomp}).

Note that \cref{cor:approxHamdecomp} requires~$D$ to be an~$(\varepsilon,n^{-2})$-robust outexpander. One can show that this condition is in fact redundant and can be omitted. Indeed, K\"{u}hn, Osthus, and Treglown~\cite{kuhn2010hamiltonian} showed that the ``reduced digraph'' of a robust outexpander inherits the robust outexpansion properties of the host graph (see~\cite[Lemma 14]{kuhn2010hamiltonian}). Thus, using \cref{lm:Chernoff} and basic properties of ``$\varepsilon$-regular pairs'', one can easily show that the following \lcnamecref{lm:epsilonrob} holds. We omit the details. 

\begin{lm}\label{lm:epsilonrob}
	Let $0<\frac{1}{n}\ll \varepsilon\ll \nu'\ll \alpha, \nu, \tau\ll 1$. Suppose~$D$ is a robust~$(\nu,\tau)$-outexpander on~$n$ vertices satisfying $\delta^0(D)\geq \alpha n$. Then,~$D$ is an~$(\varepsilon,n^{-2})$-robust~$(\nu',4\tau)$-outexpander.
\end{lm}

Thus, \cref{cor:approxHamdecomp,lm:epsilonrob} imply \cref{thm:approxHamdecomp}. As the proof of \cref{cor:approxHamdecomp} only relies on \cref{lm:approxdecomp} (which in turn makes use of \cref{cor:robpaths} as the main tool), this gives a much shorter proof than the original one.

\subsection{A remark about Conjecture~\ref{conj:all}}

\cref{conj:all,thm:annoyingT} state that any (large) tournament~$T$ can be decomposed into at most~$\texc(T)+1$ paths. This cannot be generalised to digraphs or even oriented graphs. Indeed, it is easy to see that if~$D$ is a disconnected oriented graph then more than~$\texc(D)+1$ paths may be required to decompose~$D$%
	\COMMENT{E.g.\ if each component $C$ satisfies $\texc(C)=\Delta^0(C)\geq 2$.}.
In fact, \cref{conj:all,thm:annoyingT} cannot even be generalised to strongly connected oriented graphs.

\begin{prop}
	For any~$\varepsilon>0$ and $n_0\in \mathbb{N}$, there exists a strongly connected oriented graph~$D$ on~$n\geq n_0$ vertices such that $\pn(D)\geq \texc(D)+\frac{(1-\varepsilon)n}{2}$.
\end{prop}

\begin{proof}
	Fix additional integers~$m$ and~$k$ satisfying $0<\frac{1}{m}\ll\frac{1}{k}\ll \varepsilon$ and $m\geq n_0$.
	Let $V_1, \dots, V_k$ be disjoint sets of~$2m+1$ vertices each.
	For each $i\in [k]$, let~$T_i$ be a regular tournament on~$V_i$ and $x_iy_i\in E(T_i)$. Let~$D$ be obtained from $\bigcup_{i\in [k]} T_i$ by deleting, for each~$i\in [k]$, the edge~$x_iy_i$ and adding, for each~$i\in [k]$, the edge~$x_iy_{i+1}$, where~$y_{k+1}\coloneqq y_1$.
	Observe that~$D$ is a strongly connected~$m$-regular oriented graph on $n\coloneqq k(2m+1)$ vertices. Therefore, $\texc(D)=\Delta^0(D)=m$. Moreover, note that, for each~$i\in [k]$, $\pn(D[V_i])\geq \texc(D[V_i])=m$.
	
	Let~$\cP$ be a path decomposition of~$D$ of size~$\pn(D)$. For each~$i\in [k]$, let~$\cP_i$ be the set of paths $P\in \cP$ such that $V(P)\subseteq V_i$. Then, by construction, $|\cP_i|\geq \pn(D[V_i])-2\geq m-2$. Thus, $\pn(D)=|\cP|\geq k(m-2)=\texc(D)+(k-1)m-2k\geq \texc(D)+\frac{(1-\varepsilon)n}{2}$.
\end{proof}

\COMMENT{\begin{prop}
		For any~$n_0\in \mathbb{N}$, there exists a connected oriented graph~$D$ on~$n\geq n_0$ vertices such that $\pn(D)\geq \texc(D)+\frac{n}{3}-1$.
	\end{prop}
	\begin{proof}
		Let~$m$ be a sufficiently large odd integer and let~$T$ be a regular tournament on~$m$ vertices. Let~$S$ be a star on~$\frac{m+1}{2}$ vertices, all edges ending at the central vertex~$v$ of~$S$. Suppose $V(T)\cap V(S)=\emptyset$ and let~$D$ be obtained from~$T\cup S$ by adding an edge~$uv$ for some~$u\in V(T)$. 
		Let $n\coloneqq |V(D)|= \frac{3m+1}{2}$.
		Note that $\Delta^0(D)=d_D^+(u)=d_D^-(v)=\frac{m+1}{2}$ and $\exc(D)=\exc_D^-(v)=d_D^-(v)=\frac{m+1}{2}$, so $\texc(D)=\frac{m+1}{2}$. Let~$\cP$ be a path decomposition of~$D$ of size~$\pn(D)$. Then, for each~$e\in E(S)$,~$e\in \cP$ and $V(e)\cap \{u\}=\emptyset$. Thus, $\pn(D)=|\cP|\geq d_D^+(u)+|E(S)|= \texc(D)+ \frac{m-1}{2}\geq \texc(D)+\frac{n}{3}-1$, as desired.
\end{proof}}

\onlyinsubfile{\bibliographystyle{abbrv}
	\bibliography{Bibliography/Bibliography}}